%% file: principal2006_submit.tex
\documentclass{amsart}
\usepackage{amsmath,amsfonts,amssymb,amsthm,color}
\renewcommand{\emptyset}{\varnothing}
\renewcommand{\Im}{\operatorname{Im}}
\newcommand\union{\cup}
\newcommand\cv{``$\circ$''}

\renewcommand{\epsilon}{\varepsilon}
\newcommand\G{\mathbb G}

\newcommand\ZZ{{\mathbb Z}/2{\mathbb Z}}
\newcommand\Z{\mathbb Z}

\newcommand{\integers}{{\mathbb Z}}

\newcommand\R[1]{{\mathbb R}^{#1}}
\newcommand{\reals}{{\mathbb R}}

\newcommand\C[1]{{\mathbb C}^{#1}}
\newcommand{\cx}{{\mathbb C}}

\newcommand{\cB}{{\mathcal B}}
\newcommand{\cC}{{\mathcal C}}
\newcommand{\cH}{{\mathcal H}}

\newcommand{\cQ}{{\mathcal Q}}

\newlength{\halfbls}

\newtheorem{Theorem}{Theorem}
\newtheorem{theorem}[Theorem]{Theorem}
\newtheorem*{NNTheorem}{Theorem}

\newtheorem{prop}{Proposition}
\newtheorem{proposition}[prop]{Proposition}
\newtheorem{Proposition}[prop]{Proposition}
\newtheorem*{NNProposition}{Proposition}

\newtheorem{Lemma}{Lemma}[section]
\newtheorem{lemma}[Lemma]{Lemma}

\newtheorem{cor}{Corollary}

\newtheorem{Corollary}[cor]{Corollary}

\newtheorem*{MainTheorem}{Main Theorem}

\theoremstyle{remark}

\newtheorem{Remark}{Remark}

\newtheorem*{NNRemark}{Remark}

\newtheorem{Example}{Example}

\theoremstyle{definition}

\newtheorem{Definition}{Definition}
\newtheorem{definition}[Definition]{Definition}

\begin{document}


\title[Saddle connections and Boundary of the Moduli Space]
{Multiple Saddle  Connections  on  Flat  Surfaces  and  the Principal
Boundary of the Moduli Spaces of Quadratic Differentials}

\author{Howard Masur}
\thanks{Research  of  the first author is partially supported  by
NSF grant 0244472}
\address{
Department of Mathematics,
UIC, Chicago, IL 60607-7045, USA\\
}
\email{masur@math.uic.edu}

\author{Anton Zorich}
\address{
IRMAR,
Universit\'e Rennes-1,
Campus de Beaulieu,
35042 Rennes, cedex, France
}
\email{Anton.Zorich@univ-rennes1.fr}

\begin{abstract}
We describe typical degenerations of quadratic differentials thus
describing ``generic cusps''  of  the moduli space of meromorphic
quadratic differentials with at most  simple  poles.  The part of
the  boundary  of the  moduli  space which  does  not arise  from
``generic''  degenerations  is  often   negligible   in  problems
involving information on compactification of the moduli space.

However, even for a typical degeneration  one  may  have  several
short loops on  the  Riemann surface which shrink simultaneously.
We explain this phenomenon, describe all  rigid configurations of
short  loops,  present  a  detailed  description  of  analogs  of
desingularized stable curves arising  here,  and show how one can
reconstruct  a   Riemann   surface   endowed   with  a  quadratic
differential which is close to  a  ``cusp''  by the corresponding
point at the principal boundary.
\end{abstract}
\maketitle

\setcounter{tocdepth}{1}
\tableofcontents

\section*{Introduction}

\subsection*{Saddle connections on flat surfaces}

We study flat  metrics on a  closed orientable surface  of  genus
$g$,  which   have  isolated  conical  singularities  and  linear
holonomy restricted to $\{Id,-Id\}$. If the linear holonomy group
is trivial, then the surface is referred to as a {\em translation
surface},  such  a   flat   surface  corresponds  to  an  Abelian
differential $\omega$ on a Riemann surface. If the holonomy group
is nontrivial, then such a flat surface arises from a meromorphic
quadratic differential $q$ with at most simple poles on a Riemann
surface. In  this  paper,  unless  otherwise  stated, a quadratic
differential is {\em not}  the  square of an Abelian differential
and a {\em  flat  surface} is the Riemann  surface  with the flat
metric  corresponding   to   an   Abelian   or   to  a  quadratic
differential.

It is  natural to consider  families of flat surfaces sharing the
same  combinatorial   geometry:   the   genus,   the   number  of
singularities and the cone angles at singularities. Such families
correspond to the {\it strata} $\cQ(d_1,\dots,d_m)$ in the moduli
space       of       quadratic        differentials,        where
$d_i\in\{-1,0,1,2,3,\dots\}$   stands    for   the   orders    of
singularities (simple  poles, marked points, zeroes) of quadratic
differentials.  The  collection   $\alpha=\{d_1,\dots,d_m\}$   is
called the {\it singularity data} of the stratum.

A {\em saddle connection} is a geodesic segment joining a pair of
conical singularities or a conical singularity  to itself without
any  singularities  in  its  interior.  For  the flat metrics  as
described  above,  regular  closed  geodesics  always  appear  in
families; any such family fills a  maximal  cylinder  bounded  on
each side by a closed saddle connection or by a chain of parallel
saddle connections.  Thus,  when  some  regular  closed  geodesic
becomes short the corresponding saddle connection(s) become short
as well. More  generally,  a degeneration of an  Abelian  or of a
quadratic differential  corresponds  to  collapse  of some saddle
connections.

Any  saddle   connection  on  a  flat  surface  $S\in\cQ(\alpha)$
persists under small deformations of $S$ inside $\cQ(\alpha)$. It
might happen  that any deformation  of a given flat surface which
shortens  some  specific  saddle connection necessarily  shortens
some  other  saddle   connections.   We  say  that  a  collection
$\{\gamma_1,  \dots,   \gamma_n\}$   of   saddle  connections  is
\emph{rigid}  if  any sufficiently small deformation of the  flat
surface   inside   the   stratum   preserves   the    proportions
$|\gamma_1|:|\gamma_2|: \dots:|\gamma_n|$ of the  lengths  of all
saddle connections in the collection.

\subsection*{Degeneration of Abelian differentials}

In the case of Abelian differentials  $\omega$, rigid collections
of     saddle     connections     were     studied     in     the
paper~\cite{Eskin:Masur:Zorich}.  It  was  shown that all  saddle
connections in  any  rigid  collection  are  {\it homologous}. In
particular,  they  are  all  parallel and have equal  length  and
either  all  of them  join  the same  pair  of distinct  singular
points, or they are all closed.

This  implies  that  when  the  saddle  connections  in  a  rigid
collection  are  contracted  by  a  continuous  deformation,  the
limiting  flat  surface  generically   decomposes   into  several
components  represented  by nondegenerate  flat  surfaces  $S'_1,
\dots, S'_k$, where $k$  might vary from one to the genus  of the
initial surface. Let  $\cH(\beta'_j)$  be the stratum ambient for
$S'_j$. The stratum $\cH(\beta')=\cH(\beta'_1)\sqcup \dots \sqcup
\cH(\beta'_k)$   of  disconnected   flat   surfaces   $S'_1\sqcup
\dots\sqcup S'_k$ is referred to as  a  {\it  principal  boundary
stratum} of the stratum $\cH(\beta)$. For any connected component
of any  stratum $\cH(\beta)$ the  paper~\cite{Eskin:Masur:Zorich}
describes  all  principal boundary strata; their union is  called
the  \emph{principal  boundary} of  the  corresponding  connected
component of $\cH(\beta)$.

The  paper~\cite{Eskin:Masur:Zorich} also  presents  the  inverse
construction. Consider  any flat surface $S'_1\sqcup \dots \sqcup
S'_k\in  \cH(\beta')$   in   the   {\it  principal  boundary}  of
$\cH(\beta)$;  consider  a sufficiently small value of a  complex
parameter $\delta\in\C{}$.  One  can reconstruct the flat surface
$S\in\cH(\beta)$ endowed  with  a collection of homologous saddle
connections    $\gamma_1,    \dots,     \gamma_n$    such    that
$\int_{\gamma_i}\omega=\delta$, and such that the degeneration of
$S$  that   consists   of   contracting  the  saddle  connections
$\gamma_i$ in the collection gives the  surface $S'_1\sqcup \dots
\sqcup S'_k$.  This  inverse  construction  involves several {\it
basic surgeries} of the flat structure. Given a disconnected flat
surface $S'_1\sqcup \dots \sqcup S'_k$ one applies an appropriate
surgery to each $S'_j$  producing  a surface $S_j$ with boundary.
The surgery depends on the  parameter  $\delta$:  the boundary of
each  $S_j$  is  composed  of  two  geodesic segments of  lengths
$|\delta|$;  moreover,  the  boundary  components  of  $S_j$  and
$S_{j+1}$ are compatible, which  allows  one to glue the compound
surface $S$ from the collection of surfaces with boundary.

A collection $\gamma=\{\gamma_1, \dots, \gamma_n\}$ of homologous
saddle connections determines the following data on combinatorial
geometry of the decomposition $S\setminus \gamma$:  the number of
components, their  boundary  structure,  the singularity data for
each  component,  the cyclic order in which  the  components  are
glued   to   each   other.   These  data  are  referred   to   as
\emph{configuration}   of  homologous   saddle   connections.   A
configuration   $\cC$   uniquely  determines   the  corresponding
boundary stratum $\cH(\beta'_{\cC})$.

The constructions  above explain how configurations of homologous
saddle connections  on  flat  surfaces $S\in\cH(\beta)$ determine
the ``cusps''  of  the  stratum  $\cH(\beta)$.  Consider a subset
$\cH_1^{\varepsilon}(\beta)\subset\cH(\beta)$ of surfaces of area
one having a  saddle connection shorter than $\varepsilon$. Up to
a  subset  $\cH_1^{\varepsilon,thin}(\beta)$ of  negligibly small
measure the set $\cH_1^{\varepsilon}(\beta)$  can  be represented
as a  disjoint  union  over  all  admissible configurations $\cC$
(i.e. as a union over  different  ``cusps'')  of neighborhoods of
the  corresponding  ``cusps''.  When  a  configuration  $\cC$  is
composed  from  homologous saddle  connections  joining  distinct
zeroes,  the  neighborhood of  the  corresponding  cusp  has  the
structure  of  a  fiber  bundle over the  corresponding  boundary
stratum $\cH(\beta'_{\cC})$  with  the  fiber  represented  by an
appropriate ramified cover over the Euclidean $\varepsilon$-disc.
Moreover, the  canonical  measure  in the corresponding connected
component       of       $\cH_1^{\varepsilon,thick}(\beta)      =
\cH_1^{\varepsilon}(\beta)                              \setminus
\cH_1^{\varepsilon,thin}(\beta)$ decomposes as  a product measure
of  the  canonical  measure  in  the  boundary  stratum  and  the
Euclidean measure in the fiber, see~\cite{Eskin:Masur:Zorich}.

\begin{NNRemark}
We   warn   the  reader  that  the  correspondence  between   the
compactification of the moduli space of Abelian differentials and
the Deligne---Mumford compactification of  the  underlying moduli
space  of  curves is  not  straightforward.  In  particular,  the
desingularized stable  curve  corresponding  to the limiting flat
surface generically \emph{is not} represented  as  the  union  of
corresponding Riemann  surfaces  $S'_1,  \dots, S'_k$: the stable
curve might contain more components.
\end{NNRemark}

\section{Structure of the paper and statements of theorems}
\label{s:structure:of:the:paper}

This paper concerns the study of similar phenomena in the case of
quadratic  differentials   that   are   not  squares  of  Abelian
differentials.

\subsection{\^Homologous saddle connections}

A meromorphic  quadratic  differential  $q$  with  at most simple
poles on a Riemann surface $S$  defines  a  canonical  (ramified)
double cover  $p:\hat S\to  S$ such that $p^* q  = \omega^2$ is a
square  of  an  Abelian  differential $\omega$ on $\hat  S$.  Let
$P=\{P_1,  \dots,   P_m\}\subset   S$   be   the   collection  of
singularities  (zeroes  and  simple  poles)  of  $q$;  let  $\hat
P=p^{-1}(P)$ be the set  of  their preimages under the projection
$p:\hat S \to S$.

Given  an  oriented   saddle   connection  $\gamma$  on  $S$  let
$\gamma',\gamma''$  be  its   lifts   to  the  double  cover.  If
$[\gamma']=-[\gamma'']$ as cycles in $H_1(\hat S,\hat P;\,\Z)$ we
let $[\hat\gamma]:=[\gamma']$, otherwise we define $[\hat\gamma]$
as $[\hat\gamma]:=[\gamma']-[\gamma'']$.

\begin{definition}
\label{def:homologous}
The saddle connections $\gamma_1,\gamma_2$ on a  flat surface $S$
defined by a quadratic differential $q$ are
\^homologous\footnote{The  notion ``homologous  in  the  relative
homology with local coefficients defined by  the canonical double
cover induced by a quadratic differential''  is unbearably bulky,
so we introduced an abbreviation ``\^homologous''. We stress that
the  circumflex over  the  ``h'' is quite  meaningful:  as it  is
indicated  in   the  definition,  the  corresponding  cycles  are
homologous {\it on the double cover}.}
if $[\hat\gamma_1]=[\hat\gamma_2]$  in $H_1(\hat S,\hat  P;\,\Z)$
under  an   appropriate  choice  of  orientations  of  $\gamma_1,
\gamma_2$.
\end{definition}

\begin{figure}[htb]
%
\includegraphics{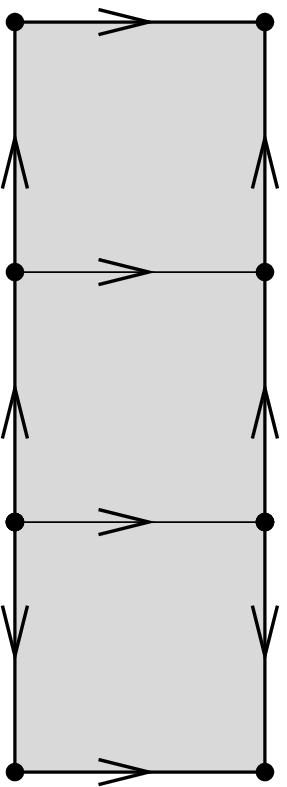}
\includegraphics{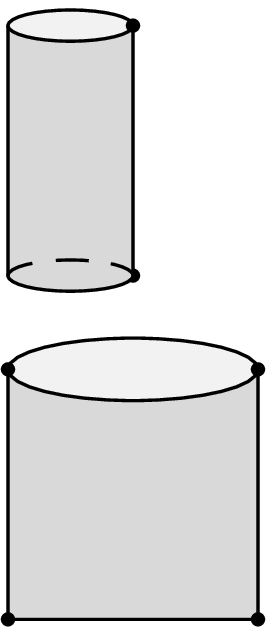}
\includegraphics{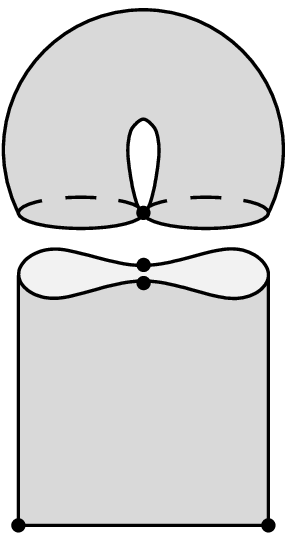}
\includegraphics{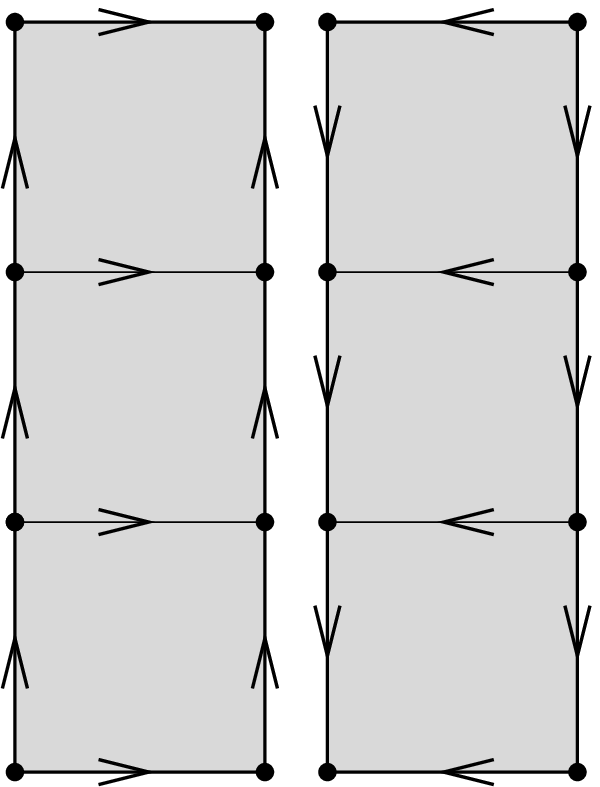}
\includegraphics{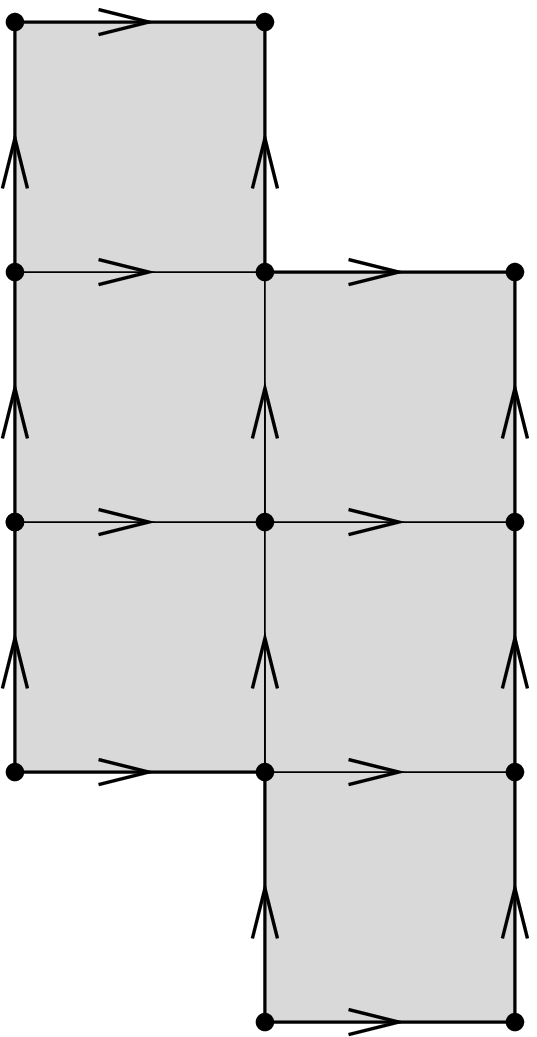}
\includegraphics{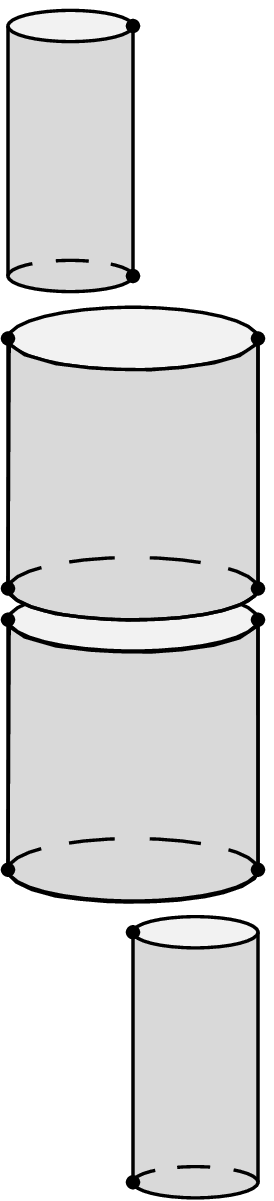}
\includegraphics{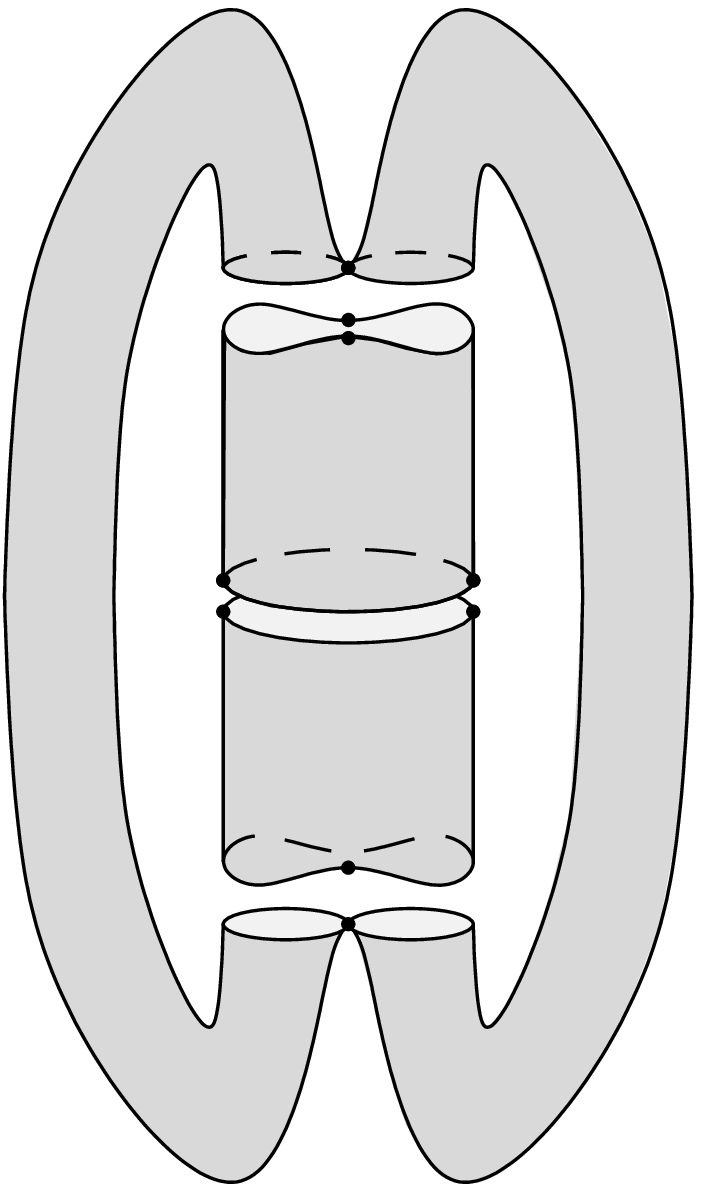}
\begin{picture}(0,0)(0,0)
\put(140,360)
{
\begin{picture}(0,-152)(0,-152) 
\put(-310,-526){\scriptsize $\gamma_3$}
\put(-332,-533){\scriptsize $P_0$}
\put(-286,-533){\scriptsize $P_0$}
\put(-318,-550){\scriptsize $\alpha$}
\put(-318,-583){\scriptsize $\beta$}
\put(-318,-615){\scriptsize $\beta$}
\put(-310,-637){\scriptsize $\gamma_3$}
\put(-332,-632){\scriptsize $P_0$}
\put(-286,-632){\scriptsize $P_0$}
\put(-297,-550){\scriptsize $\alpha$}
\put(-296,-583){\scriptsize $\delta$}
\put(-296,-615){\scriptsize $\delta$}
\put(-310,-559){\scriptsize $\gamma_2$}
\put(-332,-566){\scriptsize $P_0$}
\put(-286,-566){\scriptsize $P_0$}
\put(-310,-592){\scriptsize $\gamma_1$}
\put(-332,-599){\scriptsize $P_1$}
\put(-286,-599){\scriptsize $P_2$}
\end{picture}
\begin{picture}(4,0)(2,0) 
\put(-310,-526){\scriptsize $\gamma'_3$}
\put(-320,-550){\scriptsize $\alpha'$}
\put(-320,-583){\scriptsize $\beta'$}
\put(-321,-615){\scriptsize $-\beta''$}
\put(-310,-637){\scriptsize $\gamma_3'$}
\put(-298,-550){\scriptsize $\alpha'$}
\put(-298,-583){\scriptsize $\delta'$}
\put(-304,-615){\scriptsize $-\delta''$}
\put(-310,-559){\scriptsize $\gamma'_2$}
\put(-310,-592){\scriptsize $\gamma'_1$}
\end{picture}
\begin{picture}(-29,0)(-31,0) 
\put(-314,-526){\scriptsize $-\gamma''_3$}
\put(-321,-550){\scriptsize $-\alpha''$}
\put(-321,-583){\scriptsize $-\beta''$}
\put(-320,-615){\scriptsize $\beta'$}
\put(-314,-637){\scriptsize $-\gamma_3''$}
\put(-304,-550){\scriptsize $-\alpha''$}
\put(-304,-583){\scriptsize $-\delta''$}
\put(-297,-615){\scriptsize $\delta'$}
\put(-314,-559){\scriptsize $-\gamma''_2$}
\put(-314,-592){\scriptsize $-\gamma''_1$}
\end{picture}
\begin{picture}(-111,0)(-115,0) 
\put(-310,-526){\scriptsize $\gamma'_3$}
\put(-320,-550){\scriptsize $\alpha'$}
\put(-320,-583){\scriptsize $\beta'$}
\put(-323,-615){\scriptsize $-\beta''$}
\put(-310,-637){\scriptsize $\gamma_3'$}
\put(-298,-550){\scriptsize $\alpha'$}
\put(-298,-583){\scriptsize $\delta'$}
\put(-305,-615){\scriptsize $-\delta''$}
\put(-310,-559){\scriptsize $\gamma'_2$}
\put(-310,-592){\scriptsize $\gamma'_1$}
\end{picture}
\begin{picture}(-249,33)(-253,33) 
\put(-314,-526){\scriptsize $-\gamma''_3$}
\put(-321,-615){\scriptsize $-\alpha''$}
\put(-314,-637){\scriptsize $-\gamma_3''$}
\put(-297,-550){\scriptsize $\beta'$}
\put(-304,-583){\scriptsize $-\beta''$}
\put(-304,-615){\scriptsize $-\alpha''$}
\put(-314,-559){\scriptsize $-\gamma''_1$}
\put(-314,-592){\scriptsize $-\gamma''_2$}
\end{picture}
\begin{picture}(0,0)(-3,350) 
\put(252,-102){\scriptsize $\gamma_1$}
\put(229,-97){\scriptsize $P_1$}
\put(275,-97){\scriptsize $P_2$}
\put(229,-64){\scriptsize $P_0$}
\put(275,-64){\scriptsize $P_0$}
\put(258,-53){\scriptsize $P_0$}
\put(258,-20){\scriptsize $P_0$}
\put(257,-36){\scriptsize $\alpha$}
\put(244,-25){\scriptsize $\gamma_3$}
\put(244,-46){\scriptsize $\gamma_2$}
\put(252,-63){\scriptsize $\gamma_2$}
\put(252,-72){\scriptsize $\gamma_3$}
\put(241,-83){\scriptsize $\beta$}
\put(266,-83){\scriptsize $\delta$}
\end{picture}
\begin{picture}(0,0)(-97,350) 
\put(252,-102){\scriptsize $\gamma_1$}
\put(229,-97){\scriptsize $P_1$}
\put(275,-97){\scriptsize $P_2$}
\put(253,-39){\scriptsize $\alpha$}
\put(231,-59){\scriptsize $\gamma_2$}
\put(272,-60){\scriptsize $\gamma_3$}
\put(252,-71){\scriptsize $P_0$}
\put(241,-83){\scriptsize $\beta$}
\put(266,-83){\scriptsize $\delta$}
\end{picture}
\begin{picture}(0,0)(5,500) 
\put(244,-15){\scriptsize $\gamma'_3$}
\put(244,-48){\scriptsize $\gamma'_2$}
\put(258,-55){\scriptsize $-\gamma''_3$}
\put(244,-72){\scriptsize $\gamma'_2$}
\put(258,-88){\scriptsize $-\gamma''_1$}
\put(244,-108){\scriptsize $\gamma'_1$}
\put(258,-124){\scriptsize $-\gamma''_2$}
\put(244,-141){\scriptsize $\gamma'_3$}
\put(258,-147){\scriptsize $-\gamma''_2$}
\put(258,-180){\scriptsize $-\gamma''_3$}
\end{picture}
\begin{picture}(0,0)(-77,500) 
\put(253,-48){\scriptsize $\gamma'_2$}
\put(270,-48){\scriptsize $-\gamma''_3$}
\put(253,-72){\scriptsize $\gamma'_2$}
\put(270,-72){\scriptsize $-\gamma''_3$}
\put(261,-87){\scriptsize $-\gamma''_1$}
\put(261,-110){\scriptsize $\gamma'_1$}
\put(253,-125){\scriptsize $\gamma'_3$}
\put(270,-125){\scriptsize $-\gamma''_2$}
\put(253,-148){\scriptsize $\gamma'_3$}
\put(270,-148){\scriptsize $-\gamma''_2$}
\end{picture}
}
\end{picture}
\vspace{330bp}
\caption{
\label{fig:examples:of:homologous:sad:connections:1}
Saddle  connections  $\gamma_1,\gamma_2,\gamma_3$  on  the  torus
(above picture) are \^homologous, though $\gamma_1$  is a segment
joining distinct points and $\gamma_2$ and  $\gamma_3$ are closed
loops.
}
\end{figure}

We begin with the following example which illustrates many of the
main ideas.

\begin{Example}
\label{ex:homologous:sad:connections}
Consider three unit squares, or  rather  a  rectangle $1\times 3$
and glue a torus from it as indicated at the  top  left corner of
Figure~\ref{fig:examples:of:homologous:sad:connections:1}.
Identifying the three corresponding sides $\beta$, $\gamma_1$ and
$\delta$ of  the two bottom squares  we obtain a  ``pocket'' with
two  ``corners''  $P_1$ and  $P_2$  at the  bottom  and with  two
``corners'' $P_0$ at  the boundary on top. Identifying the points
$P_0$ we obtain a ``pocket''  with  a  ``figure-eight''  boundary
(the   bottom   fragment   of   the   top    right   picture   at
Figure~\ref{fig:examples:of:homologous:sad:connections:1}).
Identifying the sides  $\alpha$ of the remaining square we obtain
a cylinder which  we glue to the previous fragment. Topologically
the surface  thus obtained is  a torus. Metrically this torus has
three conical singularities. Two of  them  (``the  corners  $P_1,
P_2$ of  the pocket'') have cone  angle $\pi$; the  third conical
singularity $P_0$ has cone angle $4\pi$. Such a  flat torus gives
us a point in the stratum $\cQ(2,-1,-1)$.

The bottom picture illustrates the canonical double covering over
the  above  torus.   The   cycle  $\gamma'_2$  is  homologous  to
$\gamma'_3$ on the double cover and  the  cycle  $\gamma''_2$  is
homologous  to   $\gamma''_3$.   This  implies  that  the  cycles
$\hat\gamma_1$, $\hat\gamma_2$ and $\hat\gamma_3$  on  the double
cover are homologous  to the waist  curve of the  thick  cylinder
fragment   of   the  right  bottom  picture.  Thus,  the   saddle
connections   $\gamma_1$,    $\gamma_2$   and   $\gamma_3$    are
\^homologous, though  $\gamma_1$  is  a  segment joining distinct
points $P_1$ and $P_2$,  and  $\gamma_2, \gamma_3$ are the closed
loops with the base point $P_0$.
\end{Example}

It essentially  follows  from  the  definition  that \^homologous
saddle connections are parallel on $S$  and  that  their  lengths
either  coincide or  differ  by a factor  of  two. The  following
simple              statement              proved              in
appendix~\ref{ap:Long:saddle:connections}   characterizes   rigid
collections  of  saddle  connections  on  a   flat  surface  with
nontrivial linear holonomy.

\begin{Proposition}
\label{pr:rigid:configurations:hat:homologous}
Let  $S$  be  a  flat  surface  corresponding  to  a  meromorphic
quadratic  differential  $q$   with   at  most  simple  poles.  A
collection $\gamma_1,  \dots,  \gamma_n$ of saddle connections on
$S$ is rigid  if and only  if all saddle  connections  $\gamma_1,
\dots, \gamma_n$ are \^homologous.
\end{Proposition}

There  is  an  obvious  geometric test for deciding  when  saddle
connections $\gamma_1, \gamma_2$ on a translation surface $S$ are
homologous:  it   is  sufficient  to  check  whether  $S\setminus
(\gamma_1\cup\gamma_2)$   is    connected   or   not    (provided
$S\setminus\gamma_1$ and $S\setminus\gamma_2$  are connected). It
is slightly  less  obvious  to  check  whether saddle connections
$\gamma_1, \gamma_2$ on a flat surface $S$ with nontrivial linear
holonomy are \^homologous or not. In particular, a pair of closed
saddle connections might be homologous  in  the  usual sense, but
not \^homologous; a pair of  closed  saddle  connections might be
\^homologous even if one of them represents a  loop homologous to
zero,  and  the  other  does  not;  finally, a saddle  connection
joining  a  pair   of   {\it  distinct}  singularities  might  be
\^homologous to  a  saddle  connection  joining  a singularity to
itself.

Section~\ref{s:hat:homologous:saddle:connections}       describes
geometric criteria for  deciding  when two saddle connections are
\^homologous  and  what  is  the  structure   of  the  complement
$S\setminus(\gamma_1\cup\gamma_2)$.     These     criteria    are
intensively  used  in   the  remaining  part  of  the  paper.  In
particular,               we               prove               in
section~\ref{s:hat:homologous:saddle:connections}  the  following
statement.

\begin{Theorem}
\label{th:unique:trivial:holonomy}
Let  $S$  be  a  flat  surface  corresponding  to  a  meromorphic
quadratic differential $q$  with at most simple poles. Two saddle
connections $\gamma_1, \gamma_2$  on  $S$ are \^homologous if and
only  if  they have  no  interior intersections  and  one of  the
connected    components    of    the    complement    $S\setminus
(\gamma_1\cup\gamma_2)$ has trivial linear holonomy. Moreover, if
such a component exists, it is unique.
\end{Theorem}

\subsection{Graph of connected components}

A  collection   $\gamma$   of   \^homologous  saddle  connections
$\gamma=\{\gamma_1,\ldots,\gamma_n\}$  divides  $S$ into  simpler
surfaces   $S_j$   with  boundary.  We  associate  to  any   such
decomposition a  graph  $\Gamma(S,\gamma)$.  The  vertices of the
graph correspond to the connected components $S_j$ of $S\setminus
(\gamma_1\cup\dots\cup\gamma_n)$.   We    denote   the   vertices
corresponding to cylinders (if any) by small circles ``$\circ$''.
The remaining  vertices are labelled with  a ``$+$'' sign  if the
corresponding surface $S_j$ has trivial linear  holonomy and with
a  ``$-$''  sign  if it does not. We do not label the vertices of
``$\circ$''-type: it is easy to see  that  the  cylinders  always
have trivial linear holonomy.

The edges of the graph are in the  one-to-one correspondence with
the  saddle  connections  $\gamma_i$.   Each   saddle  connection
$\gamma_i$ is  on the boundary of either one  or two surfaces. If
$\gamma_i$ is on the boundary of pair of surfaces, it corresponds
to an edge joining  the  corresponding vertices. If $\gamma_i$ is
on the boundary of only  one  surface, then it corresponds to  an
edge of the graph which joins  the vertex to itself; such an edge
contributes $2$ to the valence of the vertex.

\begin{Remark}
\label{rm:dual:graph}
The   union   $\gamma=\gamma_1\cup\dots\cup\gamma_n$  of   saddle
connections can be considered  as  a graph $\gamma$ embedded into
the surface  $S$.  By  definition  $\Gamma(S,\gamma)$  is a graph
\emph{dual}  to  $\gamma$.  Namely,  $\Gamma(S,\gamma)$  can   be
realized as graph embedded into the surface $S$  in the following
way. A vertex of $\Gamma(S,\gamma)$ corresponding  to a connected
component $S_j$ of  $S\setminus\gamma$ is mapped to a point $v_j$
located  in  the  interior  of  the  corresponding  surface  with
boundary $S_j$. The  line representing the  image of an  edge  of
$\Gamma(S,\gamma)$   corresponding   to   a   saddle   connection
$\gamma_i$ has  a single transversal intersection with $\gamma_i$
in some  interior point; this  line does not intersect itself nor
any  other   such   line   nor   some   other  saddle  connection
$\gamma_{i'}$, where $i'\neq i$, in an interior point.
\end{Remark}

\begin{figure}[htb]
%
\includegraphics{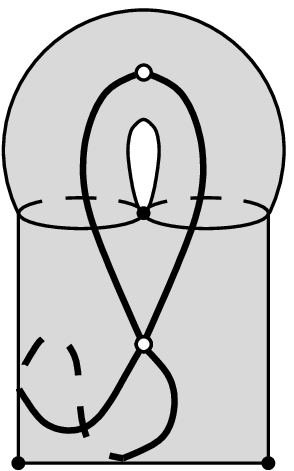}
\includegraphics{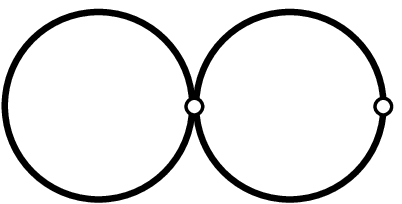}
\begin{picture}(0,0)(197,0)
\put(88,-35){$S_2$}
\put(90,-80){$S_1$}
\put(131,-14){$v_2$}
\put(138,-84){$v_1$}
\put(113,-61){$\gamma_2$}
\put(147,-61){$\gamma_3$}
\put(130,-118){$\gamma_1$}
\end{picture}
\begin{picture}(0,0)(66,0)
\put(131,-16){$v_2$}
\put(131,-75){$v_1$}
\put(100,-45){$\gamma_2$}
\put(159,-45){$\gamma_3$}
\put(131,-118){$\gamma_1$}
\end{picture}
\vspace{120bp}
\caption{
\label{fig:example:of:a:graph}
Graph $\Gamma(S,\gamma)$ of connected components
}
\end{figure}

\begin{Example}
\label{ex:graph:Gamma:for:fig:1}
Consider  the  surface   $S$   and  the  collection  $\gamma$  of
\^homologous    saddle    connections   $\{\gamma_1,    \gamma_2,
\gamma_3\}$  as  in   example~\ref{ex:homologous:sad:connections}
above                                                        (see
figure~\ref{fig:examples:of:homologous:sad:connections:1}).   The
complement $S\setminus\gamma$ has two connected components;  both
represented  by  flat  cylinders.  The  graph  $\Gamma(S,\gamma)$
contains  two  vertices, both of the ``$\circ$''-type, and  three
edges. The graph $\Gamma(S,\gamma)\subset S$ is dual to the graph
$\gamma\subset S$, see figure~\ref{fig:example:of:a:graph}.
\end{Example}

It follows from the definition of \^homologous saddle connections
that their lengths are  either the same or differ by a  factor of
two.  Having   a   collection  $\gamma$  of  \^homologous  saddle
connections $\gamma_1,  \dots,  \gamma_n$  we  can  normalize the
length  of  the  shortest  one  to  $1$. Then  the  other  saddle
connections  have lengths either  $1$  or  $2$, which endows
the edges of the graph $\Gamma$ with the weights $1$ or $2$.

The theorem below classifies all possible graphs corresponding to
nonempty collections of \^homologous saddle connections.

\begin{figure}[htb]
%
\includegraphics{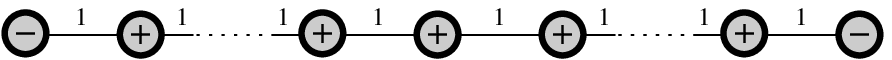}
   %
\includegraphics{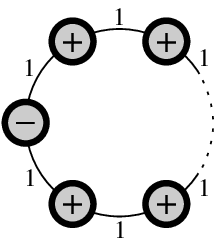}
   %
\includegraphics{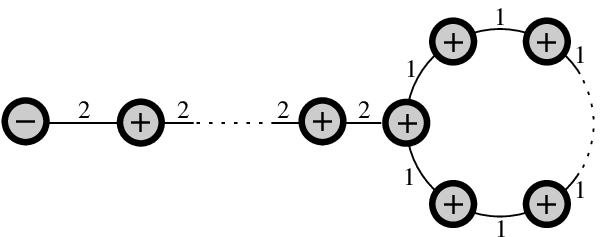}
\includegraphics{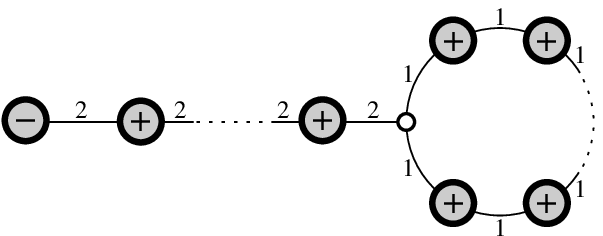}
   %
\includegraphics{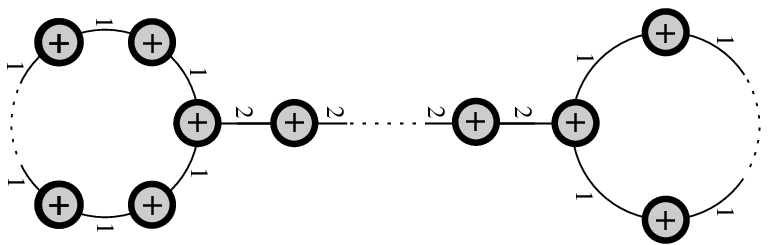}
\includegraphics{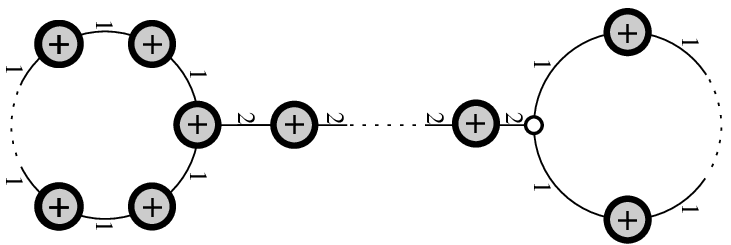}
\includegraphics{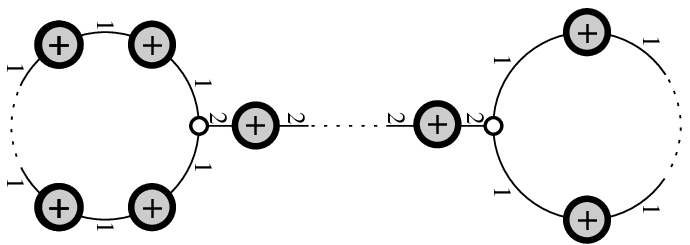}
   %
\includegraphics{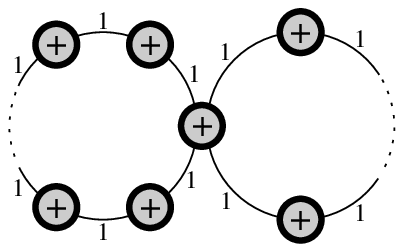}
\includegraphics{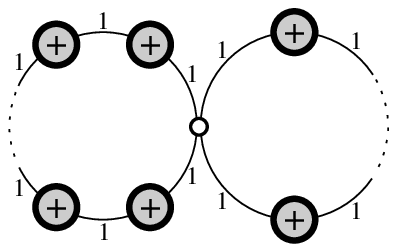}
%
%
%
\vspace{490bp} 
\begin{picture}(0,0)(0,0)
\put(-140,470.5) 
{
\begin{picture}(0,0)(0,0)
\put(38,-29){a)}
\put(285,-53){b)}
\put(132,-130){c)}
\put(132,-362){d)}
\put(132,-450){e)}
\end{picture}}
\end{picture}
\caption{
\label{fig:classification:of:graphs}
Classification of admissible graphs.
}
\end{figure}

\begin{theorem}
\label{th:graphs}
Let  $S$  be  a  flat  surface  corresponding  to  a  meromorphic
quadratic  differential  $q$  with  at  most  simple  poles;  let
$\gamma$ be  a  collection  of  \^homologous  saddle  connections
$\{\gamma_1, \dots,  \gamma_n\}$,  and  let $\Gamma(S,\gamma)$ be
the graph  of  connected  components  encoding  the decomposition
$S\setminus (\gamma_1\cup\dots\cup\gamma_n)$.

The graph $\Gamma(S,\gamma)$ either has  one  of  the basic types
listed below or  can  be  obtained from one of  these  graphs  by
placing  additional   \cv-vertices   of   valence   two   at  any
subcollection of edges subject to the  following restrictions. At
most one \cv-vertex may be placed at the same edge;  a \cv-vertex
cannot be  placed at an edge adjacent to  a \cv-vertex of valence
$3$ if this is the edge separating the graph.

The     graphs     of     basic      types,      presented     in
Figure~\ref{fig:classification:of:graphs},  are  given   by   the
following list:
\begin{itemize}
\item[a)]
An  arbitrary  (possibly  empty)  chain  of  ``$+$''-vertices  of
valence two bounded by a pair of ``$-$''-vertices of valence one;
\item[b)]
A  single loop  of  vertices of valence  two  having exactly  one
``$-$''-vertex and arbitrary number of ``$+$''-vertices (possibly
no ``$+$''-vertices at all);
\item[c)]
A single  chain and  a single loop joined at  a vertex of valence
three. The graph has  exactly  one ``$-$''-vertex of valence one;
it is  located at the  end of  the chain. The  vertex of  valence
three is either a ``$+$''-vertex, or a \cv-vertex  (vertex of the
cylinder  type).  Both the  chain,  and  the  cycle  may  have in
addition an arbitrary  number  of ``$+$''-vertices of valence two
(possibly no ``$+$''-vertices at all);
\item[d)]
Two nonintersecting  cycles joined by a  chain. The graph  has no
``$-$''-vertices. Each  of the two  cycles has a single vertex of
valence three (the one where the chain is attached to the cycle);
this vertex  is either a  ``$+$''-vertex or a \cv-vertex. If both
vertices of valence three are \cv-vertices, the chain joining two
cycles  is  nonempty:   it   has  at  least  one  ``$+$''-vertex.
Otherwise, each  of the cycles and  the chain may  have arbitrary
number  of   ``$+$''-vertices   of   valence   two  (possibly  no
``$+$''-vertices of valence two at all);
\item[e)]
``Figure-eight'' graph: two  cycles joined at a vertex of valence
four, which  is either a  ``$+$''-vertex or a \cv-vertex. All the
other vertices (if  any) are the ``$+$''-vertices of valence two.
Each  of  the  two  cycles  may  have arbitrary  number  of  such
``$+$''-vertices of  valence two (possibly no ``$+$''-vertices of
valence two at all).
\end{itemize}

Every graph listed above corresponds to some flat surface $S$ and
to some collection of saddle connections $\gamma$.
\end{theorem}

Theorem~\ref{th:graphs} is  proved in section~\ref{s:graph}  with
exception  of  the  final  statement on realizability,  which  is
proved                                                         in
sections~\ref{s:Local:Constructions}--\ref{s:Nonlocal:constructions}.

\subsection{Parities of boundary singularities}
\label{ss:Anatomy:of:decomposition:of:a:flat:surface}
In section~\ref{s:Anatomy}
we give  a detailed analysis of each connected component $S_j$
of  the  decomposition
$S\setminus\gamma$.

It  is  convenient  to  consider a closed surface  with  boundary
$S^{comp}_j$  canonically  associated  to  $S_j$  by  taking  the
natural {\it  compactification} of $S_j$. Note, that $S^{comp}_j$
need  not  be  the same as the {\it closure} of $S_j$ in $S$. For
example, if we cut a surface $S$ along a single saddle connection
$\gamma$ joining a pair of distinct  singularities  we  obtain  a
surface $S_1$ whose  compactification  is a surface with boundary
composed of a  pair of parallel  distinct geodesics of  the  same
length, while  the  closure  of  $S_1=S\setminus\gamma_1$  in $S$
coincides with $S$. The closure of $S_j$ in $S$ is  obtained from
the compactification  $S^{comp}_j$  of $S_j$ by identification of
some boundary points (if necessary), or by identification of some
boundary saddle connections (if necessary).

\subsection*{Ribbon graph}

Given a vertex  $v$  of a finite graph  $\Gamma$  consider a tree
$\Gamma_v$ obtained as a small neighborhood of $v$ in $\Gamma$ in
the natural topology  of a one-dimensional cell complex. The tree
$\Gamma_v$ together  with  the  canonical  mapping  of the graphs
$\Gamma_v\to\Gamma$ will be referred to as the \emph{boundary} of
$v$. The number of edges of $\Gamma_v$ is exactly the  valence of
$v$   (and  hence   is   at  most  $4$   for   the  graphs   from
figure~\ref{fig:classification:of:graphs}).

Suppose that the boundary of $S^{comp}_j$  has $r=r(j)$ connected
components (called for brevity \emph{boundary components}). Every
boundary  component  is  composed  of  a  closed chain of  saddle
connections $\gamma_{j_{i,1}}, \dots, \gamma_{j_{i,p(i)}}$, where
$1\le i \le  r$.  The case $p(i)=1$ is  not  excluded: a boundary
component might be  composed of a single saddle connection. The
canonical orientation of $S^{comp}_j$ determines the  orientation
of every boundary component $\cB_i$ of  $\partial S^{comp}_j$ and
hence determines the cyclic order
\begin{equation}
\label{eq:chain:i}
\to\gamma_{j_{i,1}}\to  \dots \to\gamma_{j_{i,p(i)}}\to
\end{equation}
on    every     such     chain;     by    convention    we    let
$j_{i,p(i)+1}:=j_{i,1}$. Thus,  we get a natural decomposition of
the  set of  edges  of $\Gamma_{v_j}$ into  a  disjoint union  of
subsets, each endowed with a cyclic order,
\begin{equation}
\label{eq:all:chains}
\{\to\gamma_{j_{1,1}}\to \gamma_{j_{1,2}}\to \dots \to\gamma_{j_{1,p(1)}}\to\}\
\sqcup\ \dots\ \sqcup\
\{\to\gamma_{j_{r,1}}\to  \dots \to\gamma_{j_{r,p(r)}}\to\}
\end{equation}
It  is  convenient to encode such combinatorial  structure  by  a
\emph{local ribbon  graph}  $\G_{v_j}$  which  is  defined in the
following way.

\begin{figure}[ht]
%
%
\includegraphics{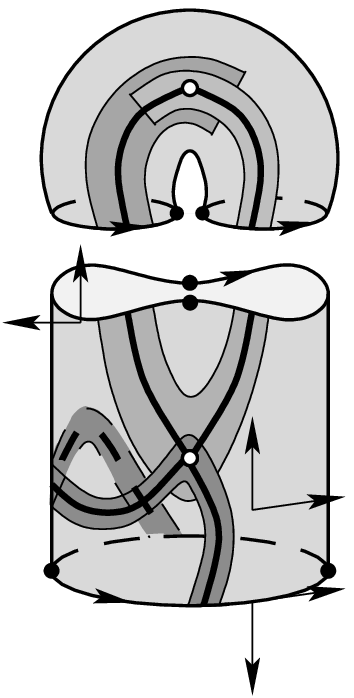}
\includegraphics{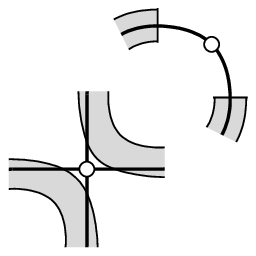}
\includegraphics{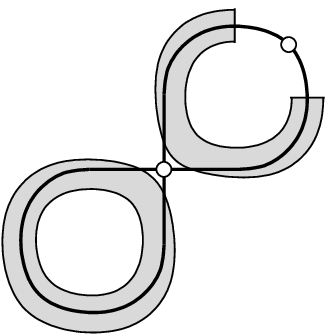}
\begin{picture}(0,0)(250,0)
\put(80,-35){$S^{comp}_2$}
\put(80,-105){$S^{comp}_1$}
\put(133,-53){$\gamma_2$}
\put(159,-53){$\gamma_3$}
\put(138,-127){$\gamma''_1$}
\put(148,-143){$\gamma'_1$}
\put(118,-54){$1$}
\put(109,-67){$2$}
\put(166,-100){$2$}
\put(174,-107){$1$}
\put(166,-150){$1$}
\put(174,-141){$2$}
\end{picture}
\begin{picture}(0,0)(130,0)
\put(133,-25){$\G_{v_2}$}
\put(111,-58){$\gamma_2$}
\put(154,-58){$\gamma_3$}
\put(109,-110){$\gamma''_1$}
\put(154,-110){$\gamma'_1$}
\put(133,-120){$\G_{v_1}$}
\end{picture}
\begin{picture}(0,0)(200,0)
\put(303,-15){$\G(S,\gamma)$}
\end{picture}
\vspace{150bp}
\caption{
\label{fig:example:of:ribbon:graph}
Compactifications   $S_1^{comp},   S^{comp}_2$    of    connected
components  of  $S\setminus\gamma$,  the associated local  ribbon
graphs  $\G_{v_1},   \G_{v_2}$   and   the  global  ribbon  graph
$\G(S,\gamma)$}
\end{figure}

Consider a realization of $\Gamma(S,\gamma)$ by an embedded graph
dual to the graph $\gamma$ in $S$ (see remark~\ref{rm:dual:graph}
above). For every vertex $v_j$  of  $\Gamma(S,\gamma)$  we get an
induced embedding $\Gamma_{v_j}\hookrightarrow S^{comp}_j$. Let a
connected   component   $\cB_i$   of  $\partial  S^{comp}_j$   be
represented by a chain~\eqref{eq:chain:i} of saddle  connections.
A  tubular  neighborhood  in  $S^{comp}_j$  of  the union of  the
corresponding      edges      $\{\gamma_{j_{i,1}}\cup       \dots
\cup\gamma_{j_{i,p(i)}}\}$  of  $\Gamma_{v_j}\subset  S^{comp}_j$
(as         in        the         left         picture         of
figure~\ref{fig:example:of:ribbon:graph}) inherits the  canonical
orientation of $S$. This orientation  induces  a  natural  cyclic
order on the  edges  $\gamma_{j_{i,1}},\dots,\gamma_{j_{i,p(i)}}$
of      $\Gamma_{v_j}$.      We     choose      the     embedding
$\Gamma_{v_j}\hookrightarrow S^{comp}_j$ in such way that turning
counterclockwise  around   $v_j$   (considered   as  a  point  of
$S_j^{comp}$)          we          see         the          edges
$\gamma_{j_{i,1}},\dots,\gamma_{j_{i,p(i)}}$ appear in the cyclic
order~\eqref{eq:chain:i}.

When  the   boundary   $\partial   S_j^{comp}$  contains  several
connected  components,   the   ribbon   graphs  corresponding  to
different components overlap at $v_j$ (as in the  left picture of
figure~\ref{fig:example:of:ribbon:graph}). However, it is easy to
make  them  disjoint  by  a  small  deformation,  subject  to  an
appropriate     choice      of     the     initial      embedding
$\Gamma_{v_j}\hookrightarrow S^{comp}_j$.  From  now  on we shall
always assume that the embedding is chosen appropriately.

\begin{Example}
\label{ex:ribbon:graph:Gamma:for:fig:1}
Consider once again the  surface  $S$ and the collection $\gamma$
of   \^homologous   saddle  connections   $\{\gamma_1,  \gamma_2,
\gamma_3\}$  as  in  example~\ref{ex:homologous:sad:connections},
see figure~\ref{fig:examples:of:homologous:sad:connections:1}. In
example~\ref{ex:graph:Gamma:for:fig:1}  we  have constructed  the
associated         graph         $\Gamma(S,\gamma)$,          see
figure~\ref{fig:example:of:a:graph}.

The complement $S\setminus\gamma$ has  two  connected components;
their compactifications $S_1^{comp},  S_2^{comp}$ are represented
by a pair of flat cylinders. Each of the two connected components
of  the   boundary   of   $S_2^{comp}$   (the   top  cylinder  in
figure~\ref{fig:example:of:ribbon:graph}) is  formed by a  single
saddle      connection,       so      we      get       $\partial
S_2^{comp}=\{\gamma_2\}\sqcup\{\gamma_3\}$.  Each  of   the   two
connected components of the boundary of  $S_1^{comp}$ (the bottom
cylinder in figure~\ref{fig:example:of:ribbon:graph})  is  formed
by  a  pair  of  saddle  connections,   so   we   get   $\partial
S_2^{comp}=\{\gamma_2\to           \gamma_3\}\sqcup\{\gamma'_1\to
\gamma''_1\}$. The orientation of the boundary components induced
by the canonical  orientation  of $S$  is  indicated in the  left
picture.

The        picture         in        the        center         of
figure~\ref{fig:example:of:ribbon:graph} shows the  corresponding
local ribbon graphs and the picture on the right shows the global
ribbon graph $\G(S,\gamma)$ for this example.
\end{Example}

\begin{figure}
%
\includegraphics{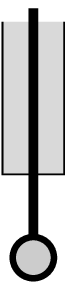}
   %
\includegraphics{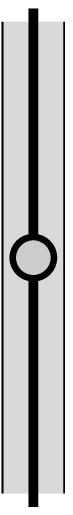}
   %
   %
\includegraphics{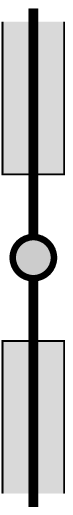}
   %
   %
\includegraphics{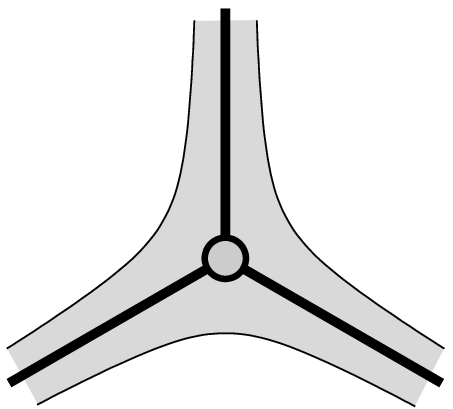}
\includegraphics{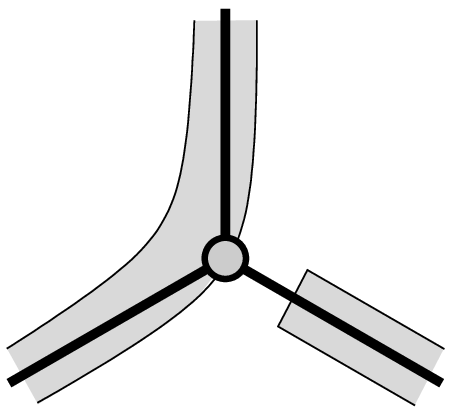}
\includegraphics{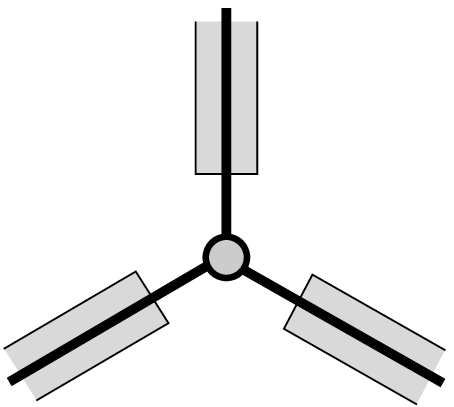}
   %
\includegraphics{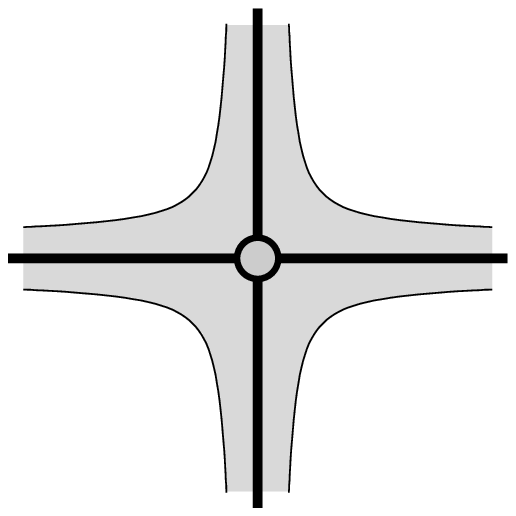}
\includegraphics{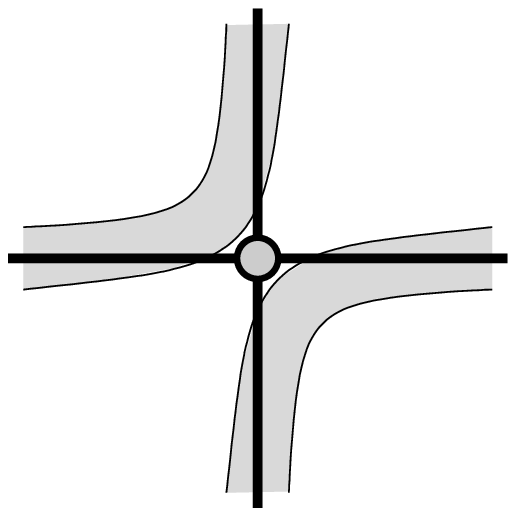}
\includegraphics{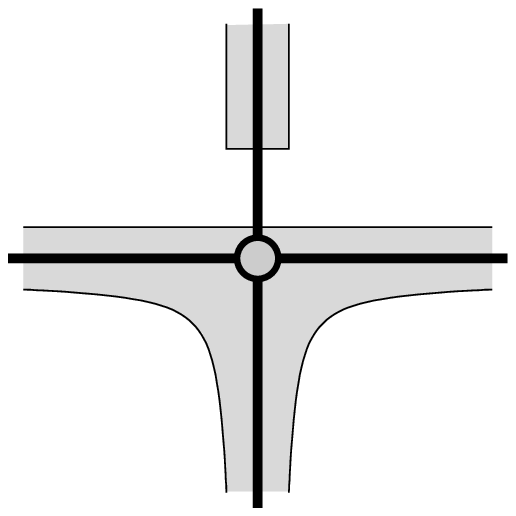}
\includegraphics{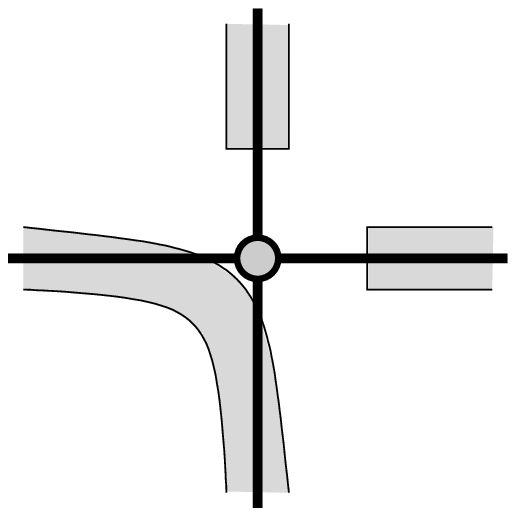}
\includegraphics{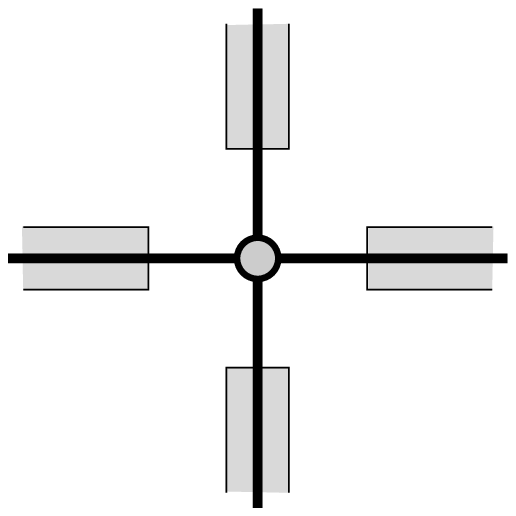}
%
%
\vspace{490bp} 
\begin{picture}(0,0)(0,-120) 
\put(-140,362)
{
\begin{picture}(0,0)(0,10)
\put(  0,   4){$\gamma_{j_{11}}$} 
\put( 72,   4){$\gamma_{j_{11}}$} 
\put(163,   4){$\gamma_{j_{12}}$}
\put(219,   4){$\gamma_{j_{11}}$} 
\put(308,   4){$\gamma_{j_{12}}$}
\end{picture}
\begin{picture}(0,-75)(0,-35)      
\put(-58, -115){$\gamma_{j_{13}}$}
\put( 20, -115){$\gamma_{j_{12}}$}
\put( -19, -175){$\gamma_{j_{11}}$}
\end{picture}
\begin{picture}(-102,-75)(-102,-35)
\put(-28, -115){$\gamma_{j_{21}}$} 
\put( 50, -115){$\gamma_{j_{22}}$}
\put( 11, -175){$\gamma_{j_{11}}$}
\end{picture}
\begin{picture}(-377,-75)(-377,-35) 
\put(-58, -115){$\gamma_{j_{31}}$}
\put( 20, -115){$\gamma_{j_{21}}$}
\put( -19, -175){$\gamma_{j_{11}}$}
\end{picture}
\begin{picture}(0,0)(0,80)   
\put(407, -175){$\gamma_{j_{14}}$}
\put(450, -135){$\gamma_{j_{13}}$}
\put(495, -175){$\gamma_{j_{12}}$}
\put(450, -216){$\gamma_{j_{11}}$}
\end{picture}
\begin{picture}(0,0)(0,80)  
\put(537, -175){$\gamma_{j_{22}}$}
\put(580, -135){$\gamma_{j_{21}}$}
\put(625, -175){$\gamma_{j_{12}}$}
\put(580, -216){$\gamma_{j_{11}}$}
\end{picture}
\begin{picture}(0,0)(0,80)  
\put(681, -175){$\gamma_{j_{23}}$}
\put(722, -135){$\gamma_{j_{11}}$} %
\put(768, -175){$\gamma_{j_{22}}$}
\put(722, -216){$\gamma_{j_{21}}$}
\end{picture}
\begin{picture}(9,0)(9,115)  
\put(407, -295){$\gamma_{j_{32}}$}
\put(450, -255){$\gamma_{j_{11}}$}
\put(495, -295){$\gamma_{j_{21}}$}
\put(450, -336){$\gamma_{j_{31}}$}
\end{picture}
\begin{picture}(16,0)(16,115) 
\put(681, -295){$\gamma_{j_{41}}$}
\put(722, -255){$\gamma_{j_{31}}$}
\put(768, -295){$\gamma_{j_{21}}$}
\put(722, -336){$\gamma_{j_{11}}$}
\end{picture}
\begin{picture}(-420,-55)(-420,-40) 
\put(-40, -70){$\{\to \gamma_{j_{11}} \to\}$}
\put( 80, -70){$\{\to \gamma_{j_{11}} \to \gamma_{j_{12}} \to \}$}
\put(215, -70){$\{\to \gamma_{j_{11}} \to\} \sqcup \{\to \gamma_{j_{21}} \to\}$}
\put(-45,-205){$\{\gamma_{j_{11}}\!\to\!\gamma_{j_{12}}\!\to\!\gamma_{j_{13}}\}$}
\put( 80,-205){$\{\gamma_{j_{11}}\!\to\} \sqcup \{\gamma_{j_{21}}\!\to\gamma_{j_{22}}\}$}
\put(225,-205){$\{\gamma_{j_{11}}\} \sqcup \{\gamma_{j_{21}}\} \sqcup \{\gamma_{j_{31}}\}$}
        %
\put(-45,-360){$\{\gamma_{j_{11}}\!\!\to\! \gamma_{j_{12}}\!\!\to\!\gamma_{j_{13}}\!\!\to\!\gamma_{j_{14}}\}$}
\put( 80,-360){$\{\gamma_{j_{11}}\!\!\to\!\gamma_{j_{12}}\}\!\sqcup\!\{\gamma_{j_{21}}\!\!\to\gamma_{j_{22}}\}$}
\put(210,-360){$\{\gamma_{j_{11}}\}\!\sqcup\!\{\gamma_{j_{21}}\!\!\to\!\gamma_{j_{22}}\!\!\to\gamma_{j_{23}}\}$}
        %
\put(-45,-515){$\{\gamma_{j_{11}}\}\sqcup\{\gamma_{j_{21}}\}\sqcup\{\gamma_{j_{31}}\!\to\!\gamma_{j_{32}}\}$}
\put(185,-515){$\{\gamma_{j_{11}}\}\sqcup\{\gamma_{j_{21}}\}\sqcup\{\gamma_{j_{31}}\}\sqcup\{\gamma_{j_{41}}\}$}
\end{picture}}
\end{picture}
\caption{
\label{fig:local:ribbon:graphs}
All local ribbon  graphs  $\G_v$ of  valences  from one to  four}
\end{figure}

For       vertices       $v$      of       valence      $1,2,3,4$
figure~\ref{fig:local:ribbon:graphs} gives a complete list of all
possible partitions of the edges  of  $\Gamma_v$  into a disjoint
union  of  subsets  endowed  with  a  cyclic  order  and  of  the
corresponding local ribbon graphs $\G_v$. Note that the canonical
orientation of $S$  induces  the counterclockwise ordering of the
edges of $\Gamma_v$.

\subsection*{Boundary singularities}
Let  $S_j$  be   a   connected  component  of  the  decomposition
$S\setminus (\gamma_1\cup\dots\cup\gamma_n)$; let $S_j^{comp}$ be
its  compactification,  and let a connected component $\cB_i$  of
the  boundary   $\partial   S_j^{comp}$   be   represented  by  a
chain~\eqref{eq:chain:i}  of   saddle  connections.  The   common
endpoint of  $\gamma_{j_i}$  and $\gamma_{j_{i+1}}$ is called the
{\it boundary  singularity}  of  $S_j^{comp}$.  Since  all saddle
connections  $\gamma_1,  \dots,  \gamma_n$   are   parallel,  the
corresponding angle between $\gamma_{j_i}$ and $\gamma_{j_{i+1}}$
is an integer  multiple  of $\pi$.  There  might be also  several
conical singularities in  the  interior of $S_j^{comp}$; they are
called {\it interior singularities}.

\begin{Definition}
\label{def:order:of:boundary:singularity}
If the  total angle at a  boundary singularity is  $(k+1)\pi$ the
{\it order of the boundary singularity} is defined to be $k$, and
the {\it parity of the boundary singularity} is defined to be the
parity of  $k$. If the total angle at  an interior singularity is
$(d+2)\pi$ the {\it order of the interior singularity} is defined
to be $d$.
\end{Definition}

The order of the interior singularity coincides with the order of
the zero  (simple pole) of  the corresponding germ of a quadratic
differential. By  convention,  boundary  singularities, and their
orders will always refer to the compactification $S^{comp}_j$.

When  $S_j$  is
represented by  a ``$+$''-vertex of the graph $\Gamma(S,\gamma)$,
we  include  the parities of the boundary  singularities  in  our
combinatorial structure  represented by the embedded local ribbon
graph $\G_{v_j}$. Let $\cB_i$ be a  connected  component  of  the
boundary     $\partial     S^{comp}_j$    constituted     by    a
chain~\eqref{eq:chain:i}   of   saddle  connections.   The  edges
$\gamma_{j_{i,1}},\dots,\gamma_{j_{i,p(i)}}$   of   the  embedded
graph   $\Gamma_{v_j}\hookrightarrow   S^{comp}_j$  subdivide   a
neighborhood  of  $v_j$  in  $S_j$ into $p(i)$ sectors.  To  each
sector bounded by a pair of  consecutive edges $\gamma_{j_{i,l}}$
and $\gamma_{j_{i,l+1}}$ we associate the  parity  of  the  order
$k_{j_{i,l}}$  of  the  corresponding  boundary  singularity   of
$S^{comp}_j$:  of  the common endpoint of the consecutive  saddle
connections $\gamma_{j_{i,l}}\to\gamma_{j_{i,l+1}}$ in $\cB_i$.

Any connected  component  $S_j$  of the decomposition $S\setminus
\{\gamma_1,   \dots,   \gamma_n\}$   determines   the   following
combinatorial data  which  we  refer  to  as the \textit{boundary
type} of  $S_j$: the structure~\eqref{eq:all:chains} of the local
ribbon          graph          $\G_{v_j}$          as          in
figure~\ref{fig:local:ribbon:graphs};        an         embedding
$\Gamma_{v_j}\hookrightarrow\Gamma(S,\gamma)$ and a collection of
parities of boundary singularities of $S_j$.

\begin{Theorem}
\label{th:all:local:ribon:graphs}
Consider  a  decomposition  of  a  flat   surface   $S$   as   in
theorem~\ref{th:graphs}. Every connected component  $S_j$  of the
decomposition  has  one  of  the  boundary   types  presented  in
figure~\ref{fig:embedded:local:ribbon:graphs}  and  all indicated
boundary types are realizable.
\end{Theorem}

The dotted lines in figure~\ref{fig:embedded:local:ribbon:graphs}
indicate pairs of edges of  a  vertex  $v\in\Gamma(S,\gamma)$  of
valence $3$ or  $4$, which  are joined  by  a loop  in the  graph
$\Gamma(S,\gamma)$                                           (see
figure~\ref{fig:classification:of:graphs}) and encode in this way
the embedding $\Gamma_{v_j}\hookrightarrow\Gamma(S,\gamma)$.

\begin{Remark}
\label{rm:indexation}
We use the following convention on  indexation  of  local  ribbon
graphs   in   figure~\ref{fig:embedded:local:ribbon:graphs}:  the
first  symbol   represents   the   type   (``$+$'',  ``$-$'',  or
``$\circ$'') of the vertex $v_j$ in the graph $\Gamma(S,\gamma)$;
the second symbol is the valence of $v_j$; the number after a dot
is the  number of boundary components  of $S_j$. An  extra letter
``$a,b,c$''  is  employed  when  it is necessary  to  distinguish
different  embedded  local ribbon graphs sharing the same  vertex
type, valence and number of boundary components.
\end{Remark}

\input{ribbon_graphs_submit.tex}

The      first      part      of      the      statement       of
theorem~\ref{th:all:local:ribon:graphs} which claims  that  every
connected component of  the decomposition has one of the boundary
types in figure~\ref{fig:embedded:local:ribbon:graphs}  is  quite
elementary;    it     is     proved     at    the     end     of
section~\ref{s:Anatomy}. The  statement
about  the  realizability of  all  boundary  types  presented  in
figure~\ref{fig:embedded:local:ribbon:graphs}   is   much    less
trivial;               it               follows              from
theorem~\ref{th:from:boundary:to:neighborhood:of:the:cusp}  which
is      proved      in       sections~\ref{s:Local:Constructions}
and~\ref{s:Nonlocal:constructions}.

\subsection{Configurations of \^homologous saddle connections}
We   formalize   the  data  on  combinatorial  geometry  of   the
decomposition           $S\setminus           \gamma$          in
definition~\ref{def:configuration} below.

\begin{Definition}
\label{def:configuration}
The    following    combinatorial    structure   is   called    a
\emph{configuration of \^homologous saddle connections}.
\begin{enumerate}
\item
A finite  graph $\Gamma$ endowed  with a labelling of each vertex
by one  of the symbols ``$+$'', ``$-$, or  ``$\circ$'', of one of
one  of  the  types  described  in  theorem~\ref{th:graphs}  (see
figure~\ref{fig:classification:of:graphs}).
\item
For any vertex $v$ of the graph $\Gamma$ an embedded ribbon graph
$\G_v$ (encoding the decomposition of $\Gamma_v$  into a disjoint
union of subsets, called \emph{boundary components}, each endowed
with a  cyclic  order; see equation~\eqref{eq:all:chains}) of one
of the types described in theorem~\ref{th:all:local:ribon:graphs}
(see figure~\ref{fig:embedded:local:ribbon:graphs}).
\item
For every  ``$+$''-vertex $v$ of $\Gamma$  and for every  pair of
consecutive    elements    $\gamma_{i,l}\to\gamma_{i,l+1}$     of
$\G_v$ (called  \emph{boundary singularities}) an  associated
parity         (even        or         odd)         as         in
figure~\ref{fig:embedded:local:ribbon:graphs}.
\item
For  every  vertex   $v$  of  $\Gamma$  and  for  every  boundary
singularity of  $\G_v$  a nonnegative integer $k_{i,l}$ (referred
to  as  the \emph{order of the boundary singularity})  satisfying
the following conditions.  The  order of the boundary singularity
has the  same parity as  the parity of the corresponding boundary
singularity when $v$  is of the  ``$+$''-type; the order  of  any
boundary singularity of any vertex  of  the  ``$\circ$''-type  is
equal to zero. The sum $D_i+2=k_{i,1}+\dots+k_{i,p(i)}$ of orders
of boundary singularities along any boundary component $\cB_i$ of
$v$  satisfies  $ D_i\ge 0$ for  a  vertex of ``$+$''-type and  $
D_i\ge -1$ for a vertex of ``$-$''-type.
\item
For every vertex $v$ of $\Gamma$  of  ``$-$''-type  an  unordered
(possibly   empty)   collection  of   integers   $\{d_1,   \dots,
d_{s(v)}\}$, where  $d_l\in\{-1,1,2,3,\dots\}$; for every  vertex
$v$ of $\Gamma$ of ``$+$''-type  an  unordered  (possibly  empty)
collection of  positive even integers $\{d_1, \dots, d_{s(v)}\}$,
where $d_l\in\{2,4,\dots\}$.  In  both cases these collections of
integers (called \emph{orders of interior singularities}) satisfy
the following  compatibility  conditions  with  the collection of
boundary singularities of $\G_v$:
$$
-4\ \le\ \big( \sum d_l + \sum D_i \big) \ \equiv\ 0\mod 4\,,
$$
where
the  first sum is  taken over all interior singularities and
the second sum is taken over all boundary components $\cB_i$
of $\G_v$.
\item
When the vertex $v$ is of the ``$-$''-type  the couple [unordered
collection  of  interior  singularities, unordered collection  of
boundary singularities] is in  addition  not allowed to belong to
the following exceptional list:

$$
\begin{array}{|c|c|}
\hline
   %
\includegraphics{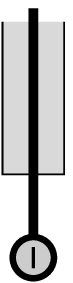}
&
[\emptyset,\{2\}];\qquad\quad
[\{-1\},\{3\}];\quad
[\{1\},\{1\}];\quad
[\{-1,1\},\{2\}]\\
\phantom{\quad\text{Type}\quad}
&
[\{1\},\{5\}];\quad
[\{3\},\{3\}];\quad
[\{1,3\},\{2\}];\qquad\quad
[\emptyset,\{6\}];\quad
[\{4\},\{2\}]
\\
\hline 
\includegraphics{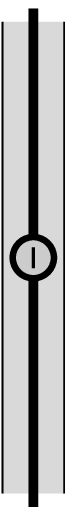}
&
[\emptyset,\{2, 0\}];\quad [\emptyset, \{1, 1\}];\qquad\quad
[\{-1\},\{0,3\}];\quad          [\{-1\},\{1,2\}]\\
   &
[\{1\},\{0,1\}];\quad
[\{1,-1\},\{0,2\}];\quad  [\{1,-1\},\{1,1\}]\\
&
[\{3,1\},\{2,0\}];\   [\{3,1\},\{1,1\}];\  [\{3\},\{3,0\}];\    [\{3\},\{2,1\}]\\
&
[\{1\},\{5,0\}];\      [\{1\},\{4,1\}];\    [\{1\},\{3,2\}];\
[\{4\},\{2,0\}];\           [\{4\},\{1,1\}]\\
&
[\emptyset,\{6,0\}];\       [\emptyset,\{5,1\}];\
[\emptyset,\{4,2\}];\       [\emptyset,\{3,3\}]\\
\hline 
\includegraphics{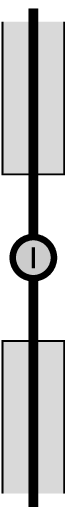}
&
[\emptyset,\{2,2\}];\
[\emptyset,\{1,3\}];\
[\{-1\},\{2,3\}];\
[\{1\},\{1,2\}];\
[\{-1,1\},\{2,2\}]
\\
&
[\emptyset,\{3,5\}];\quad
[\{1\},\{2,5\}];\quad
[\{3\},\{2,3\}];\quad
[\{1,3\},\{2,2\}]\\
&
[\emptyset,\{2,6\}];\quad
[\{4\},\{2,2\}]
\\
\hline
\end{array}
$$
\end{enumerate}
\end{Definition}

\subsection*{Singularity data corresponding to a configuration}
Any two flat surfaces  realizing  the same configuration $\cC$ of
\^homologous  saddle  connections  belong  to  the  same  stratum
$\cQ(\alpha)$ of  quadratic  differentials.  The singularity data
$\alpha$ are defined by the configuration $\cC$ as follows.

First note  that  any  configuration  $\cC$  determines a natural
\emph{global ribbon graph} $\G$ in  the  following  way.  We
have  defined  a structure of a  local  ribbon graph for a  small
neighborhood $\Gamma_v$  of  every  vertex $v\in\Gamma$. For
every vertex $v$ of $\Gamma$ we  have a ribbon going along a germ
of every edge of $\Gamma_v\subset\Gamma$ in direction from $v$ to
the center  of the edge. Note that all  local ribbon graphs carry
the canonical  orientation induced from the canonical orientation
of the  embodying plane. For every  edge of $\Gamma$  we can
extend the  ribbons from the  endpoints towards the center of the
edge and glue them together respecting the canonical orientation.
Applying this procedure to all edges of $\Gamma$ we get  a global
ribbon graph endowed with the canonical orientation.

Consider the global ribbon graph $\G$ as a surface with boundary.
The  boundary  components  of  this surface are in  a  one-to-one
correspondence with the subset of those  conical  points  of  $S$
which serve as the endpoints of the saddle connections $\gamma_i$
in the collection  $\gamma_1, \dots, \gamma_n$. The orders of the
corresponding singularities  are  calculated  as follows. For any
connected component  $(\partial\G)_m$  of  its boundary define an
integer $b_m$ as
\begin{equation}
\label{eq:order:of:a:newborn:singularity}
b_m = \sum_{\substack{\text{boundary singularities}\\
\text{which belong to } (\partial\G)_m}} (k_{i,l}+1)\quad-2
\end{equation}
The set with multiplicities $\alpha$ can be defined now as
\begin{equation}
\label{eq:singularity:data:defined:by:configuration}
\alpha = \bigg(\bigcup_{\substack{\pm\text{-vertices}\\
v_j\in\Gamma(\cC)}} \text{interior singularities of }v_j\bigg) \
\bigcup\ \bigg( \bigcup_{\substack{\text{components
}(\partial\G)_m
\\\text{of the boundary}\\\text{of } \G(\cC)}} b_m\bigg)
\end{equation}

\begin{figure}[ht]
%
\includegraphics{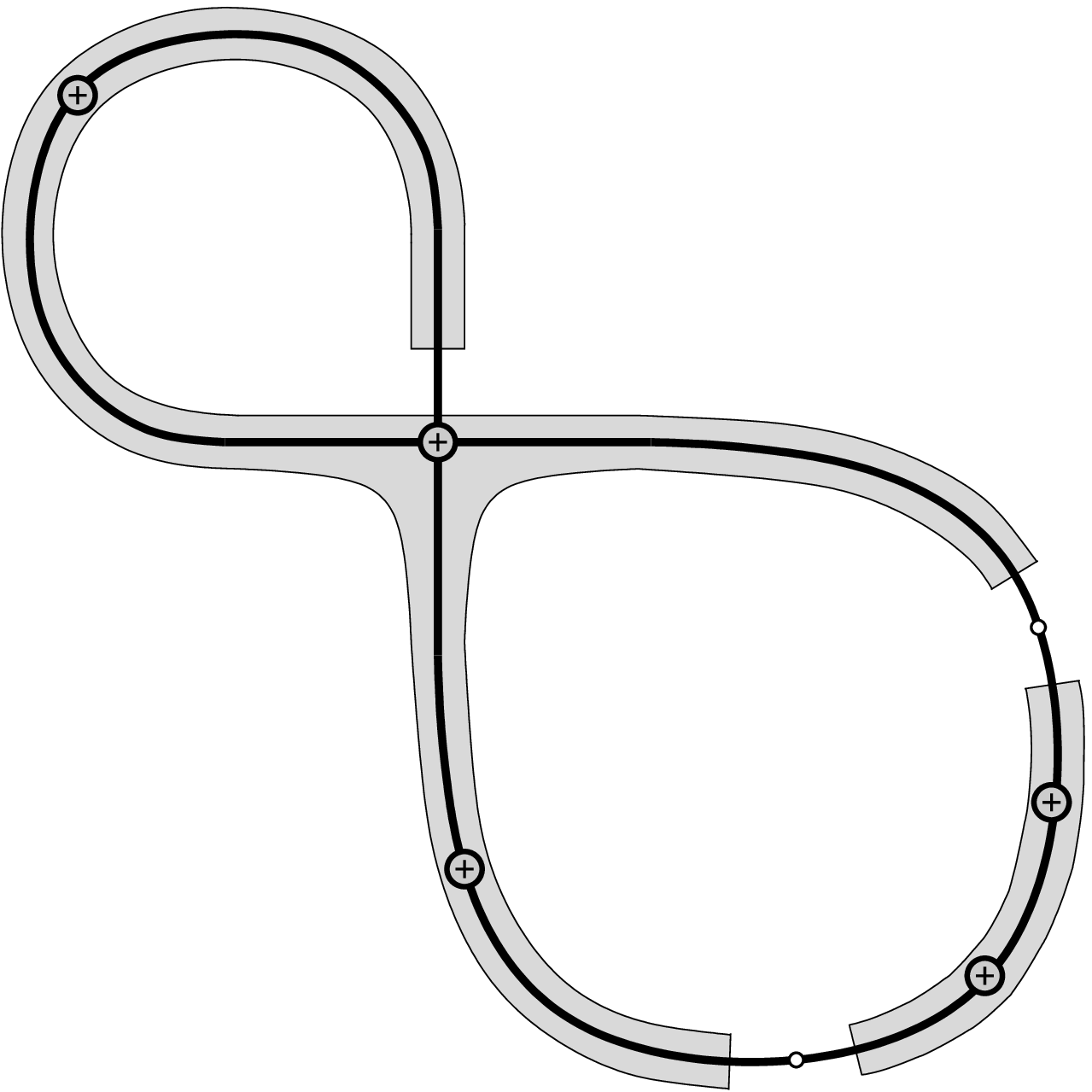}
\includegraphics{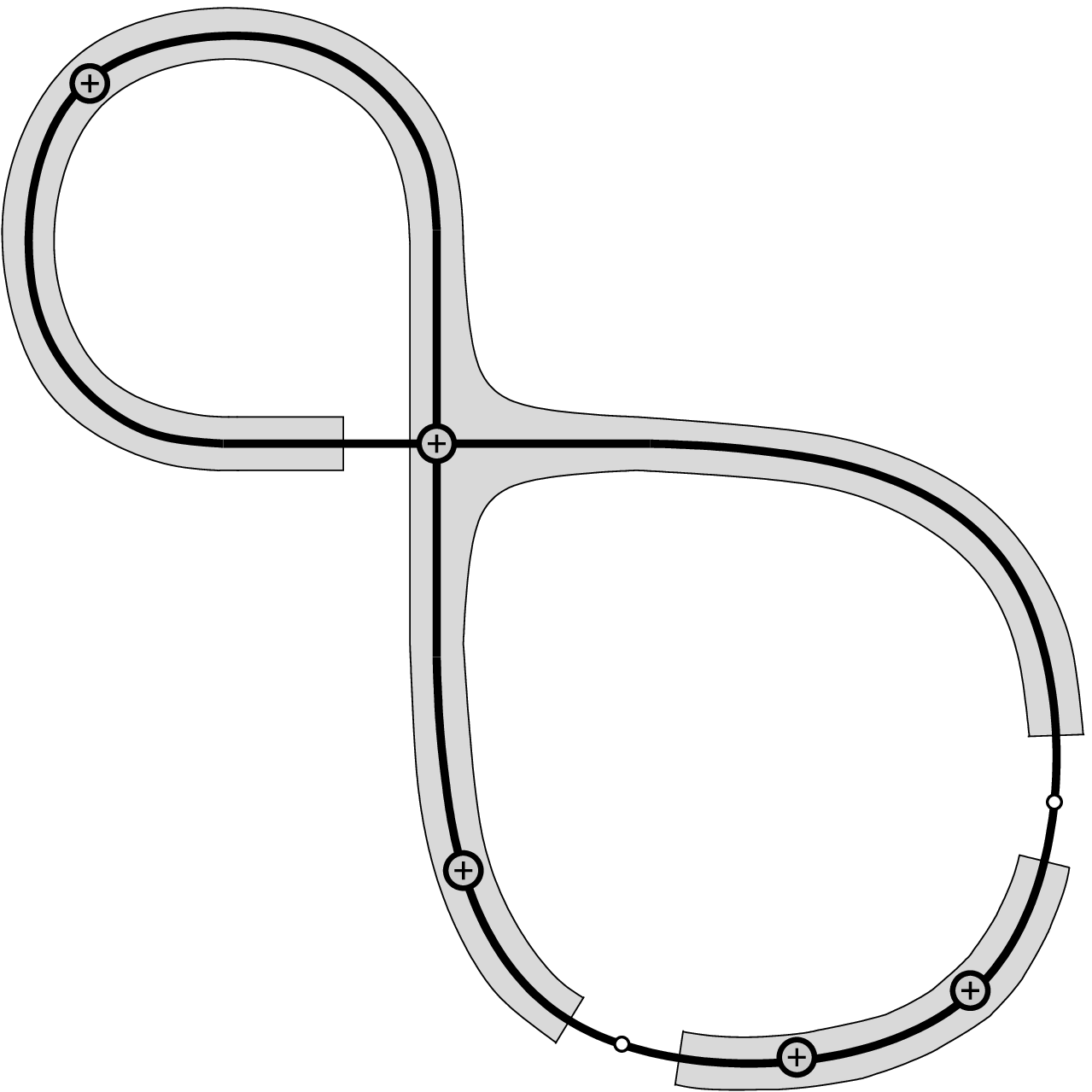}
\begin{picture}(0,0)(197,-20) 
\put(24,-32){$\scriptstyle 1$}
\put(34,-43){$\scriptstyle 3$}
\put(85,-70){$\scriptstyle 2$}
\put(92,-74){$\scriptstyle \emptyset$}
\put(72,-78){$\scriptstyle 1$}
\put(68,-97){$\scriptstyle 1$}
\put(87,-97){$\scriptstyle 0$}
\put(10,-40){$\scriptstyle \{2\}$}
\put(75,-30){$\cB_1$}
\put(60,-45){$\gamma_1$}
\put(40,-70){$\gamma_2$}
\put(85,-116){$\gamma_3$}
\put(119,-95){$\gamma_8$}
\put(100,-152){$\gamma_4$}
\put(132,-157){$\gamma_5$}
\put(145,-140){$\gamma_6$}
\put(148,-121){$\gamma_7$}
\put(57,-135){$\scriptstyle \{4,4\}$}
\put(75,-145){$\scriptstyle 5$}
\put(88,-140){$\scriptstyle 9$}
\put(154,-107){$\scriptstyle 0$}
\put(163,-114){$\scriptstyle 0$}
\put(120,-165){$\scriptstyle 0$}
\put(132,-172){$\scriptstyle 0$}
\put(154,-132){$\scriptstyle 1$}
\put(169,-134){$\scriptstyle 1$}
\put(148,-151){$\scriptstyle 1$}
\put(156,-164){$\scriptstyle 1$}
\put(171,-125){$\scriptstyle \emptyset$}
\put(163,-156){$\scriptstyle \emptyset$}
   %
\put(183,-142){$\cB_2$}
\end{picture}
   %
   %
\begin{picture}(0,0)(0,-20) 
\put(26,-32){$\scriptstyle 3$}
\put(34,-43){$\scriptstyle 1$}
\put(87,-73){$\scriptstyle 1$}
\put(87,-97){$\scriptstyle 0$}
\put(72,-92){$\scriptstyle 1$}
\put(64,-78){$\scriptstyle 2$}
\put(95,-74){$\scriptstyle \emptyset$}
\put(10,-40){$\scriptstyle \{2\}$}
\put(60,-45){$\gamma_2$}
\put(40,-70){$\gamma_1$}
\put(85,-116){$\gamma_3$}
\put(119,-95){$\gamma_8$}
\put(96,-149){$\gamma_4$}
\put(115,-156){$\gamma_5$}
\put(132,-154){$\gamma_6$}
\put(144,-139){$\gamma_7$}
\put(57,-135){$\scriptstyle \{4,4\}$}
\put(75,-145){$\scriptstyle 5$}
\put(88,-140){$\scriptstyle 9$}
\put(159,-129){$\scriptstyle 0$}
\put(165,-139){$\scriptstyle 0$}
\put(98,-170){$\scriptstyle 0$}
\put(108,-163){$\scriptstyle 0$}
\put(146,-153){$\scriptstyle 1$}
\put(155,-165){$\scriptstyle 1$}
\put(127,-161){$\scriptstyle 1$}
\put(128,-176){$\scriptstyle 1$}
\put(162,-160){$\scriptstyle \emptyset$}
\put(133,-178){$\scriptstyle \emptyset$}
\end{picture}
\vspace{160bp} 

\caption{
\label{fig:global:ribbon:graph}
An example of a configuration}
\end{figure}

\begin{Example}
\label{ex:example:of:a:configuration}
The  configurations  $\cC$  presented  in  the  left  picture  of
figure~\ref{fig:global:ribbon:graph} has $8$  saddle  connections
$\gamma=\{\gamma_1\cup\dots\cup\gamma_{8}\}$;     the     surface
$S\setminus\gamma$  decomposes  into   $7$  connected  components
$S_1\sqcup\dots\sqcup  S_7$.  Two  components are represented  by
cylinders and  thus  have  no  interior  singularities. Among the
remaining five  components  three  have no interior singularities
and   are   denoted  with  $\emptyset$,  one  has  one   interior
singularity of order  $2$, and one has two interior singularities
of orders $4$. Thus, we get
$$
\bigcup_{\substack{\pm\text{-vertices}\\
v_j\in\Gamma(\cC)}}\text{interior singularities of }v_j =
\{2,4,4\}
$$

The boundary  of the global  ribbon graph $\G$ has two components
$(\partial\G)_1$  and  $(\partial\G)_2$  which correspond to  two
conical  singularities  $P_1$  and  $P_2$  of   $S$.  The  saddle
connections $\gamma_5,  \gamma_6, \gamma_7$ join $P_2$ to itself;
the  other  saddle connections  join  $P_1$  to  itself.  Turning
counterclockwise around the point $P_l$, $l=1,2$, we see geodesic
rays parallel to  $\gamma_i$  appear in  the  same order as  they
appear when we  follow the corresponding component $\cB_l$ in the
positive direction. Denoting  by  ``$x$'' the geodesic rays which
do  not  belong   to  the  configuration  we  get  the  following
(cyclically ordered) list for the zero $P_2$:
$$
\gamma_5 x \gamma_6 x \gamma_7 \gamma_7 x \gamma_6 x \gamma_5
$$
We have  $10$ geodesic rays; this  corresponds to the  cone angle
$10\pi$ matching our  formula for the  order  $b_2$ of  the  zero
$P_2$:
$$
(0+1)+(1+1)+(1+1)+(1+1)+(1+1)+(0+1)-2 =8\,,
$$
The analogous list for $P_1$ is as follows
$$
xx\gamma_1 x \gamma_2 x \gamma_3 xxxxx \gamma_4 \gamma_4
xxxxxxxxx \gamma_3 \gamma_8 \gamma_8 x \gamma_2 xxx \gamma_1
$$
The number of consecutive ``$x$'' coincides with the order of the
corresponding          boundary          singularity         (see
definition~\ref{def:order:of:boundary:singularity}).   Thus,   at
$P_1$ we find $32$ geodesic  rays  parallel  to $\gamma_i$, which
corresponds to the  cone angle $32\pi$,  and the order  $b_1$  of
$P_1$ equals to
$$
(2+1)+(1+1)+(1+1)+(5+1)+(0+1)+(9+1)+
(0+1)+(0+1)+(1+1)+(3+1)-2\!=\!30
$$
Finally, we get the following set with multiplicities:
$$
\alpha=(2,4,4,8,30).
$$
The  surface  $S$  has  genus  $g=13$;  the  configuration  $\cC$
represents the  stratum $\cQ(2,4,4,8,30)$. Note, that the picture
on the right represents the same configuration as  the picture of
the left. \end{Example}

\begin{Remark}
The  example  above gives an idea  of  how can one construct  all
configurations         (in         the          sense          of
definition~\ref{def:configuration})   for   a    given    stratum
$\cQ(\alpha)$ of meromorphic quadratic differentials with at most
simple  poles.  This algorithm is discussed in  more  details  in
appendix~\ref{a:List:of:configurations:in:genus:2},  where  as an
illustration we present  a complete list of all configurations of
\^homologous   saddle   connections  for   holomorphic  quadratic
differentials in genus two.
\end{Remark}

\subsection*{Principal boundary}
Analogously to the case of Abelian  differentials a configuration
$\cC$   of   \^homologous   saddle  connections  determines   the
corresponding principal boundary  stratum $\cQ(\alpha'_{\cC})$ or
$\cH(\beta'_{\cC})$. Namely, to each boundary component $\cB_i$
$$\{\to\gamma_{j_{i,1}}\to  \dots \to\gamma_{j_{i,p(i)}}\to\}$$
of  every  ``$+$''  or  ``$-$''-vertex   $v_j$   of   the   graph
$\Gamma(\cC)$  (i.e.   to   each   connected   component  of  the
corresponding local ribbon graph $\G_j$) we assign a number
\begin{equation}
\label{eq:Dji}
D_{j_i}=k_{j_{i,1}}+\dots+k_{j_{i,p(i)}}-2,
\end{equation}
where $k_{j_{i,1}}, \dots, k_{j_{i,p(i)}}$ are the  orders of the
boundary singularities corresponding to this boundary  component.
By lemma~\ref{lm:sum:of:boundary:singularities:is:even} proved in
the                         beginning                          of
section~\ref{s:hat:homologous:saddle:connections}   the    number
$D_{j_i}$ is always  a nonnegative even integer whenever $v_j$ is
a ``$+$''-vertex.  To  every  ``$+$''-vertex  $v_j$  of the graph
$\Gamma(\cC)$ we assign the stratum
\begin{equation}
\label{eq:H:boundary:stratum}
\cH(\beta'_j)=\cH\left(\frac{d_1}{2}, \dots, \frac{d_{s(j)}}{2}, \frac{D_1}{2}, \dots,  \frac{D_{r(j)}}{2}\right)
\end{equation}
of  holomorphic   Abelian   differentials,   where  $d_1,  \dots,
d_{s(j)}$ are the orders of interior singularities of $v_j$. Note
that conditions (4) and (5) in definition~\ref{def:configuration}
of a configuration of \^homologous saddle  connections imply that
the entries of  $\beta'_j$  are integers  and  that their sum  is
even, so the stratum $\cH(\beta'_j)$ is nonempty.

We assign to a ``$-$''-vertex $v_j$
the  stratum
\begin{equation}
\label{eq:Q:boundary:stratum}
\cQ(\alpha'_j)=\cQ(d_1, \dots, d_{s(j)}, D_1, \dots,  D_{r(j)})
\end{equation}
of meromorphic quadratic differentials with at most simple poles,
where  $d_1,  \dots,   d_{s(j)}$   are  the  orders  of  interior
singularities   of    $v_j$.   Note   that   condition   (5)   in
definition~\ref{def:configuration}   of    a   configuration   of
\^homologous  saddle  connections  guarantees  that  the  sum  of
entries of $\alpha'_j$ defined above equals $0$ modulo $4$, while
condition             (6)             guarantees             that
$\alpha'\not\in\big\{(\emptyset,\{-1,1\},\{3,1\},\{4\}\big\}$,
which implies that the stratum $\cQ(\alpha'_j)$ is nonempty.

Given a configuration  $\cC$  we assign to every ``$\pm$''-vertex
of  the  graph $\Gamma$ the corresponding stratum. When  $\Gamma$
does   not   contain  ``$-$''   vertices   we   get   a   stratum
$\cH(\beta'_\cC)$    of    disconnected   translation    surfaces
$S'_1\sqcup\dots\sqcup   S'_k$,   where   $S'_j\in\cH(\beta'_j)$,
$j=1,\dots,k$. Otherwise  we  get a stratum $\cQ(\alpha'_\cC)$ of
disconnected  flat  surfaces $S'_1\sqcup\dots\sqcup  S'_k$, where
$S'_j\in\cH(\beta'_j)$   when   $S'_j$  is   represented   by   a
``$+$''-vertex   and   $S'_j\in\cQ(\alpha'_j)$  when   $S'_j$  is
represented by a ``$-$''-vertex. The resulting  stratum is called
the  {\it   principal  boundary  stratum}  corresponding  to  the
admissible configuration $\cC$.

\begin{Example}
\label{ex:principal:boundary}
Let  us  compute the principal boundary stratum corresponding  to
the                      configuration                       from
example~\ref{ex:example:of:a:configuration},                  see
figure~\ref{fig:global:ribbon:graph}.
The   components   represented   by    cylinders,    encoded   by
$\circ$-vertices do  not  contribute  to  the principal boundary:
they shrink and disappear.
The  vertex  $v_1$  of   valence   four  has  type  $+4.2c$,  see
figure~\ref{fig:embedded:local:ribbon:graphs}; the  corresponding
local ribbon  graph  $\G_{v_1}$  has  two  connected  components,
$r(1)=2$, which  correspond  to  two connected components $\cB_1,
\cB_2$ of  the  boundary $\partial S^{comp}_1$. The corresponding
zeroes  of  the  induced  Abelian  differential   on  $S'_1$  are
calculated   in   terms   of  $D_1=2-2=0$  and   $D_2=1+0+1-2=0$,
see~\eqref{eq:Dji}.  Since  $S^{comp}_1$  does not have  interior
singularities, the corresponding closed flat surface  $S'_1$ is a
torus    with     two    marked    points,     $S'_1\in\cH(0,0)$,
see~\eqref{eq:H:boundary:stratum}.

The remaining  four  vertices  of  $\Gamma(S,\gamma)$  have  type
$+2.1$; the  boundary  of  each  of  the corresponding components
$S_2,\dots,S_4$  is  connected. Applying  formulae~\eqref{eq:Dji}
and~\eqref{eq:H:boundary:stratum} we  get  the  following list of
surfaces $S'_j$:
$$
S'_2\in\cH\left(\frac{2}{2},\ \frac{3+1-2}{2},\right)\quad
S'_3\in\cH\left(\frac{4}{2},\ \frac{4}{2}, \frac{5+9-2}{2},\right)\quad
S'_4,
S'_5\in\cH\left(\frac{1+1-2}{2}\right)
$$
The corresponding principal boundary stratum is
$$
\cH(0,0)\sqcup\cH(1,1)\sqcup\cH(6,2,2)\sqcup\cH(0)\sqcup\cH(0)
$$

\end{Example}

\subsection*{Main Theorems}
In                           sections~\ref{s:Local:Constructions}
and~\ref{s:Nonlocal:constructions}   we   describe   some   basic
surgeries which depend continuously on a  small complex parameter
$\delta\in\C{}$ (responsible for the length and direction of
the  saddle  connections  which  form  the  boundary)  and  on an
additional discrete  parameter  having  finitely many values. The
theorem below makes a bridge  between  the  formal  combinatorial
constructions discussed  above  and  the  geometry  of the moduli
spaces  of  quadratic  differentials  and  is   proved  in  those
sections.

We  denote  by  $\cQ_1^{\varepsilon}(\alpha)\subset\cQ_1(\alpha)$
the  subset  of those flat surfaces  of  area one, which have  at
least one saddle connection of length at most $\varepsilon$.

\begin{Theorem}
\label{th:from:boundary:to:neighborhood:of:the:cusp}
For each  configuration  $\cC$ of \^homologous saddle connections
as  in  definition~\ref{def:configuration},  let $\Gamma$ be  the
graph   of   connected   components    corresponding    to   this
configuration. Let  $\cQ(\alpha'_\cC)$ (or $\cH(\beta'_\cC)$)  be
the boundary stratum corresponding to the configuration $\cC$.

For any flat surface $S'\in \cQ(\alpha'_\cC)$ (correspondingly in
$\cH(\beta'_\cC)$), and  any  sufficiently  small  value  of  the
complex parameter $\delta$, if one applies the basic surgeries to
the connected components of  $S'$  and assembles a closed surface
$S$ from the resulting surfaces  with  boundary  according to the
structure  of  the  graph  $\Gamma(\cC)$,  then  the result is  a
surface in $\cQ^{\varepsilon}(\alpha)$.
\end{Theorem}

Similar  to  the  case  of  Abelian  differentials, we denote  by
$\cQ_1^{\varepsilon,thick}(\alpha)\subset\cQ_1(\alpha)$       the
subset  of  those  flat  surfaces  of  area  one,  which  have  a
collection of \^homologous  saddle  connections of length at most
$\varepsilon$  and   no   other  short  saddle  connection.  Here
``short'' means, of length less  than  $\lambda  \epsilon^r$  for
some  parameters  $\lambda\geq  1$  and $0< r\leq 1$,  where  the
values of the parameters depend on the stratum. Then one can show
that any  surface  in  $\cQ_1^{\varepsilon,thick}(\alpha)$ can be
obtained by this construction.  We  will not prove this statement
in order not to overload this paper.

We                  put                  theorem~\ref{th:graphs},
theorem~\ref{th:all:local:ribon:graphs}                       and
theorem~\ref{th:from:boundary:to:neighborhood:of:the:cusp}
together in one  statement which may  be considered as  our  main
theorem.

We  say  that  a  collection  $\gamma$   of  \^homologous  saddle
connections $\{\gamma_1,  \dots,  \gamma_n\}$  on  a flat surface
$S\in\cQ(\alpha)$ is \emph{in  general  position} if there are no
other saddle connections on $S$ parallel to saddle connections in
the      collection      $\gamma$.      It      follows      from
proposition~\ref{pr:homologous:equiv:parallel} stated  in the end
of  section~\ref{s:hat:homologous:saddle:connections}  that   for
almost  all  flat  surfaces  in  any  stratum  any  collection of
\^homologous saddle  connections  is  in  general  position. This
implies,  that  we can always put a  collection  of  \^homologous
saddle connections  in  general  position  by  an arbitrary small
deformation of the flat surface inside the stratum.

\begin{MainTheorem}
Any  collection  $\gamma$  of  \^homologous  saddle   connections
$\{\gamma_1, \dots, \gamma_n\}$ in general  position  on  a  flat
surface  $S\in\cQ(\alpha)$  naturally   defines  a  corresponding
configuration $\cC(S,\gamma)$.

Any ``formal'' configuration of  \^homologous  saddle connections
as  in  definition~\ref{def:configuration}  corresponds  to  some
actual  collection  of  \^homologous  saddle  connections  on  an
appropriate flat surface.
\end{MainTheorem}
\begin{proof}
By   theorem~\ref{th:graphs}    any   collection   $\gamma$    of
\^homologous saddle  connections $\{\gamma_1, \dots,\!\gamma_n\}$
on a flat  surface $S\in\cQ(\alpha)$ naturally defines a graph of
connected  components  $\Gamma(S,\gamma)$  (structure   1   of  a
configuration).                   According                    to
theorem~\ref{th:all:local:ribon:graphs}, for every  vertex $v$ of
$\Gamma(S,\gamma)$ the  collection  $\gamma$ also defines a local
ribbon graph (structure  2  of a  configuration)  as well as  the
orders  $d_l$  and   $k_{i,l}$   of  all  interior  and  boundary
singularities.  By  theorem~\ref{th:all:local:ribon:graphs},  for
vertices  of  ``$+$''-type, the orders $k_{i,l}$ of the  boundary
singularities  are  compatible with  the  corresponding  parities
(structures 3 and 4 of a configuration). The lower bounds for the
sums  $D_i$  of orders  of  boundary  singularities  follow  from
lemma~\ref{lm:lower:bounds:for:D:i}. The necessary  condition  of
the  compatibility  of the orders of interior singularities  with
the orders of boundary singularities formalized as structure 5 is
proved                                                         in
lemma~\ref{lm:sum:of:boundary:singularities:minus:vertex}.    The
list of nonrealizable singularity data  for  the  vertices of the
``$-$''-types presented  in  structure  6  of  a configuration is
justified  in  lemma~\ref{lm:nonrealizable:data} at  the  end  of
section~\ref{s:Nonlocal:constructions}. This  completes the proof
of the first part of the statement.

The  realizability  of  all  formal  configurations   immediately
follows                                                      from
theorem~\ref{th:from:boundary:to:neighborhood:of:the:cusp}.
\end{proof}

\subsection*{Appendices. Long saddle connections}

In appendix~\ref{ap:Long:saddle:connections} we study collections
of \^homologous saddle connections when they  are not necessarily
short.

The     next     proposition     follows     immediately     from
definition~\ref{def:homologous} and the notion of configuration.

\begin{NNProposition}
Let $\gamma(S_0)=\{\gamma_1, \dots, \gamma_n\}$  be  a collection
of \^homologous saddle connections on a  flat  surface  $S_0$  in
$\cQ(\alpha)$.  Let  a  flat   surface   $S$  be  obtained  by  a
sufficiently   small   continuous   deformation   of   $S_0$   in
$\cQ(\alpha)$  and  $\gamma(S)$  the corresponding collection  of
saddle connections. Then all saddle connections in the collection
$\gamma(S)$     are     \^homologous.      The      configuration
$\cC(S,\gamma(S))$  defined  by  the  collection  $\gamma(S)$  of
\^homologous saddle connections on $S$ coincides with the initial
configuration $\cC(S_0,\gamma(S_0))$.
\end{NNProposition}

By definition,  a  configuration  $\cC$  of  \^homologous  saddle
connections  is   admissible  for  a  given  connected  component
$\cQ^{c}(\alpha)$  of  the  stratum  $\cQ(\alpha)$  if  there  is
\emph{at least  one}  flat surface $S_0\in\cQ^{c}(\alpha)$ and at
least one  collection $\gamma$ of \^homologous saddle connections
$\gamma=\{\gamma_1,\dots,\gamma_n\}$  on  $S_0$ realizing  $\cC$.
Consider  any  surface  $S$  in  the   same  connected  component
$\cQ^{c}(\alpha)$.   By   $N_\cC(S,L)$   denote  the  number   of
collections $\gamma=\{\gamma_1,\dots,\gamma_n\}$ of  \^homologous
saddle  connections   on  $S$  defining  the  same  configuration
$\cC(S,\gamma)=\cC$   and   such   that   $\max_{1\le   i\le   n}
|\gamma_i|\le  L$.  The  results in~\cite{Eskin:Masur} imply  the
following statement proved in the appendix.

\begin{Proposition}
\label{pr:counting}
For  almost  every flat surface $S$ in  the  connected  component
$\cQ^{c}(\alpha)$ containing $S_0$ the following limit exists
$$
\lim_{L\to +\infty} \cfrac{N_\cC(S,L)}{L^2}=c_\cC(S)
$$
and is strictly positive.  Moreover,  for almost all surfaces $S$
in    $\cQ^{c}(\alpha)$     this     limit     is    the    same,
$c_\cC(S)=const_\cC$.   (This   limit  is   called  Siegel--Veech
constant.)
\end{Proposition}

In particular,  any  admissible  configuration  is  presented  on
almost  every   flat   surface  in  the  corresponding  connected
component of the stratum by numerous  collections of \^homologous
saddle connections.

\subsection*{Final comments, open problems, applications}

The  thick  part  $\cQ_1^{\varepsilon,thick}(\alpha)$  decomposes
into a disjoint union
$$
\cQ_1^{\varepsilon,thick}(\alpha)=
\bigsqcup_{{\text{configurations }\cC}} \cQ_1^{\varepsilon}(\alpha,\cC)
$$
of  (not  necessarily  connected)  components  corresponding   to
admissible configurations; the surfaces in any  such component of
$\cQ_1^{\varepsilon}(\alpha,\cC)$  share  the same  configuration
$\cC$ of \^homologous saddle connections. Following  the lines of
the    paper~\cite{Eskin:Masur:Zorich}     one    could    extend
theorem~\ref{th:from:boundary:to:neighborhood:of:the:cusp}    and
prove that  up to a  defect of  a very small  measure, for  every
configuration  $\cC$  there is  an  integer  $M(\cC)$  such  that
$\cQ_1^{\varepsilon}(\alpha,\cC)$  is a  (ramified)  covering  of
order $M(\cC)$ over  the  following space. The  space is a  fiber
bundle   over    the   boundary   stratum    $\cQ_1(\alpha'_\cC)$
(correspondingly   $\cH_1(\beta'_\cC)$).   It  has   a  Euclidean
$\varepsilon$-disc  as  a  fiber  when  $\cC$  does  not  contain
cylinders,  and  the  space  $\cH_1^\varepsilon(0,\dots,0)$  when
$\cC$ contains cylinders (number of  marked  points  on the torus
equals  the  number of cylinders). In  both  cases it is easy  to
express the measure on $\cQ_1^{\varepsilon}(\alpha,\cC)$ in terms
of  the product  measure  on the fiber  bundle,  and compute  the
volume of $\cQ_1^{\varepsilon}(\alpha,\cC)$ in  terms  of volumes
of the strata, and  using  the Siegel---Veech formula compute the
constants $c_\cC$.

However, the evaluation of the constants $M$ (which depend on the
configuration   $\cC$)   requires  some   additional   work.   In
particular,     if     the    corresponding     surgeries    (see
theorem~\ref{th:from:boundary:to:neighborhood:of:the:cusp})   are
nonlocal  (i.e.  those,  which  use  a  path on  a  surface,  see
section~\ref{s:Nonlocal:constructions})  one  needs to  study the
dependence of  the resulting surface  on the homotopy type of the
path.  These  and   related  issues  will  be  discussed  in  the
forthcoming paper~\cite{Boissy:in:progress}.

Another  subject  which we do not  discuss  in this paper is  the
individual study of the connected components  of  the  strata  of
quadratic differentials:  different  connected  components of the
same  stratum   $\cQ(\alpha)$  have  their  individual  lists  of
admissible  configurations,  graphs,  boundary  strata,  etc.  In
particular, one can use the lists of admissible configurations to
determine the connected component  to  which a given flat surface
belongs. For example, a saddle  connection  joining  the zero and
the  simple  pole  on   any   flat  surface  from  the  component
$\cQ^{ir}(9,-1)$ has a \^homologous saddle connection joining the
zero  to  itself, while analogous saddle connections on  surfaces
from the complementary connected component $\cQ^{reg}(9,-1)$  may
have  multiplicity   one.   The  existing  invariant  called  the
\emph{Rauzy class} used to distinguish these components is rather
complicated, see~\cite{Lanneau}. Configurations  of  \^homologous
saddle connections for nonconnected strata will be studied in the
papers~\cite{Boissy:to:appear} and~\cite{Boissy:in:progress}.

Given a billiard  in  a rational polygon $\Pi$,  one  can build a
translation surface $\hat S$  from  an appropriate number $2N$ of
copies of $\Pi$ such that  geodesics  on $S$ will project to  the
billiard trajectories in $\Pi$. Taking $N$ copies instead of $2N$
one  obtains  a flat surface with $\ZZ$-holonomy  with  the  same
properties of geodesics.  In  some cases this latter construction
is  more  advantageous. In  the paper~\cite{Athreya:Eskin:Zorich}
there is the  study  of billiards  in  polygons whose angles  are
multiples of $\pi/2$.  Identifying two copies of such polygons by
their boundaries one obtains a  flat  surface  corresponding to a
meromorphic quadratic  differential  on  $\C{}P^1$  with  at most
simple poles. The  results  of this  paper  are used to  classify
closed billiard  trajectories  and  generalized  diagonals in the
paper~\cite{Athreya:Eskin:Zorich},                            see
also~\cite{Boissy:to:appear}.

\bigskip
\centerline{\sc Acknowledgements}
\medskip

Conceptually   this    paper    is    a   continuation   of   the
paper~\cite{Eskin:Masur:Zorich}. We  want  to  thank A.~Eskin for
his participation  in the early stage  of this project.  We thank
C.~Boissy  and  E.~Lanneau  for  several  helpful   conversations
concerning nonconnected  strata.  The  first  author  thanks  the
University of Rennes for  its  support and hospitality during the
preparation  of  this   paper.   The  second  author  thanks  the
University of Chicago, UIC, Max-Planck-Institut f\"ur  Mathematik
at  Bonn  and   IHES  for  hospitality  and  support  during  the
preparation of this paper.

\section{Preliminaries  on  flat  surfaces  and  on  \^homologous
saddle connections}
\label{s:hat:homologous:saddle:connections}

In this  section  of  preliminary  results  we describe geometric
criteria  for   deciding   when   two   saddle   connections  are
\^homologous  and   describe  the  structure  of  the  complement
$S\setminus(\gamma_1\cup\gamma_2)$.  The   key  result  in   this
section is proposition~\ref{pr:homologous:saddle:connections}.

In the case of  a translation surface $S$ it is obvious  that two
saddle connections $\gamma_1,\gamma_2$ are homologous if and only
if $S\setminus (\gamma_1\cup\gamma_2)$  is disconnected (provided
$S\setminus\gamma_1$ and $S\setminus\gamma_2$  are connected). It
is  less  obvious to check whether saddle connections  $\gamma_1,
\gamma_2$ on a flat  surface  $S$ with nontrivial linear holonomy
are \^homologous  or not. In  particular, a pair of closed saddle
connections  might  be  homologous  in the usual sense,  but  not
\^homologous;  a  pair of  closed  saddle  connections  might  be
\^homologous even if one of them represents a  loop homologous to
zero,  and  the  other  does  not;  finally, a saddle  connection
joining  a  pair   of   {\it  distinct}  singularities  might  be
\^homologous to  a  saddle  connection  joining  a singularity to
itself. The flat torus described in in the  introduction gives an
example         of         these          phenomena          (see
example~\ref{ex:homologous:sad:connections}                   and
figure~\ref{fig:examples:of:homologous:sad:connections:1}).

We  start  this  section  with several lemmas  establishing  some
restrictions on the orders of singularities  of  a  flat  surface
with  boundary.  By  convention  we  consider   only  those  flat
structures  which   have   linear   holonomy  in  $\{Id,  -Id\}$.
Throughout this paper we  assume  that the boundary components of
any flat surface  with boundary are  made up of  parallel  saddle
connections,  unless  otherwise noted.  We  also  assume  that  a
sufficiently small  collar neighborhood of any boundary component
is  a  topological annulus  (or,  in the  other  words, that  the
natural  compactification  of  $S\setminus \partial S$  coincides
with $S$).

\begin{lemma}
\label{lm:sum:of:boundary:singularities:is:even}
If  a  flat  surface  $S_j$  with  boundary  has  trivial  linear
holonomy,  then  the  sum  of  the   orders   of   the   boundary
singularities along each boundary component is even:
$$
k_{j_{i,1}}+  \dots + k_{j_{i,p(i)}} \equiv 0 \mod 2
$$
\end{lemma}
\begin{proof}
Take  a loop following the $i$-th boundary component
$$\{\to\gamma_{j_{i,1}}\to  \dots \to\gamma_{j_{i,p(i)}}\to\}$$
at  a  sufficiently  small  constant  distance.  Recall  that  by
definition~\ref{def:order:of:boundary:singularity}  of  the order
of  a  boundary   singularity,   the  angle  between  the  saddle
connection   $\gamma_{j_{i,l}}$   and   the   saddle   connection
$\gamma_{j_{i,l+1}}$     at      the     boundary     singularity
$\gamma_{j_{i,l}}\to\gamma_{j_{i,l+1}}$                    equals
$(k_{j_{i,l}}+1)\pi$. Thus,  the  linear holonomy around the loop
is  trivial  if  and  only  if  the  total  sum  of   the  angles
$k_{j_{i,1}}\pi+\dots+ k_{j_{i,p(i)}}\pi$ is  an integer multiple
of   $2\pi$,   or,  equivalently,  if  and  only   if   the   sum
$k_{j_{i,1}}+\dots+ k_{j_{i,p(i)}}$ of the orders of the boundary
singularities is even.
\end{proof}

\begin{lemma}
\label{lm:sum:of:boundary:singularities:minus:vertex}
Let $d_{j_l},k_{j_{i,l}}$  denote  the  orders of correspondingly
interior singularities  and  boundary  singularities  of  a  flat
surface with boundary $S_j$. Then
$$
2 r(S_j)-4\le\
\sum d_{j_l} + \sum k_{j_{i,l}} \ \equiv\ 2 r(S_j)\mod 4,
$$
where $r(S_j)$ is the number  of  boundary  components, the first
sum is  taken over all  interior singularities and the second sum
is taken over all boundary singularities.
\end{lemma}
\begin{proof}
Consider  one  more copy  of  the surface  $S_j$  taken with  the
opposite orientation. We can naturally identify  these two copies
along the common boundary. It follows  from  our  assumptions  on
$S_j$ that the resulting surface  $S$  is  a nonsingular oriented
closed flat surface  without boundary. In other words, the closed
flat  surface   $S$   corresponds   to  a  meromorphic  quadratic
differential on a Riemann surface.

Each  interior  singularity of $S_j$ of order $d_{j_l}$  produces
two distinct singularities of  $S$  of order $d_i$. Each boundary
singularity  of  $S_j$ of order $k_{j_{i,l}}$ gives  rise  to  an
interior singularity of $S$ of order  $2k_{j_{i,l}}$. The surface
$S$ has  genus  $\hat{g}=2g+r(S_j)-1$.  Now  recall  that for any
quadratic differential on a  closed  Riemann surface $S$ of genus
$\hat{g}$ the sum of orders of singularities equals $4\hat{g}-4$.
Hence,
$$
2\big(\sum_{\substack{interior\\singularities\\of\ S_j}} d_l +
\sum_{\substack{boundary\\singularities\\of\ S_j}} k_{j_{i,l}} \big) =
4\hat{g}-4=4(2g+r(S_j)-1)-4=8(g-1)+4r(S_j)\,
$$
which implies the desired relation.
\end{proof}

\begin{lemma}
\label{lm:lower:bounds:for:D:i}
The sum  $D_i+2=k_{i,1}+\dots+k_{i,p(i)}$  of  orders of boundary
singularities along some boundary component  $\cB_i$  of  a  flat
surface $S$ is equal to zero if and only if a sufficiently narrow
collar neighborhood of  $\cB_i$  in $S$  is  isometric to a  flat
cylinder.

When  $S$ has  trivial  linear holonomy and  the  sum $D_i+2$  of
orders of  boundary  singularities  along  a  boundary  component
$\cB_i$ is strictly positive, then $D_i\ge 0$.
\end{lemma}
\begin{proof}
By  definition~\ref{def:order:of:boundary:singularity}  the order
$k_{i,l}$ of  any  boundary  singularity  is  nonnegative.  Thus,
$D_i+2$  is equal  to  zero  if  and only  if the  orders  of all
boundary singularities  along  the boundary component $\cB_i$ are
equal to zero, which implies the first part of the statement.

The second statement  is an obvious  corollary of the  first  one
combined                                                     with
lemma~\ref{lm:sum:of:boundary:singularities:is:even}.
\end{proof}

\begin{lemma}
\label{lm:total:boundary:holonomy:is:trivial}
Let $\beta$ denote the total  boundary  of  a translation surface
defined by a holomorphic 1-form $\omega$. Then
\begin{itemize}
\item[1.] $\int_\beta\omega=0$.
\item[2.] $\beta$ cannot consist of a single saddle connection.
\item[3.]   If   $\beta$  is  composed  of  exactly  two   saddle
connections  $\gamma_1,\gamma_2$   then  $\gamma_1,\gamma_2$  are
parallel  and  have equal length. Moreover, the oriented  surface
obtained by isometric identification of $\gamma_1$ and $\gamma_2$
is a  translation surface (i.e. it is a  closed flat surface with
trivial linear holonomy).
\end{itemize}
\end{lemma}

\begin{proof}
Note  that  the canonical orientation of the  surface  induces  a
canonical  orientation  of the boundary $\beta$. Thus, the  first
statement is  an  immediate  consequence  of  Stokes formula. The
second  statement  follows from  the  first  since  the  holonomy
$\int_\gamma\omega$ along  a saddle connection $\gamma$ cannot be
$0$.

For    the    third     let    $\beta=\gamma_1-\gamma_2$.    Then
$\int_{\gamma_1}\omega=\int_{\gamma_2}\omega$. This implies  that
$\gamma_1,\gamma_2$ are  parallel,  have  equal  length  and that
their  directions   defined   by   the  chosen  orientations  are
compatible with  linear  holonomy.  We can isometrically identify
$\gamma_1$ either  with  $\gamma_2$ or with $-\gamma_2$. However,
the second  identification  produces  a nonorientable surface, so
$\gamma_1$ must be identified with $\gamma_2$  which implies that
the  resulting  surface  is  a  translation  surface.
\end{proof}

Let $S$ be a flat surface with boundary; let $\gamma_1, \gamma_2$
be a pair  of parallel saddle connections $\gamma_1, \gamma_2$ of
equal length at the boundary  of  $S$.  By convention, throughout
this paper we always  identify  $\gamma_1$ and $\gamma_2$ in such
way that the resulting flat surface is orientable.

Suppose that, moreover, $S$ has trivial linear holonomy.
\begin{Definition}
\label{def:flip}
We  say  that  $\gamma_1$  and  $\gamma_2$  are  identified  by a
\textit{translation} if  the  resulting  flat surface has trivial
linear holonomy; otherwise we say that  $\gamma_1$ and $\gamma_2$
are identified by a \textit{flip}.
\end{Definition}

\begin{Lemma}
\label{lm:trivial:plus:trivial:is:trivial}
Assume that a flat surface $S$ with nontrivial linear holonomy is
divided   by    a    pair    of   parallel   saddle   connections
$\gamma_1,\gamma_2$  into  two  connected components $S_1,  S_2$.
Then at  least one of  the components must have nontrivial linear
holonomy.
\end{Lemma}
\begin{proof}
If one  of the $\gamma_1,  \gamma_2$ is a closed curve homologous
to zero,  say $\gamma_1$, then $\gamma_2$ lies in  one of the two
components  of  the complement  $S\setminus\gamma_1$.  Then,  the
boundary of the other component,  say,  $S_1$  consists solely of
$\gamma_1$,   $\partial   S_1=\gamma_1$.   By   property   2   of
lemma~\ref{lm:total:boundary:holonomy:is:trivial}  the  component
$S_1$ has nontrivial linear holonomy.

Therefore, we may assume  that  $\gamma_1$ and $\gamma_2$ are not
homologous to zero so they are homologous to each other. Choosing
appropriate orientations of $\gamma_1$ and $\gamma_2$ we get
$$
\partial S_1 = \gamma_1 - \gamma_2 \qquad
\partial S_2 = -\gamma_1 + \gamma_2
$$
where  the  orientations  of  $S_1,  S_2$  are  induced  from the
canonical  orientation  of  $S$.  If  both  $S_1$ and $S_2$  have
trivial  linear  holonomy we can choose the defining  holomorphic
1-forms $\omega_1,\omega_2$ on $S_1$ and $S_2$ in such way that
$$
\int_{\gamma_1}\omega_1 = \int_{\gamma_1}\omega_2=
\int_{\gamma_2}\omega_1 = \int_{\gamma_2}\omega_2.
$$
The  latter  relations imply the compatibility of $\omega_1$  and
$\omega_2$ on $S$. Thus, the flat structure on $S$ can be defined
by     a     holomorphic    1-form     $\omega$     such     that
$\omega|_{S_1}=\omega_1,  \omega|_{S_2}=\omega_2$,  and  $S$  has
trivial linear holonomy contrary to the initial assumption.
\end{proof}

\begin{Lemma}
\label{lm:gamma:2gamma}
Any two \^homologous saddle connections $\gamma_1, \gamma_2$ on a
flat   surface    $S$   are   parallel.   When   both   relations
$[\gamma'_1]=-[\gamma''_1]$  and  $[\gamma'_2]=-[\gamma''_2]$ are
simultaneously valid  or  simultaneously  not  valid  the  saddle
connections   $\gamma_1,   \gamma_2$  have   the   same   length,
$|\gamma_1|=|\gamma_2|$.  When   one   of   the  relations,  say,
$[\gamma'_1]=-[\gamma''_1]$, is valid while the other one is not,
$[\gamma'_2]\neq-[\gamma''_2]$, the lengths differ by the  factor
of two, $|\gamma_1|=2|\gamma_2|$.
\end{Lemma}
\begin{proof}
The     proof     is     a    straightforward    corollary     of
definition~\ref{def:homologous}  and  the fact that the length
of a saddle  connection $\delta$ on the translation surface $\hat
S$  is   defined   as  $|\delta|=|\int_\delta  \omega|$  and  its
direction is defined by the  argument  of  $\int_\delta  \omega$.
\end{proof}

\begin{Lemma}
\label{lm:antiinvariant:cycle}
Let $\gamma$ be a saddle connection on a flat surface  $S$ having
nontrivial  linear   holonomy.   The   following  properties  are
equivalent
\begin{itemize}
\item[(a)]
$[\gamma']=-[\gamma'']$ in $H_1(\hat{S},\hat P;\,\Z)$;
\item[(b)]
$\hat  S\setminus  (\gamma'\cup\gamma'')$ contains  two connected
components;
\item[(c)]
$[\gamma]=0$ in $H_1(S,P;\,\Z)$.
\end{itemize}
\end{Lemma}
\begin{proof}

\noindent(a)$\Rightarrow$(c).   Consider   the    map    $p_\ast:
H_1(\hat{S},\hat   P;\,\Z)\to   H_1(S,P;\,\Z)$  induced   by  the
covering  $p$.   By  definition  of  $\gamma',\gamma''$  we  have
$[\gamma]=p_\ast[\gamma']=p_\ast[\gamma'']$.      Thus,      when
$[\gamma']=-[\gamma'']$   we   get    $[\gamma]=-[\gamma]$,    so
$[\gamma]=0$.

(c)$\Rightarrow$(b).   Since   $[\gamma]=0$,   $S\setminus\gamma$
contains two connected components $S_1, S_2$, such that $\partial
S_1=\gamma$,  $\partial  S_2   =   -\gamma$.  By  property  2  of
lemma~\ref{lm:total:boundary:holonomy:is:trivial} both $S_1$  and
$S_2$  have  nontrivial linear holonomy, which implies that  both
$\hat  S_1=p^{-1}(S_1)$,  $\hat S_2=p^{-1}(S_2)$,  are connected.
Thus, $\hat S\setminus  (\gamma'\cup\gamma'')=\hat S_1\sqcup \hat
S_2$ contains two connected components.

(b)$\Rightarrow$(a).   Since   $\gamma'$   and   $\gamma''$   are
symmetric, the two connected components $\hat  S',  \hat  S'$  of
$\hat S\setminus (\gamma'\cup\gamma'')$ are  also  symmetric with
respect to the involution. This restricts the possible situations
to the following three (up to an interchange  of the superscripts
of $\hat S_1', \hat S_1''$ if necessary):

--- either $\partial  \hat  S'$ is  composed  of two copies  of
$\gamma'$ and $\partial \hat S''$ from two copies of $\gamma''$;

--- or $\partial \hat  S'  = \gamma'-\gamma''$ and $\partial \hat
S'' = \gamma''-\gamma'$;

--- or $\partial \hat  S'  = \gamma'+\gamma''$ and $\partial \hat
S'' = -\gamma''-\gamma'$.

The first situation  implies that $\hat S$ contains two connected
components  which   contradicts  the  assumptions  that  $S$  has
nontrivial  linear   holonomy.  Hence,  the  first  situation  is
excluded.  The  second  situation   implies   that  isometrically
identifying the  boundary  components $\gamma'$ and $\gamma''$ of
$S'$ we obtain a flat surface  isometric to $S$. By property 3 of
lemma~\ref{lm:total:boundary:holonomy:is:trivial}   this    again
implies  that  $S$ has trivial linear holonomy which  contradicts
the  assumptions.  This  case  is  also  excluded.  In  the  only
remaining case  we  have  $\partial  \hat  S' = \gamma'+\gamma''$
which implies $[\gamma']=-[\gamma'']$.
\end{proof}

The   next   proposition   is   the   key   to   the   proofs  of
theorems~\ref{th:unique:trivial:holonomy} and~\ref{th:graphs}. We
do  not  assume that the saddle connections  in  the  proposition
below are parallel.

\begin{proposition}
\label{pr:homologous:saddle:connections}
Two saddle connections $\gamma_1,\gamma_2$ on a  flat surface $S$
are  \^homologous  if  and  only   if   they   have  no  interior
intersections and one of the following holds
\begin{enumerate}
\item The  union  $\gamma_1\cup\gamma_2$  does  not  separate the
surface          $S$          and         the          complement
$S\setminus\{\gamma_1\cup\gamma_2\}$ has trivial linear holonomy.
    \noindent
(In this case $|\gamma_1|=|\gamma_2|$; all combinations: loop ---
loop, loop --- segment, segment --- segment are possible.)
\item  The  union  $\gamma_1\cup\gamma_2$ separates $S$;  neither
$\gamma_1$ nor  $\gamma_2$  by  itself  separates; the complement
$S\setminus\{\gamma_1\cup\gamma_2\}$     has     two    connected
components, one of them has  trivial  linear  holonomy, the other
--- nontrivial.
    \noindent
(In this case $|\gamma_1|=|\gamma_2|$; the saddle connections are
either     two      segments      or      two     loops.)
\item One of $\gamma_1, \gamma_2$, say, $\gamma_1$ separates $S$,
the     other      one      does      not;     the     complement
$S\setminus\{\gamma_1\cup\gamma_2\}$     has     two    connected
components, one of them has  trivial  linear  holonomy, the other
one, whose boundary is $\gamma_1$, has nontrivial holonomy.
    \noindent
(In  this  case  $|\gamma_1|=2|\gamma_2|$; the separating  saddle
connection $\gamma_1$ is a loop, $\gamma_2$ might be a segment or
a loop.)
\item Both $\gamma_1$ and $\gamma_2$ separate $S$; the complement
$S\setminus\{\gamma_1\cup\gamma_2\}$    has    three    connected
components; two of  which  have nontrivial linear holonomy, while
the remaining  one, whose boundary is $\gamma_1\cup\gamma_2$, has
trivial linear holonomy.
    \noindent
(In  this  case  $|\gamma_1|=|\gamma_2|$;  both  $\gamma_1$   and
$\gamma_2$ are loops.)
\end{enumerate}
\end{proposition}
\begin{proof}
According  to   lemma~\ref{lm:gamma:2gamma}  \^homologous  saddle
connections are  parallel. If two \^homologous saddle connections
$\gamma_1$ and  $\gamma_2$ have a  common point, this point is an
endpoint for both $\gamma_1$ and $\gamma_2$. Thus, from now on we
can   assume   that   $\gamma_1,\gamma_2$   have   no    interior
intersections.

Two  saddle  connections   $\gamma_1,\gamma_2$  without interior
intersections subdivide  a  flat  surface  $S$  in  one  of  the
following ways:
\begin{itemize}
\item[(i)]
The union  $\gamma_1\cup\gamma_2$  does  not separate the surface
$S$.
\item[(ii)]
The   union   $\gamma_1\cup\gamma_2$   separates   $S$;   neither
$\gamma_1$ nor $\gamma_2$ by itself separates.
\item[(iii)]
One of $\gamma_1,  \gamma_2$,  say, $\gamma_1$ separates $S$, the
other one does not.
\item[(iv)]
Both $\gamma_1$ and $\gamma_2$ separate $S$.
\end{itemize}

For each  of these cases  we prove that the additional assumption
that $\gamma_1$ and $\gamma_2$ are \^homologous  is equivalent to
the corresponding additional assumptions  (1)--(4)  on triviality
of linear holonomy  of the corresponding components. In each case
we                use                lemmas~\ref{lm:gamma:2gamma}
and~\ref{lm:antiinvariant:cycle} to determine  the  corresponding
relation between  the  lengths  $|\gamma_1|$ and $|\gamma_2|$. We
combine this  information with lemma~\ref{lm:antiinvariant:cycle}
(when appropriate) to  prove that one of $\gamma_1, \gamma_2$ (or
both $\gamma_1$  and $\gamma_2$) is  a closed cycle. The proof of
realizability  of  combinations loop --- loop, loop ---  segment,
segment         ---          segment         indicated         in
proposition~\ref{pr:homologous:saddle:connections} is left to the
reader        as        an       exercise.        Note       that
example~\ref{ex:homologous:sad:connections}    already     proves
realizability   of   combinations   loop---segment  in  (1)   and
loop---loop in (2). The remaining  combinations  can  be found in
sections~\ref{s:Local:Constructions}
and~\ref{s:Nonlocal:constructions}.

Let $X\subseteq S$ be a subset  of $S$. By $\hat X$ we denote the
preimage $\hat X=  p^{-1}(X)$. Let $S_j$ be a connected component
of  $S\setminus  (\gamma_1\cup\gamma_2)$. We  use  the  following
obvious  criterion:  $S_j$ has nontrivial linear holonomy if  and
only if the preimage $\hat S_j$  is connected. Now let us pass to
consideration of cases (i)--(iv).

(i)   In   this   case  $S\setminus  (\gamma_1\cup\gamma_2)$   is
connected; denote it by $S_1$. Neither  of $\gamma_1,\gamma_2$ is
homologous to zero, so $[\hat\gamma_1]=[\gamma_1']-[\gamma_1'']$,
and        $[\hat\gamma_2]=[\gamma_2']-[\gamma_2'']$         (see
lemma~\ref{lm:antiinvariant:cycle}).                           By
lemma~\ref{lm:gamma:2gamma} when  such $\gamma_1$ and  $\gamma_2$
are \^homologous, we have $|\gamma_1|=|\gamma_2|$.

If  the   saddle   connections   $\gamma_1$  and  $\gamma_2$  are
\^homologous   then   the  cycle   $[\gamma_1']-[\gamma_1'']$  is
homologous   (in   the   usual   sense)    to    one    of    the
$\pm([\gamma_2']-[\gamma_2''])$    which    means   that    $\hat
S_1=\hat{S}\setminus
(\gamma_1'\cup\gamma''_1\cup\gamma'_2\cup\gamma''_2)$   is    not
connected. Hence, by the above criterion $S_1$ has trivial linear
holonomy.

Suppose  now  that  $S_1=S\setminus  (\gamma_1\cup\gamma_2)$  has
trivial  linear  holonomy. Then  $\hat  S_1$  has  two  connected
components  $\hat  S_1'$ and $\hat S_1''$. Note,  that  the  flat
surface $S$ has nontrivial linear  holonomy.  Hence,  it  follows
from                 property                 3                in
lemma~\ref{lm:total:boundary:holonomy:is:trivial}    that    both
$S\setminus\gamma_1$  and  $S\setminus\gamma_2$  have  nontrivial
linear holonomy.  This implies that there  exist a pair  of loops
$\rho_1, \rho_2$ on $S$ such  that  $\rho_i$  and $\gamma_i$ have
single    transversal   intersection,    $i=1,2$;    such    that
$\rho_1\cap\gamma_2=\emptyset$,   $\rho_2\cap\gamma_1=\emptyset$;
and  such  that   holonomy   along  each  $\rho_i$,  $i=1,2$,  is
nontrivial. Interchanging  the  superscripts  of $\hat S_1', \hat
S_1''$ if necessary,  we may assume  that $\gamma'_1$ is  on  the
boundary of $\hat S_1'$. Since  $\rho_1$  has  nontrivial  linear
holonomy, the  lift  $\rho'_1\subset\hat S_1'$ of $\rho$ starting
at $\gamma'_1$ is not closed and hence it  ends on $-\gamma''_1$.
This  implies  that both $\gamma'_1$ and $-\gamma''_1$ belong  to
the   boundary    of    $\hat    S_1'$.   Since   $S_1=S\setminus
(\gamma_1\cup\gamma_2)$  is  connected,  at  least  one  of  both
$\pm\gamma'_2$ and  $\pm\gamma''_2$  belongs  to  the boundary of
$\hat S_1'$.  Applying the same argument  as above and  using the
obvious symmetry between $\hat S_1'$ and $\hat S_1''$ we conclude
that  under  an appropriate choice of orientations of  $\gamma_1$
and $\gamma_2$ one has
$$
\partial \hat S_1' = \gamma'_1 - \gamma''_1 - \gamma'_2 + \gamma''_2
$$
which is equivalent to
$$
[\gamma'_1] - [\gamma''_1] = [\gamma'_2] - [\gamma''_2]
$$
and hence, $\gamma_1$ and $\gamma_2$ are \^homologous.

(ii) In this case $\gamma_1$ and $\gamma_2$ are homologous in the
usual  sense,  and   not   homologous  to  zero;  the  complement
$S\setminus\{\gamma_1\cup\gamma_2\}$ has two connected components
$S_1,  S_2$.  This  implies  that either both of  $\gamma_1$  and
$\gamma_2$ are segments, or both are closed cycles. Since neither
of    $\gamma_1,    \gamma_2$    is    homologous    to     zero,
lemma~\ref{lm:antiinvariant:cycle}          claims           that
$[\gamma'_i]\neq-[\gamma''_i]$   for   $i=1,2$.  Thus,   if  such
$\gamma_1,     \gamma_2$     are      \^homologous     we     get
$|\gamma_1|=|\gamma_2|$ by lemma~\ref{lm:gamma:2gamma}.

By lemma~\ref{lm:trivial:plus:trivial:is:trivial} at least one of
two components, say, $S_1$ has nontrivial  linear holonomy. Under
an appropriate  choice  of orientations of $\gamma_1,\gamma_2$ we
have $\partial S_1 =\gamma_1-\gamma_2$, which implies
$$
\partial \hat S_1 =(\gamma'_1+\gamma''_1)-(\gamma'_2+\gamma''_2).
$$

Since $[\gamma'_i]\neq-[\gamma''_i]$, for  $i=1,2$, the condition
that $\gamma_1$ and $\gamma_2$ are \^homologous  is equivalent to
one                of                the                relations
$([\gamma'_1]-[\gamma''_1])=\pm([\gamma'_2]-[\gamma''_2])$.
Together with the  above equation on  $\partial \hat S_1$  it  is
equivalent to one of the following systems
$$
\begin{cases}
[\gamma'_1]=[\gamma'_2]\\
[\gamma''_1]=[\gamma''_2]
\end{cases}
\qquad
\begin{cases}
[\gamma'_1]=[\gamma''_2]\\
[\gamma''_1]=[\gamma'_2]
\end{cases}
$$
The  systems  might be identified by interchange of  superscripts
of, say, $\gamma_2'$  and $\gamma_2''$, thus we can consider just
the first one.

Since        by       the        second        property        of
lemma~\ref{lm:total:boundary:holonomy:is:trivial}   neither    of
$[\gamma'_i], [\gamma''_i]$,  $i=1,2$, is homologous to zero, the
latter system  is valid if and only if  cutting $\partial \hat S$
by    any    of   two    pairs    $[\gamma'_1],[\gamma'_2]$    or
$[\gamma''_1],[\gamma''_2]$  of  saddle connections  we  get  two
connected components. Since $\hat S_1$ is connected the latter is
true if and only if $\hat S_2$ contains two connected components.
By the criterion  formulated  above this is true  if  and only if
$S_2$ has trivial linear holonomy. The equivalence is proved.

(iii) In this  case $\gamma_1$ is  a closed cycle  homologous  to
zero, while $\gamma_2$ is not  homologous  to  zero. This implies
that the complement  $S\setminus\{\gamma_1\cup\gamma_2\}$ has two
connected       components       $S_1,      S_2$.       Combining
lemma~\ref{lm:antiinvariant:cycle}                           with
lemma~\ref{lm:gamma:2gamma} we  see  that  if such $\gamma_1$ and
$\gamma_2$ are \^homologous, we have $|\gamma_1|=2|\gamma_2|$.

Choose  the  orientation of  $\gamma_1$  and  enumeration  of  the
components in such way that
$$
\partial S_1 = \gamma_1 \qquad
\partial S_2 = -\gamma_1 + \gamma_2 - \gamma_2
$$
Then
$$
\partial \hat S_1 = \gamma'_1 + \gamma''_1\qquad
\partial \hat S_2 = -\gamma'_1 -\gamma''_1 + \gamma'_2 - \gamma'_2 + \gamma''_2 - \gamma''_2
$$
Note that $S_1$ has nontrivial linear holonomy (see property 2 of
lemma~\ref{lm:total:boundary:holonomy:is:trivial}) so $\hat  S_1$
is connected. If $\gamma_1$ and $\gamma_2$ are \^homologous, then
$[\gamma'_1]=\pm ([\gamma'_2]-[\gamma''_2])$. This  implies  that
$\hat S\setminus(\gamma'_1\cup\gamma'_2\cup\gamma''_2)$  contains
at    least    two    connected    components.    Since     $\hat
S\setminus(\gamma'_1\cup\gamma'_2\cup\gamma''_2)=            \hat
S_1\cup\gamma''_1\cup \hat S_2$ where $\hat S_1$ is connected and
$\gamma_1''$ connects $\hat S_1$ and  $\hat  S_2$,  this  implies
that $\hat S_2$  is  not be  connected.  Hence $S_2$ has  trivial
linear holonomy.

Conversely,     consider     the    connected     component    of
$S\setminus\gamma_1$ containing $\gamma_2$; denote it by  $\tilde
S_2$.                Property                2                 of
lemma~\ref{lm:total:boundary:holonomy:is:trivial}  implies   that
$\tilde  S_2$   has   nontrivial   linear   holonomy.  Note  that
$S_2=\tilde S_2\setminus\gamma_2$. Thus, when  $S_2$  has trivial
linear holonomy,  there exist a a  closed path $\rho$  on $\tilde
S_2$  transversally  intersecting  $\gamma_2$ such that  holonomy
along  $\rho$  is nontrivial.  Since  $S_2$  has  trivial  linear
holonomy, $\hat S_2$ has two  connected  components  $\hat  S_2',
\hat S_2''$.  Changing  if  necessary  the  superscripts of $\hat
S_2',  \hat  S_2''$ we  may  assume that  $\gamma_2'$  is on  the
boundary  of  $\hat  S'_2$.  Since the holonomy along  $\rho$  is
nontrivial,  following  the  lift  of  $\rho$   which  starts  at
$\gamma_2'$ and  goes inside $\hat  S'_2$ the path $\rho$ ends at
$-\gamma''_2$,  which  shows  that $\gamma_2'$ and  $-\gamma''_2$
make part of  the  boundary of  the  same component $\hat  S'_2$.
Taking into  consideration  the symmetry between components $\hat
S_2', \hat  S_2''$  and  choosing  an  appropriate orientation of
$\gamma_1$ we get
$$
\partial \hat S'_2 = -\gamma'_1 + \gamma'_2 - \gamma''_2
\qquad
\partial \hat S''_2 = -\gamma''_1 - \gamma'_2 + \gamma''_2
$$
which implies that $\gamma_1$ and $\gamma_2$ are \^homologous.

(iv)
In this case  the complement $S\setminus\{\gamma_1\cup\gamma_2\}$
has three  connected  components.  Both $\gamma_1$ and $\gamma_2$
are homologous to zero, so they are represented by closed cycles.
This also  implies  that $[\hat \gamma_i]=\gamma'_i$, $i=1,2$. If
such  $\gamma_1$   and   $\gamma_2$  are  \^homologous,  we  have
$|\gamma_1|=|\gamma_2|$ (see lemma~\ref{lm:gamma:2gamma}).

Denote         the         connected        components         of
$S\setminus(\gamma_1\cup\gamma_2)$  in  such way  that  under  an
appropriate  choice  of  orientations of $\gamma_1,\gamma_2$  one
gets
$$
\partial S_1 = \gamma_1 \qquad
\partial S_2 = -\gamma_2 \qquad
\partial S_3 = -\gamma_1+\gamma_2
$$
By                 property                 2                  of
lemma~\ref{lm:total:boundary:holonomy:is:trivial} the  components
$S_1$ and $S_2$ have nontrivial  linear  holonomy,  so $\hat S_1$
and $\hat S_2$ are connected. We get
$$
\partial \hat S_1 = \gamma'_1 + \gamma''_1 \qquad
\partial \hat S_2 = -\gamma'_2 - \gamma''_2\qquad
\partial \hat S_3 = -\gamma'_1-\gamma''_1 + \gamma'_2 + \gamma''_2
$$

If   $\gamma_1$    and    $\gamma_2$    are   \^homologous   then
$[\gamma'_1]=\pm[\gamma'_2]$ which implies that cutting $\hat S$
by $\gamma'_1, \gamma'_2$ we get two  connected components, which
means  that  $\hat S_3$  is  not connected  and  hence $S_3$  has
trivial linear holonomy.

Conversely, if $S_3$ has trivial  linear  holonomy  then  $\hat
S_3$  contains  two connected components $\hat S'_3, \hat  S''_3$
which (under appropriate enumeration) have boundaries
$$
\partial \hat S'_3 = -\gamma'_1 + \gamma'_2 \qquad
\partial \hat S''_3 = -\gamma''_1+ \gamma''_2
$$
which implies  that  $\gamma_1$  and  $\gamma_2$  are homologous.
Proposition~\ref{pr:homologous:saddle:connections}   is   proved.
\end{proof}

Theorem~\ref{th:unique:trivial:holonomy}       follows       from
proposition~\ref{pr:homologous:saddle:connections}.

\begin{proof}[Proof of theorem~\ref{th:unique:trivial:holonomy}]
Cutting  $S$  by  $\gamma_1,  \gamma_2$   we   get   one  of  the
decompositions          (i)--(iv).          According          to
proposition~\ref{pr:homologous:saddle:connections},           the
additional assumptions (1)--(4)  on  the triviality of the linear
holonomy  of   the  corresponding  component  are  necessary  and
sufficient   conditions   for   $\gamma_1,   \gamma_2$   to    be
\^homologous.  It  remains to  note  that in  each  of the  cases
(1)--(4)  there  is  a  unique  component   with  trivial  linear
holonomy.
\end{proof}

The  following   criterion  is  analogous  to  the  corresponding
statement   in~\cite{Eskin:Masur:Zorich}.   It   is   proved   in
appendix~\ref{ap:Long:saddle:connections}, where the  notion of a
measure on each stratum is discussed.

\begin{Proposition}
\label{pr:homologous:equiv:parallel}
For  almost  every  flat  surface  in  any  stratum,  two  saddle
connections are parallel if and only if they are \^homologous.
\end{Proposition}

\section{Graph of connected components}
\label{s:graph}

In this  section we  give the proof that every  graph is given by
the list  in  theorem~\ref{th:graphs}.  Denote  by  $\dot{S}$ the
surface  $S$  punctured  at  the singularities. Any  closed  loop
$\rho$ on  $\dot{S}$ can be homotoped to have  a finite number of
transverse intersections  with  the  saddle  connections from the
collection   $\gamma=\{\gamma_1\ldots\gamma_n\}$.   It  naturally
induces a path $\rho_\ast$ on  the  graph  $\Gamma(S,\gamma)$  by
recording the surfaces $S_j$ intersected by $\rho$. It is easy to
see that paths $\rho\sim\rho'$ homotopic on the punctured surface
$\dot{S}$ define  paths  $\rho_\ast\sim  \rho'_\ast$ homotopic on
the graph.  Mark a point $x\in\dot{S}\setminus\{\gamma_i\} $; let
$v(x)$   be    the    corresponding    vertex    of   the   graph
$\Gamma(S,\gamma)$.    We    get     a    natural    homomorphism
$\pi_1(\dot{S},x)\to\pi_1\big(\Gamma(S,\gamma),v(x)\big)$.

Any finite connected graph can be retracted by a deformation to a
bouquet  of circles  (possibly  to a point).  We  can choose  the
retraction in such way that $v(x)$ retracts to the base  point of
the bouquet  of circles. We can  consider the bouquet  of circles
$\operatorname{B}$  as   a   graph   obtained   from   the  graph
$\Gamma(S,\gamma)$    by    identifying     some    subtree    of
$\Gamma(S,\gamma)$  to  a single  vertex  of  $\operatorname{B}$.
Thus, some edges of $\Gamma(S,\gamma)$ remain nondegenerate under
the retraction, and some edges collapse to a point.

Now we  can prove  the lemma which is the  main technical tool in
the proof of theorem~\ref{th:graphs}.

\begin{Lemma}
\label{lm:nontriv:lin:holon:along:loops:of:the:graph}
Let   $\alpha\subset\Gamma(S,\gamma)$   be   a  closed  path   on
$\Gamma(S,\gamma)$ realized as a subgraph of  $\Gamma(S,\gamma)$.
If under some retraction of  $\Gamma(S,\gamma)$  to  a bouquet of
circles,  $\alpha$  retracts to one of the  circles,  then  there
exists a closed path $\rho$  on  the  punctured surface $\dot{S}$
such that $\rho_\ast=\alpha$ and the linear holonomy along $\rho$
is nontrivial.
\end{Lemma}
\begin{proof}
We suppose that a  retraction  of $\Gamma(S,\gamma)$ to a bouquet
of circles  is fixed. We  start with consideration of the general
case, when the bouquet of circles contains at least two circles.

Let   $\gamma_1$   and  $\gamma_2$  be  a  pair   of   edges   of
$\Gamma(S,\gamma)$, which remain  nondegenerate under retraction,
such that $\gamma_1\in\alpha$  and $\gamma_2\not\in\alpha$ (since
the  bouquet  contains at  least  two  circles,  such  $\gamma_2$
exists).  Cutting  $\Gamma(S,\gamma)$  by  these edges we  get  a
connected graph. Equivalently,  cutting the surface $S$ by a pair
of \^homologous saddle connections $\gamma_1, \gamma_2$  we get a
connected   surface    $S_{(1,2)}=S\setminus(\gamma_\cup\gamma_2$
which, by proposition~\ref{pr:homologous:saddle:connections}, has
trivial  linear  holonomy.  By construction $\partial  S_{(1,2)}=
\gamma_1\cup-\gamma_1\cup\gamma_2\cup-\gamma_2$. Gluing back  the
boundary components  $\gamma_1$ and $-\gamma_1$ of $S_{(1,2)}$ we
get a surface $S_{(2)}=S\setminus\gamma_2$  which  has nontrivial
linear                        holonomy                         by
lemma~\ref{lm:total:boundary:holonomy:is:trivial}.   Thus,    the
boundary  components  $\gamma_1$  and  $-\gamma_1$  of  the  {\it
translation} surface $S_{(1,2)}=S\setminus(\gamma_1\cup\gamma_2)$
are  identified  by a flip (see definition~\ref{def:flip} in  the
previous section).

Consider any path  $\rho$ in $S$ such that $\rho_\ast=\alpha$ and
such  that   $\rho$  has  unique  transversal  intersection  with
$\gamma_1$. By construction $\rho$ it gives a nonclosed connected
path  $\rho'$  on   $S_{(1,2)}$  joining  a  pair  of  points  on
$P_+\in\gamma_1$ and $P_-\in-\gamma_1$ corresponding to the  same
point  $P\in\gamma_1$  on $S$  upon  gluing  of  $\gamma_1$  with
$-\gamma_1$. Since $\gamma_1$ and $-\gamma_1$ are identified by a
flip, we see that the linear holonomy along $\rho$ is nontrivial.

It remains to  consider  the case,  when  the bouquet of  circles
corresponding  to  the  graph  $\Gamma(S,\gamma)$  has  a  single
circle.               It               follows               from
proposition~\ref{pr:homologous:saddle:connections} that the graph
cannot  be  just a single loop composed  of  ``$+$''-vertices  of
valence two and of ``cylinder  vertices''  of  valence two. Thus,
either  $\Gamma(S,\gamma)$  is  a  loop composed of  vertices  of
valence two with some ``$-$''-vertices, or there is  at least one
nontrivial subtree with a vertex on the base loop.

In  the  first  case choose any  path  $\rho'$  on  $S$ such that
$\rho'_\ast=\alpha$.  If  the  linear  holonomy  along  the  path
$\rho'$ is nontrivial, we  choose  $\rho:=\rho'$ and the lemma is
proved. If the holonomy is trivial, we can compose $\rho'$ with a
closed  path  $\rho''$, such that $\rho''$ is contained  entirely
inside  some  $S^-_j$, and such that the  linear  holonomy  along
$\rho''$   is   nontrivial.   Since  $\rho''\subset  S^-_j$   the
projection    $\rho''_\ast$    is    a   trivial   path.    Thus,
$(\rho'\cdot\rho'')_\ast   =   \rho'_\ast   =  \alpha$,  and   by
construction the  linear  holonomy  along $(\rho'\cdot\rho'')$ is
nontrivial.   The    required    path    $\rho$   is   given   by
$\rho'\cdot\rho''$.

In the second case the subtree necessarily has a vertex  $S_j$ of
valence                one,                which               by
Lemma~\ref{lm:total:boundary:holonomy:is:trivial}      is       a
``$-$''-vertex.  Denote  by $\gamma_1$ the edge adjacent to  this
vertex of  valence one; we denote  by the same  symbol $\gamma_1$
the corresponding saddle connection in  $S$.  Consider  any  path
$\rho'$    on    $S$   such    that    $\rho'_\ast=\alpha$.    If
$\alpha\subset\Gamma(S,\gamma)$ passes through  $S^-_j$, we apply
the same  construction as in  the previous case. If $\alpha$ does
not pass  through  $S^-_j$  then $\gamma_1\cap \alpha=\emptyset$,
and any path $\rho\in S$ such that $\rho_\ast=\alpha$ has trivial
intersection with  the  saddle connection $\gamma_1$. Choose some
edge $\gamma_2\in\alpha$ which is nondegenerate under retraction.
Cutting  $S$  by the  pair  of  \^homologous  saddle  connections
$\gamma_1, \gamma_2$ we get two connected components: a connected
surface $S_{(1,2)}$ and  a  surface $S^-_j$ (corresponding to the
vertex of  valence one). The closed path $\rho$  on $S$ becomes a
nonclosed connected  path  on  $S_{(1,2)}$  joining  the boundary
components      $\gamma_2'$      and       $\gamma_2''$.       By
proposition~\ref{pr:homologous:saddle:connections}  the   surface
$S_{(1,2)}$  has   trivial   linear   holonomy.  By  construction
$\partial  S_{(1,2)}=  \gamma_1\cup\gamma_2\cup-\gamma_2$.   Glue
back  the  boundary  components  $\gamma_2$  and  $-\gamma_2$  of
$S_{(1,2)}$. We get  a surface $S_{(1)}$ which coincides with one
of the two components of the  initial surface $S$ cut by a single
saddle connection $\gamma_1$. Since $\partial S_{(1)} = \gamma_1$
,  by  lemma~\ref   {lm:total:boundary:holonomy:is:trivial}   the
surface $S_{(1)}$  has  nontrivial  linear  holonomy.  Thus,  the
boundary  components  $\gamma_2$  and  $-\gamma_2$  of  the  {\it
translation} surface $S_{(1,2)}$ were identified by  a flip which
implies that the linear holonomy along $\rho$ is nontrivial.
\end{proof}

\begin{Lemma}
\label{lm:subgraph:with:minus:vertex:has:nontrivial:holonomy}
Consider a connected subgraph $\Upsilon$  of  the  initial  graph
$\Gamma(S,\gamma)$. If  $\Upsilon$  has  a  vertex  labelled with
``$-$''  or if it  is  not  a  tree, the  surface  with  boundary
$S_\Upsilon$ corresponding to this subgraph has nontrivial linear
holonomy.
\end{Lemma}
\begin{proof}
If  the  subgraph  has  some  vertex  labelled with ``$-$'',  the
corresponding surface $S^-_j$  has  a closed path with nontrivial
linear holonomy. The bigger surface  $S_\Upsilon$  has  the  same
path, so it also has nontrivial linear holonomy.  If the subgraph
is not  a tree, then  it has  a loop which  is not  homotopically
trivial.                                                       By
Lemma~\ref{lm:nontriv:lin:holon:along:loops:of:the:graph}   there
is a closed path $\rho$ on  $S_\Upsilon$  corresponding  to  this
loop such that $\rho$ has nontrivial linear holonomy.
\end{proof}

\begin{proof}[Proof of theorem~\ref{th:graphs}. (Necessity)]
First we note that a cylinder has trivial linear holonomy,  so by
lemma~\ref{lm:total:boundary:holonomy:is:trivial}  a   \cv-vertex
cannot have valence $1$.

If the valence  of a  \cv-vertex is  two,  then each  of the  two
boundary components  of  the  corresponding cylinder represents a
single saddle connection. Hence, two \cv-vertices  of valence two
cannot have a common edge,  otherwise  the  pair of corresponding
cylinders would be  identified along a boundary component of each
which  would  result  in  a  longer  cylinder  contradicting  the
assumption that each cylinder is maximal.

Now  note  that  the  bouquet  of  circles  to  which  the  graph
$\Gamma(S,\gamma)$ is  retracted  contains  at  most two circles.
Otherwise there would be edges  $\gamma_1$  and  $\gamma_2$  such
that  $\Gamma(S,\gamma)\setminus(\gamma_1\cup\gamma_2)$  would be
connected  but   not   simply   connected.   Thus,  according  to
lemma~\ref{lm:subgraph:with:minus:vertex:has:nontrivial:holonomy}
the  surface  $S_{(1,2)}=S\setminus(\gamma_1\cup\gamma_2)$  would
have    nontrivial    linear    holonomy,    which    contradicts
proposition~\ref{pr:homologous:saddle:connections}.

{\bf Two loops.} Suppose that the  bouquet  of  circles  contains
exactly  two  circles.  Cut  them  by  some edges $\gamma_1$  and
$\gamma_2$ which correspond to different circles  of the bouquet.
By     proposition~\ref{pr:homologous:saddle:connections}     the
resulting     surface    has     trivial     linear     holonomy.
Lemma~\ref{lm:subgraph:with:minus:vertex:has:nontrivial:holonomy}
implies that  the surface and therefore  the graph does  not have
any ``$-$''-vertices, in particular, no vertices  of valence $1$.
Since the Euler characteristic of $S^1\vee S^1$ equals to $-1$ we
get
\begin{multline*}
-1=\chi(S^1\vee S^1)=\chi(\Gamma(S,\gamma))=
-\frac{1}{2}\cdot(\text{number of vertices of valence 3})-\\
-\frac{2}{2}\cdot(\text{number of vertices of valence 4})
-\frac{3}{2}\cdot(\text{number of vertices of valence 5})-\dots
\end{multline*}
which  means  that either $\Gamma(S,\gamma)$ has two vertices  of
valence $3$ while  all the other  vertices have valence  $2$,  or
$\Gamma(S,\gamma)$ has a single  vertex  of valence $4$ while all
the other vertices  have  valence $2$.  All  graphs of this  type
except one are  in the list of theorem~\ref{th:graphs}, see types
d and  e. The type  which we  have to rule  out is  schematically
presented at figure~\ref{fig:banned:graph}.

\begin{figure}[htb]
%
\includegraphics{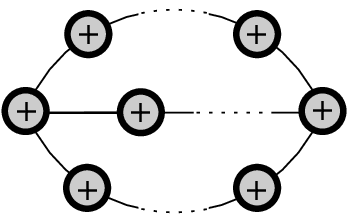}
\vspace{80bp} 
\caption{
\label{fig:banned:graph}
A graph of this type is not realizable as $\Gamma(S,\gamma)$.
 }
\end{figure}

We prove  by contradiction that this  graph is not  realizable as
$\Gamma(S,\gamma)$. Let $S^+_n$ be a vertex of valence three; let
$\gamma_1,\gamma_2,\gamma_3$ be the edges adjacent to it. Cutting
$\Gamma(S,\gamma)$ by  any  pair  of  distinct  edges  $\gamma_i,
\gamma_j$, $i=1,2,3$, we  still get a connected graph. This means
that    no    pair    of    \^homologous    saddle    connections
$\gamma_i\union\gamma_j$,  $i=1,2,3$, separates  $S$.  Hence,  by
proposition~\ref{pr:homologous:saddle:connections}  the   lengths
$|\gamma_i|$,  $i=1,2,3$   are   equal  and  all  $\gamma_i$  are
parallel.

Let  $\omega$   be  the  holomorphic  $1$-form  representing  the
translation structure on $S^+_n$. Under an  appropriate choice of
orientations  of $\gamma_1,\gamma_2,\gamma_3$  we  get  $\partial
S^+_n    =     \gamma_1\cup\gamma_2\cup\gamma_3$,    and    hence
$\int_{\gamma_1}\omega       +       \int_{\gamma_2}\omega      +
\int_{\gamma_3}\omega =  0$. On the  other hand the fact that all
the lengths $|\gamma_i|$, $i=1,2,3$ are equal  and all $\gamma_i$
are    parallel    implies    that    $\int_{\gamma_1}\omega    =
\pm\int_{\gamma_2}\omega = \pm\int_{\gamma_3}\omega$.  These  two
relations are incompatible, which is a contradiction.

{\bf One  loop.} If $\Gamma(S,\gamma)$ is  a loop of  vertices of
valence                two                then                 by
proposition~\ref{pr:homologous:saddle:connections}   it   has  at
least one  ``$-$''-vertex  $S_i^-$. Let $\gamma_1$ and $\gamma_2$
be  the  edges  of  $\Gamma(S,\gamma)$ adjacent to  $S_i^-$.  The
complement $S\setminus(\gamma_1\cup\gamma_2)$ has  two  connected
components: $S_i^-$  and  $S\setminus  S_i^-$.  Since $S_i^-$ has
nontrivial             linear             holonomy,            by
proposition~\ref{pr:homologous:saddle:connections}    the    flat
surface  $S\setminus  S^-_i$  has  trivial  linear  holonomy.  It
follows                         now                          from
Lemma~\ref{lm:subgraph:with:minus:vertex:has:nontrivial:holonomy}
that $S\setminus S^-_i$ is a chain of ``$+$''-vertices of valence
two and of \cv-vertices of valence two.    Thus, in this case the
graph   $\Gamma(S,\gamma)$    is    of    the    type    b,   see
theorem~\ref{th:graphs}                                       and
figure~\ref{fig:classification:of:graphs}.

If the  graph has a nontrivial subtree attached  to the base loop
then any  such subtree necessarily has  a vertex of  valence one,
which by lemma~\ref{lm:total:boundary:holonomy:is:trivial}  is  a
``$-$''-vertex. Let us show that $\Gamma(S,\gamma)$ can have only
one ``$-$''-vertex of valence one. Suppose  that  there  are  two
vertices $S^-_i$ and $S^-_j$ of valence one; denote by $\gamma_1$
and $\gamma_2$ the edges of $\Gamma(S,\gamma)$  adjacent to these
vertices. Cutting  $S$  by  $\gamma_1$  and  $\gamma_2$ we obtain
three     connected    components:     $S^-_i,     S^-_j$     and
$S_{(1,2)}:=S\setminus(S^-_i\union S^-_j)$.  Since the first  two
surfaces  have   nontrivial  linear  holonomy,  it  follows  from
proposition~\ref{pr:homologous:saddle:connections}           that
$S_{(1,2)}$  has  trivial linear holonomy. But by assumption  the
graph $\Upsilon$  corresponding  to  $S_{(1,2)}$ has a nontrivial
loop,                            so                            by
Lemma~\ref{lm:nontriv:lin:holon:along:loops:of:the:graph}     the
flat surface $S_{(1,2)}$ has nontrivial linear homology, which is
a contradiction.

Thus, the graph has the structure of a union of  a  circle with a
segment attached to  a  circle. The  graph  has single vertex  of
valence  three,  a  single  vertex of valence one  and  arbitrary
number of vertices of  valence  two. Choosing an appropriate pair
of     edges      $\gamma_1,     \gamma_2$     and      combining
proposition~\ref{pr:homologous:saddle:connections}           with
Lemma~\ref{lm:nontriv:lin:holon:along:loops:of:the:graph}     and
Lemma~\ref{lm:subgraph:with:minus:vertex:has:nontrivial:holonomy}
we see that the only ``$-$''-vertex of the graph is the vertex of
valence one located at  the free end of the segment. This  is the
graph   of   the   type   c   in   the   list    of   graphs   in
theorem~\ref{th:graphs}.

{\bf A tree.}
In  this  case $\Gamma(S,\gamma)$ has at least  two  vertices  of
valence one which  are  therefore of ``$-$''-type. Let $\gamma_1$
and $\gamma_2$ be the edges  of  $\Gamma(S,\gamma)$  adjacent  to
this pair of vertices  $S^-_i$,  $S^-_j$. Cutting the surface $S$
by $\gamma_1,\gamma_2$ we get three connected components $S^-_i$,
$S^-_j$,      and      $S\setminus(S^-_i\cup     S^-_j)$.      By
proposition~\ref{pr:homologous:saddle:connections} the  component
$S\setminus(S^-_i\cup S^-_j)$ has trivial linear holonomy.  Thus,
by
Lemma~\ref{lm:subgraph:with:minus:vertex:has:nontrivial:holonomy}
it does  not  have any ``$-$''-vertices. Since $\Gamma(S,\gamma)$
is  a  tree  it  means  that  $\Gamma(S,\gamma)$ is  a  chain  of
``$+$''-vertices of valence two bounded at the ends by a  pair of
``$-$''-vertices  of  valence  one.  This  is  the  graph  a from
theorem~\ref{th:graphs}                 (see                 also
figure~\ref{fig:classification:of:graphs}).

Two \cv-vertices of valence  two  cannot be neighbors. It remains
to prove that a  \cv-vertex of valence $3$ cannot be joined  by a
separating edge to a \cv-vertex of valence $2$. If that  were the
case then  on one boundary  component of the cylinder there would
be a marked point. If  this  boundary component were joined to  a
\cv-vertex of valence $2$ it would produce a ``fake singularity''
on $S$.

We  have  proved  that  all  graphs  must  be  of  the   type  in
theorem~\ref{th:graphs}.  The   fact  that  the  weights  are  as
described follows from the next lemma.

\begin{Lemma}
\label{lm:lengths:of:saddle:connections}
The  type  of  the  graph $\Gamma(\cC)$ uniquely  determines  the
distribution of unsigned weights $1$ and $2$ on the edges  of the
graph;    the    corresponding    weights   are   presented    in
Figure~\ref{fig:classification:of:graphs}.
\end{Lemma}
\begin{proof}
For every vertex representing a  component  with  trivial  linear
holonomy we can choose  signs for the weights $1$ and $2$  on the
edges adjacent to the vertex. The sum of these signed  weights is
zero. This immediately implies that the globally defined unsigned
weights on both edges adjacent to a valence two ``$+$''-vertex or
to a valence two \cv-vertex  are  the same. This in turn  implies
that all the weights  on the graphs of types a) and  b) coincide,
and hence are marked by $1$.

The remaining graphs do  not  have ``$-$''-vertices. The edges of
any vertex of valence three are  weighted by $1, 1$ and $2$. This
implies that the weights of the graphs of types c) and d) are
as in figure~\ref{fig:classification:of:graphs}.

Let $\gamma_1, \gamma_2$ be a pair of edges adjacent to a valence
four vertex, and belonging to  two  different  loops. The surface
cut   along   these   saddle   connections   is   connected.   By
proposition~\ref{pr:homologous:saddle:connections},
$|\gamma_1|=|\gamma_2|$, and hence the  corresponding  edges have
the  same  weight.   Since   all  the   edges   in  a  chain   of
``$+$''-vertices or \cv-vertices of valence  $2$  have  the  same
weight,  we  see  that  all  edges of a graph of type e)  are
weighted by $1$.
\end{proof}

This   completes   the   proof   of   the   necessity   part   of
theorem~\ref{th:graphs}.

\end{proof}

\section{Parities of boundary singularities}
\label{s:Anatomy}

In    this    section   we   prove   the   necessity   part    of
theorem~\ref{th:all:local:ribon:graphs} which  says that for  any
decomposition of a flat surface $S$ as in theorem~\ref{th:graphs}
every connected component $S_j$ has one  of  the  boundary  types
presented in figure~\ref{fig:embedded:local:ribbon:graphs}.

theorem~\ref{th:graphs}                                       and
figure~\ref{fig:classification:of:graphs}  give  the   types   of
graphs $\Gamma$; figure~\ref{fig:local:ribbon:graphs}  gives  the
list of all abstract local ribbon graphs of valences from  one to
four.  Basically,  what  remains   to   check  is  that  for  any
``$+$''-vertex $v$ of $\Gamma$ an embedding  $\G_v\hookrightarrow
\Gamma$ of the local ribbon graph $\G_v$ into  the graph $\Gamma$
uniquely determines  the  parities of the boundary singularities,
and    that     these     parities     are    exactly    as    in
figure~\ref{fig:embedded:local:ribbon:graphs}.

\subsection*{Signs of the weights}
Given a  collection  $\gamma$  of \^homologous saddle connections
$\gamma_1,  \dots,  \gamma_n$  on  a  flat  surface  $S$  we have
assigned weights $1$ and $2$  to  saddle  connections  $\gamma_i$
(see  the  paragraph preceding  theorem~\ref{th:graphs}  for  the
definition          of          the          weights          and
figure~\ref{fig:classification:of:graphs} for the distribution of
the  weights  in  $\Gamma$).  If a connected component  $S_j$  of
$S\setminus\gamma$ has trivial linear holonomy  (i.e.  if  it  is
represented by  a ``$+$'' or by a \cv-vertex  of $\Gamma$) we may
assign signs  $\pm$ to the weights  of the saddle  connections on
the boundary of $S_j$. The canonical orientation of $S_j$ induces
the canonical orientation of the  boundary  $\partial  S_j$.  Let
$\omega$  be   a   holomorphic  $1$-form  representing  the  flat
structure on $S_j$ normalized so that
$$
\int_{\gamma_i}\omega =\text{weight of }\gamma_i,
$$
for some saddle connection  $\gamma_i$  on the boundary of $S_j$.
Then for  the other saddle connections  on $\partial S_j$  we get
$\int_{\gamma_{i'}}\omega =\pm 1$  or $\int_{\gamma_{i'}}\omega =
\pm 2$ (see also the tables in section~\ref{ss:tables}).

There  is   an  ambiguity  in  the   choice  of  signs:   we  may
simultaneously change the signs of  all  weights  to the opposite
ones. This corresponds to choosing $-\omega$ instead of $\omega$.

\begin{Lemma}
\label{lm:weights:define:parity}
Consider      two      consecutive       saddle       connections
$\gamma_{j_{i,l}}\to\gamma_{j_{i,l+1}}$  on  the  same  boundary
component   $\cB_i$   of  $\partial  S_j$.  The  parity  of   the
corresponding boundary  singularity  is  even  if  the weights of
$\gamma_{j_{i,l}}$ and $\gamma_{j_{i,l+1}}$ have the same  signs,
and   odd    if    the    weights   of   $\gamma_{j_{i,l}}$   and
$\gamma_{j_{i,l+1}}$ have opposite signs.
\end{Lemma}
\begin{proof}
The holomorphic 1-form $\omega$ chosen above  defines an oriented
horizontal  foliation   on   $S_j$:   the   kernel  foliation  of
$\Im(\omega)$. The  above  normalization of $\omega$ implies that
any  saddle  connection   at   the  boundary  $\partial  S_j$  is
horizontal. The  weight of a saddle connection $\gamma_{j_{i,l}}$
on  the  boundary of  $S_j$  is positive  if  the orientation  of
$\gamma_{j_{i,l}}$ induced  from  the orientation of the boundary
matches the orientation of the foliation and negative  if it does
not.

The cone angle between two  incoming  or  two outgoing separatrix
rays (in the  sense  of the orientation of  the  foliation) is an
even multiple of $\pi$ and the cone angle between an incoming and
an outgoing  separatrix ray (in the  sense of the  orientation of
the foliation) is an odd multiple of $\pi$. The statement  of the
lemma               now               follows                from
definition~\ref{def:order:of:boundary:singularity}  of  the order
of a boundary singularity.
\end{proof}

Consider now a  particular case when  $S_j$ is represented  by  a
vertex $v_j$ of valence four  of  the  graph  $\Gamma(S,\gamma)$.
Four edges of $\Gamma_{v_j}$ are distributed into two pairs: each
pair bounds one of the two loops of the graph $\Gamma(S,\gamma)$,
see figure~\ref{fig:classification:of:graphs}.

\begin{Lemma}
\label{lm:signs:of:weights:for:valence:4}
The weights of saddle connections on the boundary  of a component
$S_j$ represented by a vertex of valence four have same  signs if
they bound the same loop in $\Gamma(S,\gamma)$ and opposite signs
otherwise.
\end{Lemma}
\begin{proof}
From  lemma~\ref{lm:lengths:of:saddle:connections}  we know  that
the absolute values of weights of all edges of $\Gamma_{v_j}$ for
a  vertex $v_j$ of  valence  four  are  equal to  one  (see  also
figure~\ref{fig:classification:of:graphs}).  Hence,  it   follows
from Stokes theorem that we have two edges of weight $+1$ and two
edges of weight $-1$ in $\Gamma_{v_j}$. We want to show  that the
weights of  a pair of edges  of $\Gamma_{v_j}$ bounding  the same
loop in $\Gamma$ have the same signs.

Let $\gamma_1, \gamma_2\in \Gamma_{v_j}$ bound \textit{different}
loops in  $\Gamma$. Cutting $S$  by $\gamma_1, \gamma_2$ we get a
connected flat surface $S_{12}$.  Using  the same notation as in
the  proof  of proposition~\ref{pr:homologous:saddle:connections}
we get
$$
\partial S_{12}=\gamma'_1\cup -\gamma''_1\cup\gamma_2'\cup-\gamma''_2
$$
By theorem~\ref{th:unique:trivial:holonomy}  the surface $S_{12}$
has  trivial  linear  holonomy.  Hence,  we  can extend the  form
$\omega$ to $S_{12}$  which  enables us  to  assign signs to  the
weights of saddle connections $\gamma'_1, \gamma''_1,  \gamma'_2,
\gamma''_2$ on the boundary $\partial S_{12}$ of $S_{12}$.

The                last                statement               of
lemma~\ref{lm:total:boundary:holonomy:is:trivial}  implies   that
gluing the initial  closed surface $S$ from $S_{12}$ the boundary
component $\gamma'_1$ is glued to  $-\gamma''_1$  by  a flip (see
definition~\ref{def:flip}                                      in
section~\ref{s:hat:homologous:saddle:connections}).     Similarly
$\gamma'_2$  is  glued  to  $-\gamma''_2$ by a flip.  Hence,  the
weights of $\gamma'_1$ and  of  $\gamma''_1$ have the same signs,
and the weights of $\gamma'_2$ and of $\gamma''_2$  have the same
signs.

This  completes  the  proof of the  lemma  in  the  case when the
corresponding loop contains no  vertices  at all. An induction on
the number of vertices in the loop completes the proof in general
case.
\end{proof}

\begin{Lemma}
\label{lm:signed:lengths:of:saddle:connections}
For  any  ``$+$''-vertex  or  \cv-vertex   $v$   of   the   graph
$\Gamma(S,\gamma)$ the type of the graph  uniquely determines the
distribution of  signed weights $\pm 1$ and $\pm  2$ on the edges
of $\Gamma_v$ (up to simultaneous  interchange  of  all signes to
the opposite ones).
\end{Lemma}
\begin{proof}
By Stokes theorem the sum of weights of all saddle connections of
$\Gamma_v$  is   equal   to   zero.   Taking  into  consideration
lemma~\ref{lm:lengths:of:saddle:connections}      (see       also
figure~\ref{fig:classification:of:graphs}) this implies that when
the vertex $v$  has  valence  $2$, the weights of  the  edges  of
$\Gamma_v$ are $+1,-1$; when $v$ has valence $3$, the weights are
$+1,+1,-2$;  when  $v$  has  valence   $4$,   the   weights   are
$+1,+1,-1,-1$.          Moreover,          according           to
lemma~\ref{lm:signs:of:weights:for:valence:4}  the  weights   of
edges   of   $\Gamma_v$   which   bound   the    same   loop   in
$\Gamma(S,\gamma)$ coincide. \end{proof}


Now  we  are  ready  to  prove  the following proposition,  which
corresponds       to      the       necessity       part       of
theorem~\ref{th:all:local:ribon:graphs}. (The sufficiency part of
theorem~\ref{th:all:local:ribon:graphs} immediately follows  from
theorem~\ref{th:from:boundary:to:neighborhood:of:the:cusp} proved
in the next section).

\begin{Proposition}
For   any   decomposition   of   a   flat   surface  $S$  as   in
theorem~\ref{th:graphs} every connected component  $S_j$  has one
of       the       boundary        types       presented       in
figure~\ref{fig:embedded:local:ribbon:graphs}.
\end{Proposition}
\begin{proof}
The   necessity   part   of  theorem~\ref{th:graphs}  proved   in
section~\ref{s:graph} claims that the graph $\Gamma(S,\gamma)$ of
the   decomposition   has   one   of  the  types   presented   in
figure~\ref{fig:classification:of:graphs}. Note  that for ``$+$''
and                                             ``$-$''-vertices,
figure~\ref{fig:embedded:local:ribbon:graphs}           describes
\textit{all} possible  embeddings of abstract local ribbon graphs
$\G_v$  that are  given  in  figure~\ref{fig:local:ribbon:graphs}
into           graphs           $\Gamma$          as           in
figure~\ref{fig:classification:of:graphs}. We use dotted lines to
indicate the pairs  of  edges bounding  cycles  in the graphs  in
figure~\ref{fig:embedded:local:ribbon:graphs};  dotted  lines are
not  indicated  in  symmetric  situations.  Since  there  are  no
restrictions  on   the  parities  of  boundary  singularities  of
``$-$''-vertices this complets the proof for ``$-$''-vertices.

Any \cv-vertex  $S^{comp}_j$  corresponds  to  a  flat  cylinder.
Hence,  it  has  exactly  two distinct boundary  components.  The
boundary  singularities  on each of the components correspond  to
marked points, so  the  order of  any  boundary singularity of  a
\cv-vertex is zero.  By lemma~\ref{lm:weights:define:parity} this
implies that all  edges of $\Gamma_{v_j}$ which correspond to the
same boundary component of the cylinder $S_j$ have weights of the
same        sign.         Taking        into        consideration
lemma~\ref{lm:signed:lengths:of:saddle:connections}   these   two
conditions restrict the possible structures of  an embedded local
ribbon     graph     $\G_v\hookrightarrow\Gamma(S,\gamma)$    for
\cv-vertices  to  structures $\circ 2.2$, $\circ 3.2$ and  $\circ
4.2$ in figure~\ref{fig:embedded:local:ribbon:graphs}.

By  lemma~\ref{lm:signed:lengths:of:saddle:connections}  for  any
``$+$''-vertex of  $\Gamma(S,\gamma)$  we know the signed weights
of the edges of $\Gamma_v$ (up to simultaneous interchange of all
signs to the  opposite ones). For ``$+$''-vertices of valence two
and   three   this   distribution   follows   immediately    from
figure~\ref{fig:classification:of:graphs}    and    from   Stokes
theorem; for  ``$+$''-vertices  of valence four this distribution
is  described  by  lemma~\ref{lm:signs:of:weights:for:valence:4}.
Hence,   using    lemma~\ref{lm:weights:define:parity}   we   can
determine the  parities  of  all  boundary  singularities for any
embedded  local  ribbon graph  $\G_v\hookrightarrow  \Gamma$.  It
remains  to  check  that  for all possible embeddings  listed  in
figure~\ref{fig:embedded:local:ribbon:graphs}  the  parities  are
the ones listed. This is an easy exercise.
\end{proof}


\begin{Corollary}
\label{cr:unique:choice:of:parities}
Given any  abstract  graph $\Gamma$ as in theorem~\ref{th:graphs}
(see figure~\ref{fig:classification:of:graphs}), any  ``$+$''  or
``$-$''-vertex $v_j$ of  $\Gamma$, any choice of the structure of
a  local  ribbon  graph  $\G_{v_j}$  on  $\Gamma_{v_j}$  and  any
embedding $\G_{v_j}\hookrightarrow  \Gamma$,  one can find a flat
surface $S$  and  a  collection  $\gamma$  of \^homologous saddle
connections on  it  such  that $\Gamma(S,\gamma)=\Gamma$ and such
that the  boundary type of  the component $S_j$ is represented by
the  chosen  embedded   ribbon   graph.  Moreover,  if  $v_j$  is
represented by a ``$+$''-vertex of $\Gamma$, then the parities of
boundary singularities of $S_j$ are completely  determined by the
choice of the embedded ribbon graph.

Conversely,   given    an   abstract   graph   $\Gamma$   as   in
theorem~\ref{th:graphs}                                      (see
figure~\ref{fig:classification:of:graphs}),   a    ``$+$''-vertex
$v_j$ of $\Gamma$, an abstract local ribbon graph $\G_{v_j}$, and
a choice  of the parities of  boundary singularities as  given in
figure~\ref{fig:embedded:local:ribbon:graphs} there  is a unique
way (up to a symmetry  of  the ribbon graph $\G_{v_j}$) to  embed
the  local  ribbon  graph  with  marked  parities into the  graph
$\Gamma$. This unique  way is expressed  by the dotted  lines  in
figure~\ref{fig:embedded:local:ribbon:graphs}.
\end{Corollary}
\begin{proof}
For ``$+$'' and ``$-$''-vertices $v_j$ all possible embeddings of
local  ribbon  graphs as  in figure~\ref{fig:local:ribbon:graphs}
into        the        graphs        $\Gamma$        as        in
figure~\ref{fig:classification:of:graphs}   are   represented  in
figure~\ref{fig:embedded:local:ribbon:graphs}.  Thus,  the  first
statement  follows  from theorem~\ref{th:all:local:ribon:graphs}.
The     second     statement     immediately     follows     from
theorem~\ref{th:all:local:ribon:graphs}       combined       with
lemmas~\ref{lm:weights:define:parity}
and~\ref{lm:signed:lengths:of:saddle:connections}.
\end{proof}

\section{Neighborhood of the principal boundary: local constructions}
\label{s:Local:Constructions}

In this section and in  the  next one we construct surfaces  with
boundaries   representing   all   boundary   types   listed    in
figure~\ref{fig:embedded:local:ribbon:graphs}. We first prove the
key  proposition  below. Combining it with some elementary  extra
arguments                        we                         prove
theorem~\ref{th:from:boundary:to:neighborhood:of:the:cusp}  (and,
hence,      the     missing      realizibility      parts      of
theorems~\ref{th:graphs} and~\ref{th:all:local:ribon:graphs}).

\begin{Proposition}
\label{pr:realizability:of:all:vertices}
Consider      any       configuration      $\cC$      as       in
definition~\ref{def:configuration}, and any vertex  $v_j$  of the
graph  $\Gamma(\cC)$.  Let $S'_j$ be any flat  surface  from  the
component $\cQ(\alpha'_j)$ (or $\cH(\beta'_j)$) of the  principal
boundary stratum $\cQ(\alpha')$  (or $\cH(\beta')$) corresponding
to  $v_j$.  Choose  any  sufficiently  small  value of a  complex
parameter $\delta$ (depending on $S'_j$).

Applying to $S'_j$ an appropriate  basic  surgery  (depending  on
$\delta$)  as  described  below  one  gets  a surface $S_j$  with
boundary,  such  that  the   boundary   type  of  $S_j$  and  the
collections   of   interior   singularities   and   of   boundary
singularities of $S_j$  are represented by the local ribbon graph
$\G_{v_j}$  and  by  the  corresponding  structures   $\{d_l\}_j,
\{k_{i,l}\}_j$ of the configuration $\cC$.
\end{Proposition}

Recall that the principal boundary  stratum  corresponding  to  a
``$+$''-vertex is of type $\cH(\beta'_j)$; the principal boundary
stratum  corresponding   to   a   ``$-$''-vertex   is   of   type
$\cQ(\alpha'_j)$. The singularity data $\beta'_j, \alpha'_j$  are
defined        by         equations~\eqref{eq:H:boundary:stratum}
and~\eqref{eq:Q:boundary:stratum}  correspondingly.  Unlike   the
initial singularity  data $\alpha$ the collections $\beta'_j$ and
$\alpha'_j$  might  contain  entries ``$0$'' representing  marked
points of the surface $S'_j$.

Though  the   principal   boundary  stratum  corresponding  to  a
\cv-vertex                       is                        empty,
proposition~\ref{pr:realizability:of:all:vertices}     is     not
meaningless (though very simple) even for such vertices. We leave
the construction of  surfaces  $S_j$ with boundary realizing each
of $\circ 2.2, \circ 3.2, \circ 4.2$-boundary types to the reader
as an elementary exercise.

We split proposition~\ref{pr:realizability:of:all:vertices}  into
a   collection   of    propositions~\ref{pr:local:constructions},
\ref{pr:parallelogram:constructions}  and~\ref{pr:minus:2:2}.  To
avoid excessive repetitions  we  abbreviate the statements of the
corresponding propositions; they should be read  as the statement
of proposition~\ref{pr:realizability:of:all:vertices} applied  to
vertices of specified types.

Part of the surgeries (namely,  ``breaking  up a zero by a  local
construction'' and  a  ``parallelogram  construction'') are taken
from  the  paper~\cite{Eskin:Masur:Zorich}.  For   the   sake  of
completeness we present their  outline  in the current paper. For
more   details   we  address   the   reader   to   the   original
paper~\cite{Eskin:Masur:Zorich}.

\input{ribbon_versus_holes_submit.tex}

\subsection{Local constructions}
We reserve  the word ``degree'' for  the zeroes of  {\it Abelian}
differentials. A zero of  degree  $l$ has cone angle $\pi(2l+2)$.
We reserve  the word ``order''  for the zeroes of {\it quadratic}
differentials. A zero of order $m$  has  cone  angle  $\pi(m+2)$.
Recall that  a boundary singularity of  order $k$ has  cone angle
$\pi(k+1)$.

We distinguish two kinds of surgeries. The surgeries of the first
type  are  purely local: they do  not  change the flat metric  on
$S'_j$  outside a  small  neighborhood of one  or  two points  on
$S'_j$. The surgeries  of the second  type depend on  a  nonlocal
construction. In the remaining  part  of this section we describe
local surgeries.

\begin{Proposition}
\label{pr:local:constructions}
Every  surface   with  boundary  type  $+2.1$,  $+3.1$,  $+4.1a$,
$+4.1b$, $+4.2a$ is realizable by a local construction.
\end{Proposition}

We   use   the   indexation   of   the   boundary  types  as   in
figure~\ref{fig:embedded:local:ribbon:graphs}       and        in
remark~\ref{rm:indexation}                                     in
section~\ref{ss:Anatomy:of:decomposition:of:a:flat:surface}.

\begin{proof}
The principle boundary stratum corresponding to a vertex $v_j$ of
a ``$+$''-type  has  type  $\cH(\beta'_j)$.  The singularity data
$\beta'_j$  is  given by  equation~\eqref{eq:H:boundary:stratum},
namely
$$
\beta'_j=\{d_1/2, \dots, d_{s(j)}/2,\ D_1/2, \dots,  D_{r(j)}/2\},
$$
where  $d_1,  \dots,   d_{s(j)}$   are  the  orders  of  interior
singularities, and $D_1, \dots, D_{r(j)}$ are  expressed in terms
of     the     orders    of     boundary     singularities     by
formula~\eqref{eq:Dji}.     Conditions     4     and     5     in
definition~\ref{def:configuration} of a  configuration  guarantee
that  all  the entries  of  $\beta'_j$  are  nonnegative  integer
numbers,  and  that the  total  sum  of  these  numbers  is even.
According to~\cite{Masur:Smillie:realizability} this implies that
the stratum $\cH(\beta'_j)$ is nonempty.

Consider any surface $S'_j$ in $\cH(\beta'_j)$. Denote the length
of the shortest saddle connection on $S'_j$ by $4\varepsilon$. We
shall apply a  surgery to $S'_j$, which would continuously depend
on a  small complex parameter considered as a  vector $\vec v$ in
$\R{2}\simeq\C{}$.  It  is  convenient  to  change  slightly  the
notations and to denote by  $\delta$  the  \textit{norm} of $\vec
v$. We always assume that $\delta<\varepsilon$. Our surgery would
not affect interior singularities of $S'_j$.

We  provide all  the details  of  the proof  in the  case of  the
boundary type $+2.1$  and  we point  out  the differences in  the
other cases.

\subsection*{Boundary type +2.1}
In this case the boundary has  single  component,  $r(j)=1$,  and
$D_1=k_{1,1}+k_{1,2}-2$, where $k_{1,1}, k_{1,2}$ are the  orders
of the two  boundary  singularities of $\G_{v_j}$. Both $k_{1,1},
k_{1,2}$      are      odd       positive      integers,      see
figure~\ref{fig:embedded:local:ribbon:graphs}.

Let $P$ be the zero of  $S'_j$ of degree $m$, where $m=D_1/2$. We
can represent $m$ as the sum $m=m'+m''$, where $m'=(k_{1,1}-1)/2$
and  $m''=(k_{1,2}-1)/2$.   Consider  a  metric  disc  of  radius
$\varepsilon$ centered at $P$. By the choice of $\varepsilon$ the
disc  does   not   contain   any   other   singularities  and  is
isometrically embedded into $S'_j$. It can be glued from $2(m+1)$
copies of standard metric half-discs of radius $\varepsilon$; see
the picture  at  the  top of figure~\ref{fig:breaking:up:a:zero}.
Let $\vec{v}\in\R{2}$ be a vector of length $\delta<\varepsilon$.
Following~\cite{Eskin:Masur:Zorich}  we  may  {\it break up}  the
zero $P$ of degree $m$ into  a pair of zeroes of degrees $m'$ and
$m''$ joined by a  single  saddle connection with affine holonomy
$\vec{v}$.

\begin{figure}[ht]
%
\includegraphics{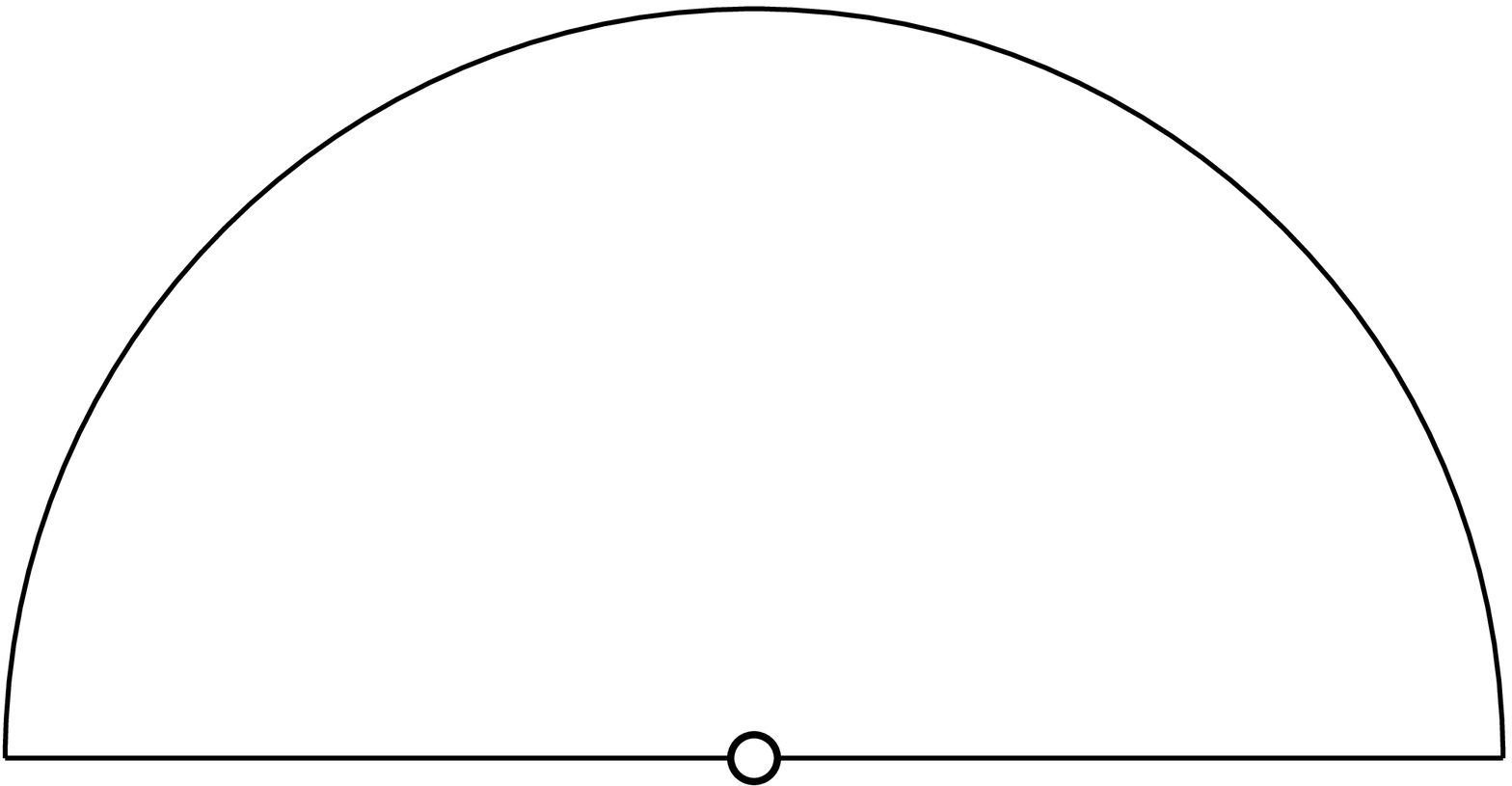}
\includegraphics{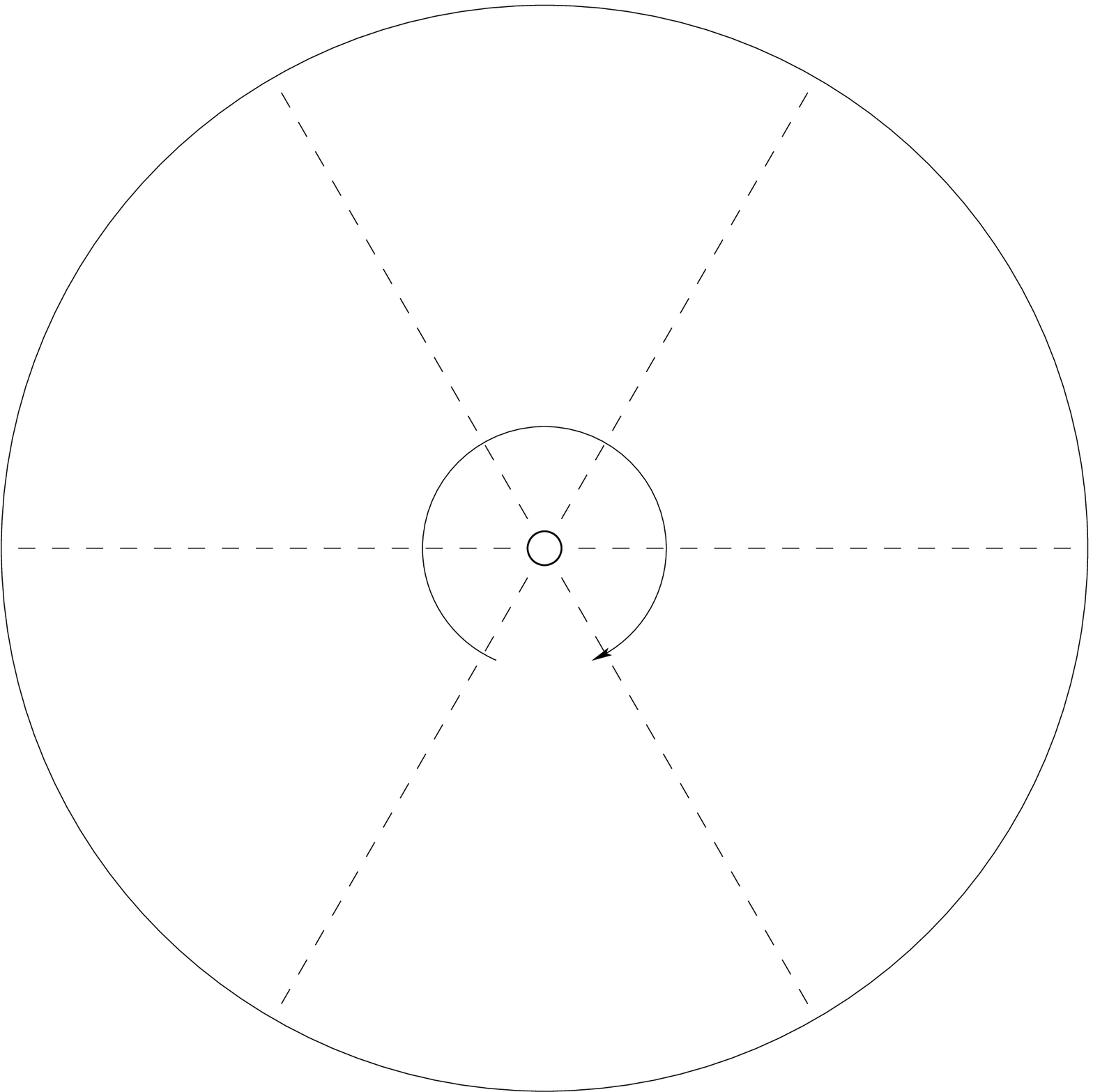}
\includegraphics{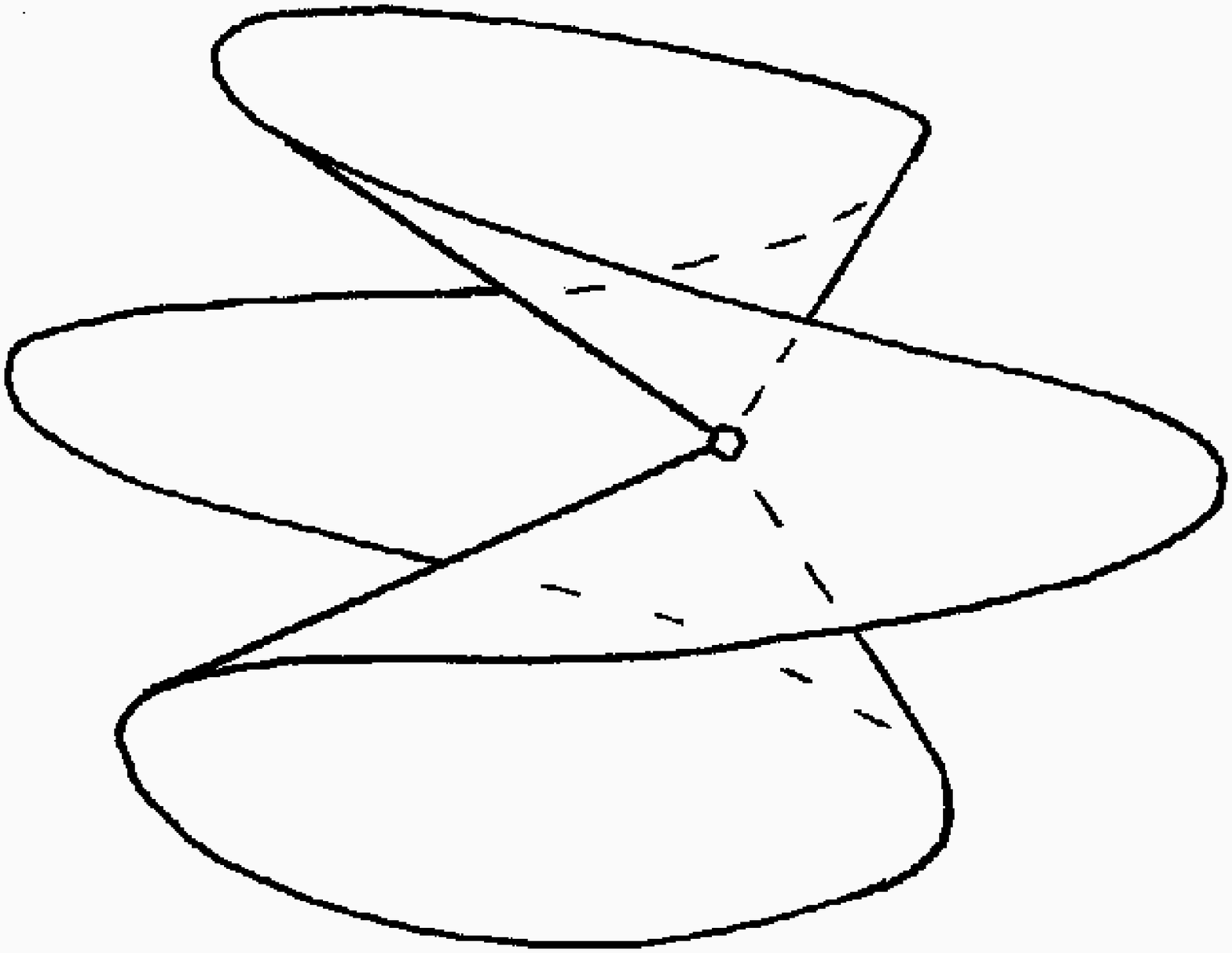}
%
%
\includegraphics{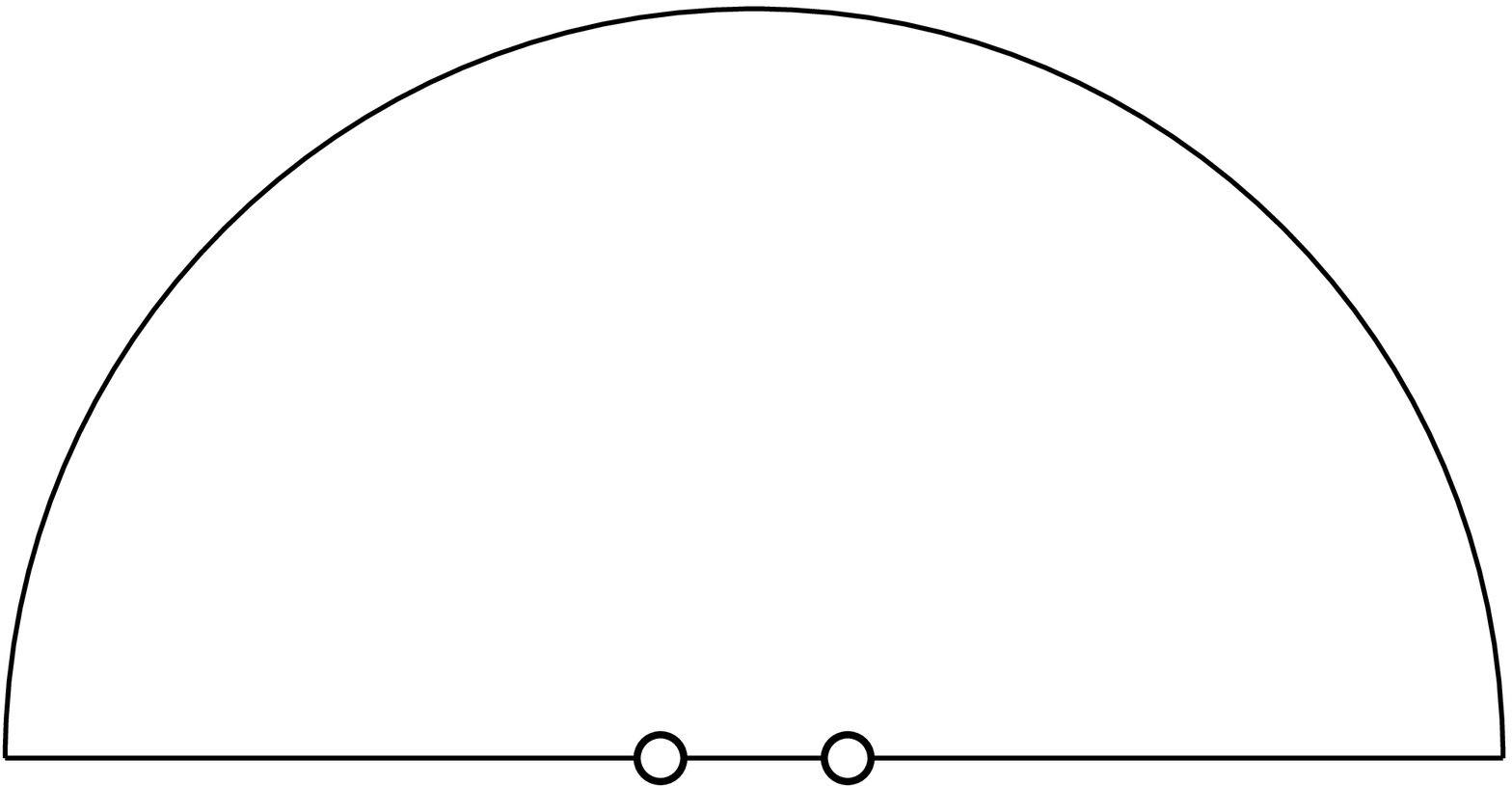}
\includegraphics{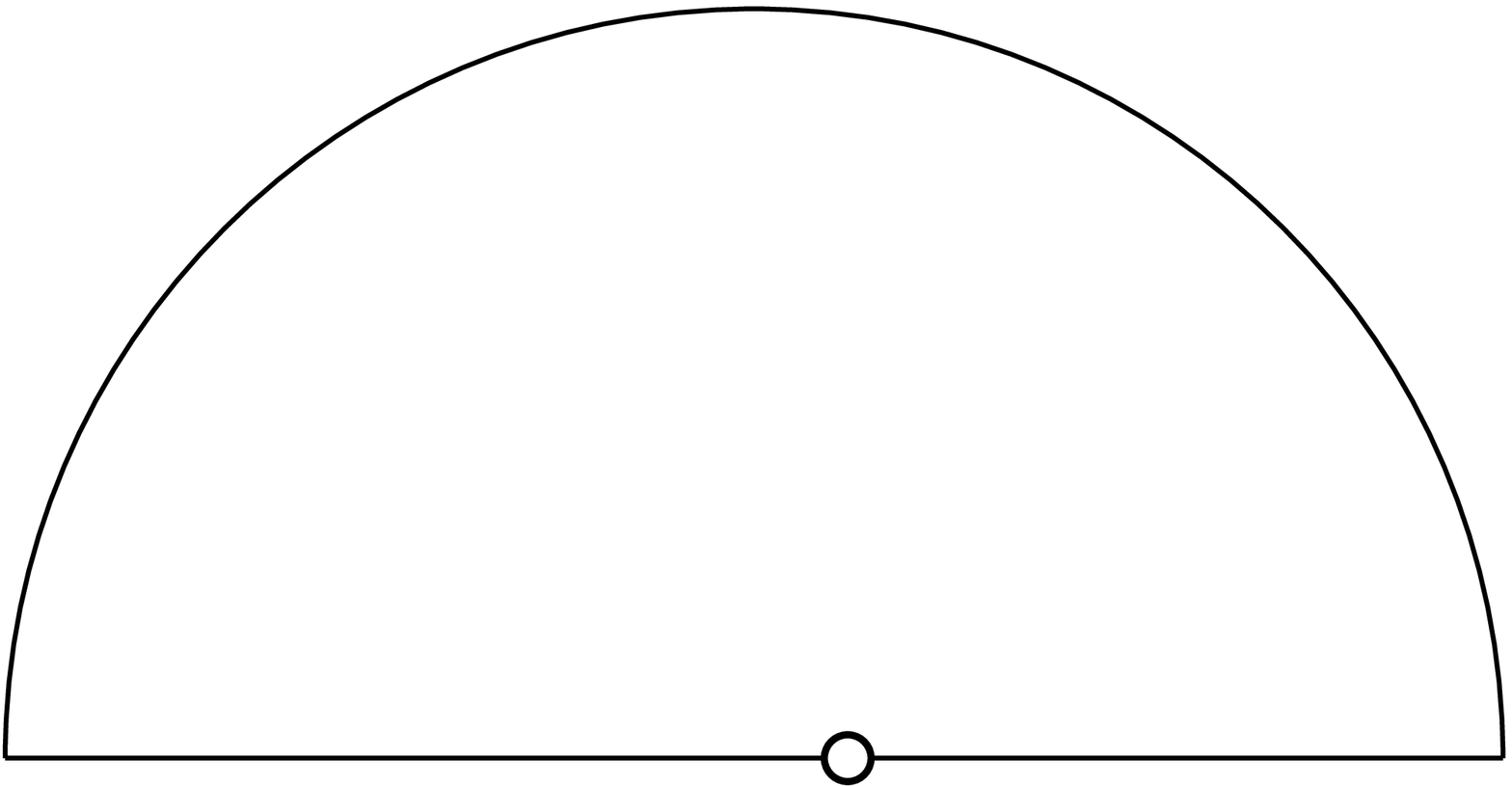}
\includegraphics{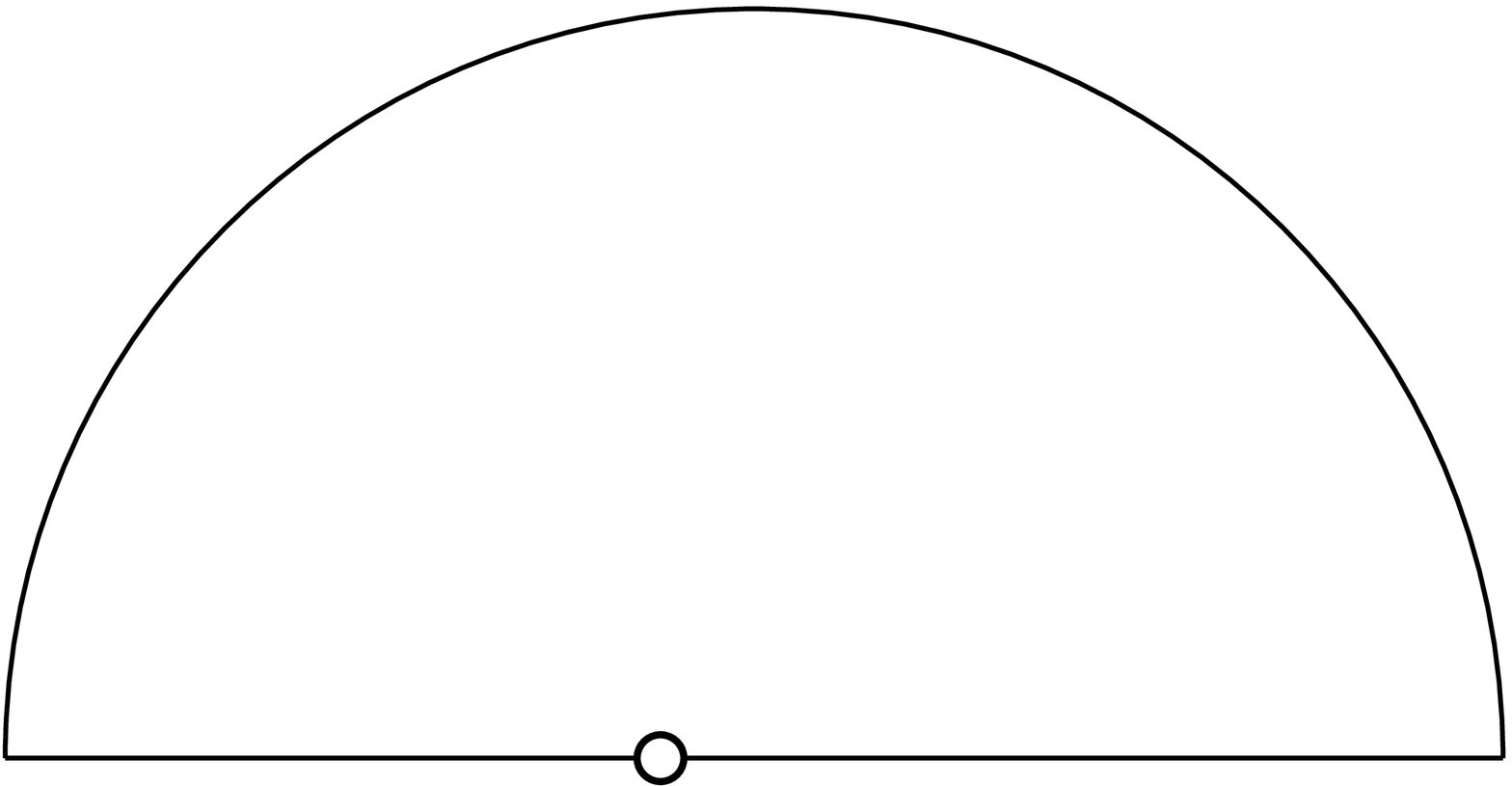}
\includegraphics{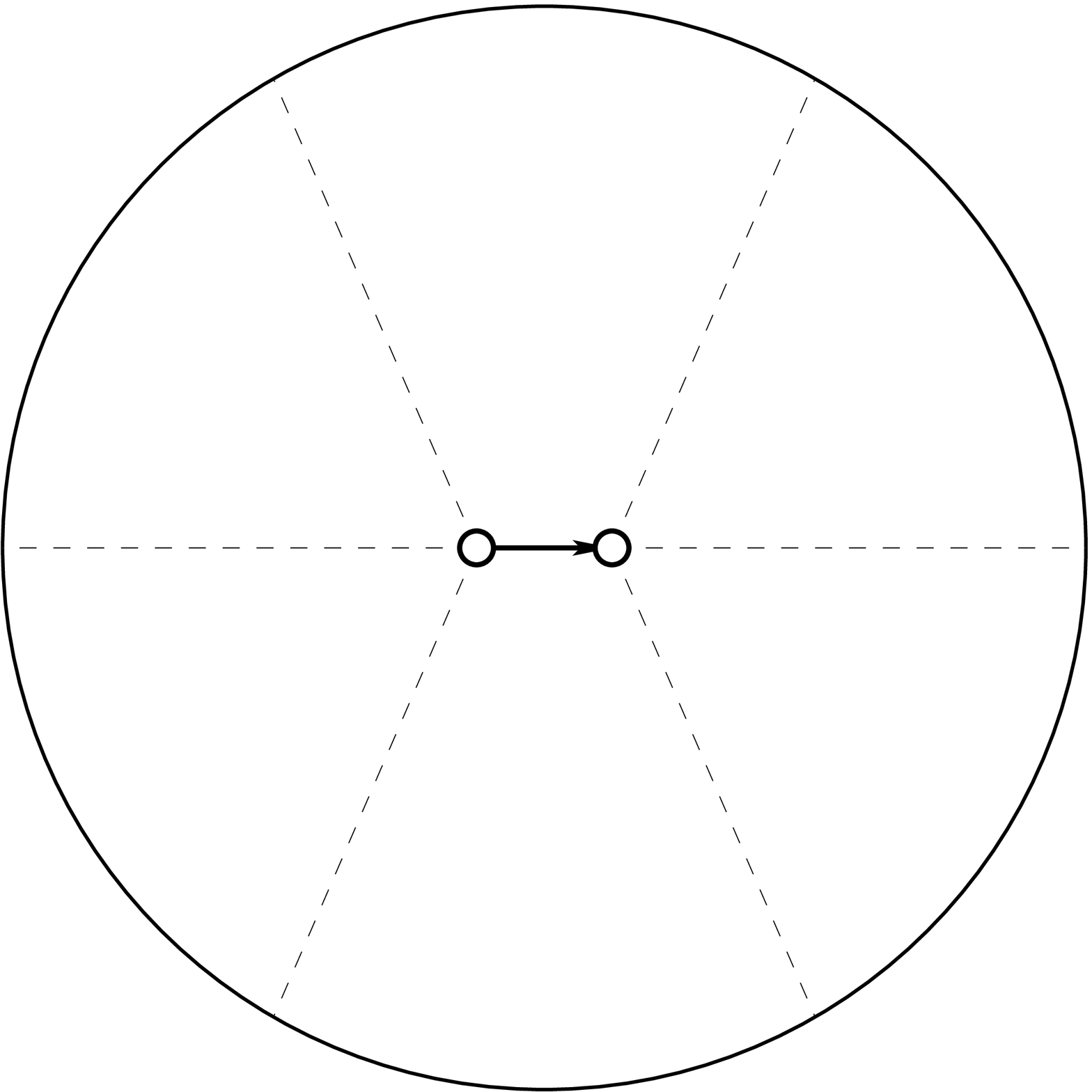}
\includegraphics{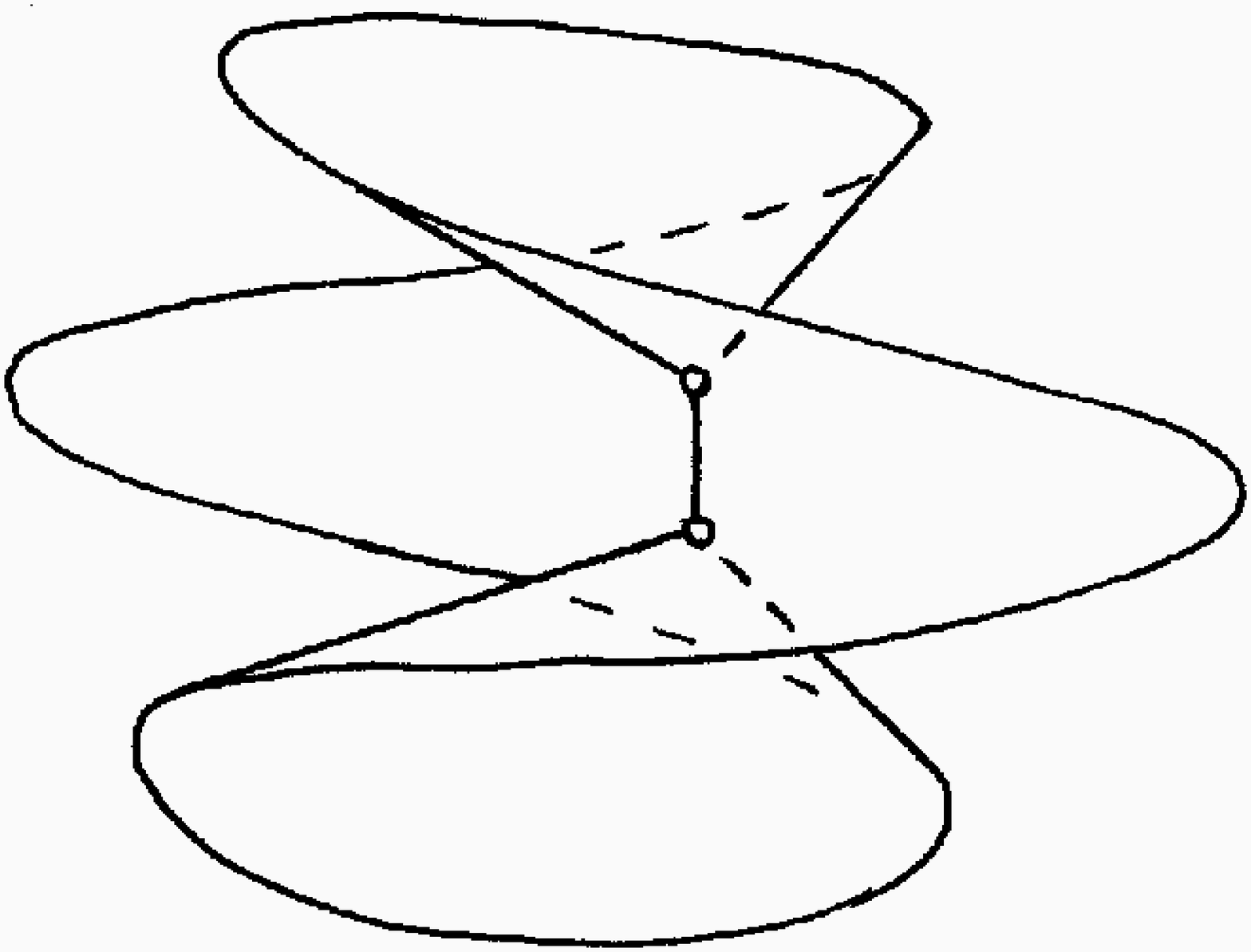}
%
%
\begin{picture}(0,0)(0,-20)
\put(10,-10)
 {\begin{picture}(0,0)(0,0)
 \put(-154,-55){$\scriptstyle \varepsilon$}
 \put(-125,-55){$\scriptstyle \varepsilon$}
 \put(-38,-66){$\scriptstyle 6\pi$}
 \end{picture}}
\put(10,-30)
 {\begin{picture}(0,0)(0,0)
 \put(-145,-148){$\scriptstyle 2\delta$}
 \put(-160,-197){$\scriptstyle \varepsilon+\delta$}
 \put(-130,-197){$\scriptstyle \varepsilon-\delta$}
 \put(-160,-247){$\scriptstyle \varepsilon-\delta$}
 \put(-130,-247){$\scriptstyle \varepsilon+\delta$}
 \put(-39,-197){$\scriptstyle 2\delta$}
 \put(-65,-197){$\scriptstyle \varepsilon+\delta$}
 \put(-19,-197){$\scriptstyle \varepsilon+\delta$}
 \put(-64,-173){$\scriptstyle \varepsilon-\delta$}
 \put(-19,-173){$\scriptstyle \varepsilon-\delta$}
 \put(-65,-217){$\scriptstyle \varepsilon-\delta$}
 \put(-17,-217){$\scriptstyle \varepsilon-\delta$}
\end{picture}}
\end{picture}
\vspace{260bp} 
\caption{
\label{fig:breaking:up:a:zero}
Breaking up a zero into two zeroes
(after~\cite{Eskin:Masur:Zorich})}
\end{figure}

We  do  this  by changing the  way  of  gluing  the half-discs as
indicated       on       the        bottom       picture       of
figure~\ref{fig:breaking:up:a:zero}. As patterns we still use the
standard metric half-discs, but  move  the marked points on their
diameters. Two special half-discs  have  two marked points on the
diameter  at  distance  $\delta$  from  the  center. Each of  the
remaining $2m$ half-discs has  a  single marked point at distance
$\delta$ from the center. We  alternate  the  half-discs with the
marked point  moved to  the right and to the  left of the center.
The picture  shows  that  all  the  lengths along identifications
match; gluing the half-discs we obtain a topological  disc with a
flat metric. Now  the flat metric has two cone-type singularities
with cone angles  $2\pi(m'+1)$  and $2\pi(m''+1)$. Here $2m'$ and
$2m''$ are the numbers of half-discs with one  marked point glued
in between the distinguished  pair  of half-discs with two marked
points.

The case when one of $m', m''$ (or both of them) is equal to zero
is  not  excluded,  in  this case the  corresponding  ``newborn''
singularity is just a marked point.

Note that a  small annular neighborhood  of the boundary  of  the
initial  disc   is   isometric   to   the  corresponding  annular
neighborhood of  the boundary of  the deformed disc. Thus, we can
glue the deformed disc back into the surface. Gluing back  we can
turn it by any angle $\varphi$, where $0\le\varphi< 2\pi(m+1)$ in
such  way  that  the  newborn  saddle  connection  will  have the
prescribed affine holonomy $\vec v$.

Making a slit  along  the resulting  saddle  connection we get  a
surface $S_j$  with  boundary  having  prescribed  boundary  type
$+2.1$, a pair of boundary  singularities  of  prescribed  orders
$k_{1,1}, k_{1,2}$, and a collection of interior singularities of
prescribed orders  $d_1,  \dots, d_{s(j)}$ (see the corresponding
entry in the table in section~\ref{ss:tables}). We have completed
the  proof of  proposition~\ref{pr:local:constructions}  for  the
boundary type $+2.1$.

\subsection*{Boundary types +3.1 and +4.1a}
Boundary type  $+3.1$ can be considered  as a particular  case of
boundary  type  $+4.1a$. To see this compare  the  surfaces  with
boundary representing  the  corresponding  ribbon graphs (see the
appropriate  entries  in  the table in  section~\ref{ss:tables}).
Marking a  point in the  middle of the saddle connection labelled
by ``$+2$'' on the boundary of  the surface of type $+3.1$ we get
a  surface  with  boundary  type  $+4.1a$,   where  the  boundary
singularity joining the pair of  edges  labelled  by ``$+1$'' has
order $0$.

Consider  a  local  ribbon  graph of type $+4.1a$,  a  collection
$\{2m_1, \dots, 2m_n\}$ of orders of interior singularities and a
collection  $\{2a_1+1,2a_2,2a_3+1,2a_4\}$   of  orders  of   four
boundary                    singularities                    (see
figure~\ref{fig:embedded:local:ribbon:graphs}      for      their
parities). The  singularity  data $\beta'_j$ of the corresponding
component $\cH(\beta'_j)$  of  the principal boundary stratum has
the  form  $\beta'=\{m_1,  \dots,  m_n,  a_1+a_2+a_3+a_4\}$,  see
equations~\eqref{eq:Dji} and~\eqref{eq:H:boundary:stratum}.

\begin{figure}[ht]
%
\includegraphics{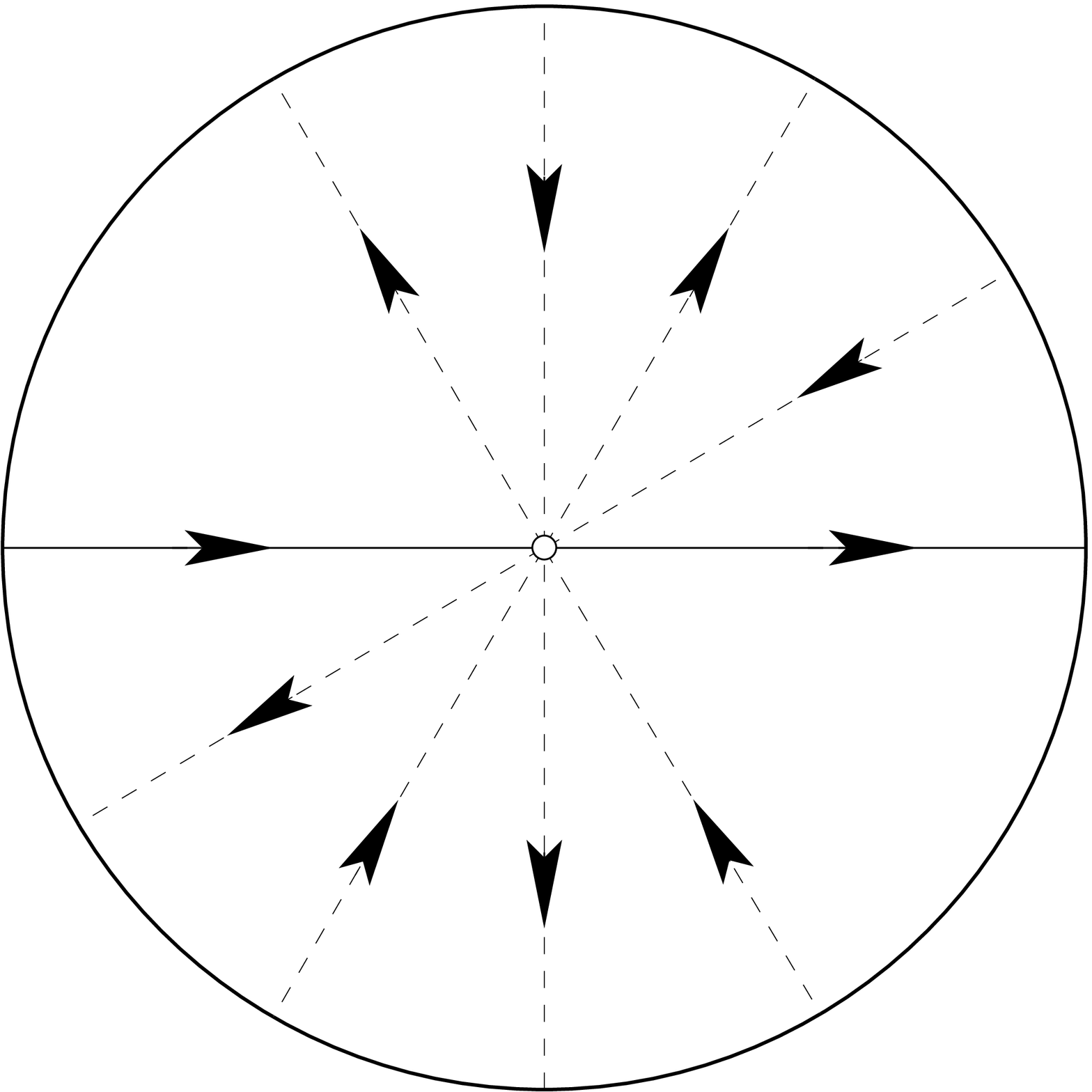}
\includegraphics{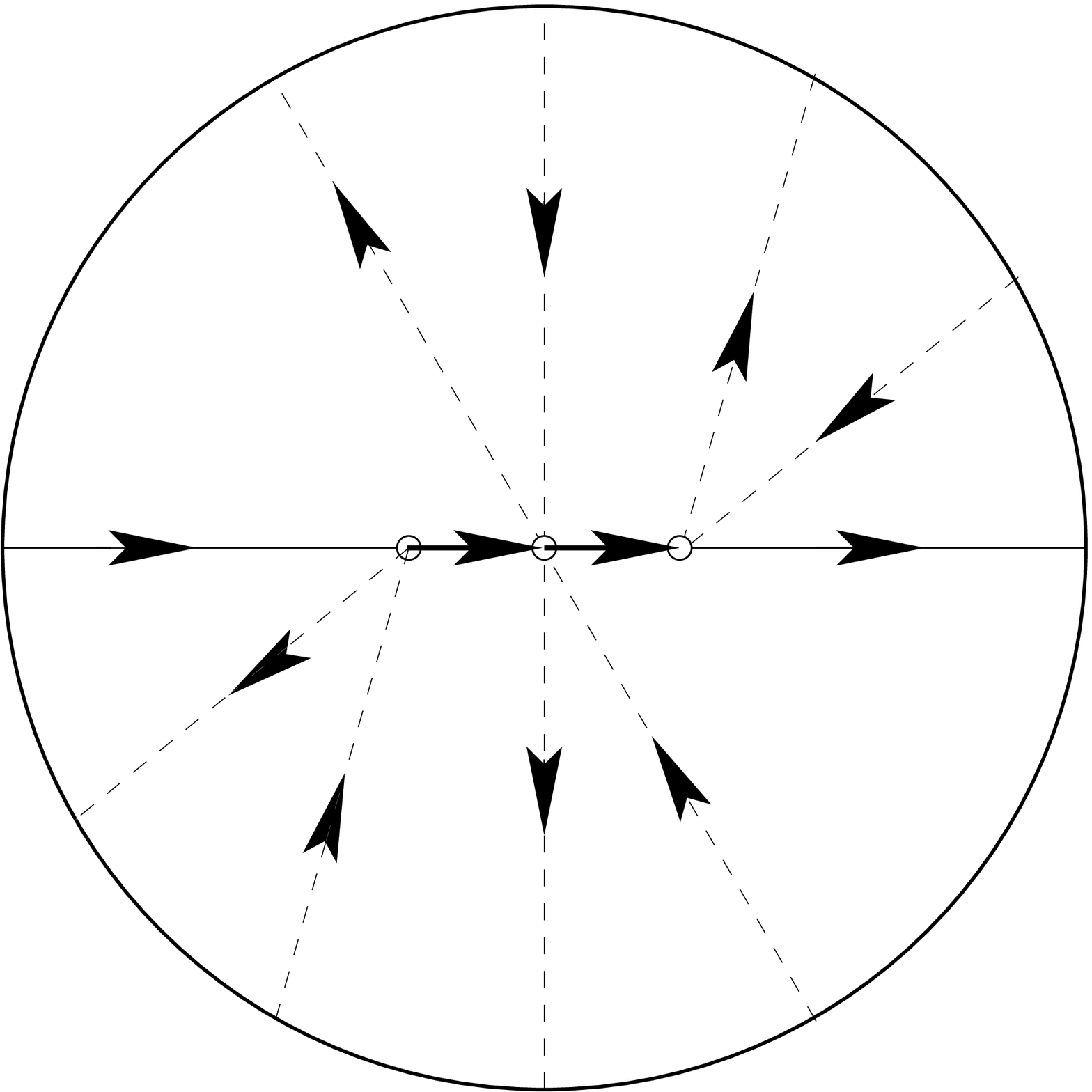}
\includegraphics{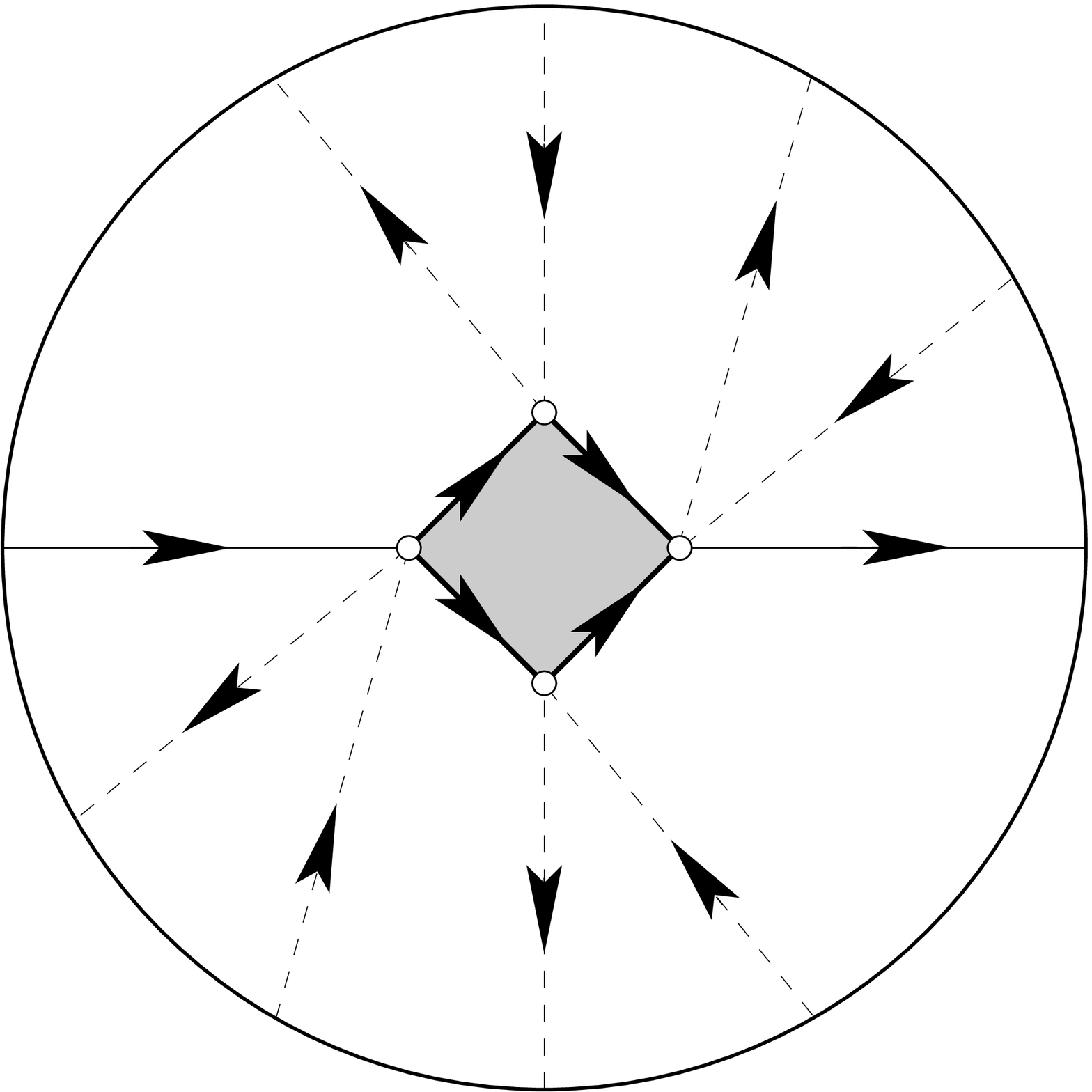}
%
%
\begin{picture}(0,0)(0,-150)
\put(10,-6)
 {\begin{picture}(0,0)(0,0)
 \put(-98,-199){$\scriptstyle \varepsilon$}
 \put(20,-199){$\scriptstyle \varepsilon-\delta$}
 \put(28,-185){$\scriptstyle \varepsilon+\delta$}
 \put(3,-160){$\scriptstyle \varepsilon-\delta$}
 \put(3,-209){$\scriptstyle P_3$} 
 \put(-9,-199){$\scriptstyle P$}
 \put(-8,-173){$\scriptstyle \varepsilon$}
 \put(-29,-178){$\scriptstyle \varepsilon$}
 \put(-52,-199){$\scriptstyle \varepsilon-\delta$}
 \put(-4,-209){$\scriptstyle \delta$}
 \put(-20,-209){$\scriptstyle \delta$}
 \put(-28,-199){$\scriptstyle P_1$}
 \put(-29,-229){$\scriptstyle \varepsilon-\delta$}
 \put(-58,-217){$\scriptstyle \varepsilon+\delta$}
 \put(-8,-236){$\scriptstyle \varepsilon$}
 \put(10,-232){$\scriptstyle \varepsilon$}
 \put(125,-209){$\scriptstyle P_3$} 
 \put(89,-198){$\scriptstyle P_1$}
 \put(112,-185){$\scriptstyle P_4$}
 \put(112,-221){$\scriptstyle P_2$}
\end{picture}}
\end{picture}
\vspace{115bp} %
\caption{
\label{fig:4:1a}
Breaking  up   a   zero  into  three  ones  and
performing a surgery we get a surface of type ``$+4.1a$''}
\end{figure}

   %

Choose an  Abelian  differential  $S'_j\in\cH(\beta')$. As before
denote the length of the shortest saddle connection  on $S'_j$ by
$4\varepsilon$. This time  we split the distinguished zero $P$ of
degree $a_1+a_2+a_3+a_4$ into {\it three} zeroes $P_1,P,P_3$ such
that the zero $P_1$ of degree  $a_1$ is joined to the zero $P$ of
degree $a_2+a_4$ by  a  saddle connection,  and  the zero $P$  of
degree $a_2+a_4$ is joined to the zero $P_3$ of degree $a_3$ by a
saddle  connection,  see figure~\ref{fig:4:1a}.  The  two  saddle
connections have the same holonomy vector $\vec{v}$. We assume as
before that $\|\vec{v}\|=\delta<\epsilon$. We then cut along both
saddle  connections  and detach the zero $P$  into  two  boundary
singularities  $P_2,   P_4$   of   orders   $2a_2$   and   $2a_4$
correspondingly, getting a surface $S_j$ with boundary of desired
geometric combinatorial type (see the corresponding  entry in the
table in section~\ref{ss:tables}).

\subsection*{Boundary type +4.1b}

Consider  a  local  ribbon  graphs  of   type   $+4.1b$   and   a
corresponding collection  $\{2m_1,  \dots,  2m_n\}$  of orders of
interior       singularities       and        a        collection
$\{2a_1+1,2a_2+1,2a_3+1,2a_4+1\}$  of  orders  of  four  boundary
singularities.

The singularity  data  $\beta'_j$  of the corresponding component
$\cH(\beta'_j)$  of  the principal boundary stratum has the  form
$\beta'=\{m_1,    \dots,     m_n,    a_1+a_2+a_3+a_4+1\}$,    see
equations~\eqref{eq:Dji} and~\eqref{eq:H:boundary:stratum}.

Choose an  Abelian  differential $S'_j\in \cH(\beta'_j)$; let $P$
be a zero of $S'_j$ of degree $(a_1+a_2+a_3+a_4+1)$. We split $P$
into three zeroes $P_1,P,P_3$ such that the zero  $P_1$ of degree
$a_1$ is joined to the zero $P$ of degree $a_2+a_4+1$ by a saddle
connection with a holonomy vector $\vec{v}$ and the zero $P_3$ of
degree $a_3$ is also joined  to  $P$ by a saddle connection  with
the same  holonomy  vector  $\vec{v}$, see figure~\ref{fig:4:1b}.
Note that the new saddle  connections  are  oriented  differently
than  in  the  previous  case.  We  then cut  along  both  saddle
connections  and  detach  the  zero   $P$   into   two   boundary
singularities  $P_2,   P_4$   of  orders  $2a_2+1$  and  $2a_4+1$
correspondingly.

\begin{figure}[htb]
%
\includegraphics{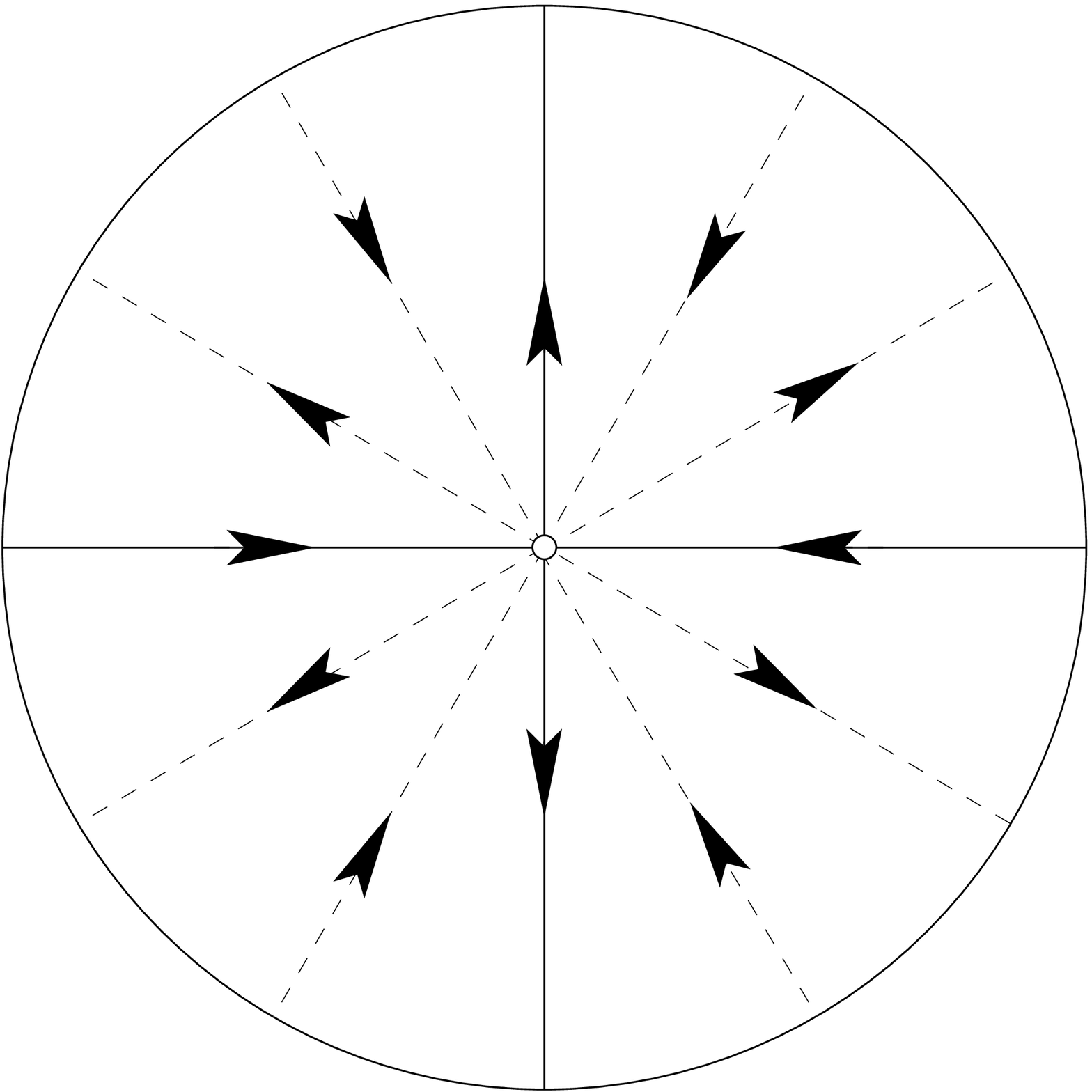}
\includegraphics{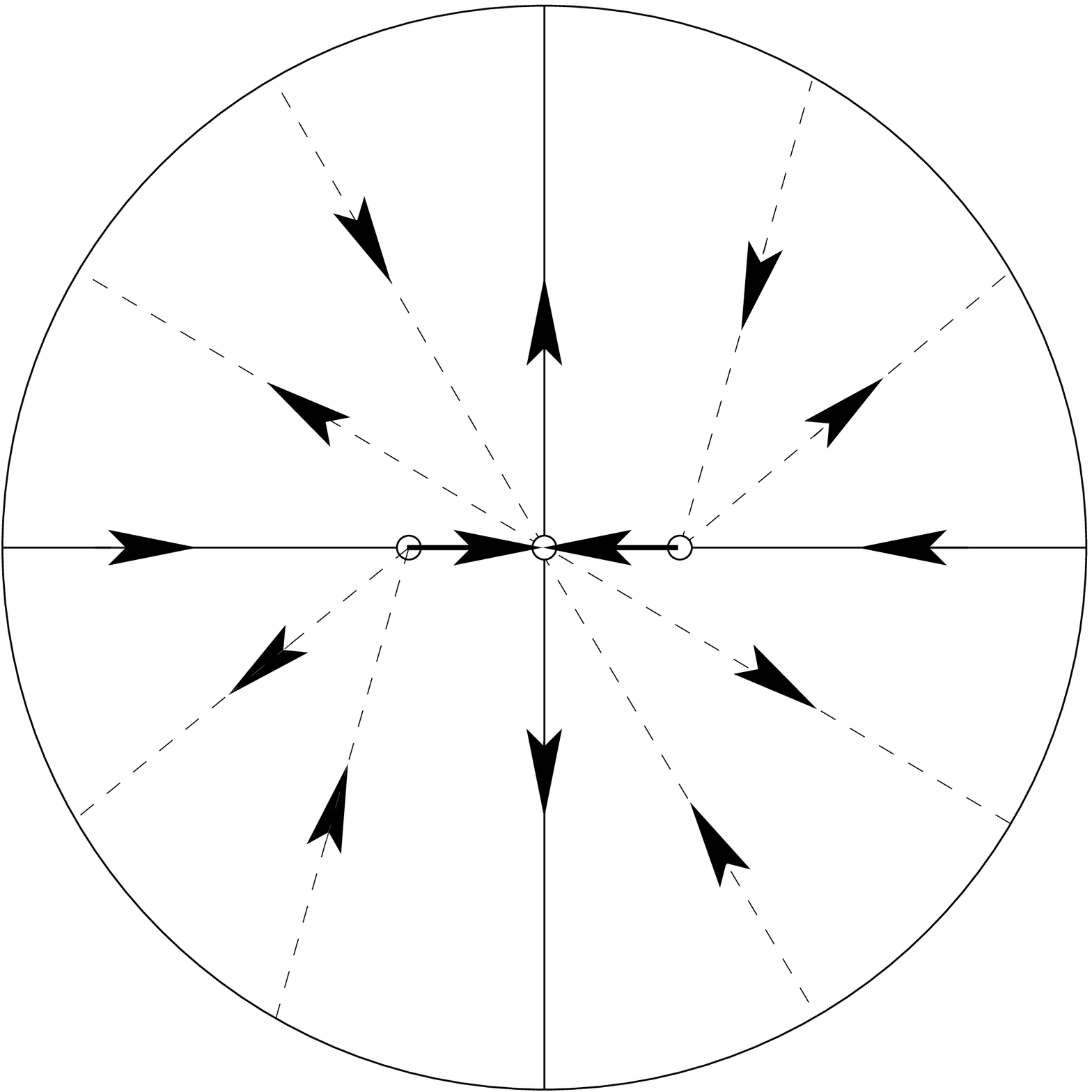}
\includegraphics{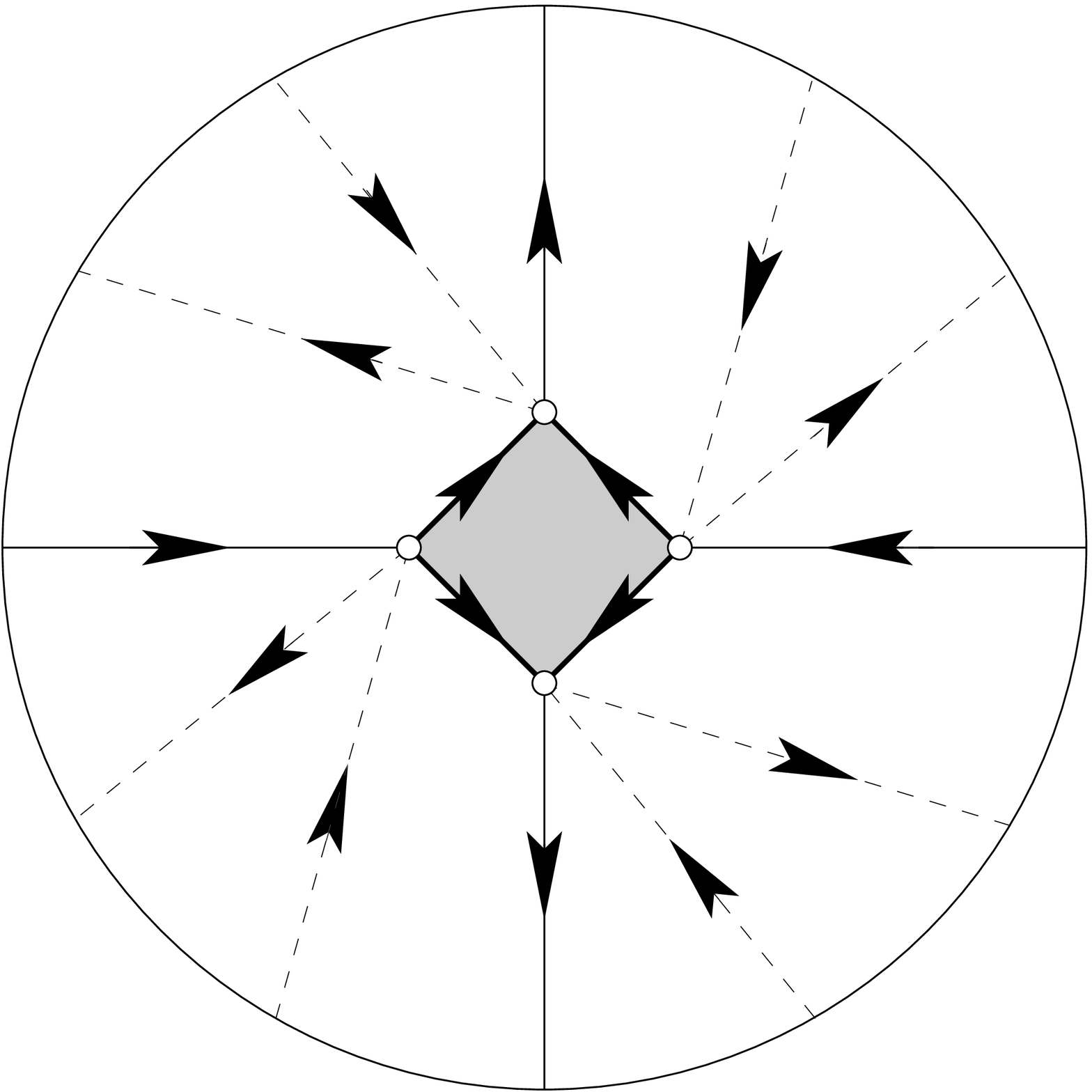}
%
%
\begin{picture}(0,0)(0,-150)
\put(10,-6)
 {\begin{picture}(0,0)(0,0)
 \put(-104,-199){$\scriptstyle \varepsilon$}
 \put(-48,-199){$\scriptstyle \varepsilon$}
 \put(-2,-198){$\scriptstyle \delta$} 
 \put(-20,-209){$\scriptstyle \delta$}
 \put(20,-199){$\scriptstyle \varepsilon$}
 \put(-50,-176){$\scriptstyle \varepsilon-\delta$}
 \put(-31,-162){$\scriptstyle \varepsilon+\delta$}
 \put(-9,-160){$\scriptstyle \varepsilon-\delta$}
 \put(17,-162){$\scriptstyle \varepsilon$}
 \put(32,-182){$\scriptstyle \varepsilon$}
 \put(-48,-217){$\scriptstyle \varepsilon$}
 \put(-30,-229){$\scriptstyle \varepsilon$}
 \put(-9,-242){$\scriptstyle \varepsilon-\delta$}
 \put(7,-229){$\scriptstyle \varepsilon+\delta$}
 \put(21,-218){$\scriptstyle \varepsilon-\delta$}
 \put(125,-209){$\scriptstyle P_3$}
 \put(89,-198){$\scriptstyle P_1$}
 \put(112,-185){$\scriptstyle P_4$}
 \put(102,-221){$\scriptstyle P_2$}
\end{picture}}
\end{picture}
\vspace{115bp} %
\caption{
\label{fig:4:1b}
Breaking up a zero into  three  ones and performing a surgery  we
get a surface of type ``$+4.1b$''}
\end{figure}

   %

By construction the resulting surface $S_j$ with boundary has the
desired boundary type ``$+4.1b$'' (see the corresponding entry in
the  table  in  section~\ref{ss:tables}),   and   collections  of
interior and boundary singularities of prescribed orders.

\subsection*{Boundary type +4.2a}

Consider a local  ribbon graphs of type $+4.2a$ and corresponding
collection  $\{2m_1,   \dots,   2m_n\}$  of  orders  of  interior
singularities  and  collections   $\{2a'+1,2a''+1\}$,   $\{2b'+1,
2b''+1\}$ of orders of  two  pairs of boundary singularities (see
figure~\ref{fig:embedded:local:ribbon:graphs}      for      their
parities). The  singularity  data $\beta'_j$ of the corresponding
component $\cH(\beta'_j)$  of  the principal boundary stratum has
the  form  $\beta'=\{m_1,  \dots,   m_n,   a'+a'',b'+b''\}$,  see
equations~\eqref{eq:Dji} and~\eqref{eq:H:boundary:stratum}.

Choose an Abelian differential $S'_j\in \cH(\beta'_j)$; let $P_1$
be  a  zero  of $S'_j$ of degree $a'+a''$; let $P_2$ be a zero of
degree $b'+b''$. As  in  the  case $+2.1$ we break  each  of  the
distinguished zeroes $P_1, P_2$ into a pair of  zeroes of degrees
$a',a''$  and  $b',b''$ correspondingly. We apply the surgery  in
such way that each of  the  two corresponding pairs of zeroes  is
joined by a saddle  connection  with a holonomy vector $\vec{v}$.
We  then cut  open  the modified flat  surface  along the  saddle
connections. As a result we get a surface $S_j$ with  boundary of
the desired  boundary  type  ``$+4.2a$''  (see  the corresponding
entry  in   the   table  in  section~\ref{ss:tables}),  and  with
collections of  interior and boundary singularities of prescribed
orders.

Proposition~\ref{pr:local:constructions} is proved.
\end{proof}

\section{Neighborhood of the principal boundary: nonlocal constructions}
\label{s:Nonlocal:constructions}

Recall that a direction $\pm v\in\R{2}\setminus\{0\}$  determines
a corresponding line field on the flat surface and a foliation in
direction $v$. The foliation is orientable if and only if $S$ has
trivial linear holonomy. Such an  auxiliary  direction  $v$ is an
element of all our  constructions.  We are creating surfaces with
boundary  from  closed flat surfaces; the direction  $v$  is  the
direction of parallel  geodesic  segments which form the boundary
components.

An  interior  singularity  $P$  of order $d$ has  $d+2$  adjacent
separatrix  rays  (or just  separatrices)  of  the  foliation  in
direction  $v$.  They divide a disc of  small  radius  $\epsilon$
centered at  $P$ into $d+2$ sectors,  each with cone  angle $\pi$
(see the  top  part  of figure~\ref{fig:breaking:up:a:zero} which
represents a singularity  of order $4$).  When $P$ is  a  regular
point (a  marked point) we  still have two such adjacent sectors,
each having cone angle $\pi$. When  $P$ is a simple pole, we have
a   single    separatrix    adjacent    to    $P$;   cutting   an
$\varepsilon$-neighborhood of $P$ by this  separatrix  we  get  a
single  sector.  When we speak about ``sectors''  adjacent  to  a
singularity  we always  mean  the sectors bounded  by  a pair  of
neighboring separatrices of the foliation in direction $v$.

When a flat surface has  trivial  linear  holonomy, the foliation
parallel to $v$ is oriented by the choice of direction $\vec{v}$.
The separatrix rays adjacent to any point $P$ inherit the natural
orientation: incoming and outgoing rays alternate with respect to
the natural  cyclic order on the  collection of rays  adjacent to
$P$.  The  sectors  adjacent  to  any  singularity  $P$  are also
naturally divided into two classes: the ones which are located to
the right of the  corresponding  oriented separatrix rays and the
ones which are located to the  left. We shall refer to them as to
the ``right'' and to the ``left'' sectors correspondingly.

In all nonlocal  constructions  we shall  use  a surgery along  a
smooth  path   without   self-intersections  joining  a  pair  of
singularities of  a  compact  flat  surface  (sometimes joining a
singularity  to  itself).  This  path $\rho$ (two paths  in  some
constructions)  will  be always chosen to be  transverse  to  the
direction  $v$  (and  hence,  transverse  to   the  foliation  in
direction  $v$);  in  particular,  $\rho$  never  passes  through
singularities. We shall often call such path a ``transversal''.

The following theorem from~\cite{Hubbard:Masur}  gives  us a  key
instrument for all  nonlocal constructions:

\begin{NNTheorem}[Hubbard---Masur]
Consider  a  closed  flat  surface  $S$  with  nontrivial  linear
holonomy, a pair  of points $P_1,P_2$  on $S$, a  direction  $\pm
v\in\R{2}\setminus\{0\}$  and  a  pair   of   sectors  $\Sigma_i$
adjacent to the corresponding points $P_i$, $i=1,2$.

For  any such  data  there exist a  transversal  $\rho$ with  the
endpoints at $P_1$ and $P_2$ which leaves $P_1$ in $\Sigma_1$ and
arrives at $P_2$  in $\Sigma_2$. The  case when $P_1$  and  $P_2$
coincide, or even  when $\Sigma_1$ and $\Sigma_2$ coincide is not
excluded.

If $S$  has trivial linear  holonomy the statement above is valid
under  additional  assumption  that  one  of  the  sectors  is  a
``right'' sector, and the other one is a ``left'' sector.
\end{NNTheorem}
  %

\subsection{Parallelogram construction}
\label{ss:Parallelogram construction}
In this  section  we  extend  the  ``parallelogram construction''
from~\cite{Eskin:Masur:Zorich} to  flat surfaces with  nontrivial
linear holonomy. For more details (including  restrictions on the
choice   of   parameter   $\delta$   in  terms  of   the   length
$4\varepsilon$ of  the  shortest  saddle  connection  on $S$, and
generalization   of   the   ``parallelogram   construction''   to
piecewise-transverse paths) we address the reader to the original
paper~\cite{Eskin:Masur:Zorich}    and    to   the    forthcoming
paper~\cite{Boissy:in:progress}.

Consider a transversal  $\rho$ as in  the theorem above.  In  the
construction  below,  we assume that if $S$  has  trivial  linear
holonomy,  then  $P_1\neq P_2$.  If  $S$  has  nontrivial  linear
holonomy, then we allow $P_1= P_2$ unless $P_1$  is a singularity
of order $-1$. If $P_1=P_2$ we allow $\Sigma_1=\Sigma_2$.

Fix the orientation of $\rho$ from $P_1$ to $P_2$. Since the path
is    smooth,    it    has   well-defined   tangent    directions
$\vec{u}_1=\dot\rho|_{P_1}$  and  $\vec{u}_2=\dot\rho|_{P_2}$  at
the  endpoints.  If the surface has trivial  linear  holonomy  we
assume  that  the frame  $\{\vec{u_i},\vec{v}\}$  represents  the
canonical  orientation  (upon interchanging,  if  necessary,  the
ordering of $P_1, P_2$).

For some  interior  point $P\in\rho$ let $\vec{u}=\dot\rho|_P$ be
the vector tangent  to $\rho$. Chose a vector $\vec{v}\in T_P(S)$
at  $P$  parallel to $v$ such that  the  frame  $\vec{u},\vec{v}$
represents  the  canonical orientation of the surface. Perform  a
parallel transport of $\vec{v}$ along $\rho$  to  all  points  of
$\rho$.

For a sufficiently small $\delta>0$ and any positive $s\le\delta$
we can construct  a {\it parallel  shift} $\rho_s$ of  $\rho$  in
direction  $\vec{v}$   at   the   distance   $s$.   Suppose  that
$\Sigma_1,\Sigma_2$  do  not coincide nor are adjacent. Then  for
any  $0\le   s_1<   s_2\le   \delta$   the  corresponding  shifts
$\rho_{s_1}$ and $\rho_{s_2}$ do not  intersect  and  do not have
self-intersections. If  $\Sigma_1$ and $\Sigma_2$ coincide or are
adjacent,  the  same  is  true  upon  an  appropriate  choice  of
orientation of $\rho$. .

\begin{figure}[ht]
%
\includegraphics{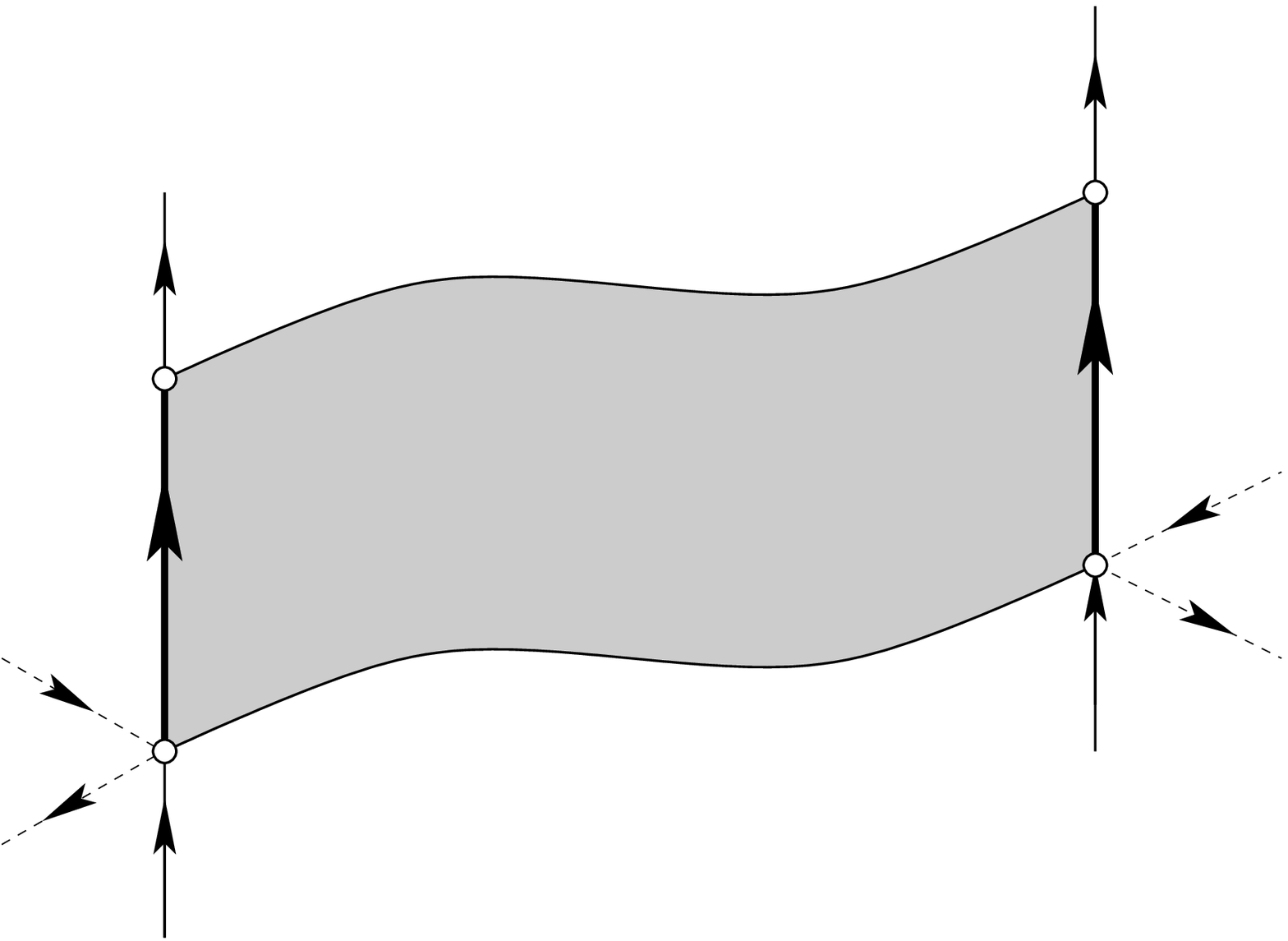}
\includegraphics{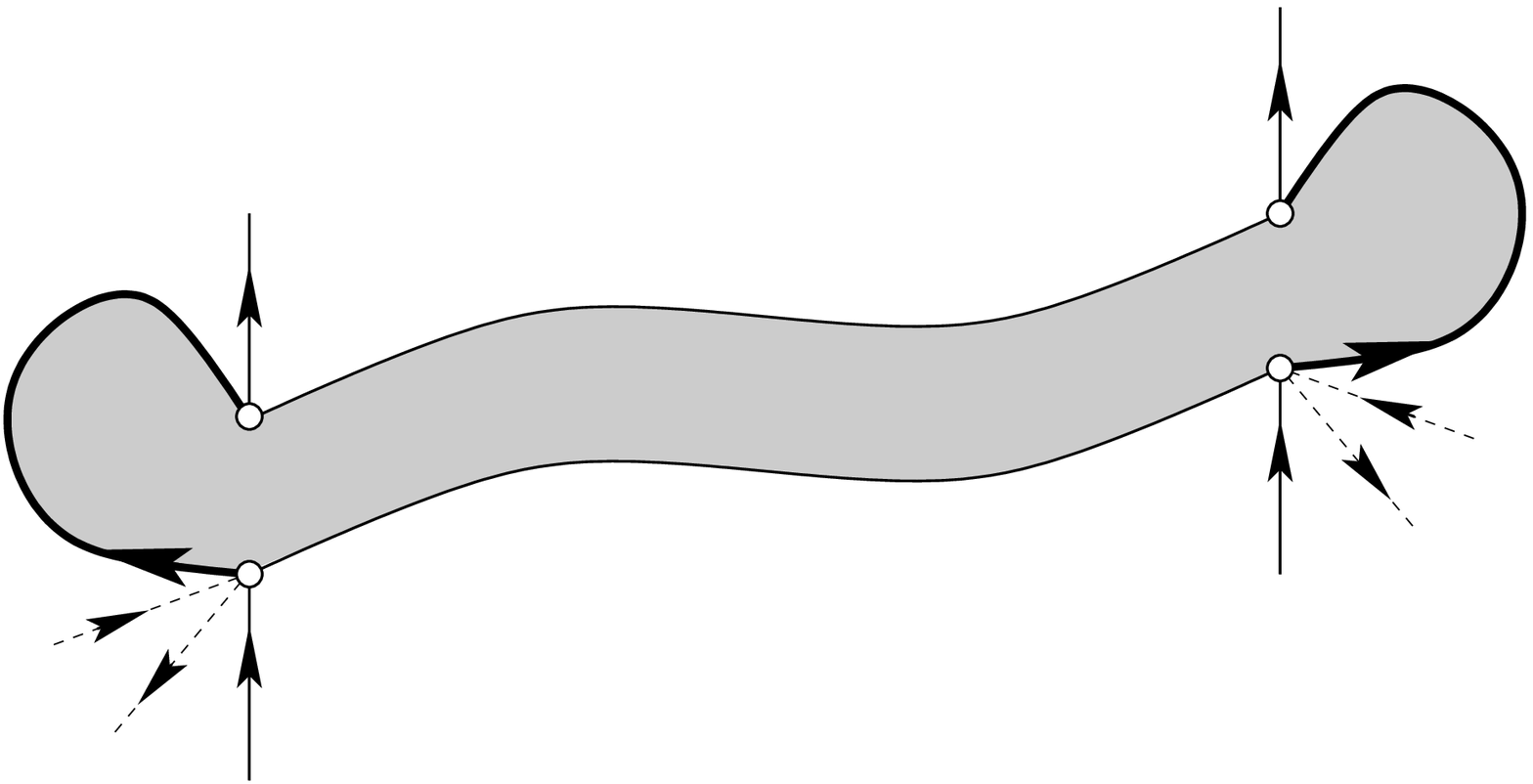}
\includegraphics{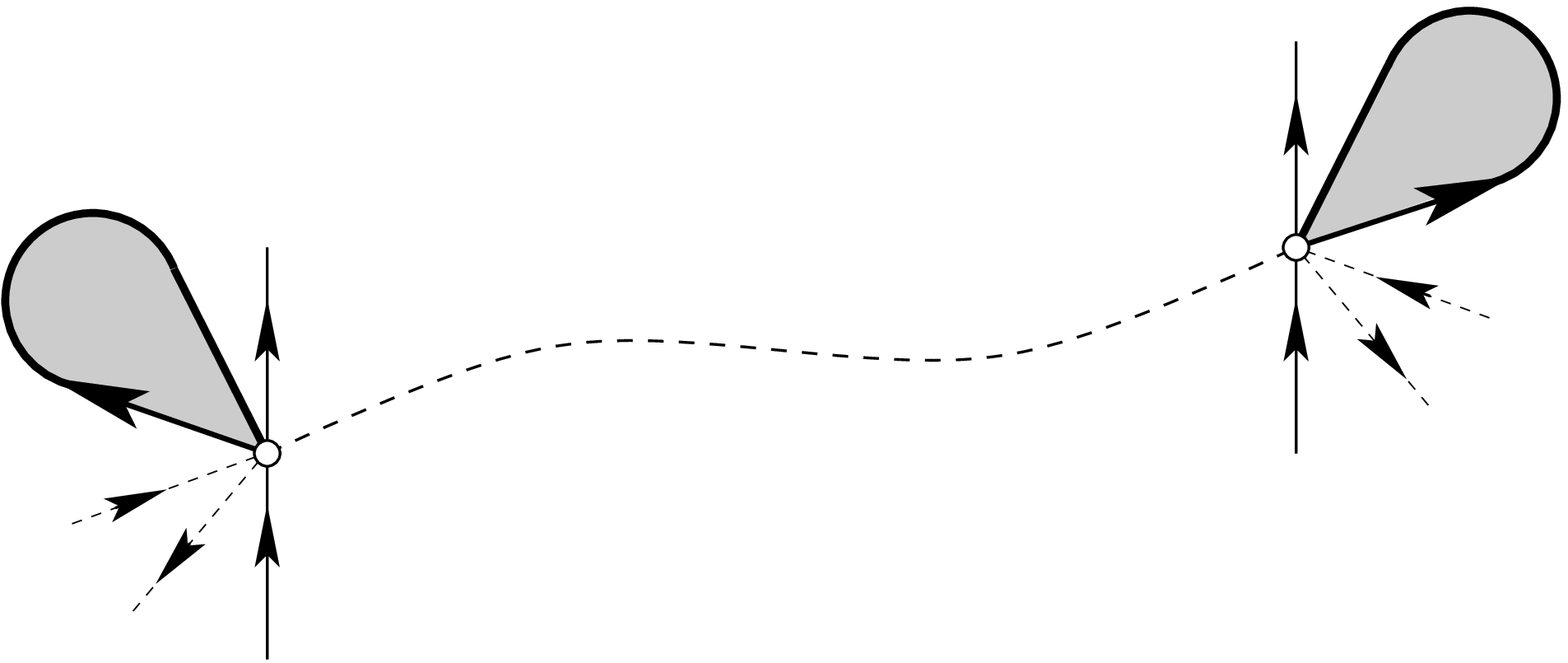}
%
%
\begin{picture}(0,0)(0,-150)
\put(10,-6)
 {\begin{picture}(0,0)(0,0)
 \put(-114,-152){$\scriptstyle P_2'$}
 \put(-114,-183){$\scriptstyle P_2$}
 \put(-173,-206){$\scriptstyle P_1$}
 \put(-173,-162){$\scriptstyle P_1'$}
 \put(-142,-159){$\scriptstyle \rho'$}
 \put(-142,-198){$\scriptstyle \rho$}
 \put(-142,-180){$\scriptstyle \Omega$}
 \put(11,-171){$\scriptstyle P_2'$}
 \put(12,-195){$\scriptstyle P_2$}
 \put(-49,-206){$\scriptstyle P_1$}
 \put(-49,-181){$\scriptstyle P_1'$}
 \put(-20,-177){$\scriptstyle \rho'$}
 \put(-20,-198){$\scriptstyle \rho$}
 \put(133,-182){$\scriptstyle P_2$}
 \put(74,-205){$\scriptstyle P_1$}
 \put(105,-198){$\scriptstyle \rho$}
\end{picture}}
\end{picture}
\vspace{70bp} %
\caption{
\label{fig:paral} ``Parallelogram construction''}
\end{figure}

Let $\gamma_i$  be a segment of  the separatrix ray  in direction
$\vec{v}$ at $P_i$ of length  $\delta$.  Even  when $P_1=P_2$ (in
particular,     when     $\Sigma_1=\Sigma_2$)    the     segments
$\gamma_1\neq\gamma_2$  are  well-defined. The  interior  of  the
domain  $\Omega$  bounded  by  $\rho,\rho',\gamma_1,\gamma_2$  is
homeomorphic  to  an  open  disc  and  can  be  thought of  as  a
``curvilinear parallelogram'', see figure~\ref{fig:paral}.

Remove $\Omega$ from  $S$ and identify  $\rho$ and $\rho'$  by  a
parallel translation. When $P_1\neq P_2$, as  a  result  of  this
surgery we get a surface with two boundary components each with a
single singular point, see figure~\ref{fig:paral}. When $P_1=P_2$
we detach the resulting  boundary  singularity into two getting a
surface  with  a single  boundary  component  with  two  boundary
singularities, see  figure~\ref{fig:m:2:1:nonlocal}. We refer  to
this surgery as to the {\em parallelogram construction.}

If  the  parallelogram  construction  is  applied  to  a  pair of
distinct points $P_1\neq  P_2$,  let $D_i$  be  the order of  the
corresponding singularity $P_i\in  S$,  $i=1,2$. In the case when
$P_1= P_2$, let  $\pi(a_1+1)$ be the angle between $\gamma_1$ and
$\gamma_2$  counted   in   the   positive   direction,   and  let
$\pi(a_2+1)$ be  the  angle  between  $\gamma_1$  and  $\gamma_2$
counted  in  the negative direction. By construction $a_i\ge  0$,
$i=1,2$.  The  order of  the  singularity  $P$  in  this  case is
$a_1+a_2$.

In this  notation, the orders  of the boundary singularities of a
surface obtained by a parallelogram  construction  are  equal  to
$\{D_1+2\}, \{D_2+2\}$, when $P_1\neq P_2$ and to $\{a_1,a_2+2\}$
when $P_1=P_2$. To see this, when $P_1\neq P_2$  it is sufficient
to  observe  figure~\ref{fig:paral}; in the remaining case it  is
sufficient to observe Figure~\ref{fig:m:2:1:nonlocal}.

\subsection{Nonlocal surgeries}
\label{ss:nonlocal:surgeries}

The remaining constructions are a combination of one of the local
constructions  described   in   the   previous   section  with  a
parallelogram construction. The parameters $\epsilon, \delta$ are
chosen as before.

\begin{proposition}
\label{pr:parallelogram:constructions}
Every surface of  any of boundary types $+2.2$, $+3.2a$, $+4.2b$,
$+3.2b$, $+4.2c$,  $+4.3a$, $+3.3$, $+4.3b$, $+4.4$ is realizable
by a combination of a local  construction  with  a  parallelogram
construction.
\end{proposition}
\begin{proof}
Applying the  same arguments as in the beginning  of the proof of
proposition~\ref{pr:local:constructions}   we   check  that   the
singularity        data         $\beta'$        defined        by
equation~\eqref{eq:H:boundary:stratum} from formal  combinatorial
data $\big(\G_{v_j}$, $\{d_1, \dots,  d_s\}$,  $\{k_{1,1}, \dots,
k_{r,p(r)}\big)$                      as                       in
proposition~\ref{pr:parallelogram:constructions}   represents   a
nonempty stratum  $\cH(\beta'_j)$.  Having  a closed flat surface
$S'_j\in\cH(\beta'_j)$ we now  need  to construct a surface $S_j$
with   boundary   realizing  the   initial   combinatorial   data
$\big(\G_{v_j}$,  $\{d_1,   \dots,  d_s\}$,  $\{k_{1,1},   \dots,
k_{r,p(r)}\big)$.

\subsection*{Boundary type +2.2}

We  begin  with boundary type $+2.2$. All interior  singularities
$\{2m_1, \dots, 2m_s\}$ have positive  even  orders;  each of the
two boundary  components  contains a single boundary singularity.
The boundary singularities  also  have positive even orders $2m',
2m''$   (see   figure~\ref{fig:embedded:local:ribbon:graphs}  and
condition   (4)  of   definition~\ref{def:configuration}   of   a
configuration),   so  in   this   case   $\beta'=\{m_1,\dots,m_s,
m'-1,m''-1\}$.

Choosing   an   Abelian  differential   $S'_j\in\cH(\beta')$  and
performing  the  parallelogram  construction  at  the  zeroes  of
degrees $m'-1,m''-1$  (see  figure~\ref{fig:paral}) we get a flat
surface  $S_j$  with   boundary   of  type  ``$+2.2$''  (see  the
corresponding  entry  in  the table in  section~\ref{ss:tables}),
having collections of  interior  and of boundary singularities of
prescribed orders.

\subsection*{Boundary types +3.2a and +4.2b}

Boundary type  $+3.2a$ can be  considered as a particular case of
$+4.2b$ when one of the boundary singularities has order $0$ (see
the appropriate entries in the table in section~\ref{ss:tables}).

Consider  a  ribbon graph of type $+4.2b$.  Let  $\{2m_1,  \dots,
2m_s\}$ be  a  collection  of  orders  of interior singularities.
According  to  figure~\ref{fig:embedded:local:ribbon:graphs}  the
orders of all boundary  singularities  are even for boundary type
$+4.2b$;  denote   by   $2a_1,2a_2+2$   the  orders  of  boundary
singularities corresponding  to  the first boundary component and
by  $2a_3,   2a_4+2$   the   orders   of  boundary  singularities
corresponding  to  the second  component.  By  condition  (4)  of
definition~\ref{def:configuration} of a configuration the numbers
$a_i$  are   nonnegative   integers  for  $i=1,\dots,4$.  We  get
$\beta'=\{m_1, \dots, m_s, a_1+a_2, a_3+a_4\}$.

Choose  a  flat surface $S'_j\in\cH(\beta'_j)$. Choose a pair  of
separatrices $\gamma_1,\gamma_2$ in  direction $\vec{v}$ adjacent
to the first zero. Choose $\gamma_1$ to be an outgoing separatrix
and $\gamma_2$ to  be incoming separatrix  in such way  that  the
angle from the separatrix ray  $\gamma_1$  to  the separatrix ray
$\gamma_2$ in  the  clockwise  direction  is  $(2a_1+1)\pi$.  Let
$\Sigma_1$ be the sector adjacent to $\gamma_1$ counterclockwise;
let $\Sigma_2$ be the sector adjacent to $\gamma_2$ clockwise.

Similarly, choose  a  pair of separatrices $\gamma_3,\gamma_3$ in
direction $\vec{v}$ adjacent to  the  zero of degree $a_3+a_4$ in
such way that $\gamma_3$ is outgoing, $\gamma_4$ is incoming; the
counterclockwise   angle   from   $\gamma_3$  to  $\gamma_4$   is
$(2a_3+1)\pi$.  Let   $\Sigma_3$   be   the  sector  adjacent  to
$\gamma_3$  clockwise;  let $\Sigma_4$ be the sector adjacent  to
$\gamma_4$ counterclockwise.

Join $\Sigma_3$ to $\Sigma_1$ by  a  transversal  $\rho_1$;  join
$\Sigma_2$ to $\Sigma_4$  by  a transversal $\rho_2$. If $\rho_1$
intersects $\rho_2$ we  can  resolve the intersections to achieve
nonintersecting transversals.

Suppose that  in resolving the  intersections (if any) we did not
change the correspondence between the sectors  and $\rho_1$ still
joins $\Sigma_3$ to  $\Sigma_1$  and $\rho_2$ joins $\Sigma_2$ to
$\Sigma_4$.  Choosing   some   small   $\delta$   we   can  apply
parallelogram construction  to  the  transversal $\rho_1$ and the
direction  $\vec{v}$  and  to  the transversal $\rho_2$  and  the
direction $-\vec{v}$, see Figure~\ref{fig:42b}.

In the remaining  case after the resolution of intersections, the
correspondence between  the  sectors  changed and the transversal
$\rho_1$  joins   sector   $\Sigma_4$  to  $\Sigma_1$  while  the
transversal $\rho_2$ joins $\Sigma_2$ to $\Sigma_3$. In this case
we deform  the transversals slightly  in such way that they still
do not intersect and $\rho_2$  lands  on the ray $\gamma_3$ at  a
distance $\delta$ from the  zero  (in the same sector $\Sigma_3$)
and  $\rho_1$ starts at  a  point  on  the ray  $\gamma_4$  at  a
distance $\delta$ from  the zero (in the same sector $\Sigma_4$).
We can  construct  two  ``curvilinear parallelograms'' $\Omega_1,
\Omega_2$ (see  figure~\ref{fig:42b})  which do not intersect, so
we can proceed as above.

\begin{figure}[ht]
%
\includegraphics{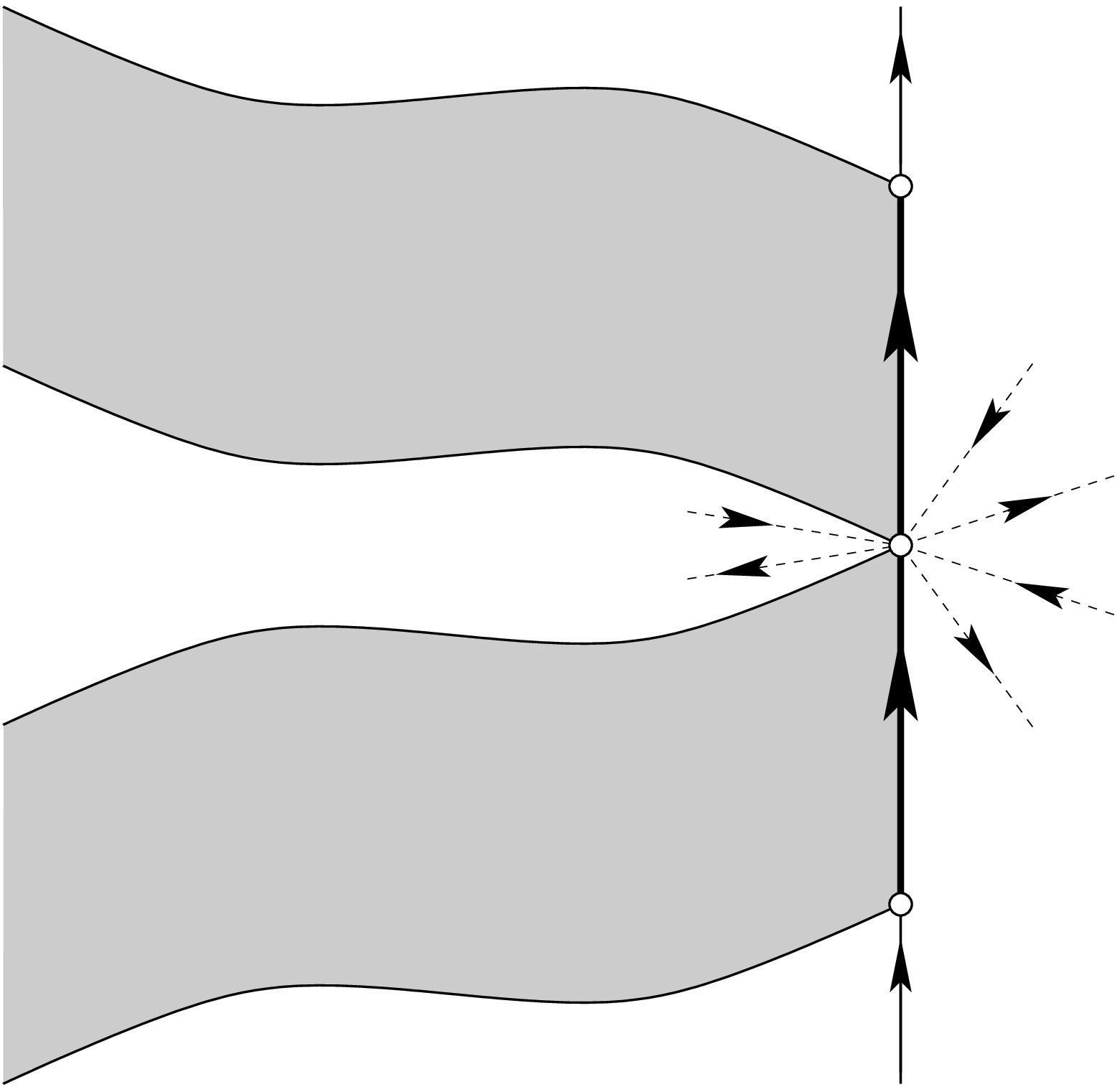}
\includegraphics{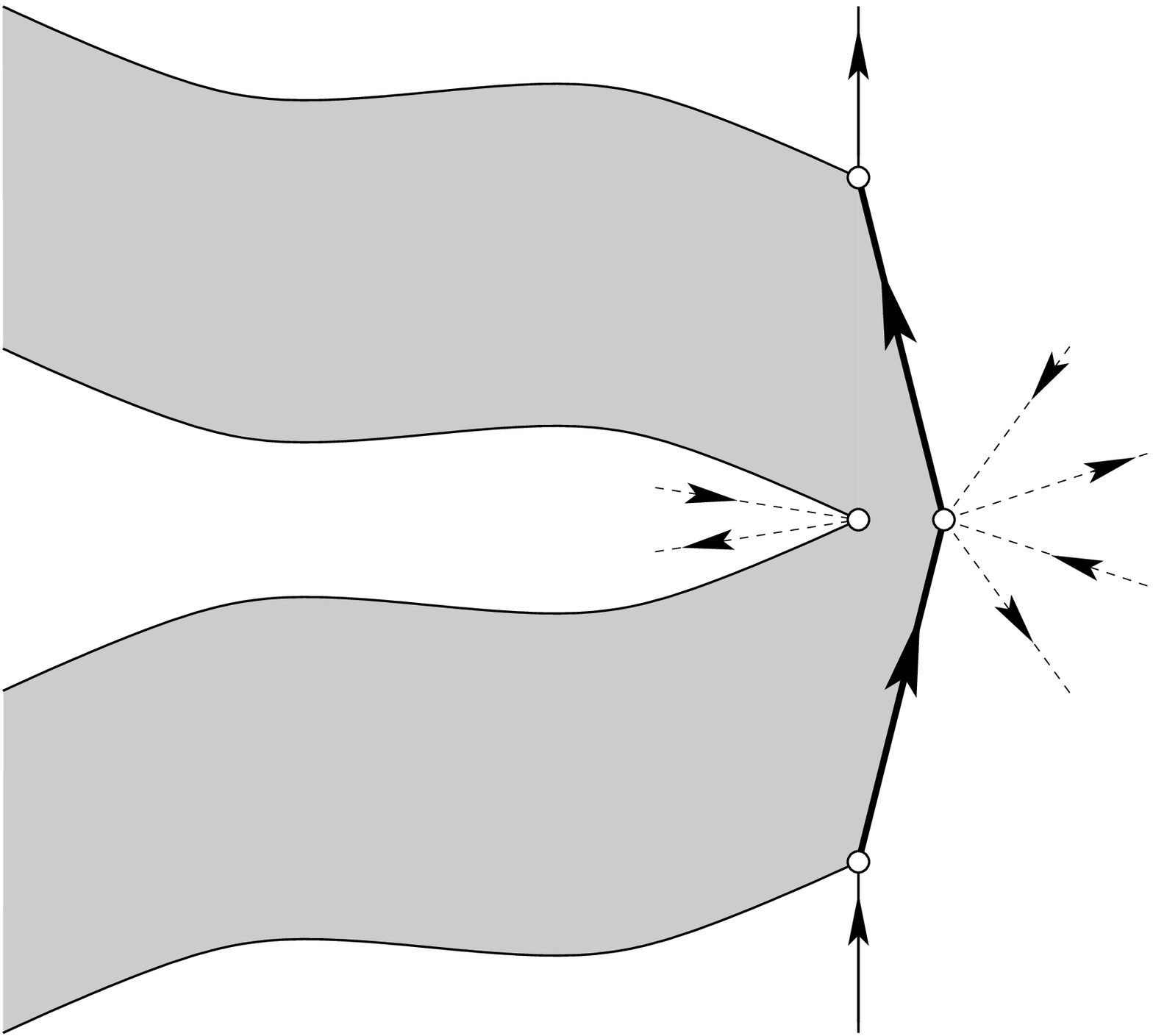}
\includegraphics{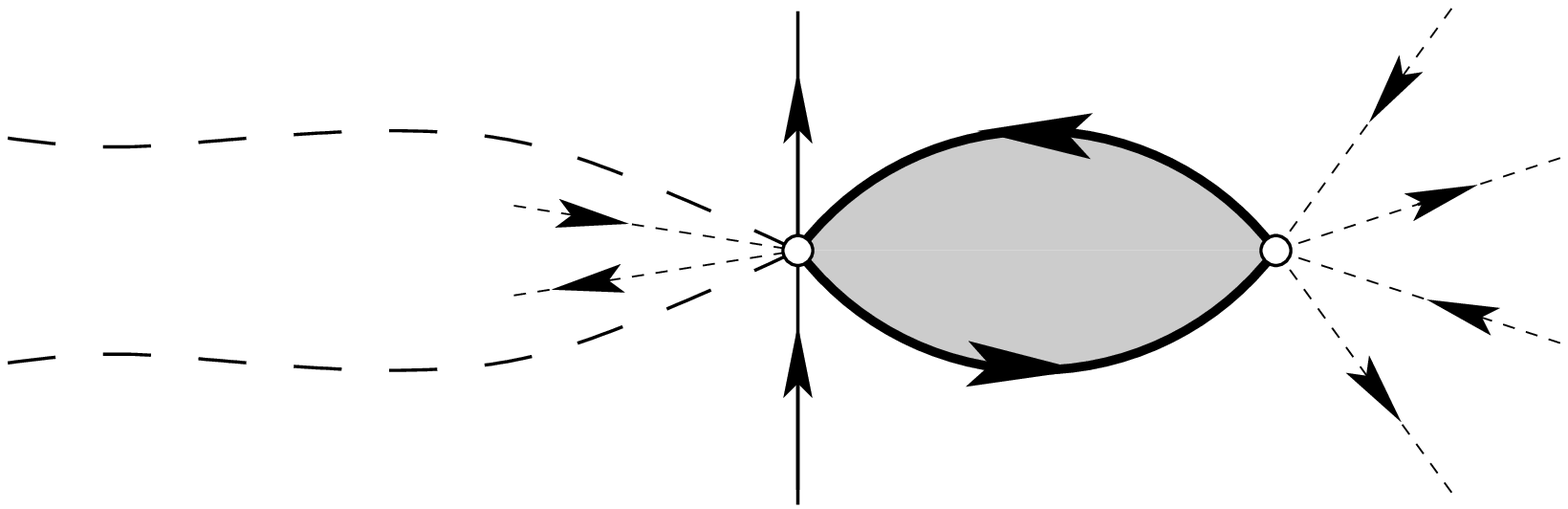}
%
%
\begin{picture}(0,0)(0,-165) 
\put(10,-6)
 {\begin{picture}(0,0)(0,0)
 \put(-114,-186){$\scriptstyle P'$}
 \put(-114,-200){$\scriptstyle \gamma_1$}
 \put(-126,-210){$\scriptstyle \Sigma_1$}
 \put(-114,-218){$\scriptstyle P$}
 \put(-126,-227){$\scriptstyle \Sigma_2$}
 \put(-114,-231){$\scriptstyle \gamma_2$}
 \put(-114,-245){$\scriptstyle P''$}
 \put(-165,-187){$\scriptstyle \rho_1'$}
 \put(-165,-205){$\scriptstyle \rho_1$}
 \put(-165,-228){$\scriptstyle \rho_2$}
 \put(-165,-247){$\scriptstyle \rho_2'$}
 \put(7,-186){$\scriptstyle P'$}
 \put(1,-223){$\scriptstyle P_1$}
 \put(12,-213){$\scriptstyle P_2$}
 \put(7,-245){$\scriptstyle P''$}
 \put(-45,-187){$\scriptstyle \rho_1'$}
 \put(-45,-205){$\scriptstyle \rho_1$}
 \put(-45,-228){$\scriptstyle \rho_2$}
 \put(-45,-247){$\scriptstyle \rho_2'$}
 \put(110,-218){$\scriptstyle P_1$}
 \put(130,-218){$\scriptstyle P_2$}
 \put(60,-204){$\scriptstyle \rho_1$}
 \put(60,-229){$\scriptstyle \rho_2$}
\end{picture}}
\end{picture}
\vspace{115bp} 
\caption{
\label{fig:42b} A pair of simultaneous parallelogram constructions.}
\end{figure}

Detaching  each  of the resulting singularities into pairs  $P_1,
P_2$ and $P_3, P_4$ (see figure~\ref{fig:42b}) we get the desired
surface  $S_j$  with  boundary  of  type   ``$+4.2b$''  (see  the
corresponding entry in the table in section~\ref{ss:tables}), and
with  prescribe   collections   of   interior   and  of  boundary
singularities.

\begin{Remark}
Recall that if we  identify the opposite sides of each hole  of a
surface constructed above we obtain a closed surface  with a pair
of even order zeroes simultaneously broken up into a pair  of odd
order zeroes.
\end{Remark}

\subsection*{Boundary types +3.2b and +4.2c}
Boundary type  $+3.2b$ can be  considered as a particular case of
boundary  type  $+4.2c$. To see this compare  the  surfaces  with
boundary representing  the  corresponding  ribbon graphs (see the
appropriate  entries  in  the table in  section~\ref{ss:tables}).
Marking a  point in the  middle of the saddle connection labelled
by ``$+2$'' on the boundary of the surface of type $+3.2b$ we get
a  surface  of  boundary  type  $+4.2c$  with  the  corresponding
boundary singularity of order $0$.

Consider  a  ribbon graph of type $+4.2c$.  Let  $\{2m_1,  \dots,
2m_s\}$ be a  collection of orders of interior singularities. Let
$2a_1+1,2a_2,2a_3+1$ be  the orders of the boundary singularities
on the boundary component composed from three saddle connections;
let $2a_4+2$  be the order  of the single boundary singularity on
the       complementary       boundary       component,       see
figure~\ref{fig:embedded:local:ribbon:graphs}.  By  condition (4)
of  definition~\ref{def:configuration}  of  a  configuration  the
numbers $a_i$ are nonnegative integers for  $i=1,\dots,4$. We get
$\beta'=\{m_1, \dots, m_s, a_1+a_2+a_3, a_4\}$.

Choose a flat surface $S'_j\in\cH(\beta'_j)$. Let $P$ be the zero
of degree  $a_1+a_2+a_3$, and $P_4$  be the zero of degree $a_4$.
Choose a separatrix $\gamma_1$ in direction $\vec{v}$ adjacent to
$P_4$ and a separatrix $\gamma_2$ in direction $\vec{v}$ adjacent
to  $P$.  Let  $\Sigma_1$  be  the  ``right'' sector adjacent  to
$\gamma_1$; let $\Sigma_2$ be the  ``left''  sector  adjacent  to
$\gamma_2$,   see  figure~\ref{fig:42c}.   Join   $\Sigma_1$   to
$\Sigma_2$ by a separatrix $\rho$. Choose  $\delta$ small enough,
so that the intersection of $\rho$  with an $\delta$-neighborhood
of  $P$  contains  a  single  connected  component  contained  in
$\Sigma_2$.

\begin{figure}[ht]
%
\includegraphics{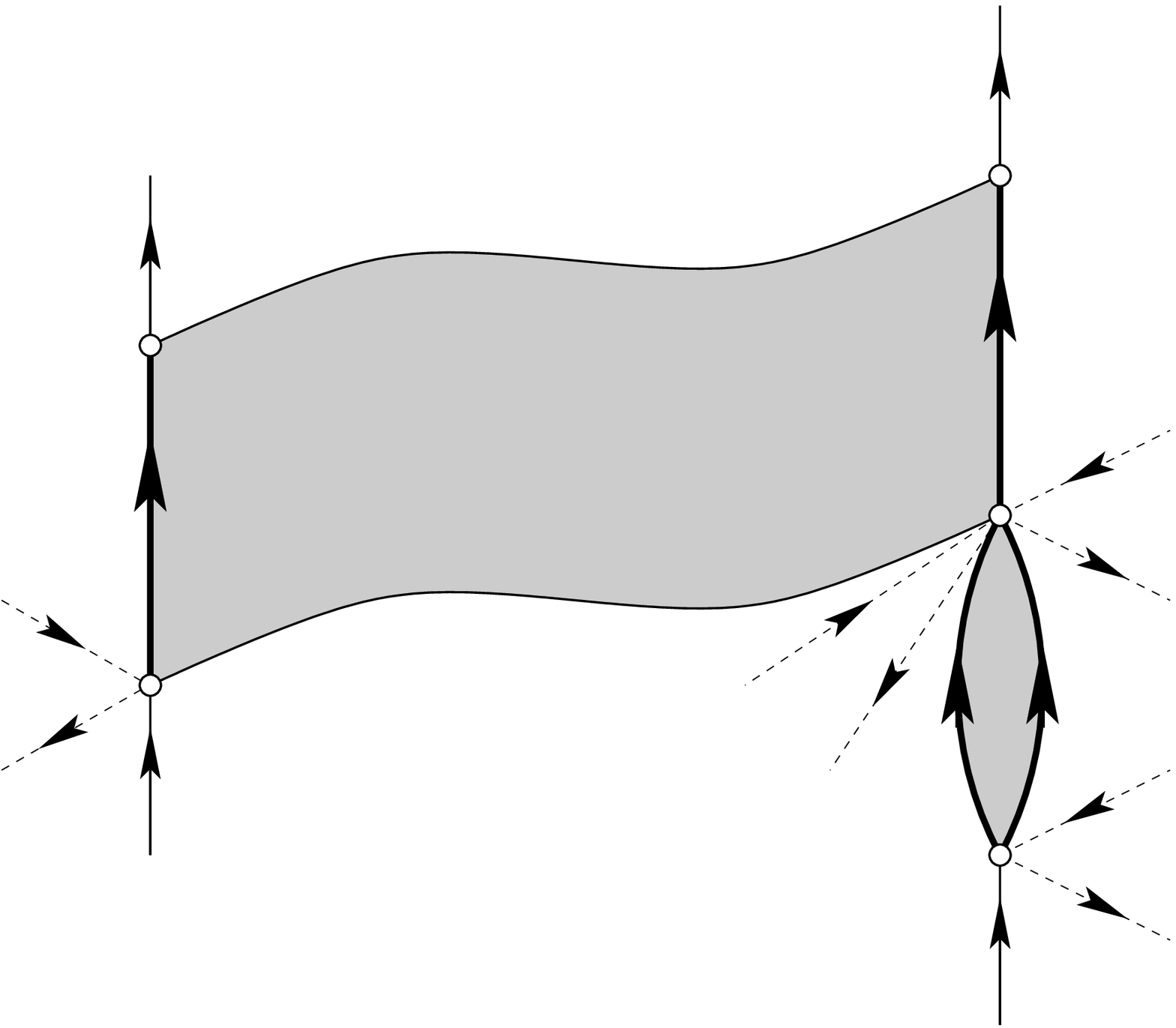}
\includegraphics{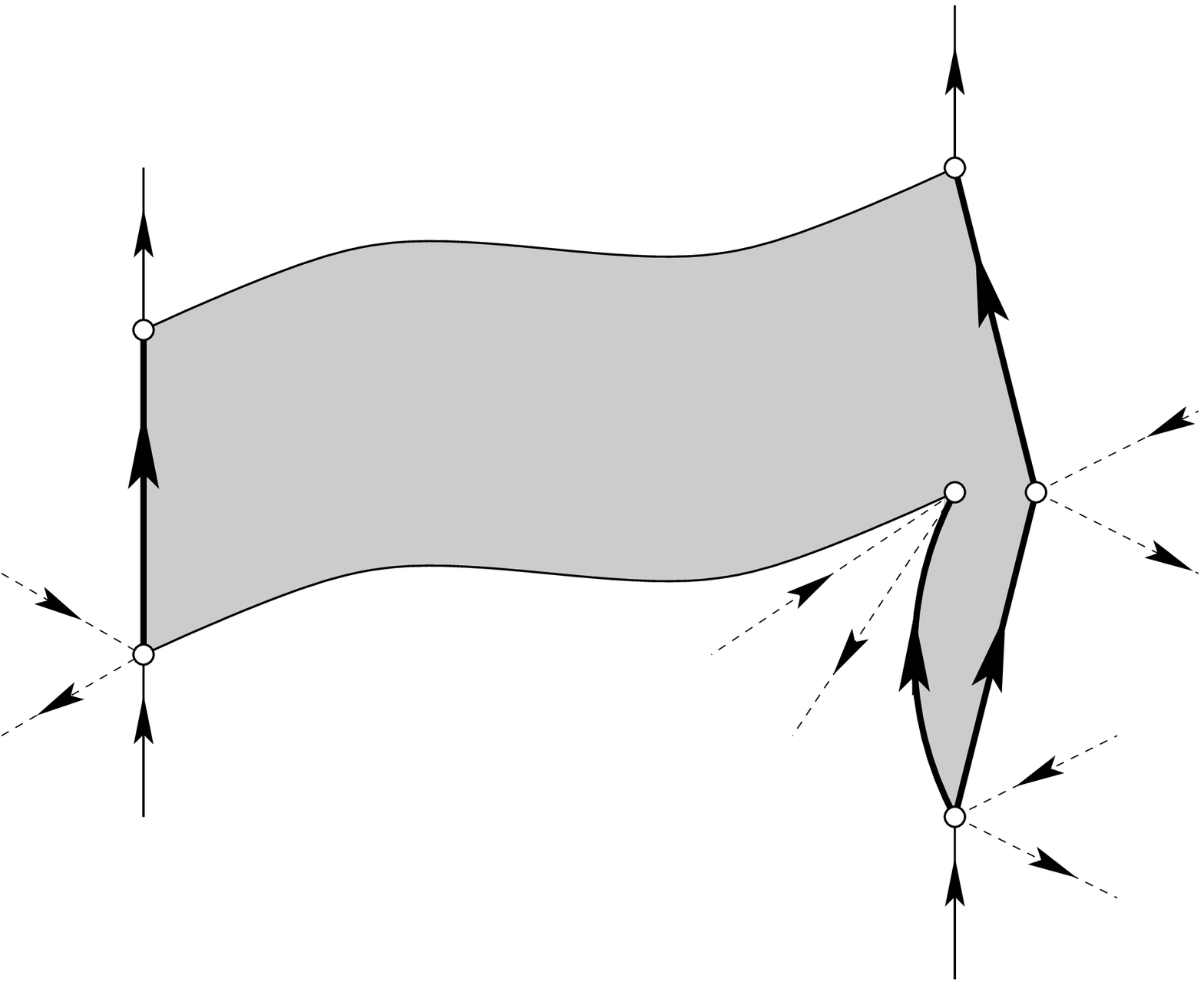}
\includegraphics{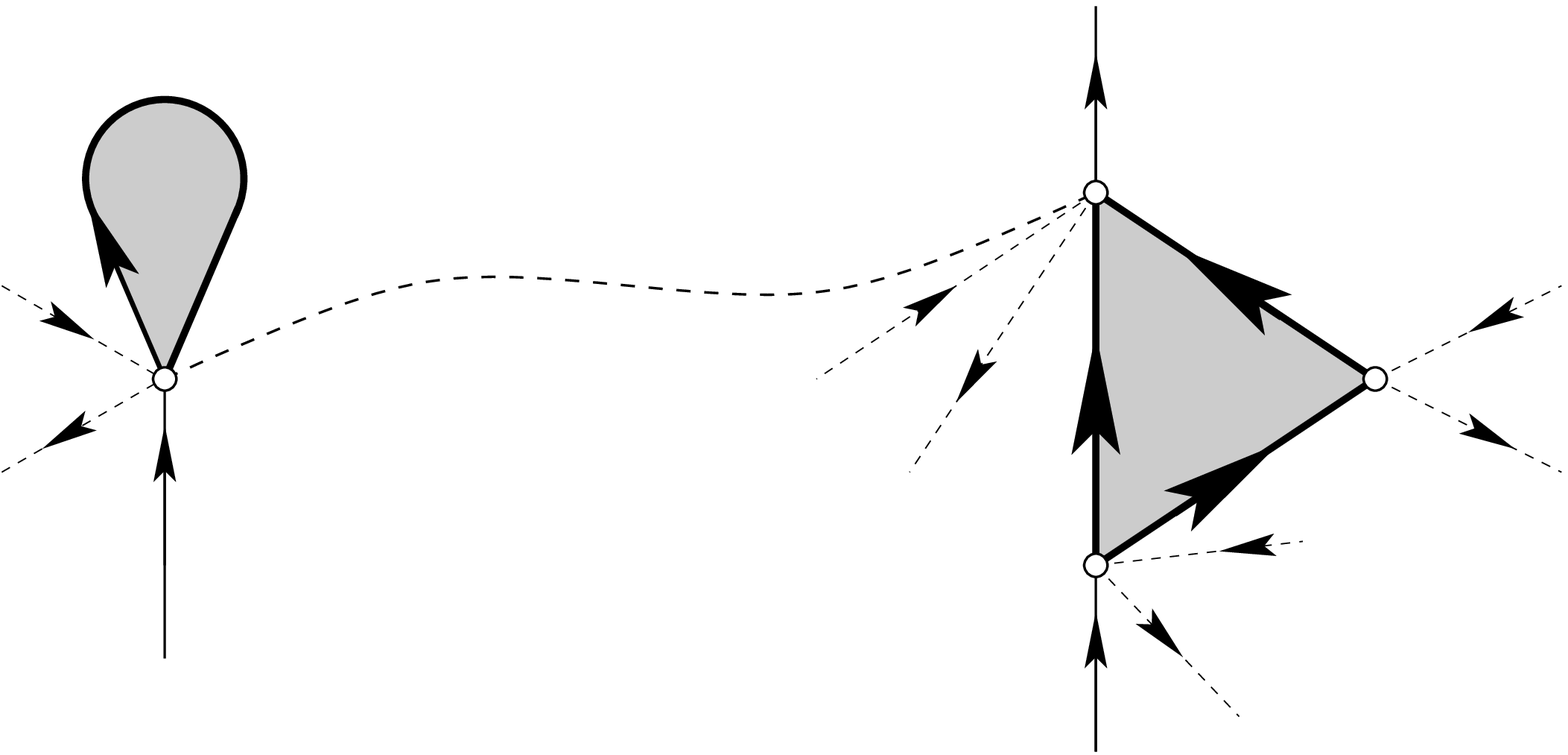}
%
%
\begin{picture}(0,0)(0,-165) 
\put(10,-6)
 {\begin{picture}(0,0)(0,0)
 \put(-173,-236){$\scriptstyle P_4$}
 \put(-113,-187){$\scriptstyle P'$}
 \put(-112,-216){$\scriptstyle P$}
 \put(-113,-249){$\scriptstyle P_1$}
 \put(-152,-190){$\scriptstyle \rho'$}
 \put(-152,-228){$\scriptstyle \rho$}
 \put(-53,-236){$\scriptstyle P_4$}
 \put(7,-187){$\scriptstyle P'$}
 \put(8,-216){$\scriptstyle P''$}
 \put(26,-216){$\scriptstyle P_2$}
 \put(7,-249){$\scriptstyle P_1$}
 \put(-32,-190){$\scriptstyle \rho'$}
 \put(-32,-228){$\scriptstyle \rho$}
 \put(67,-236){$\scriptstyle P_4$}
 \put(127,-214){$\scriptstyle P_3$}
 \put(157,-227){$\scriptstyle P_2$}
 \put(127,-249){$\scriptstyle P_1$}
 \put(88,-228){$\scriptstyle \rho$}
\end{picture}}
\end{picture}
\vspace{115bp} 

\caption{
\label{fig:42c}
Parallelogram construction  combined  with  breaking  up  a  zero
creates a surface of type $+4.2c$}
\end{figure}

Choose a  separatrix $\gamma_3$ at $P$  such that the  angle from
$\gamma_2$  to  $\gamma_3$ (in  the  counterclockwise  direction)
equals  $\pi(2a_3+1)$.  Break the zero $P$ along $\gamma_3$  into
two  zeroes  $P$   and  $P_1$  of  degrees  $a_2+a_3$  and  $a_1$
correspondingly  joined  by  a  saddle  connection  in  direction
$\vec{v}$  of  length  $\delta$  and  perform  the  parallelogram
construction along  $\rho$  (strictly  speaking  to a transversal
naturally  corresponding to  $\rho$),  see  figure~\ref{fig:42c}.
Detaching  $P$  into two  points  we  obtain  a  surface  of type
``$+4.2c$'' with the desired singularity data.
   %

\subsection*{Boundary type +4.3a}
Let  $\{2a_1+2\}$,  $\{2a_2+2\}$,  $\{2a_3+1,  2a_4+1\}$  be  the
orders of  the  boundary singularities naturally distributed into
the       corresponding       boundary      components,       see
figure~\ref{fig:embedded:local:ribbon:graphs}.  By  condition (4)
of  definition~\ref{def:configuration}  of  a  configuration  the
numbers $a_i$ are nonnegative integers for  $i=1,\dots,4$. We get
$\beta'=\{m_1, \dots, m_s, a_1, a_2, a_3+a_4\}$.

Surfaces with boundary  of this type  are obtained by  a  trivial
combination of a  parallelogram construction applied to a pair of
distinct zeroes  of degrees $a_1$ and $a_2$ and  by breaking up a
zero of  degree $a_3+a_4$ into two  zeroes of degrees  $a_3, a_4$
with a subsequent slit along the resulting saddle connection.

\subsection*{Boundary types +3.3 and +4.3b}
The boundary  type $+3.3$ can  be considered as a particular case
of  the  boundary type  $+4.3b$  when  the  appropriate  boundary
singularity has order $0$.

Surfaces of type $+4.3b$  can  be constructed in complete analogy
with surfaces of type $+4.2b$ (see figure~\ref{fig:42b}) with the
only  difference  that  now  we  choose  sectors  $\Sigma_3$  and
$\Sigma_4$ at two distinct points.

\subsection*{Boundary type +4.4}
In this  case the orders  of boundary singularities have the form
$\{2a_1+2, 2a_2+2,2a_3+2,2a_4+2\}$; according to condition (4) of
definition~\ref{def:configuration} of a configuration all numbers
$a_i$  are  nonnegative  integers.  Thus,  in  this  case  we get
$\beta'=\{m_1, \dots, m_s,  a_1, a_2, a_3, a_4\}$. To construct a
desired surface with boundary of type $+4.4$ it  is sufficient to
apply a pair of independent parallelogram constructions.

Proposition~\ref{pr:parallelogram:constructions} is proved.
\end{proof}

\subsection{Surfaces with boundary of ``--'' type}
\label{ss:minus:type}

To           complete            the           proof           of
proposition~\ref{pr:realizability:of:all:vertices} it remains  to
construct surfaces with boundary realizing any combinatorial data
$\big(\G_{v_j}$,  $\{d_1,   \dots,  d_s\}$,  $\{k_{1,1},   \dots,
k_{r,p(r)}\big)$      satisfying      conditions     2--6      of
definition~\ref{def:configuration}   for   local  ribbon   graphs
$\G_{v_j}$ of ``$-$''-types.

\begin{proposition}
\label{pr:minus:2:2}
Combinatorial data  representing  boundary types $-2.2, -1.1$ and
$-2.1$ are  realizable  by  appropriate surfaces with boundaries.
\end{proposition}
(See  initial  proposition~\ref{pr:realizability:of:all:vertices}
for the detailed formulation.)
\begin{proof}
The component of the principle boundary  stratum corresponding to
a vertex $v_j$ of ``$-$''-type  has  type  $\cQ(\alpha'_j)$.  The
singularity      data      $\alpha'_j$      is      given      by
equation~\eqref{eq:Q:boundary:stratum}, namely
$$
\alpha'_j=\{d_1, \dots, d_{s(j)},\ D_1, \dots,  D_{r(j)}\},
$$
where  $d_1,  \dots,   d_{s(j)}$   are  the  orders  of  interior
singularities, and $D_1, \dots, D_{r(j)}$ are  expressed in terms
of     the     orders    of     boundary     singularities     by
formula~\eqref{eq:Dji}.     Conditions     4     and     5     in
definition~\ref{def:configuration} of a  configuration  guarantee
that  all  the  entries  of  $\alpha'_j$   are   from   the   set
$\{-1,0,1,2,\dots  \}$,  that  the  total sum of the  entries  of
$\alpha'_j$ is divisible by $4$ and that this sum is greater than
or    equal     to     $-4$.    Moreover,    condition    6    in
definition~\eqref{def:configuration}  implies  that   $\alpha'_j$
neither   belongs    to    the    exceptional   list   given   by
equation~\eqref{eq:empty:strata} below, nor can be obtained  from
an entry of this list by adding additional  elements ``$0$'' (see
lemma~\ref{lm:nonrealizable:data:combinatrorics}   in   the  next
section).     According     to     the     results     of     the
paper~\cite{Masur:Smillie:realizability}  this  implies that  the
stratum $\cQ(\alpha'_j)$ is a nonempty.

Consider any flat  surface $S'_j$ in $\cQ(\alpha'_j)$. We use the
same conventions on  parameters $\delta, \epsilon$, and $v$ as in
the proof of  proposition~\ref{pr:local:constructions}.  Applying
an appropriate surgery to the closed surface $S'_j$  we are going
to construct a  surface $S_j$ with boundary realizing the initial
combinatorial  data   $\big(\G_{v_j}$,  $\{d_1,  \dots,   d_s\}$,
$\{k_{1,1}, \dots, k_{r,p(r)}\big)$.

\subsection*{Boundary type --2.2}
Boundary type $-2.2$ is constructed in complete analogy to $+2.2$
by  a  parallelogram  construction.  Each  of  the  two  boundary
components contains  a  single boundary singularity. The boundary
singularities have  strictly  positive  orders $k_{1,1}, k_{2,1}$
(see    inequality    on    $D_i$    in    condition    (4)    of
definition~\ref{def:configuration}  of a  configuration),  so  in
this case $\alpha'=\{d_1,\dots,d_s, k_{1,1}-2, k_{2,1}-2\}$.

Choosing   a   quadratic  differential   $S'\in\cQ(\alpha')$  and
performing the parallelogram construction at the zeroes of orders
$k_{1,1}-2, k_{2,1}-2$ (see figure~\ref{fig:paral}) we get a flat
surface  $S_j$   with   boundary   of   type  ``$-2.2$'',  having
collections  of   interior   and  of  boundary  singularities  of
prescribed orders.

\subsection*{Boundary types --1.1 and --2.1}
Note  next  that boundary  type  $-1.1$ can  be  considered as  a
particular  case of  boundary  type $-2.1$ when  one  of the  two
boundary singularities has order $0$.

Consider a ribbon graph of type $-2.1$. Let $\{d_1, \dots, d_s\}$
be  the   orders   of  interior  singularities,  let  $\{k_{1,1},
k_{1,2}\}$ be the orders of boundary  singularities. By condition
(4) of  definition~\ref{def:configuration}  we  have $D_1\ge -1$,
where  $D_1=k_{1,1}+k_{1,2}-2$, which  implies  that  nonnegative
integers $k_{1,1},  k_{1,2}$  cannot  be  simultaneously equal to
zero.  Thus,  we  may  assume   that   $k_{1,1}\ge   1$.  We  get
$\alpha'_j=\{d_1,\dots,d_s, k_{1,1}+k_{1,2}-2\}$.

Consider a flat surface $S'_j$  in  $\cQ(\alpha'_j)$.  When  both
$k_{1,1},k_{1,2}$  are odd  we  can break up  the  zero of  order
$k_{1,1}+k_{1,1}-2$ into a  pair  of zeroes of orders $k_{1,1}-1$
and $k_{1,2}-1$ as in figure~\ref{fig:breaking:up:a:zero}.

When one  of $k_{1,1},k_{1,2}$ is odd and another  one is even we
can break up the zero of order $k_{1,1}+k_{1,1}-2$ into a pair of
zeroes  of  orders  $k_{1,1}-1$  and  $k_{1,2}-1$  by  a  similar
construction, see  figure~\ref{fig:minus:2:1:local}. (Recall that
by convention  a ``zero  of order $-1$'' is a  simple pole of the
corresponding meromorphic quadratic  differential.) Cutting along
the  saddle  connection  we  obtain the desired surface  of  type
$-2.1$  with   prescribed   orders   of   interior  and  boundary
singularities.

\begin{figure}[ht]
%
\includegraphics{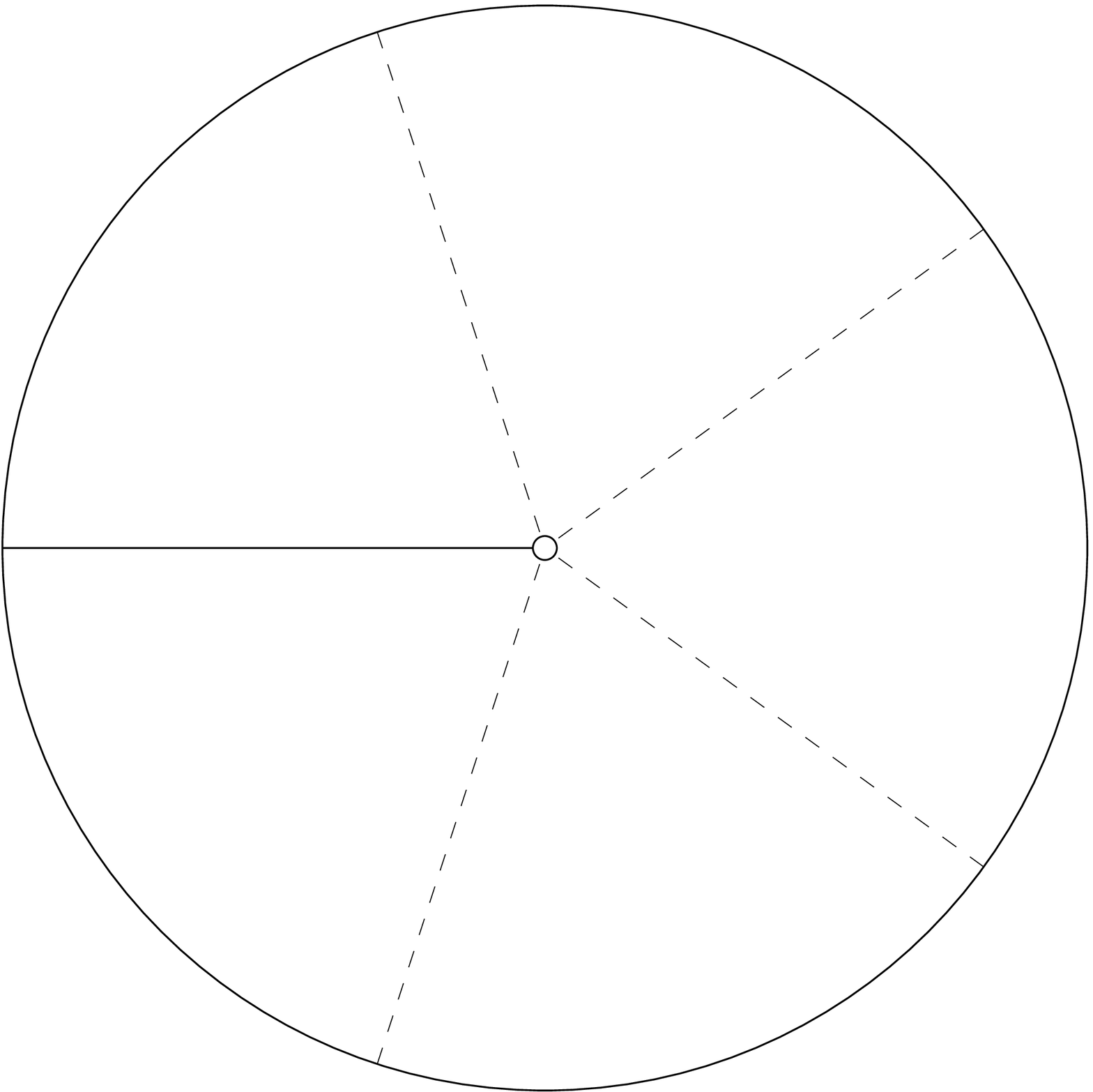}
\includegraphics{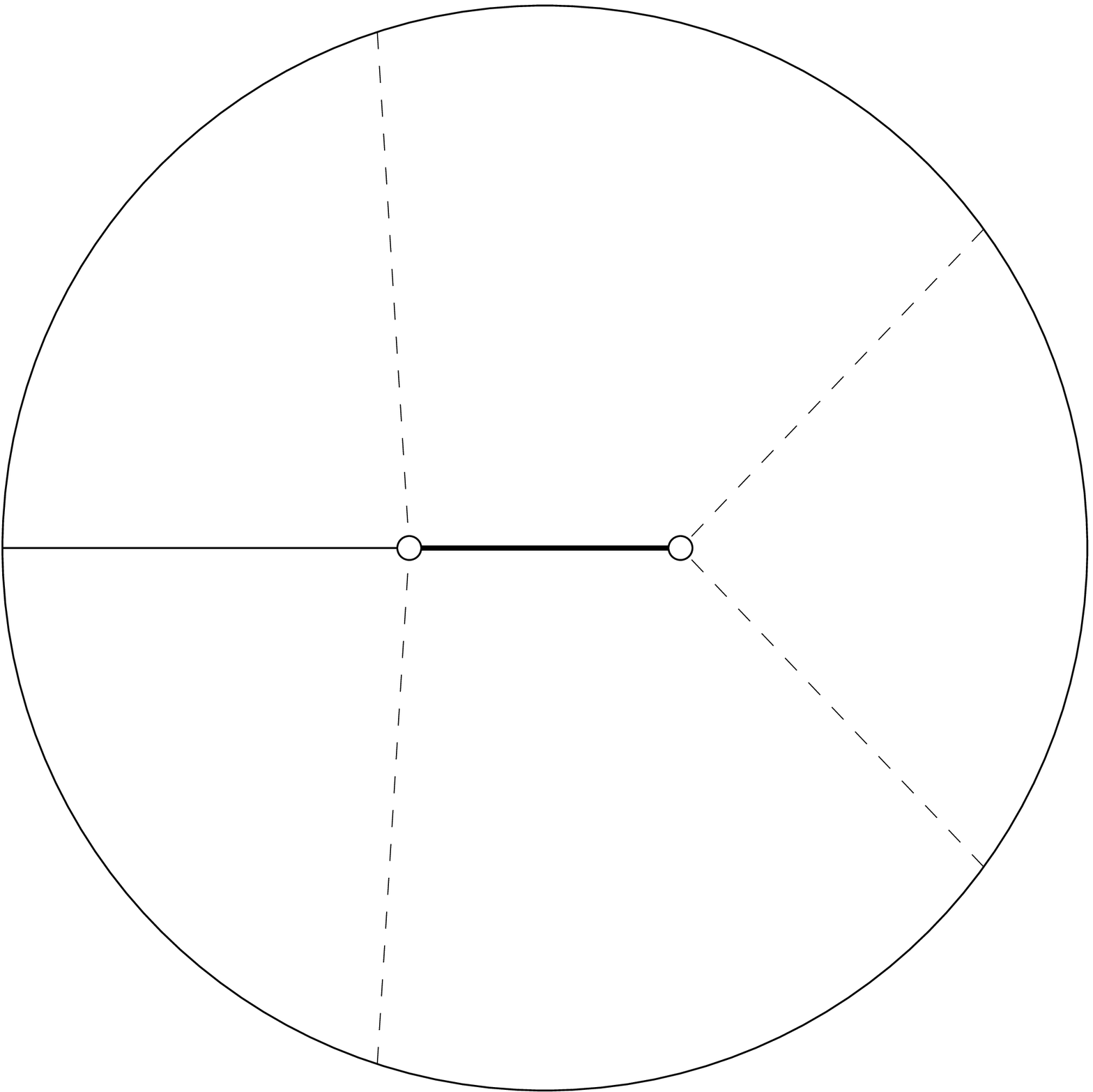}
\includegraphics{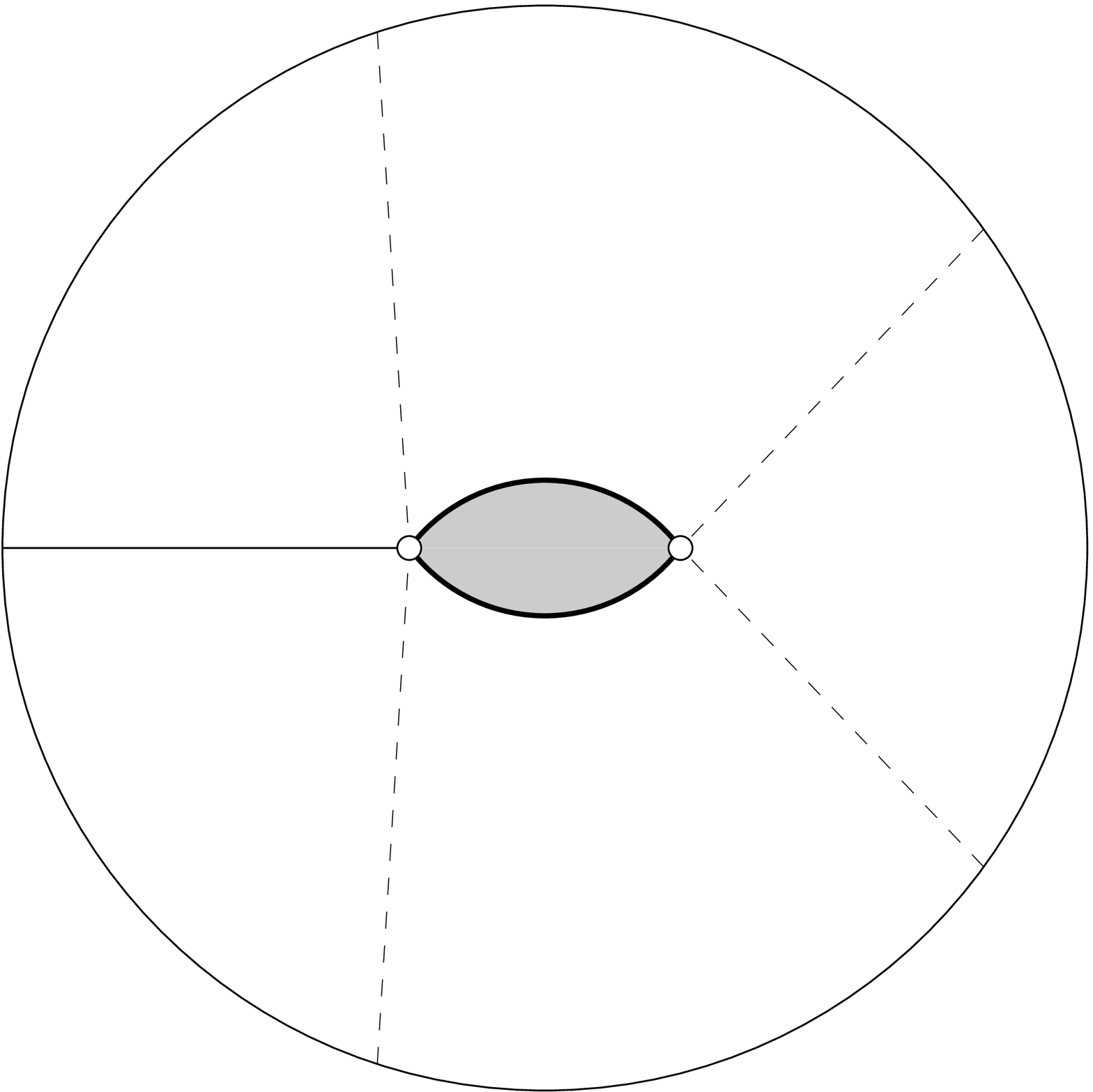}
%
%
\begin{picture}(0,0)(0,-140) 
\put(10,-6)
 {\begin{picture}(0,0)(0,0)
 \put(-160,-199){$\scriptstyle \varepsilon$}
 \put(-5,-199){$\scriptstyle P_2$}
 \put(-25,-175){$\scriptstyle \varepsilon-\delta$}
 \put(-52,-199){$\scriptstyle \varepsilon+\delta$}
 \put(-10,-209){$\scriptstyle \delta$}
 \put(-23,-199){$\scriptstyle P_1$}
 \put(-25,-229){$\scriptstyle \varepsilon-\delta$}
 \put(20,-189){$\scriptstyle \varepsilon$}
 \put(20,-216){$\scriptstyle \varepsilon$}
 \put(126,-199){$\scriptstyle P_2$}
 \put(87,-199){$\scriptstyle P_1$}
\end{picture}}
\end{picture}
\vspace{135bp} 
\caption{
\label{fig:minus:2:1:local}
Breaking up a zero of odd order into two zeroes and splitting the
saddle connection we get a surface of type $-2.1$}
\end{figure}

When both $k_{1,1},k_{1,2}$ are even, in fact $k_{1,1}\geq 2$ and
$k_{1,1}+k_{1,2}-2\ge   0$.   Let  $P$  be  the  zero  of   order
$k_{1,1}+k_{1,2}-2$  of  the quadratic  differential representing
the flat surface $S'_j$. Choose a pair of separatrices $\gamma_1,
\gamma_2$  in  such   way  that  the  angle  from  $\gamma_1$  to
$\gamma_2$  counted  counterclockwise  is  $\pi(k_{1,2}+1)$.  Let
$\Sigma_1$ be the  sector adjacent to $\gamma_1$ in the clockwise
direction and $\Sigma_2$  be the sector adjacent to $\gamma_2$ in
the   counterclockwise   direction.   Apply   the   parallelogram
construction to  $\Sigma_1,  \Sigma_2$  and  detach  $P$ into two
singularities $P_1, P_2$  (see  figure~\ref{fig:m:2:1:nonlocal}).
The orders of the boundary singularities of the resulting surface
$S_j$ are $k_{1,1}$ and $k_{1,2}$.
\begin{figure}[ht]
%
\includegraphics{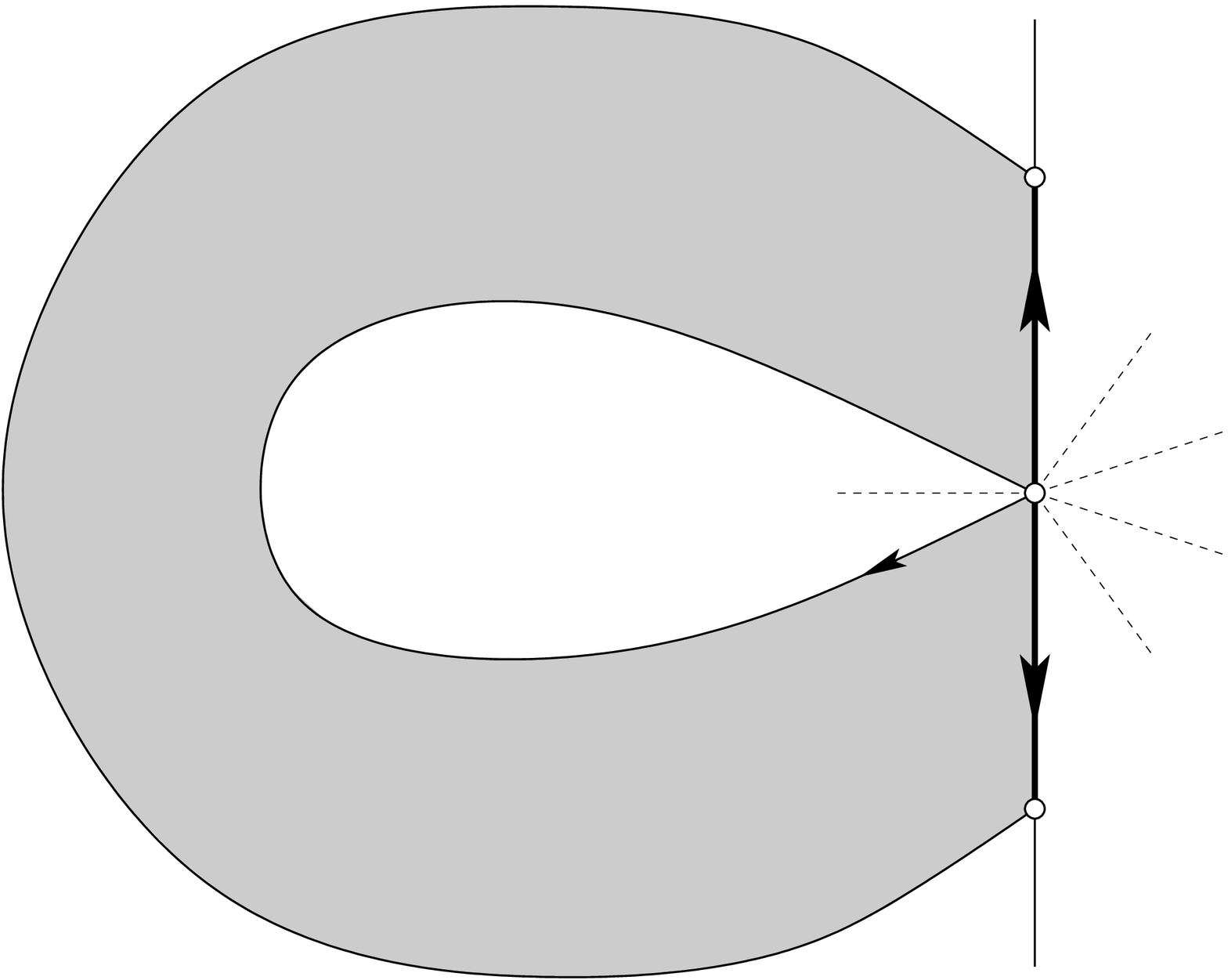}
\includegraphics{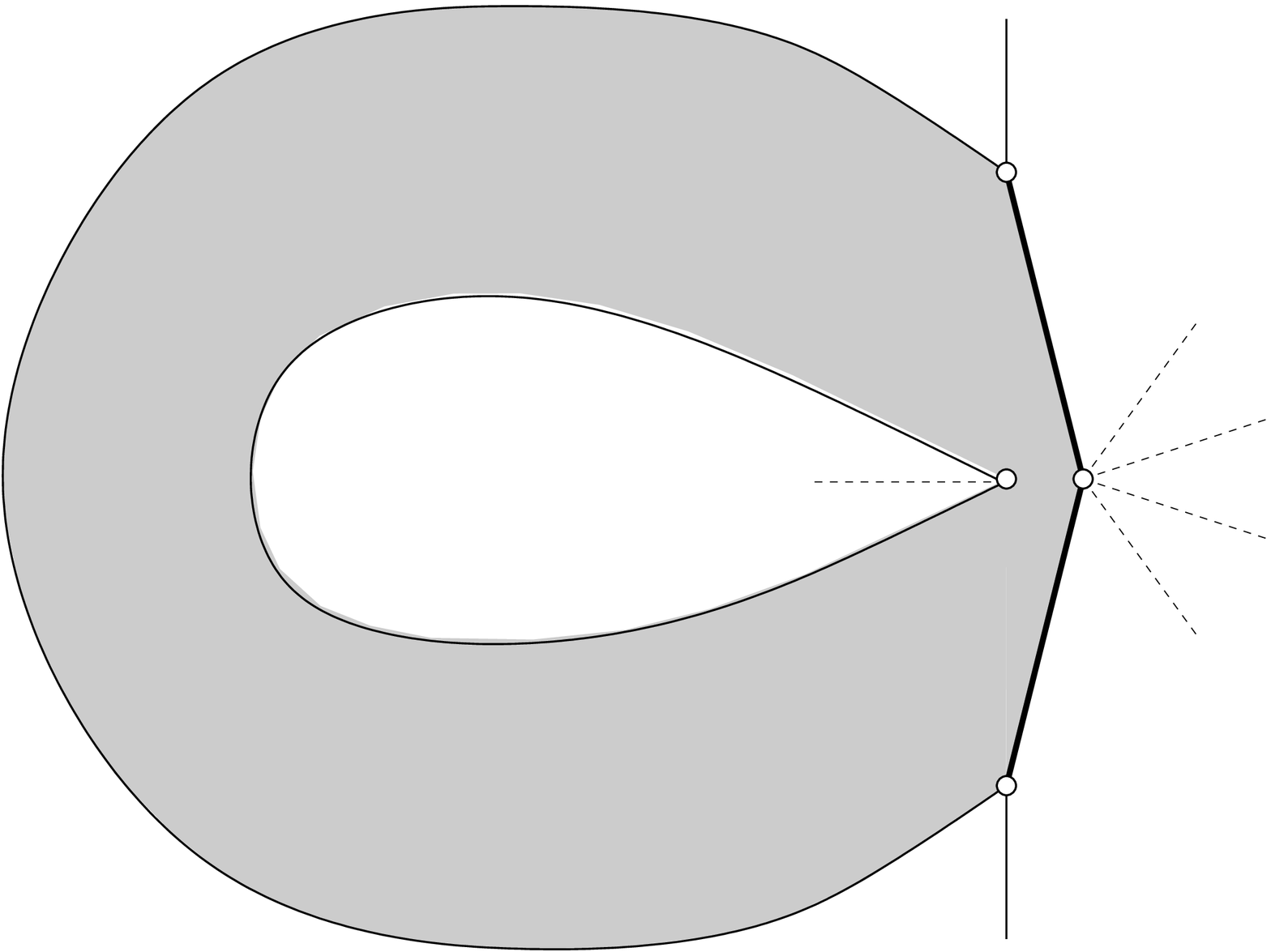}
\includegraphics{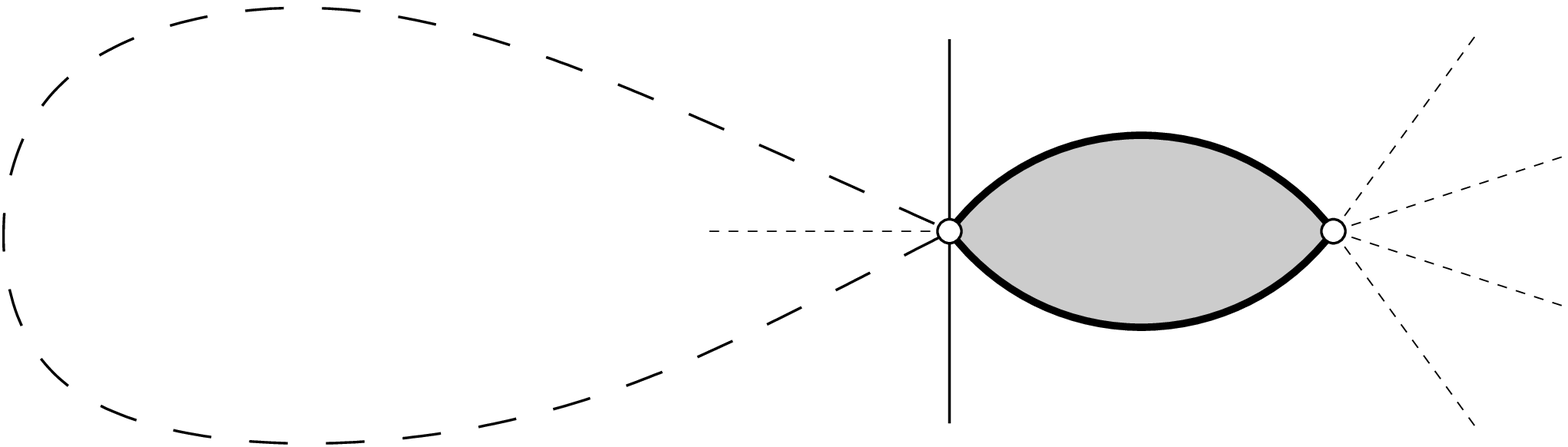}
%
%
\begin{picture}(0,0)(0,-150)
\put(28,12) 
 {\begin{picture}(0,0)(0,0)
 \put(-114,-200){$\scriptstyle \gamma_2$}
 \put(-127,-208){$\scriptstyle \Sigma_2$}
 \put(-127,-227){$\scriptstyle \Sigma_1$}
 \put(-114,-231){$\scriptstyle \gamma_1$}
\end{picture}}
\end{picture}
\vspace{115bp} 

\caption{
\label{fig:m:2:1:nonlocal}
Applying the parallelogram construction  to  a pair of sectors of
the same zero we get the missing surfaces of type $-2.1$}
\end{figure}

We have completed the proof of proposition~\ref{pr:minus:2:2}
and, thus, the proof of
proposition~\ref{pr:realizability:of:all:vertices}.
\end{proof}

Now         we         are         ready         to         prove
theorem~\ref{th:from:boundary:to:neighborhood:of:the:cusp}.  Note
that   theorem~\ref{th:from:boundary:to:neighborhood:of:the:cusp}
immediately   implies   the  missing   realizability   parts   of
theorems~\ref{th:graphs} and~\ref{th:all:local:ribon:graphs}.

\begin{proof}[Proof of theorem~\ref{th:from:boundary:to:neighborhood:of:the:cusp}]
Consider  a  configuration  $\cC$  (in  the  sense of the  formal
combinatorial      definition~\ref{def:configuration}).       Let
$\cQ(\alpha'_{\cC})$ (resp. $\cH(\beta'_{\cC}$)  be the principal
boundary stratum  corresponding  to  the configuration $\cC$. Let
$S'$   be    a    (possibly   nonconnected)   flat   surface   in
$\cQ(\alpha'_{\cC})$   (resp.   $\cH(\beta'_{\cC}$).   To   every
connected component of $S'$ apply the appropriate surgeries as in
sections~\ref{s:Local:Constructions}
and~\ref{s:Nonlocal:constructions}  realizing  the  corresponding
local ribbon graphs. We apply  the  surgeries in such a way  that
the saddle connections on the boundary of each surface $S_j$ are,
say, horizontal, and  have length proportional to their weight in
$\Gamma$ with  coefficient  $\delta$.  For  every  \cv-vertex  of
$\Gamma$ consider an appropriate flat  cylinder,  with  the  same
requirement for the boundary. Now we glue a compound surface from
the  components  $S_j$ as prescribed by the  graph  $\Gamma$.  By
construction the  result is a closed  surface $S$ endowed  with a
flat metric with linear holonomy restricted to $\{Id, -Id\}$.

By construction each  flat surface $S_j$ with boundary is endowed
with the canonical  orientation.  By definition the global ribbon
graph  $\G(\cC)$   is  endowed  with  the  canonical  orientation
compatible with the  canonical  orientation of the embedded local
ribbon graphs. This implies that the resulting closed surface $S$
inherits the canonical orientation. By  construction  $S$  has  a
collection  of  saddle connections  $\gamma_1,  \dots,  \gamma_n$
realizing the configuration $\cC$.

It remains to prove that  $S$  is nonsingular, i.e. that it  does
not  have  any  double  (triple,  ...)  points. Suppose it  does.
Detaching them we get  a  nonsingular closed flat surface $\tilde
S$. By  construction $\tilde S$  still has a collection of saddle
connections   $\gamma_1,    \dots,   \gamma_n$   realizing    the
configuration $\cC$,  which  means  that  assembling  the initial
surface $S$ we have performed some superfluous identifications of
several points of $\tilde S$.
\end{proof}

\subsection{Nonrealizable collections of singularities}
\label{ss:nonrealizable:configurations}

It was proved in~\cite{Masur:Smillie:realizability} that for  the
following exceptional list $\{\alpha''_1, \dots, \alpha''_4\}$ of
singularity data
\begin{equation}
\label{eq:empty:strata}
\big\{\emptyset,\ \{1,-1\},\ \{3,1\},\ \{4\}\big\}
\end{equation}
the four corresponding strata $\cQ(\alpha''_j)$ are  empty. It is
clear, that completing any of  these  lists  with entries ``$0$''
(which stand for  marked points) we  also get an  empty  stratum.
This      gives      rise      to      restriction      6      in
definition~\ref{def:configuration}  of  a configuration  which we
justify in this section.

Let $\G_{v_j}$ be a local  ribbon  graph of one of types  $-1.1$,
$-2.1$, $-2.2$ and let $\{d_1, \dots,  d_s\}$, $\{k_{1,1}, \dots,
k_{r,p(r)}\}$  be  a couple of unordered collections of  integers
satisfying        conditions        4       and        5       of
definition~\ref{def:configuration}  of  a configuration.  (In our
formal combinatorial definition they represent orders of interior
and  of  boundary singularities of a virtual  flat  surface  with
boundary.)     Applying     formally     equations~\eqref{eq:Dji}
and~\eqref{eq:Q:boundary:stratum} (which evaluate the singularity
data of  the  corresponding  component  of  the virtual principal
boundary  stratum)  to  our  combinatorial  data   we  obtain  an
unordered  collection   $\alpha'_j$   of   integers.  Consider  a
collection $\alpha''_j$ obtained from $\alpha'_j$ by omitting all
entries ``$0$'' (if any).

\begin{Lemma}
\label{lm:nonrealizable:data:combinatrorics}
The   collection   $\alpha''_j$  belongs   to   the   exceptional
list~\eqref{eq:empty:strata}  if  and  only if the  combinatorial
data $\big(\G_{v_j}$, $\{d_1, \dots,  d_s\}$,  $\{k_{1,1}, \dots,
k_{r,p(r)}\big)$  as   above   belongs   to   the   list  (6)  in
definition~\ref{def:configuration} of a configuration.
\end{Lemma}
\begin{proof}
The proof of the lemma is an exercise in elementary combinatorics.
\end{proof}

Haven justified  the  combinatorial obstructions we complete this
section with the corresponding geometric lemma.

\begin{Lemma}
\label{lm:nonrealizable:data}
Let $S_j$ be a flat surface with boundary of one of types $-1.1$,
$-2.1$,  $-2.2$.  Assume  that  $S_j$  does  not have any  saddle
connections parallel to  the  boundary different from those which
belong to the boundary. Then the corresponding combinatorial data
$\big(\G_{v_j}$,    $\{$unordered    collection    of    interior
singularities$\}$,   $\{$unordered   collection    of    boundary
singularities$\}\big)$ does  not  belong  to the exceptional list
(6) in definition~\ref{def:configuration} of a configuration.
\end{Lemma}
\begin{proof}
We use the following strategy to prove the lemma. If some surface
$S_j$ with boundary would define  an  entry from the list (6)  in
definition~\ref{def:configuration} we  would shrink the  boundary
of $S_j$ to get as  a  limit a nondegenerate surface $S'_j$  from
the  corresponding  component  $\cQ(\alpha'_j)$ of the  principal
boundary                    stratum.                     However,
lemma~\ref{lm:nonrealizable:data:combinatrorics}   implies   that
such $\cQ(\alpha'_j)$ is empty, which leads to a contradiction.

To complete the proof we  need  to describe how can one  ``shrink
the boundary'' of a flat surface. First note,  that boundary type
``$-1.1$'' can be considered as a  particular  case  of  boundary
type  ``$-2.1$''  when  the  order   of   one   of  the  boundary
singularities  is  equal to zero (see the corresponding  surfaces
with boundary in the table in section~\ref{ss:tables}).

Having  a  surface  $S_j$  of  type  $-2.1$ we can  isometrically
identify the pair  of boundary components  to get a  closed  flat
surface $S$.  The  corresponding  singularity  data  $\alpha$  of
$S\in\cQ(\alpha)$ is expressed  in  terms of the singularity data
of $S_j$ as follows:
$$
\alpha=\{d_1, \dots, d_{s(j)}, k_{1,1}-1, k_{1,2}-1\}
$$
This  implies  that  the  couples  $[\{d_1,  \dots,  d_{s(j)}\}$,
$\{k_{1,1}, k_{1,2}\}]$ of collections of orders  of interior and
of boundary singularities in the list below
$$ \begin{array}{c}
\emptyset, \{1, 1\}\qquad
\emptyset, \{2, 0\}\\
\{1,-1\}, \{1,1\};\qquad
\{1\}, \{1,0\};\qquad
\{-1\},\{2,1\}\\
\{3,1\},\{1,1\};\qquad
\{3\},\{2,1\};\qquad
\{1\},   \{4,1\}\\
\{4\},   \{1,1\};\qquad
\emptyset, \{5,1\};\qquad
\emptyset, \{4,2\}.
\end{array}
$$
are not realizable  by any surface  $S_j$ with boundary  of  type
$-2.1$, for in these  cases we would get a flat surface  $S$ from
an empty stratum, see equation~\eqref{eq:empty:strata}.

In the remaining  cases we get a closed surface $S\in\cQ(\alpha)$
with a distinguished pair of singularities $P_0, P_1$ joined by a
distinguished saddle  connection  $\gamma$. By assumptions of the
lemma this saddle connection is not parallel to  any other saddle
connection on $S$. This implies that deforming, if necessary, $S$
and  then  applying an appropriate element of $SL(2,\reals)$  the
surface $S$ can be continuously deformed  inside $\cQ(\alpha)$ to
a surface $\tilde S$ with a single {\it  short} saddle connection
$\tilde\gamma$  and  with no other short saddle connections.  The
deformation might be  performed in such  a way that  the  conical
singularities  $\tilde  P_0, \tilde P_1$ serving as endpoints  of
$\tilde\gamma$ would have the same cone angles as $P_0$ and $P_1$
correspondingly. But then  we  would apply an appropriate surgery
inverse        to       the        one        presented        in
figures~\ref{fig:breaking:up:a:zero},   \ref{fig:minus:2:1:local}
or~\ref{fig:m:2:1:nonlocal} to coalesce the corresponding pair of
zeroes into one. This would give  a  nondegenerate  flat  surface
$S'$.  Forgetting,  if necessary, the resulting marked points  on
$S'$ we  get  $S'\in\cQ(\alpha'')$,  where  $\alpha''$  is in the
list~\eqref{eq:empty:strata}                                 (see
lemma~\ref{lm:nonrealizable:data:combinatrorics}).   The   latter
leads to a contradiction since these strata are empty.

The  proof in  the  case of boundary  type  $-2.2$ is  completely
analogous. Lemma~\ref{lm:nonrealizable:data} is proved.
\end{proof}

\appendix
\section{Long saddle connections}
\label{ap:Long:saddle:connections}

We recall  the  definition  of the natural $GL(2;\R{})$-invariant
measure in the stratum $\cQ(\alpha)$. Let $\hat P = p^{-1}(P)$ be
the  collection  of  preimages  of  the  singularities of a  flat
surface $S\in\cQ(\alpha)$. Let $H_1^-(\hat S,\hat P;\,\Z)$ be the
subgroup in the  relative homology group  of $\hat S$,  odd  with
respect to  the  involution  $\tau$.  Similarly,  let $H^1_-(\hat
S,\hat P;\,\C{})$ be  the subspace in the relative cohomology odd
with  respect  to  the  involution  $\tau$  (i.e.  the  invariant
subspace corresponding  to  the  eigenvalue  $-1$  of the induced
linear  involution   $\tau^\ast:  H^1(\hat  S,\hat   P;\,\C{})\to
H^1(\hat S,\hat P;\,\C{})$). We can choose a basis in $H_1^-(\hat
S,\hat P; \Z)$ obtained as  lifts  $\hat\gamma_i$,  $i=1,  \dots,
\dim_{\C{}}\cQ(\alpha)$, of a collection of saddle connections on
$S$.  For  any  surface  near  $S$  the  affine  holonomy vectors
$\int_{\hat\gamma}\omega$   serve  as   local   coordinates   for
$\cQ(\alpha)$. We define a measure $d\nu(S)$  on $\cQ(\alpha)$ as
Lebesgue measure defined by these coordinates, normalized so that
the volume of a fundamental domain of the integer lattice in
$$
H^1_{-}(\hat S,\hat P;\integers\oplus i\integers)
\subset H_-^1(\hat S,\hat P;\cx)
$$
is equal to  one.

\begin{NNRemark}
Note that the  Abelian differential $\omega$  on $\hat S$  has  a
regular point at the  preimage  $P_i'\in p^{-1}(P_i)$ of a simple
pole $P_i$ of the quadratic differential $q$ on $S$. Consider the
set $\tilde P\subseteq \hat P$ obtained by removing these regular
points.  It  is  easy  to  see  that  the  canonical homomorphism
$H_-^1(\hat S,\hat P;\cx)\to H_-^1(\hat S,\tilde P;\cx)$  induced
by  the  inclusion  $\tilde  P\subseteq  \hat  P$ is actually  an
isomorphism. Thus, it does not matter which of two sets  $\hat P,
\tilde P$ is used to define the coordinate charts.
\end{NNRemark}

\begin{proof}[Proof of
proposition~\ref{pr:counting}]
Let $\cC$ be  an  admissible configuration of \^homologous saddle
connections.  Let  $\gamma=\{\gamma_1, \dots,  \gamma_n\}$  be  a
collection of \^homologous saddle connections on the flat surface
$S_0$  representing  configuration  $\cC$.  Choose  some   saddle
connection $\gamma_i$ corresponding to  an  edge of weight $1$ of
the  graph  $\Gamma(S,\gamma)$;  such  edge  always  exists,  see
figure~\ref{fig:classification:of:graphs}.  We  associate to  the
collection  $\gamma$  a pair of vectors $\pm \vec{v}(\gamma)  \in
\R{2}$  setting  $v=\int_{\gamma_i}\omega\in\cx\cong\R{2}$.   For
every surface $S$ in the same connected component we consider the
discrete subset $V_\cC(S)$ by taking the union $V_\cC(S)=\cup \pm
v(\gamma)$   over   all  collections   of   \^homologous   saddle
connections $\gamma$ realizing $\cC$.

It is easy to see that the set $V_\cC(S)$ satisfies axioms $(A)$,
$(B)$,              $(C_\mu)$              in~\cite{Eskin:Masur}.
Proposition~\ref{pr:counting}  now   follows  from  the   general
results          in~\cite{Eskin:Masur}          and          from
theorem~\ref{th:from:boundary:to:neighborhood:of:the:cusp}  which
implies that the Siegel--Veech constant $const_\cC$ is nonzero.
\end{proof}

\begin{proof}[Proof of
proposition~\ref{pr:homologous:equiv:parallel}]
If saddle  connections  $\gamma_1$  and  $\gamma_2$ are parallel,
then  $\int_{\hat\gamma_1}\omega=  r \int_{\hat\gamma_2}  \omega$
for $r$ real. If  $\gamma_1$  and $\gamma_2$ are not \^homologous
then  the  homology  classes  of  the  lifts  $\hat\gamma_1$  and
$\hat\gamma_2$ are independent in $H_1^{-}(\hat S,\hat  P;\,\Z)$.
Then the above equation  holds only for a set of measure  zero in
$H^1_{-}(\hat S,\hat P;\,\cx)$. Taking a countable  union of sets
of measure zero corresponding to possible  pairs  of  cycles  and
different coordinate  charts,  we  see  that  two non\^homologous
saddle connections on $S$  are parallel only for a set of  $S$ of
measure zero.
\end{proof}

\begin{proof}[Proof of
proposition~\ref{pr:rigid:configurations:hat:homologous}]
Suppose  that   there   are  two  saddle  connections  $\gamma_1,
\gamma_2$ in the  collection which are not \^homologous. Then the
corresponding     periods     $\int_{\hat\gamma_1}\omega$     and
$\int_{\hat\gamma_1}\omega$   correspond   to   two   independent
coordinates in a  small neighborhood of the initial flat surface,
and hence they can be  deformed  independently.  Since the length
$|\gamma|$       equals       $|\int_{\hat\gamma}\omega|$      or
$1/2|\int_{\hat\gamma}\omega|$ (depending on  whether $\gamma$ is
homologous  to  zero  or  not),  we  conclude  that  a collection
containing  two  non\^homologous  saddle  connections  cannot  be
rigid.

The       necessity        of       the       condition        in
proposition~\ref{pr:rigid:configurations:hat:homologous}       is
proved.      Sufficiency       immediately      follows      from
lemma~\ref{lm:gamma:2gamma}  which  says  that  the  lengths   of
\^homologous saddle connections  are either the same or differ by
a factor of two.
\end{proof}

\section{List of configurations in genus 2}
\label{a:List:of:configurations:in:genus:2}

Using                         definition~\ref{def:configuration},
theorem~\ref{th:from:boundary:to:neighborhood:of:the:cusp}    and
corollary~\ref{cr:unique:choice:of:parities},    and    following
examples~\ref{ex:example:of:a:configuration}
and~\ref{ex:principal:boundary}                                in
section~\ref{s:structure:of:the:paper}   one   can  construct   a
complete  list   of   configurations   for   any   given  stratum
$\cQ(\alpha)$.  In  this  section  we present an outline  of  the
algorithm and  list  all configurations for holomorphic quadratic
differentials in genus 2.

There are two natural parameters measuring  ``complexity''     of
singularity data $\alpha=\{d_1, \dots, d_m\}$: the genus $g$ of a
flat surface  $S$ in $\cQ(\alpha)$ and  the number $N$  of simple
poles on $S$  (i.e.  the number of conical  points  with the cone
angle $\pi$). Having a configuration $\cC$  denote  by  $N'$  the
number of interior  singularities  of order $-1$ corresponding to
this  configuration  and  by  $g'_1, \dots, g'_k$ the  genera  of
surfaces $S'_1,  \dots,  S'_k$  corresponding  to  the  principal
boundary  $\cQ(\alpha'_\cC)$  (correspondingly  $\cH(\beta'_\cC)$
when $\cC$ does  not have ``$-$''-vertices).  It is easy  to  see
that  the  number  of simple poles  on  $S$  (i.e.  the number of
entries ``$-1$'' of $\alpha$) might vary from $N'$ to $N'+4$, and
that  the  genus  $g$  might  vary  from  $\sum_{j=1}^k  g'_j$ to
$\sum_{j=1}^k  g'_j  +  2$  (see~\cite{Boissy:to:appear}  for  an
explicit  expression  of  $g(S)$  in  terms  of genera $g'_j$  of
components and of a structure of the global  ribbon graph). Thus,
having fixed the  upper  bounds for $g$ and  $N$,  we confine the
list of corresponding configurations to a finite one.

A  naive  algorithm of enumeration of all  configurations  for  a
given  stratum  $\cQ(\alpha)$ can be represented as follows.  Let
$g=g(\alpha)$ be the genus corresponding to  the singularity data
$\alpha$,
$$
d_1+\dots+d_m=4g(\alpha)-4
$$
Consider  complete  lists  of   (possibly   disconnected)  strata
$\cH(\beta')$ of genera $g-2$, $g-1$, $g$. These lists are finite
and  can  be  easily  constructed.  Consider  complete  lists  of
(possibly  disconnected)  strata  $\cQ(\alpha')$ of genera  $g-2,
g-1$, $g$ such that $\alpha'$ contains from $N-4$  to $N$ entries
``$-1$''  and   at  most  two  connected  components  $\alpha'_i,
\alpha'_j$   representing   strata  of   quadratic  differentials
$\cQ(\alpha'_i)$,  $\cQ(\alpha'_j)$  (the   remaining   connected
components are represented by strata of holomorphic differentials
$\cH(\alpha'_l)$). These lists  are also finite and can be easily
constructed. Add the empty set  to  these lists when $0\le g  \le
2$.

For  every  entry   $\alpha'=\alpha'_1\sqcup\dots\sqcup\alpha'_k$
(correspondingly $\beta'$) as above consider all possible ways to
organize the  set $\{\alpha'_1,\dots,\alpha'_k\}$ into one of the
graphs as  in figure~\ref{fig:classification:of:graphs}, in  such
way that  vertices  corresponding to the strata $\cH(\alpha'_j)$,
$\cH(\beta'_j)$ have ``$+$''-type, and vertices corresponding  to
the strata  $\cQ(\alpha'_j)$ have ``$-$''-type. Using these basic
graphs,  construct  all  possible   ``extended''   graphs  adding
vertices of the \cv-type as described in theorem~\ref{th:graphs}.

For every  vertex of every graph  as above consider  all possible
structures   of   an  embedded   local   ribbon   graph   as   in
figure~\ref{fig:embedded:local:ribbon:graphs}.

At    the     current    stage    we    have    already    chosen
$\alpha'=\{\alpha'_1,\dots,\alpha'_k\}$          (correspondingly
$\beta'$),    the    graph    $\Gamma$,    the    bijection    of
$\{\alpha'_1,\dots,\alpha'_k\}$   (correspondingly   $\{\beta'_1,
\dots,  \beta'_k\}$)  with  the  set  of   vertices  of  $\Gamma$
compatible with  the  structure  of ``$+$'' and ``$-$''-vertices,
and the structure  of  a local ribbon graph  for  every vertex of
$\Gamma$. Now for  every local ribbon graph $\G_j$ representing a
``$+$'' or ``$-$''-vertex $S_j$ consider  all  possible  ways  to
arrange  orders   of   interior  singularities  and  of  boundary
singularities  of  $S_j$ in  a  way  compatible  with  conditions
(3)--(6)   of   definition~\ref{def:configuration}    and    with
equation~\eqref{eq:H:boundary:stratum}   for   the  corresponding
singularity        data        $\beta'_j$        (correspondingly
equation~\eqref{eq:Q:boundary:stratum} for the  singularity  data
$\alpha'_j$). By ``compatibility'' with
equations~\eqref{eq:H:boundary:stratum}--\eqref{eq:Q:boundary:stratum}
we mean that  singularity data computed by these equations should
produce   $\beta'_j$   (correspondingly   $\alpha'_j$)   possibly
completed with several (from $1$ to $r_j$) entries ``$0$'' (where
$r_j$ is  the number of  connected components of the local ribbon
graph $\G_j$).

From  the  resulting lists of configurations extract those  which
correspond to the required singularity data $\alpha$.

Certainly this algorithm is  not  very efficient for large values
of $g$  or $N$. Nevertheless, for  strata in small  genera having
reasonable  number  of  simple  poles, it works quite   well
(especially being slightly optimized using specific properties of
given data  $\alpha$).

As an example  we present a  complete list of  configurations  of
\^homologous   saddle   connections  for   holomorphic  quadratic
differentials  in  genus 2. We are grateful  to  Alex~Eskin,  who
helped us to test completeness of this list.

$$
\begin{array}{|c|c|c|}
\hline
&&\\
[-\halfbls]
\cQ(2,2) & \cQ(2,1,1) & \cQ(1,1,1,1) \\
[-\halfbls]&&\\
\hline&&\\
%
%
& 
\includegraphics{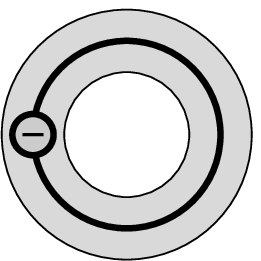}
\begin{picture}(100,50)(0,0) %
\put(-150,80)
{
\begin{picture}(0,0)(0,0)
\put(172,-53){$\scriptstyle 2$}
\put(194,-53){$\scriptstyle 2$}
\put(150,-53){$\scriptstyle \{2\}$}
\end{picture}}
\end{picture}
& 
\includegraphics{gen2_21.eps}
\begin{picture}(100,50)(0,0) %
\put(-150,80)
{
\begin{picture}(0,0)(0,0)
\put(172,-53){$\scriptstyle 2$}
\put(194,-53){$\scriptstyle 2$}
\put(145,-53){$\scriptstyle \{1,1\}$}
\end{picture}}
\end{picture}
\\\hline&&\\ 
\includegraphics{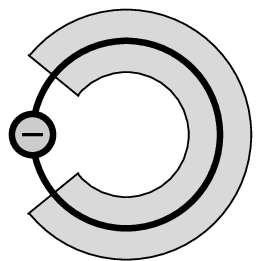}
\begin{picture}(100,50)(0,0) %
\put(-150,80)
{
\begin{picture}(0,0)(0,0)
\put(190,-45){$\scriptstyle 1$}
\put(190,-61){$\scriptstyle 1$}
\put(160,-53){$\scriptstyle \{2\}$}
\end{picture}}
\end{picture}
&  
\includegraphics{gen2_22.eps}
\begin{picture}(100,50)(0,0) %
\put(-150,80)
{
\begin{picture}(0,0)(0,0)
\put(190,-45){$\scriptstyle 1$}
\put(190,-61){$\scriptstyle 1$}
\put(150,-53){$\scriptstyle \{1,1\}$}
\end{picture}}
\end{picture}
&  
\\ \hline
%
&  
\includegraphics{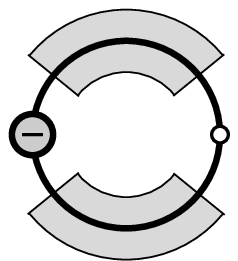}
\begin{picture}(100,50)(0,0) %
\put(-150,80)
{
\begin{picture}(0,0)(0,0)
\put(190,-45){$\scriptstyle 1$}
\put(190,-61){$\scriptstyle 1$}
\put(160,-53){$\scriptstyle \{2\}$}
\put(221,-45){$\scriptstyle 0$}
\put(221,-61){$\scriptstyle 0$}
\end{picture}}
\end{picture}
&  
\includegraphics{gen2_22_o.eps}
\begin{picture}(100,63)(0,0) %
\put(-150,80)
{
\begin{picture}(0,0)(0,0)
\put(190,-45){$\scriptstyle 1$}
\put(190,-61){$\scriptstyle 1$}
\put(150,-53){$\scriptstyle \{1,1\}$}
\put(221,-45){$\scriptstyle 0$}
\put(221,-61){$\scriptstyle 0$}
\end{picture}}
\end{picture}
\\ \hline
\includegraphics{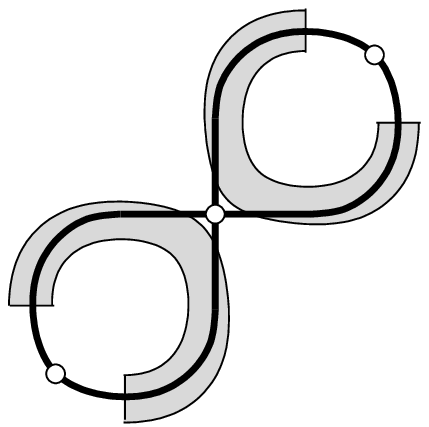}
\begin{picture}(110,100)(0,0) 
\put(-150,80)
{
\begin{picture}(0,0)(0,0)
\put(189,-46){$\scriptstyle 0$}
\put(196,-39){$\scriptstyle 0$}
\put(207,-27){$\scriptstyle 0$}
\put(214,-20){$\scriptstyle 0$}

\put(156,-58){$\scriptstyle 0$}
\put(177,-79){$\scriptstyle 0$}

\put(226,13){$\scriptstyle 0$}
\put(247,-8){$\scriptstyle 0$}
\end{picture}}
\end{picture}
&
\includegraphics{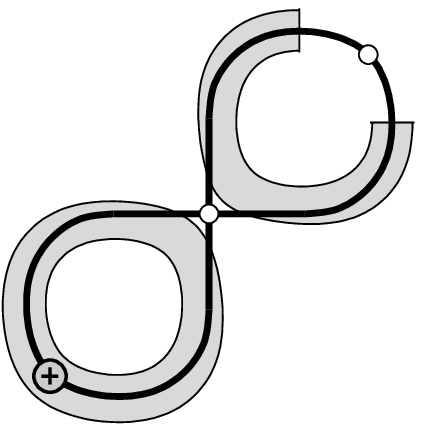}
\begin{picture}(110,100)(0,0) 
\put(-150,80)
{
\begin{picture}(0,0)(1,1)
\put(189,-46){$\scriptstyle 0$}
\put(196,-39){$\scriptstyle 0$}
\put(207,-27){$\scriptstyle 0$}
\put(214,-20){$\scriptstyle 0$}

\put(149,-71){$\scriptstyle \emptyset$}

\put(158,-77){$\scriptstyle 1$}
\put(172,-63){$\scriptstyle 1$}

\put(226,13){$\scriptstyle 0$}
\put(247,-8){$\scriptstyle 0$}
\end{picture}}
\end{picture}
& 
\includegraphics{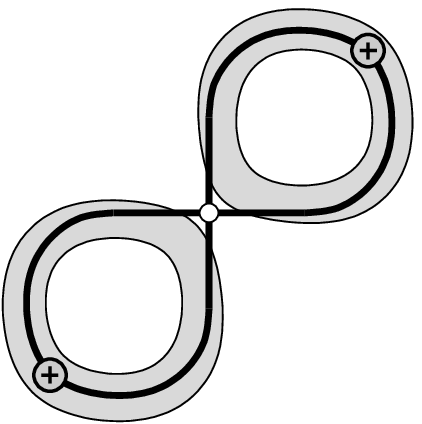}
\begin{picture}(110,100)(0,0) 
\put(-150,80)
{
\begin{picture}(0,0)(1,1)
\put(189,-46){$\scriptstyle 0$}
\put(196,-39){$\scriptstyle 0$}
\put(207,-27){$\scriptstyle 0$}
\put(214,-20){$\scriptstyle 0$}

\put(149,-71){$\scriptstyle \emptyset$}

\put(158,-77){$\scriptstyle 1$}
\put(172,-63){$\scriptstyle 1$}

\put(232,-2){$\scriptstyle 1$}
\put(246,12){$\scriptstyle 1$}

\put(253,7){$\scriptstyle \emptyset$}

\end{picture}}
\end{picture}
\\\hline
%
\includegraphics{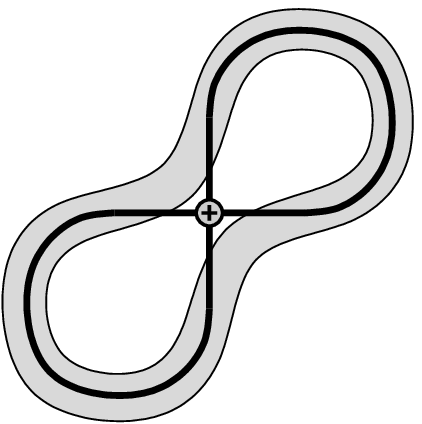}
\begin{picture}(110,100)(0,0) 
\put(-150,80)
{
\begin{picture}(0,0)(1,1)
\put(207,-25){$\scriptstyle \emptyset$}

\put(189,-18){$\scriptstyle 1$}
\put(196,-25){$\scriptstyle 1$}
\put(207,-39){$\scriptstyle 1$}
\put(214,-46){$\scriptstyle 1$}

\end{picture}}
\end{picture}
& 
\includegraphics{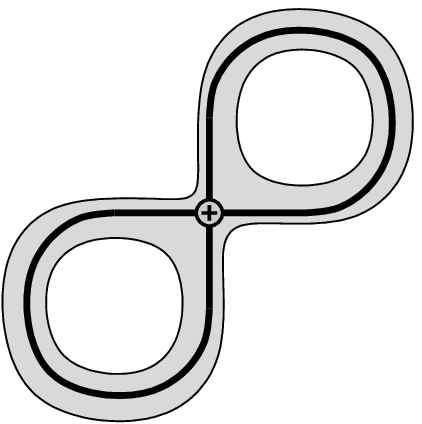}
\begin{picture}(110,100)(0,0) 
\put(-150,80)
{
\begin{picture}(0,0)(1,1)
\put(214,-46){$\scriptstyle \emptyset$}

\put(190,-45){$\scriptstyle 2$}
\put(213,-19){$\scriptstyle 2$}

\put(208,-40){$\scriptstyle 1$}
\put(195,-25){$\scriptstyle 1$}

\end{picture}}
\end{picture}
& 
\includegraphics{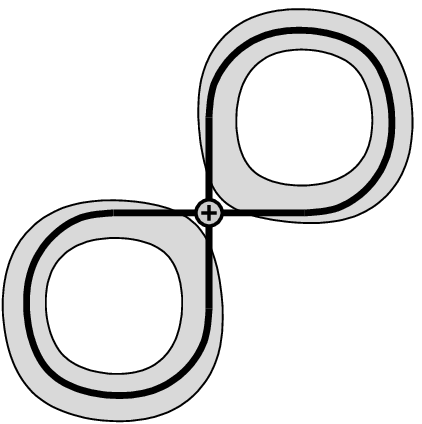}
\begin{picture}(0,90)(0,0) %
\put(-200,80)
{
\begin{picture}(0,0)(6,1) 
\put(189,-46){$\scriptstyle 2$}
\put(196,-39){$\scriptstyle 2$}
\put(207,-27){$\scriptstyle 2$}
\put(214,-20){$\scriptstyle 2$}

\put(205,-37){$\scriptstyle \emptyset$}

\end{picture}}
\end{picture}
\\ \hline
\end{array}
$$




\end{document}

%% file: ribbon_graphs_submit.tex
\begin{figure}

  %
\includegraphics{valence1_min.eps}
\includegraphics{valence21_min.eps}
\includegraphics{valence22_min.eps}
   %
\includegraphics{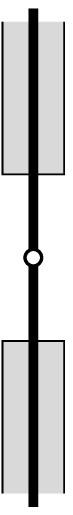}
\includegraphics{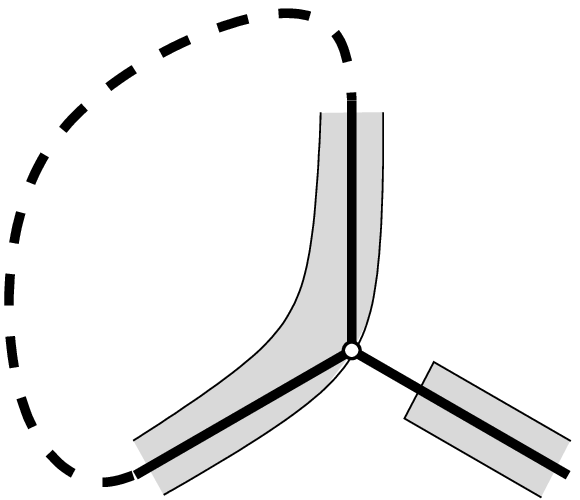}
\includegraphics{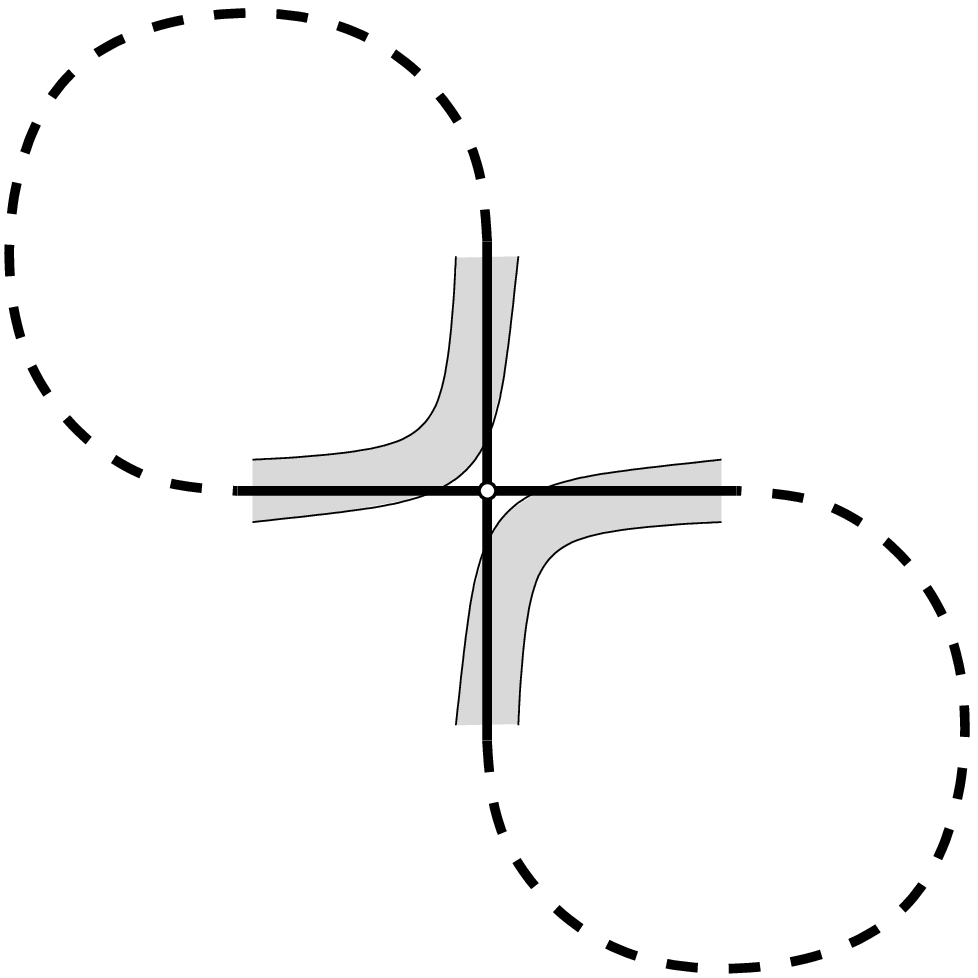}
   %
\includegraphics{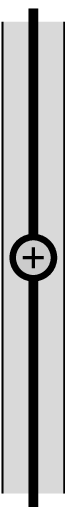}
\includegraphics{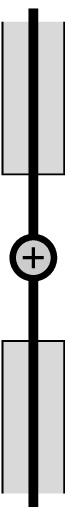}
   %
\includegraphics{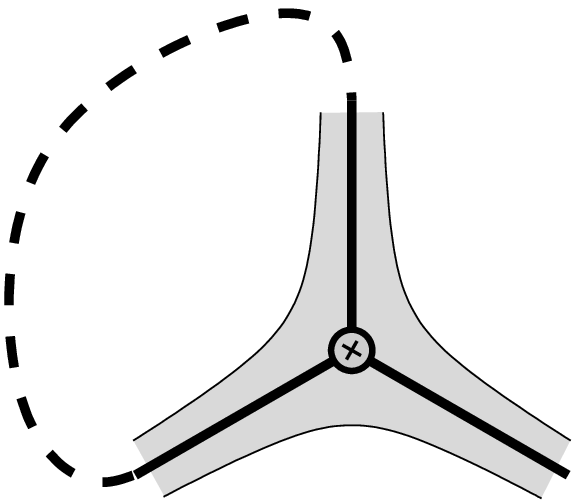}
\includegraphics{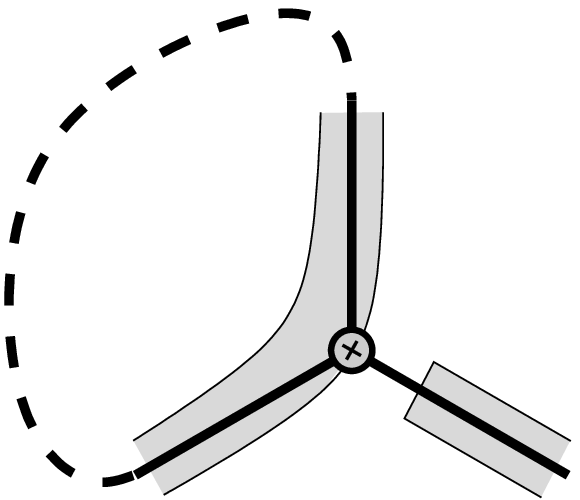}
\includegraphics{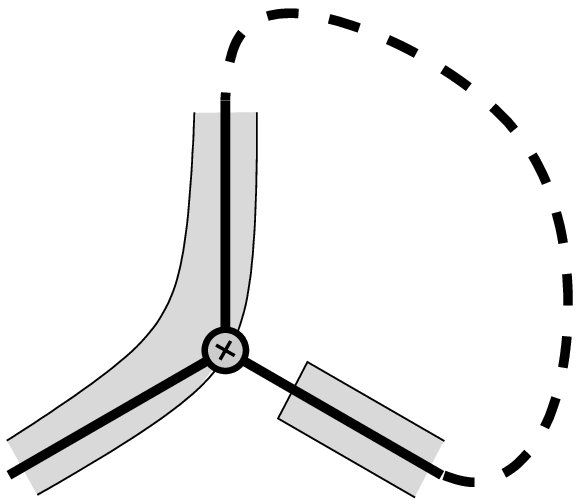}
\includegraphics{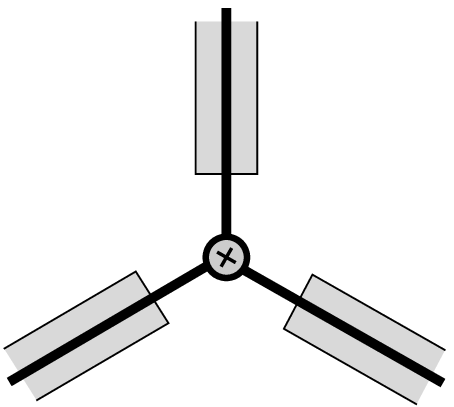}
   %
\includegraphics{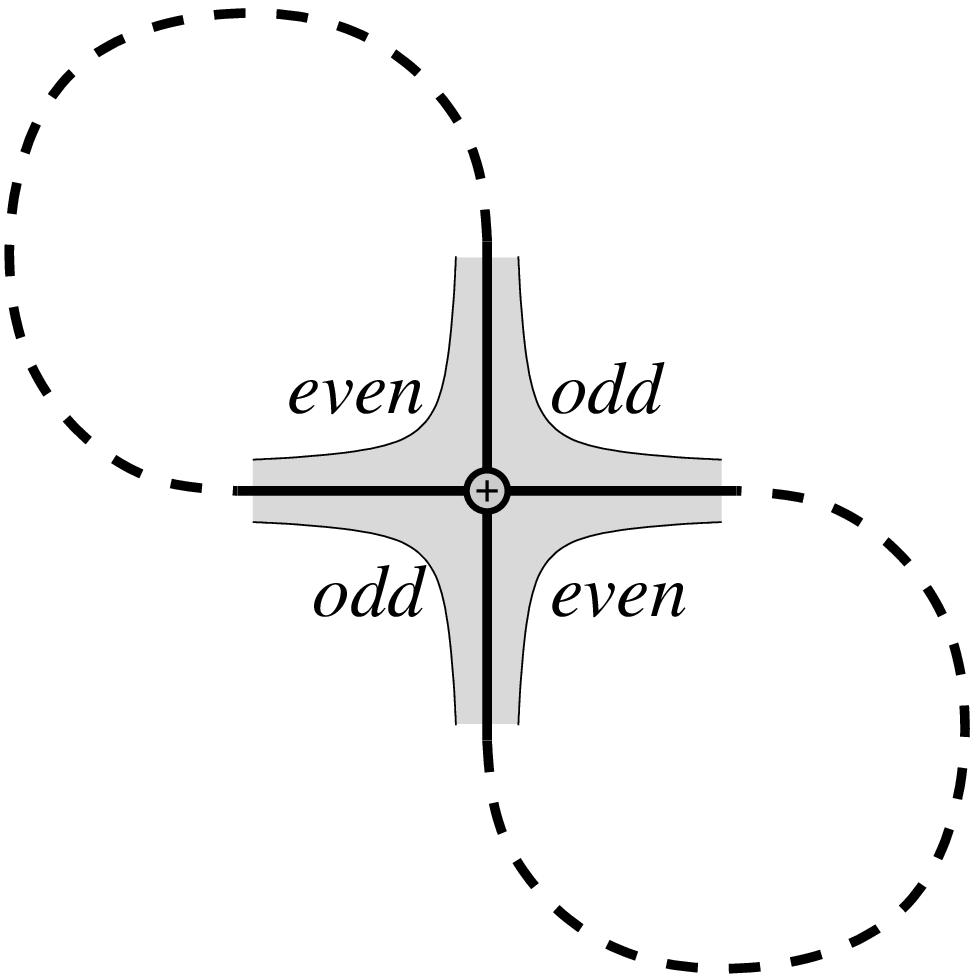}
\includegraphics{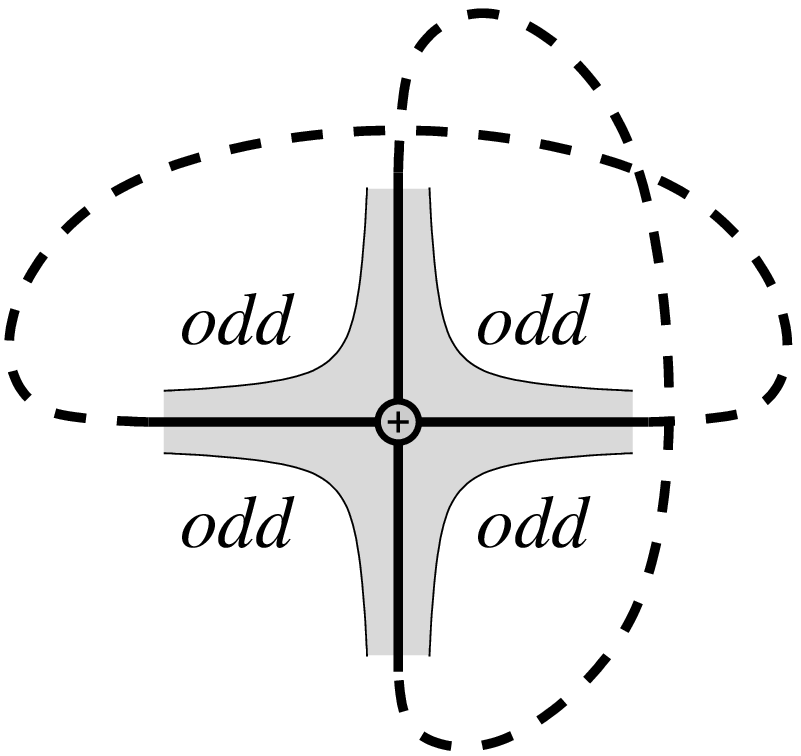}
\includegraphics{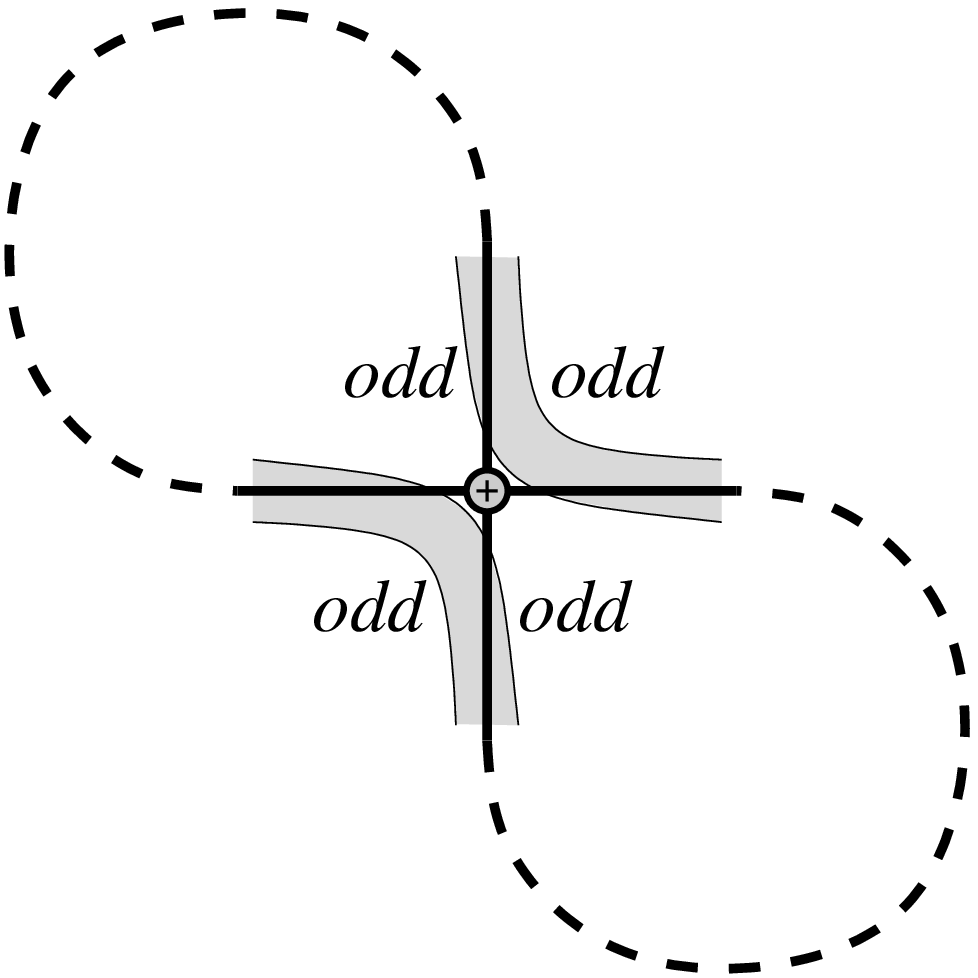}
\includegraphics{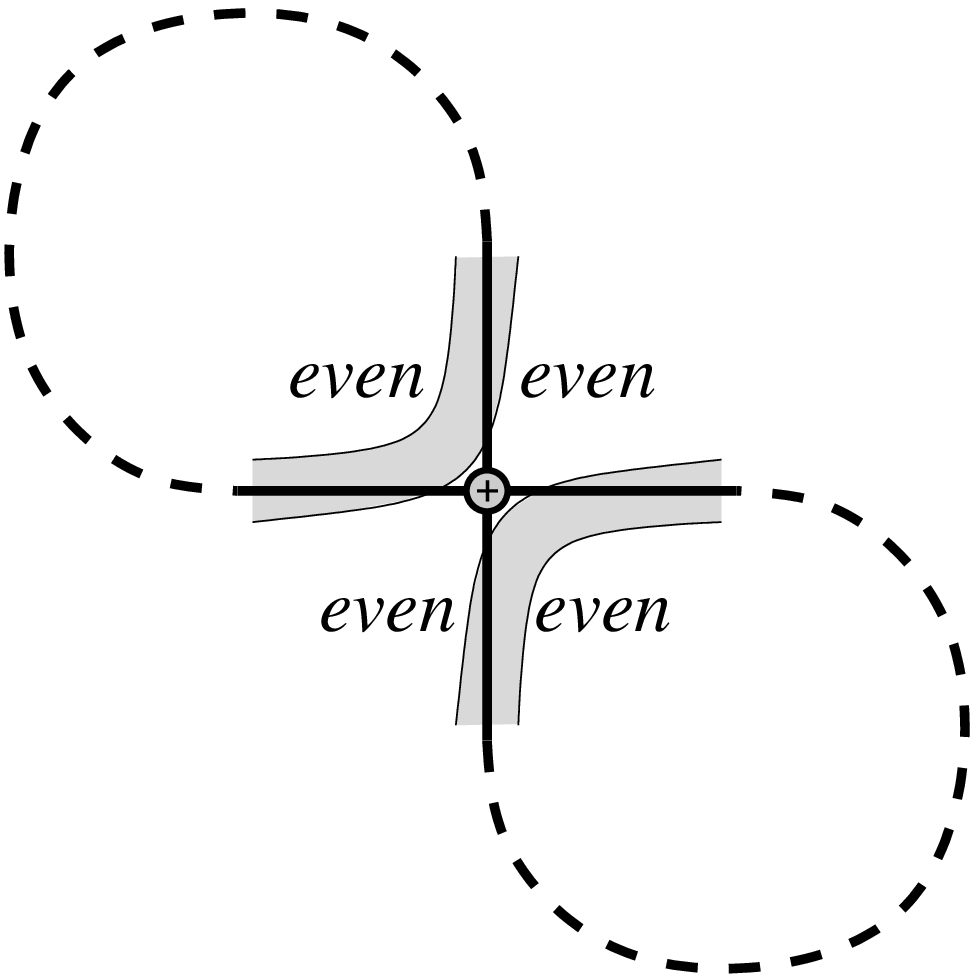}
   %
   %
\includegraphics{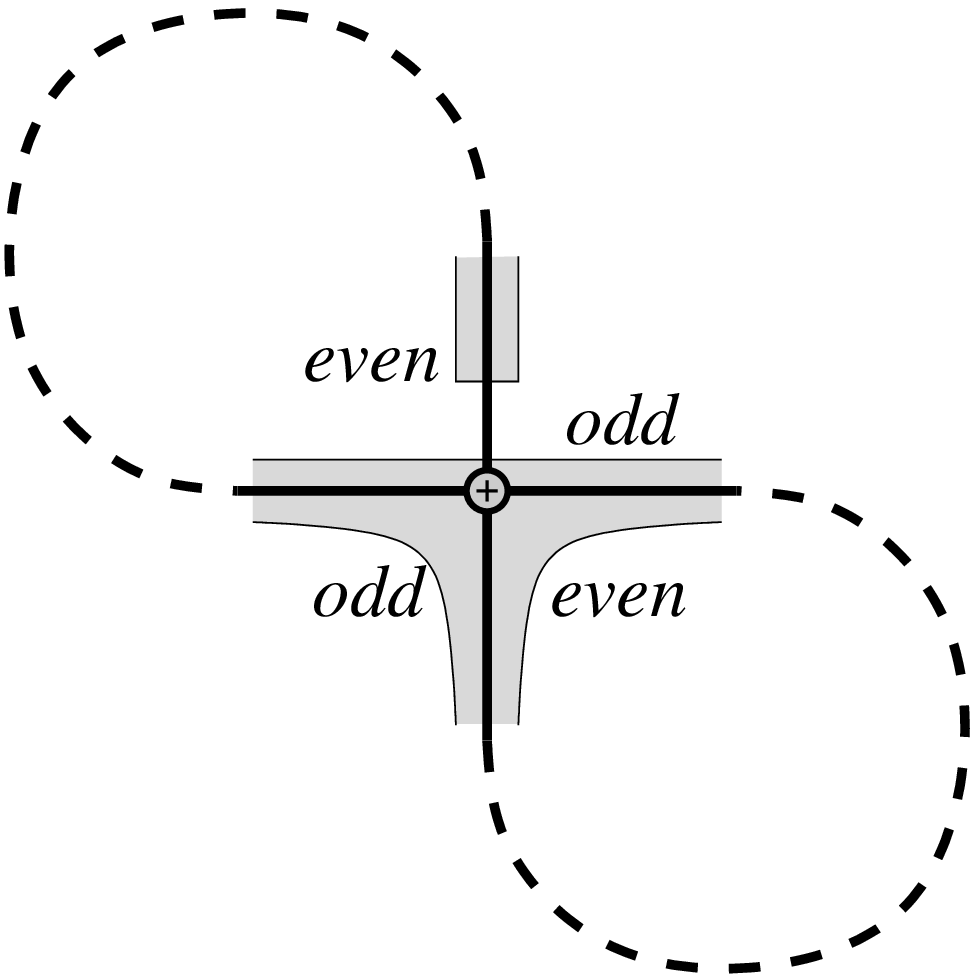}
\includegraphics{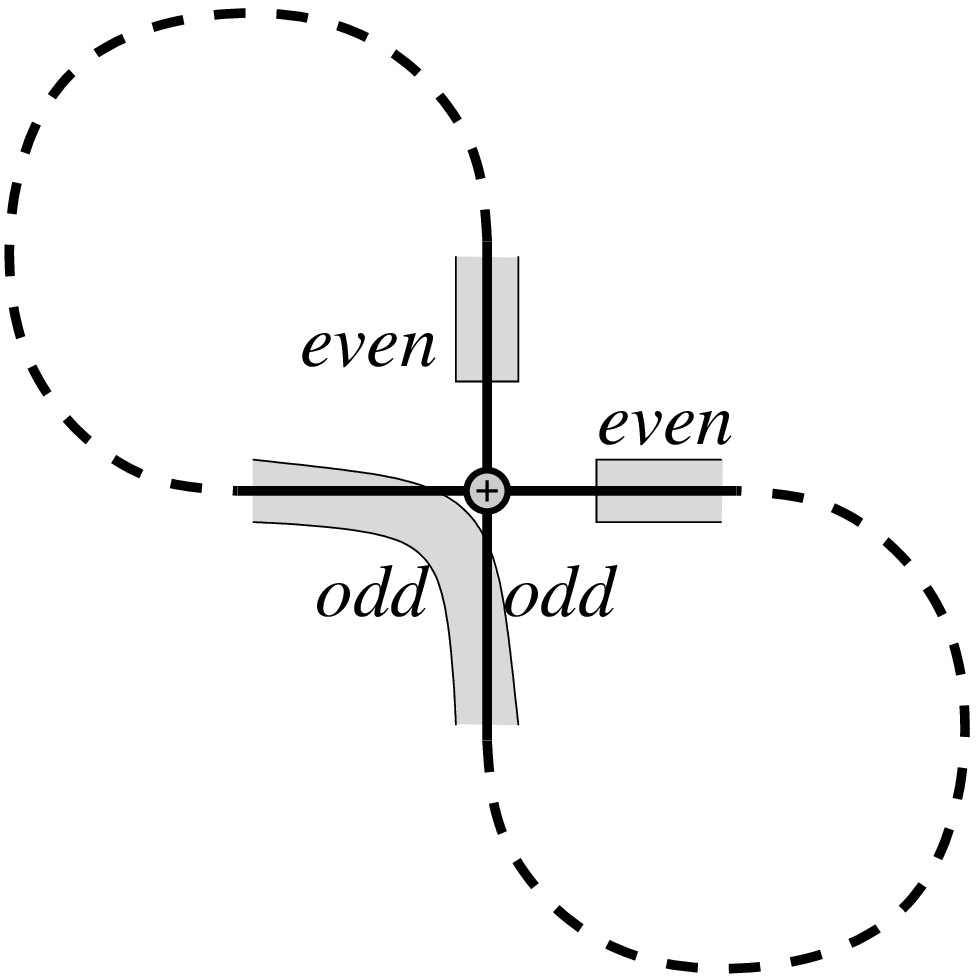}
\includegraphics{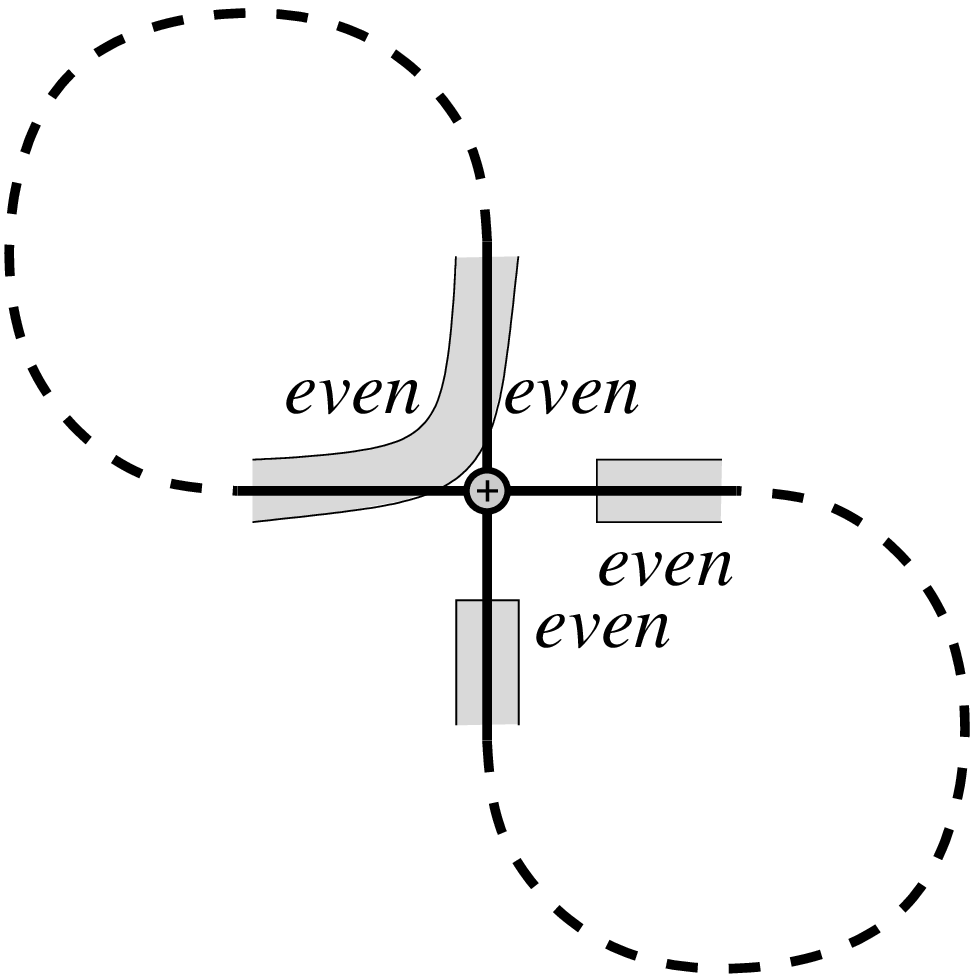}
\includegraphics{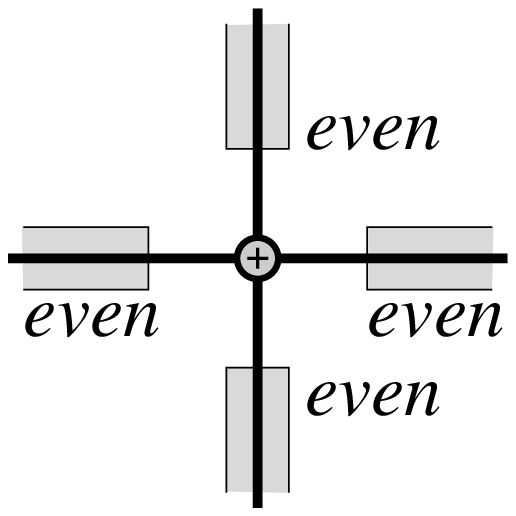}
%
   %
\vspace{523bp}
\begin{picture}(0,0)(160,-500) 
\put(-2,0){$-1.1$}
\put(11,22){\scriptsize\textit{any}}

\put(125,0){$-2.1$}
  \put(150,22){\scriptsize\textit{any}}
  \put(150,8){\scriptsize\textit{any}}

\put(278,0){$-2.2$}
  \put(286,22){\scriptsize\textit{any}}
  \put(315,22){\scriptsize\textit{any}}
\put(-2,-70){$\circ2.2$}
  \put(19,-47){$\scriptstyle 0$}
  \put(30,-47){$\scriptstyle 0$}

\put( 125,-70){$\circ3.2$}
  \put(149.5,-40){$\scriptstyle 0$} 
  \put(149.5,-47){$\scriptstyle 0$} 
  \put(155,-61){$\scriptstyle 0$}

\put( 250,-70){$\circ4.2$}
  \put(265,-44){$\scriptstyle 0$}
  \put(270,-49){$\scriptstyle 0$}
  \put(277,-58){$\scriptstyle 0$}
  \put(282,-63){$\scriptstyle 0$}
\put(-2,-125){$+2.1$}
  \put(19,-103){\scriptsize\textit{odd}}
  \put(19,-117){\scriptsize\textit{odd}}
\put(125,-125){$+2.2$}
  \put(131,-103){\scriptsize\textit{even}}
  \put(164,-103){\scriptsize\textit{even}}
\put( -2,-228){$+3.1$}
  \put(3,-173){\scriptsize\textit{even}}
  \put(-10,-192){\scriptsize\textit{odd}}
  \put( 17,-192){\scriptsize\textit{odd}}
\put( 98,-228){$+3.2a$}
  \put(103,-173){\scriptsize\textit{even}}
  \put( 90,-187){\scriptsize\textit{even}}
  \put(115,-195){\scriptsize\textit{even}}
\put(193,-228){$+3.2b$}
  \put(199,-173){\scriptsize\textit{odd}}
  \put(189,-189){\scriptsize\textit{odd}}
  \put(209,-196){\scriptsize\textit{even}}
\put(293,-228){$+3.3$}
  \put(315,-186){\scriptsize\textit{even}}
  \put(280,-186){\scriptsize\textit{even}}
  \put(309,-196){\scriptsize\textit{even}}
  \begin{picture}(0,0)(0,0) 
  \put(0,0){
  \put( 10,-365){$+4.1a$}
  \put(100,-365){$+4.1b$}
  \put(190,-365){$+4.2a$}
  \put(275,-365){$+4.2b$}
  \put( 10,-485){$+4.2c$}
  \put(100,-485){$+4.3a$}
  \put(190,-485){$+4.3b$}
  \put(275,-485){$+4.4$}
  }
  \end{picture}
\end{picture}
\caption{
\label{fig:embedded:local:ribbon:graphs}
Classification of embedded local ribbon graphs
}
\end{figure}

%% file: ribbon_versus_holes_submit.tex

\subsection{Surfaces with boundary versus local ribbon graphs}
\label{ss:tables}
Figure~\ref{fig:embedded:local:ribbon:graphs},
theorem~\ref{th:all:local:ribon:graphs}                       and
proposition~\ref{pr:realizability:of:all:vertices} are formulated
in terms of local ribbon graphs. By definition they correspond to
flat surfaces with boundary. The tables below explicitly describe
this correspondence.

$$ 
\begin{array}{|c|c|c|}
\hline 
&&\\
[-\halfbls]
-1.1 & -2.1 & -2.2 \\
   %
\hline && \\
&&\\  && \\
&&\\
\special{                     
psfile=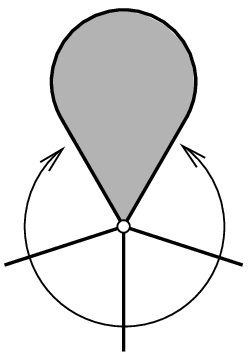
hscale=50 %
vscale=50 %
angle=0
hoffset=25  
voffset=-10 
}
\begin{picture}(90,0)(-10,60) %
\put(-150,80) {
\begin{picture}(0,0)(0,0)
\put(173,-34){\scriptsize\it any}
\end{picture}}
\end{picture}
& 
\special{                   
psfile=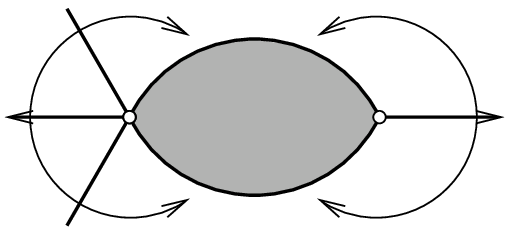 hscale=50 vscale=50 angle=0
hoffset=25  
voffset=0   
}
\begin{picture}(120,0)(-10,0)
\put(-60,-100) %
{
\begin{picture}(0,0)(0,0)
\put(54,116){\scriptsize\it any} \put(147,116){\scriptsize\it
any}
\end{picture}}
\end{picture}
& 
\special{                     
psfile=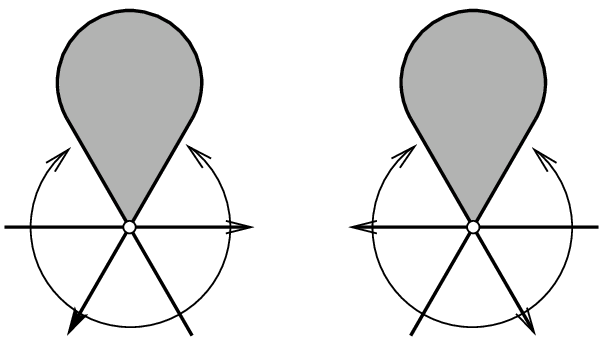 hscale=50 vscale=50 angle=0
hoffset=5    
voffset=-10  
}
\begin{picture}(100,0)(0,10) %
\put(-220,50) %
{
\begin{picture}(0,0)(-10,0)
\put(223,-56){\scriptsize\it any} \put(272,-56){\scriptsize\it
any}
\end{picture}}
\end{picture}
\\
&&\\
[-\halfbls]&&\\
\hline && \\ 
   %
&&\\
\includegraphics{valence1_min.eps}
\begin{picture}(0,0)(0,0)
\put(-15,-6){\scriptsize\it any}
\end{picture}
&
\includegraphics{valence21_min.eps}
\begin{picture}(0,0)(5,0)
\put(-8,13){\scriptsize\it any} \put(-8,-6){\scriptsize\it any}
\end{picture}&
\includegraphics{valence22_min.eps}
\begin{picture}(0,0)(-1,0)
\put(-30,-6){\scriptsize\it any} \put(10,-6){\scriptsize\it any}
\end{picture}
\\
&&\\
   %
\hline
\hline \vspace{3truept}
&&\\
[-\halfbls]
\circ2.2 & \circ3.2 & \circ4.2  \\
   %
\hline && \\
&&\\  
&&\\
\special{                  
psfile=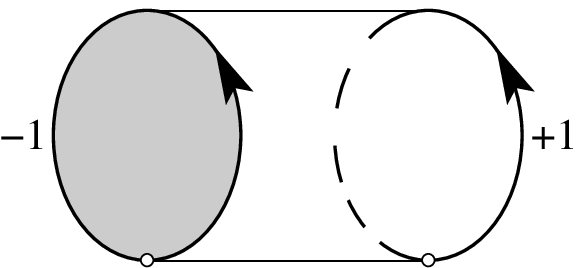 hscale=50 vscale=50 angle=0
hoffset=-40 %
voffset=-25 %
}
& 
\special{                   
psfile=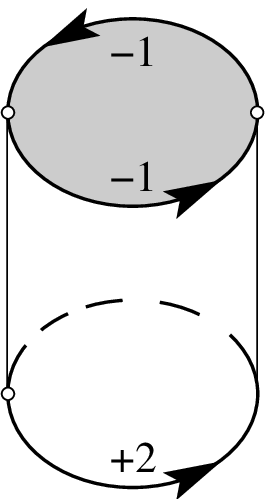 
hscale=50 vscale=50 angle=0
hoffset=-25 
voffset=-32 
}
& 
\special{                     
psfile=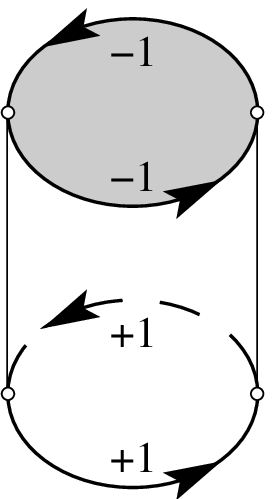 hscale=50 vscale=50 angle=0
hoffset=-20 
voffset=-32 %
}
\\
&&\\
   %
&& \\  && \\
\hline && \\ 
&&\\ &&\\
[-\halfbls]&&\\
\includegraphics{valence22o.eps}
\begin{picture}(0,0)(0,0)
  \put(-8,6){$\scriptstyle 0$}
  \put(7.5,6){$\scriptstyle 0$}
\end{picture}
&
\includegraphics{valence32o.eps}
\begin{picture}(0,0)(157,-49) 
  \put(149.5,-42){$\scriptstyle 0$} 
  \put(149.5,-47.5){$\scriptstyle 0$} 
  \put(155,-60){$\scriptstyle 0$}
\end{picture}
&
\includegraphics{valence42b_o.eps}
\begin{picture}(0,0)(275,-50) 
  \put(266,-45){$\scriptstyle 0$}
  \put(270,-49){$\scriptstyle 0$}
  \put(276.5,-57.5){$\scriptstyle 0$}
  \put(280.5,-61.5){$\scriptstyle 0$}
\end{picture}
\\
&&\\
&&\\
&&\\
&&\\
   %
   %
\hline
\end{array}
$$

$$ 
\begin{array}{|c|c|}
\hline 
&\\
[-\halfbls]
+2.1 & +2.2 \\
   %
\hline & \\
[-\halfbls]&\\  & \\ & \\
&\\
\special{                     
psfile=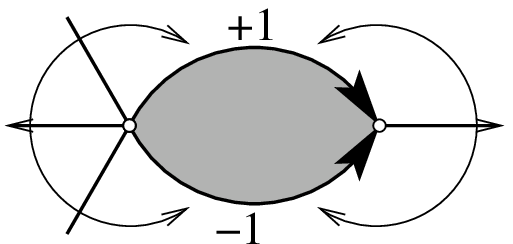 hscale=50 vscale=50 angle=0
hoffset=30  
voffset=0   
}
\begin{picture}(130,0)(-17,0)
\put(-60,-100) 
{
\begin{picture}(0,0)(0,0)
\put(54,116){\scriptsize\it odd} \put(147,116){\scriptsize\it
odd}
\end{picture}}
\end{picture}
& 
\special{                   
psfile=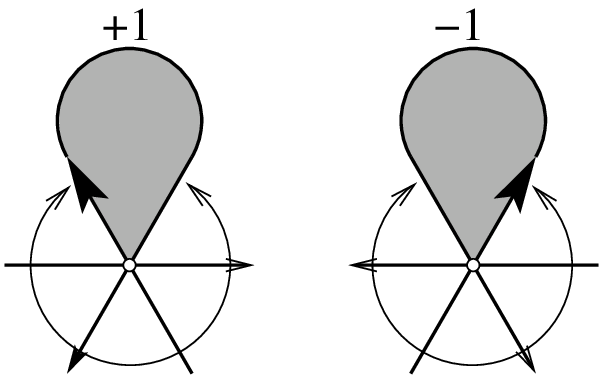 hscale=50 vscale=50 angle=0
hoffset=15  
voffset=0  
}
\begin{picture}(120,0)(-20,0) %
\put(-220,50) %
{
\begin{picture}(0,0)(0,0)
\put(223,-56){\scriptsize\it even} \put(272,-56){\scriptsize\it
even}
\end{picture}}
\end{picture}
\\
&\\
[-\halfbls]&\\
\hline & \\ 
&\\
   %
\includegraphics{valence21_plus.eps}
\begin{picture}(0,0)(0,0)
\put(-8,13){\scriptsize\it odd} \put(-8,-6){\scriptsize\it odd}
\end{picture}
&
\includegraphics{valence22_plus.eps}
\begin{picture}(0,0)(0,0)
\put(-30,-6){\scriptsize\it even} \put(10,-6){\scriptsize\it
even}
\end{picture}
\\
&\\
[-\halfbls]&\\
\hline
\end{array}
$$

$$ 
\begin{array}{|c|c|c|c|}
\hline\vspace{-1.8truept}&&&\\
[-\halfbls]
 +3.1 & +3.2a & + 3.2b & + 3.3 \\
   %
\hline &&& \\
[-\halfbls]&&&\\  &&& \\
&&&\\
\special{                     
psfile=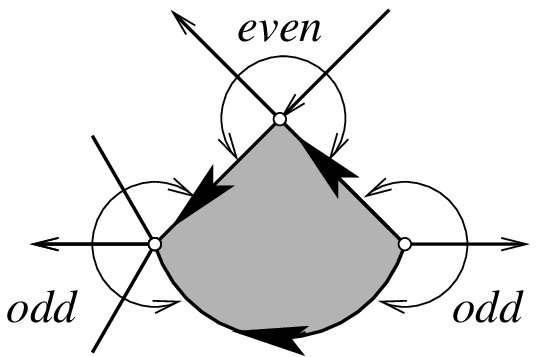
hscale=50 %
vscale=50 %
angle=0
hoffset=-2  
voffset=-20 %
}
\begin{picture}(76,0)(-38,-25) %
\put(-120,0) {
\begin{picture}(0,0)(0,0)
\put(99,-16){\tiny -1} \put(130,-16){\tiny -1}
\put(117,-50){\tiny +2}
\end{picture}}
\end{picture}
& 
\special{                   
psfile=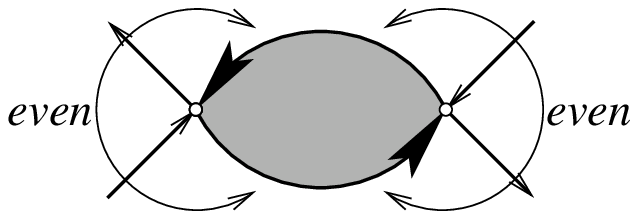 hscale=50 vscale=50 angle=0
hoffset=-5  
voffset=15 
}
\includegraphics{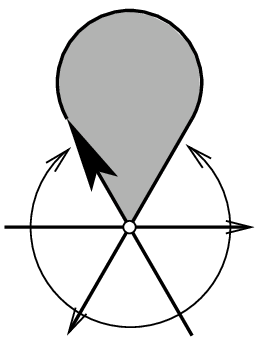}
\begin{picture}(83,0)(-35,-25) 
\put(-120,0) %
{
\begin{picture}(0,0)(-3,0)
\put(115,18){\tiny -1} \put(115,-12){\tiny -1}
\put(117,-68){\tiny +2}

\put(133,-24){\scriptsize\textit{even}}   %
\end{picture}}
\end{picture}
& 
\special{                     
psfile=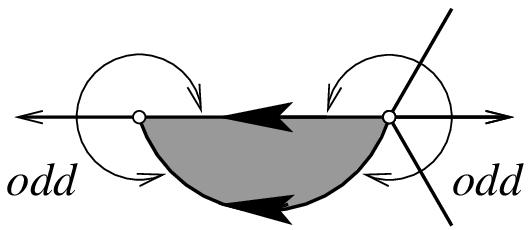 
hscale=50 vscale=50 angle=0
hoffset=2  %
voffset=15 %
}
\includegraphics{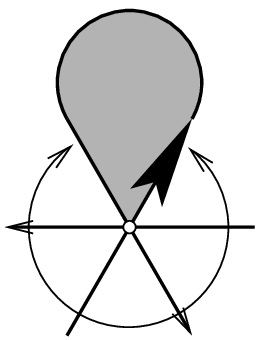}
\begin{picture}(83,0)(-33,-25) 
\put(-120,0) %
{
\begin{picture}(0,0)(-2,0)
\put(114,10){\tiny -1} \put(116,-14){\tiny +2}
\put(115,-65){\tiny -1}

\put(113,-29){\scriptsize\textit{even}}
\end{picture}}
\end{picture}
& 
\special{                       
psfile=holeright3.eps hscale=30 vscale=30 angle=0
hoffset=1 
voffset=5 
}
\includegraphics{holeright3.eps}
\includegraphics{holeleft3.eps}
\begin{picture}(78,0)(-34,-25) 
\put(-120,0) %
{
\begin{picture}(0,0)(0,0)
\put(90,12){\tiny -1} \put(139,12){\tiny -1} \put(116,-65){\tiny +2}

\put(87,-24){\scriptsize\textit{even}}
\put(137,-24){\scriptsize\textit{even}}
\put(110,-12){\scriptsize\textit{even}}
\end{picture}}
\end{picture}
\\
&&&\\
[-\halfbls]&&&\\
&&& \\ &&& \\
\hline &&& \\ 
   %
&&&\\
\includegraphics{valence31_plus.eps}
\begin{picture}(0,0)(10,-180) 
  \put(3,-173){\scriptsize\textit{even}}
  \put(-9,-192){\scriptsize\textit{odd}}
  \put( 17,-192){\scriptsize\textit{odd}}
\end{picture}
&
\includegraphics{valence32a.eps}
\begin{picture}(0,0)(109,-180)
  \put(102,-173.5){\scriptsize\textit{even}}
  \put( 89,-186){\scriptsize\textit{even}}
  \put(114,-193){\scriptsize\textit{even}}
\end{picture}
&
\includegraphics{valence32b.eps}
\begin{picture}(0,0)(205,-185)
  \put(199.5,-174){\scriptsize\textit{odd}}
  \put(188,-188){\scriptsize\textit{odd}}
  \put(209,-193){\scriptsize\textit{even}}
\end{picture}
&
\includegraphics{valence33_plus.eps}
\begin{picture}(0,0)(304,-186)
  \put(314,-186){\scriptsize\textit{even}}
  \put(280,-186){\scriptsize\textit{even}}
  \put(308,-195){\scriptsize\textit{even}}
\end{picture}
\\
&&&\\
&&&\\
   %
\hline\hline
\vspace{-1.8truept}&&&\\
[-\halfbls]
+4.1a & +4.1b & +4.2a & +4.2b \\
   %
\hline &&& \\
[-\halfbls]&&&\\  
&&&\\
\special{                  
psfile=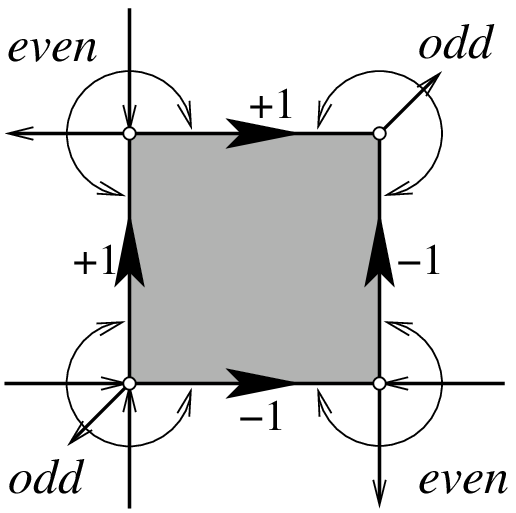 
hscale=45 vscale=45 angle=0
hoffset=-33 
voffset=-30 %
}
& 
\special{                   
psfile=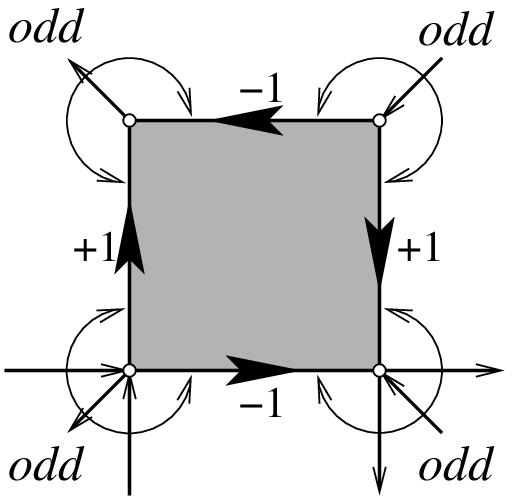 
hscale=45 vscale=45 angle=0
hoffset=-33 
voffset=-30 
}
& 
\special{                     
psfile=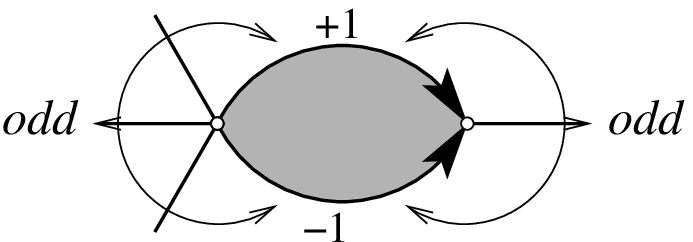 hscale=45 vscale=45 angle=0
hoffset=-45 
voffset=-32 %
} \includegraphics{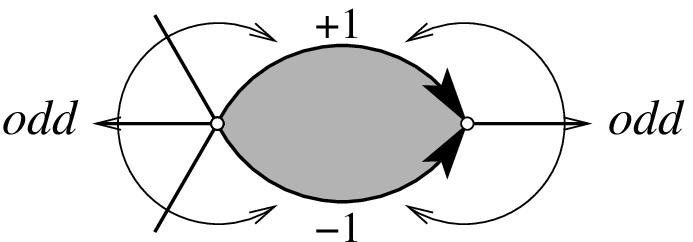}
& 
\special{                       
psfile=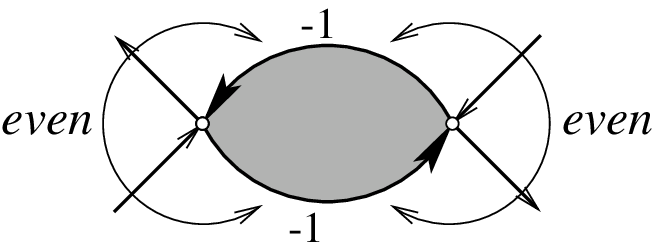 hscale=45 vscale=45 angle=0
hoffset=-42   
voffset=-32   
} \includegraphics{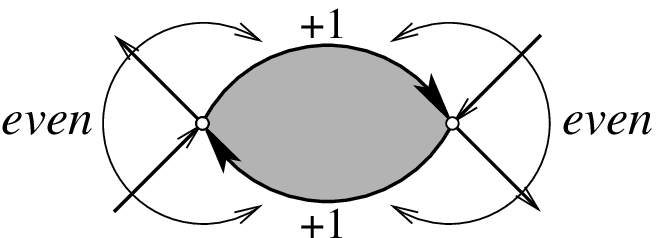}
\\
&&&\\
[-\halfbls]&&&\\
&&& \\ 
\hline &&& \\ 
&&&\\ &&&\\
\includegraphics{valence41a.eps} &
\includegraphics{valence41b.eps} &
\includegraphics{valence42a.eps} &
\includegraphics{valence42b.eps}
\\
&&&\\
&&&\\
&&&\\
[-\halfbls]&&&\\
   %
\hline\hline
\vspace{-2truept}&&&\\
[-\halfbls]
+4.2c & +4.3a & +4.3b & +4.4 \\
   %
\hline &&& \\
[-\halfbls]&&&\\  
&&&\\&&&\\&&&\\
\special{                    
psfile=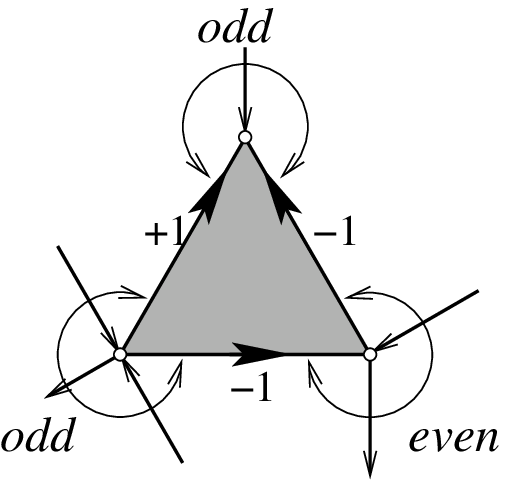       
hscale=45 vscale=45 angle=0
hoffset=-32 %
voffset=-44 %
}
\includegraphics{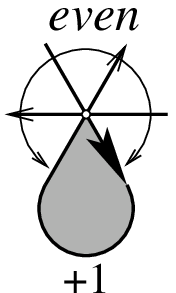}
& 
\special{                   
psfile=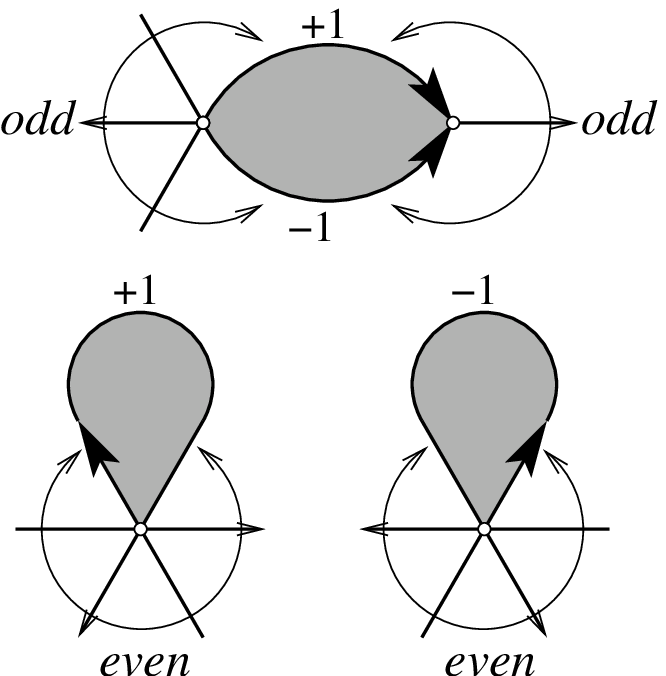 hscale=45 vscale=45 angle=0
hoffset=-44 
voffset=-36 %
}
& 
\special{                   
psfile=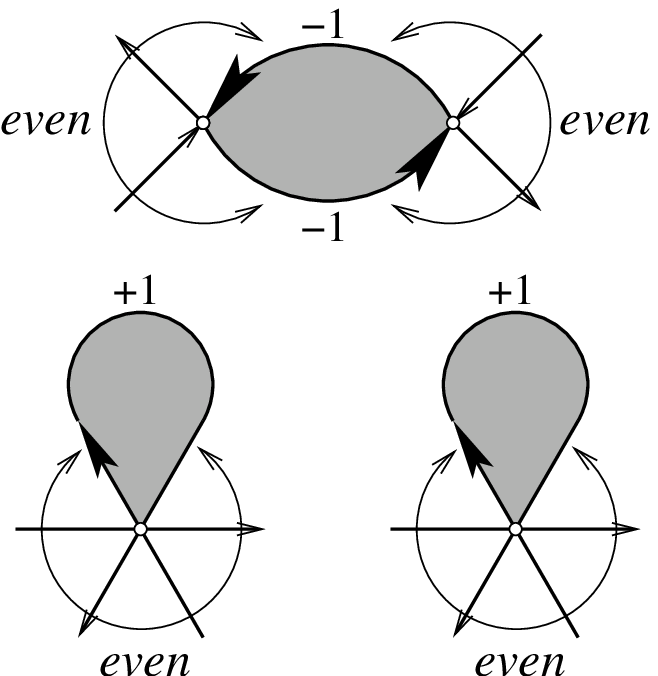 hscale=45 vscale=45 angle=0
hoffset=-42 %
voffset=-36 %
}
& 
\special{                      
psfile=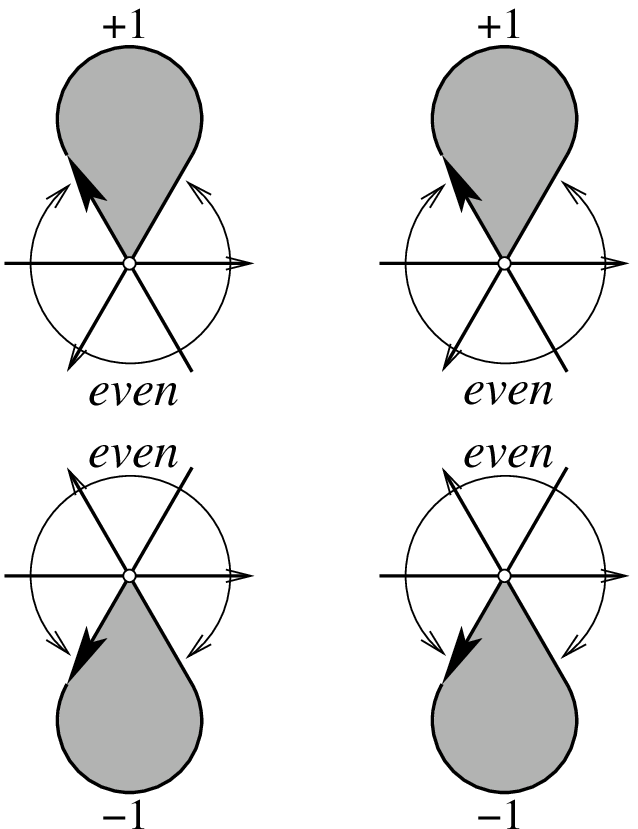 hscale=45 vscale=45 angle=0
hoffset=-40   %
voffset=-45   
}
\\
&&&\\&&&\\
[-\halfbls]&&&\\
&&& \\ 
\hline &&& \\ 
&&&\\ &&&\\
\includegraphics{valence42c.eps} &
\includegraphics{valence43a.eps} &
\includegraphics{valence43b.eps} &
\includegraphics{valence44_plus.eps}
\\
&&&\\
&&&\\
&&&\\
[-\halfbls]&&&\\
\hline
\end{array}
$$

By convention we  orient the boundary  of a surface  as  follows.
Choose an orthogonal frame  $(\vec  n, \vec\tau)$ equivalent to a
canonical frame in such way that $\vec n$ is the external normal,
and $\vec\tau$  is a vector tangent  to the boundary.  The vector
$\vec\tau$  defines   the   orientation   of  the  boundary  (see
figure~\ref{fig:example:of:ribbon:graph}).

In those pictures  above,  which represent surfaces with boundary
(and not the ribbon graphs), the shadowed regions represent small
holes  inside  a  flat  surface.   The   same   remark   concerns
figures~\ref{fig:4:1a}--\ref{fig:m:2:1:nonlocal}.

Our convention on  orientation implies that the boundaries of the
small holes  are oriented \textit{clockwise}. The same convention
implies  that  the  edges  of  the  graphs  $\G_v$  are  oriented
\textit{counterclockwise}.

Choosing a line element  in the tangent space to some point  of a
flat surface  one can transport  this line element to the tangent
space at any other point.  The  resulting  distribution defines a
foliation. For the surfaces with boundary under consideration the
foliation  can  be   chosen   parallel  to  all  boundary  saddle
connections.

When   the   flat  surface  has  trivial  linear  holonomy,   the
corresponding foliation is  orientable.  The arrows on the saddle
connections  in  the pictures above represent the orientation  of
the foliation and  \textit{not}  the canonical orientation of the
saddle connections induced from the canonical  orientation of the
boundary.     We     also     use     this     convention     for
figures~\ref{fig:breaking:up:a:zero}--\ref{fig:42c}.